\documentclass[reqno]{amsart}

\usepackage[english]{babel}
\usepackage{csquotes}
\usepackage[style=alphabetic,
            isbn=false,
            doi=false,
            url=false,
            backend=bibtex,
            maxnames=5, 
            maxalphanames=5
           ]{biblatex}
\bibliography{sources}
\usepackage{amsmath,amsfonts,amsthm,mathtools, amssymb}
\usepackage{xcolor}
\usepackage{tikz-cd}
\usepackage{hyperref}
\usepackage{aliascnt}
\usepackage{xparse}
\usepackage{tensor}
\usepackage{caption}
\usepackage{ mathrsfs }
\usepackage{tabto}
\usepackage{imakeidx}
\usepackage[normalem]{ulem}
\usepackage{comment}
\allowdisplaybreaks

 \RequirePackage{filecontents}

\indexsetup{level=\subsection*,toclevel=\paragraph,headers={ \MakeUppercase{Euler sequence on strata}}{\MakeUppercase{Euler sequence on strata}}, noclearpage,firstpagestyle=headings}
\makeindex[name=graph,title=Graphs and levels,columns=1,columnseprule,options=-s Idx.ist]

\usepackage{tikz}
\usetikzlibrary{matrix,arrows,calc}
\usetikzlibrary{positioning}
\usetikzlibrary{decorations.pathreplacing,decorations.markings,decorations.pathmorphing}
\usetikzlibrary{positioning,arrows,patterns}
\usetikzlibrary{cd}
\usetikzlibrary{intersections}

\tikzset{
	every loop/.style={very thick},
	comp/.style={circle,fill,black,inner sep=0pt,minimum size=5pt},
	order bottom left/.style={pos=.05,left,font=\tiny},
	order top left/.style={pos=.9,left,font=\tiny},
	order bottom right/.style={pos=.05,right,font=\tiny},
	order top right/.style={pos=.9,right,font=\tiny},
	order node dis/.style={text width=.75cm},
	circled number/.style={circle, draw, inner sep=0pt, minimum size=12pt},
	below left with distance/.style={below left,text height=10pt},
    below right with distance/.style={below right,text height=10pt}
	}

\makeatletter
\ifcsname phantomsection\endcsname
    \newcommand*{\@gobblenexttocentry}[9]{}
\else
    \newcommand*{\@gobblenexttocentry}[4]{}
\fi
\newcommand*{\addsubsection}{%
    \addtocontents{toc}{\protect\@gobblenexttocentry}%
    \subsection*}
\makeatother

\DeclareFieldFormat[article]{title}{{\it #1}}
\DeclareFieldFormat{journaltitle}{{\rm #1}}
\renewbibmacro{in:}{%
  \ifentrytype{article}{}{\printtext{\bibstring{in}\intitlepunct}}}
\DeclareFieldFormat[incollection]{title}{{\it #1}}
\DeclareFieldFormat{journaltitle}{{\rm #1}}

\begin{document}

\def\subsectionautorefname{Section}
\def\subsubsectionautorefname{Section}
\def\sectionautorefname{Section}
\def\equationautorefname~#1\null{(#1)\null}

\newcommand{\mynewtheorem}[4]{
  \if\relax\detokenize{#3}\relax 
    \if\relax\detokenize{#4}\relax 
      \newtheorem{#1}{#2}
    \else
      \newtheorem{#1}{#2}[#4]
    \fi
  \else
    \newaliascnt{#1}{#3}
    \newtheorem{#1}[#1]{#2}
    \aliascntresetthe{#1}
  \fi
  \expandafter\def\csname #1autorefname\endcsname{#2}
}

\mynewtheorem{theorem}{Theorem}{}{section}
\mynewtheorem{lemma}{Lemma}{theorem}{}
\mynewtheorem{rem}{Remark}{lemma}{}
\mynewtheorem{prop}{Proposition}{lemma}{}
\mynewtheorem{cor}{Corollary}{lemma}{}
\mynewtheorem{add}{Addendum}{lemma}{}
\mynewtheorem{definition}{Definition}{lemma}{}
\mynewtheorem{question}{Question}{lemma}{}
\mynewtheorem{assumption}{Assumption}{lemma}{}
\mynewtheorem{example}{Example}{lemma}{}


\def\defbb#1{\expandafter\def\csname b#1\endcsname{\mathbb{#1}}}
\def\defcal#1{\expandafter\def\csname c#1\endcsname{\mathcal{#1}}}
\def\deffrak#1{\expandafter\def\csname frak#1\endcsname{\mathfrak{#1}}}
\def\defop#1{\expandafter\def\csname#1\endcsname{\operatorname{#1}}}
\def\defbf#1{\expandafter\def\csname b#1\endcsname{\mathbf{#1}}}

\makeatletter
\def\defcals#1{\@defcals#1\@nil}
\def\@defcals#1{\ifx#1\@nil\else\defcal{#1}\expandafter\@defcals\fi}
\def\deffraks#1{\@deffraks#1\@nil}
\def\@deffraks#1{\ifx#1\@nil\else\deffrak{#1}\expandafter\@deffraks\fi}
\def\defbbs#1{\@defbbs#1\@nil}
\def\@defbbs#1{\ifx#1\@nil\else\defbb{#1}\expandafter\@defbbs\fi}
\def\defbfs#1{\@defbfs#1\@nil}
\def\@defbfs#1{\ifx#1\@nil\else\defbf{#1}\expandafter\@defbfs\fi}
\def\defops#1{\@defops#1,\@nil}
\def\@defops#1,#2\@nil{\if\relax#1\relax\else\defop{#1}\fi\if\relax#2\relax\else\expandafter\@defops#2\@nil\fi}
\makeatother

\defbbs{ZHQCNPALRVW}
\defcals{ABOPQMNXYLTRAEHZKCFI}
\deffraks{apijklgmnopqueRC}
\defops{IVC, PGL,SL,mod,Spec,Re,Gal,Tr,End,GL,Hom,PSL,H,div,Aut,rk,Mod,R,T,Tr,Mat,Vol,MV,Res,Hur, vol,Z,diag,Hyp,hyp,hl,ord,Im,ev,U,dev,c,CH,fin,pr,Pic,lcm,ch,td,LG,id,Sym,Aut,Log,tw,irr,discrep,BN,NF,NC,age,hor,lev,ram,NH,av,app,mid}
\defbfs{cuvzwp} 

\def\ep{\varepsilon}
\def\ve{\varepsilon}
\def\abs#1{\lvert#1\rvert}
\def\dd{\mathrm{d}}
\def\WP{\mathrm{WP}}
\def\inj{\hookrightarrow}
\def\eq{=}

\def\i{\mathrm{i}}
\def\e{\mathrm{e}}
\def\st{\mathrm{st}}
\def\ct{\mathrm{ct}}
\def\rel{\mathrm{rel}}
\def\odd{\mathrm{odd}}
\def\even{\mathrm{even}}

\def\uC{\underline{\bC}}
\def\ol{\overline}
  
\def\Vrel{\bV^{\mathrm{rel}}}
\def\Wrel{\bW^{\mathrm{rel}}}
\def\twolev{\mathrm{LG_1(B)}}

\def\be{\begin{equation}}   \def\ee{\end{equation}}     \def\bes{\begin{equation*}}    \def\ees{\end{equation*}}
\def\ba{\be\begin{aligned}} \def\ea{\end{aligned}\ee}   \def\bas{\bes\begin{aligned}}  \def\eas{\end{aligned}\ees}
\def\={\;=\;}  \def\+{\,+\,} \def\m{\,-\,}

\newcommand*{\proj}{\mathbb{P}}
\newcommand{\IVCst}[1][\mu]{{\mathcal{IVC}}({#1})}
\newcommand{\barmoduli}[1][g]{{\overline{\mathcal M}}_{#1}}
\newcommand{\moduli}[1][g]{{\mathcal M}_{#1}}
\newcommand{\omoduli}[1][g]{{\Omega\mathcal M}_{#1}}
\newcommand{\modulin}[1][g,n]{{\mathcal M}_{#1}}
\newcommand{\omodulin}[1][g,n]{{\Omega\mathcal M}_{#1}}
\newcommand{\zomoduli}[1][]{{\mathcal H}_{#1}}
\newcommand{\barzomoduli}[1][]{{\overline{\mathcal H}_{#1}}}
\newcommand{\pomoduli}[1][g]{{\proj\Omega\mathcal M}_{#1}}
\newcommand{\pomodulin}[1][g,n]{{\proj\Omega\mathcal M}_{#1}}
\newcommand{\pobarmoduli}[1][g]{{\proj\Omega\overline{\mathcal M}}_{#1}}
\newcommand{\pobarmodulin}[1][g,n]{{\proj\Omega\overline{\mathcal M}}_{#1}}
\newcommand{\potmoduli}[1][g]{\proj\Omega\tilde{\mathcal{M}}_{#1}}
\newcommand{\obarmoduli}[1][g]{{\Omega\overline{\mathcal M}}_{#1}}
\newcommand{\obarmodulio}[1][g]{{\Omega\overline{\mathcal M}}_{#1}^{0}}
\newcommand{\otmoduli}[1][g]{\Omega\tilde{\mathcal{M}}_{#1}}
\newcommand{\pom}[1][g]{\proj\Omega{\mathcal M}_{#1}}
\newcommand{\pobarm}[1][g]{\proj\Omega\overline{\mathcal M}_{#1}}
\newcommand{\pobarmn}[1][g,n]{\proj\Omega\overline{\mathcal M}_{#1}}
\newcommand{\princbound}{\partial\mathcal{H}}
\newcommand{\omoduliinc}[2][g,n]{{\Omega\mathcal M}_{#1}^{{\rm inc}}(#2)}
\newcommand{\obarmoduliinc}[2][g,n]{{\Omega\overline{\mathcal M}}_{#1}^{{\rm inc}}(#2)}
\newcommand{\pobarmoduliinc}[2][g,n]{{\proj\Omega\overline{\mathcal M}}_{#1}^{{\rm inc}}(#2)}
\newcommand{\otildemoduliinc}[2][g,n]{{\Omega\widetilde{\mathcal M}}_{#1}^{{\rm inc}}(#2)}
\newcommand{\potildemoduliinc}[2][g,n]{{\proj\Omega\widetilde{\mathcal M}}_{#1}^{{\rm inc}}(#2)}
\newcommand{\omoduliincp}[2][g,\lbrace n \rbrace]{{\Omega\mathcal M}_{#1}^{{\rm inc}}(#2)}
\newcommand{\obarmoduliincp}[2][g,\lbrace n \rbrace]{{\Omega\overline{\mathcal M}}_{#1}^{{\rm inc}}(#2)}
\newcommand{\obarmodulin}[1][g,n]{{\Omega\overline{\mathcal M}}_{#1}}
\newcommand{\LTH}[1][g,n]{{K \overline{\mathcal M}}_{#1}}
\newcommand{\PLS}[1][g,n]{{\bP\Xi \mathcal M}_{#1}}

\DeclareDocumentCommand{\LMS}{ O{\mu} O{g,n} O{}}{\Xi\overline{\mathcal{M}}^{#3}_{#2}(#1)}
\DeclareDocumentCommand{\Romod}{ O{\mu} O{g,n} O{}}{\Omega\mathcal{M}^{#3}_{#2}(#1)}

\newcommand*{\Tw}[1][\Lambda]{\mathrm{Tw}_{#1}}  
\newcommand*{\sTw}[1][\Lambda]{\mathrm{Tw}_{#1}^s}  

\newcommand{\HH}{{\mathbb H}}
\newcommand{\MM}{{\mathbb M}}
\newcommand{\bbC}{{\mathbb C}}

\newcommand{\bfa}{{\bf a}}
\newcommand{\bfb}{{\bf b}}
\newcommand{\bfc}{{\bf c}}
\newcommand{\bfd}{{\bf d}}
\newcommand{\bfe}{{\bf e}}
\newcommand{\bff}{{\bf f}}
\newcommand{\bfg}{{\bf g}}
\newcommand{\bfh}{{\bf h}}
\newcommand{\bfm}{{\bf m}}
\newcommand{\bfn}{{\bf n}}
\newcommand{\bfp}{{\bf p}}
\newcommand{\bfq}{{\bf q}}
\newcommand{\bft}{{\bf t}}
\newcommand{\bfP}{{\bf P}}
\newcommand{\bfR}{{\bf R}}
\newcommand{\bfU}{{\bf U}}
\newcommand{\bfu}{{\bf u}}
\newcommand{\bfx}{{\bf x}}
\newcommand{\bfz}{{\bf z}}

\newcommand{\bfl}{{\boldsymbol{\ell}}}
\newcommand{\bfmu}{{\boldsymbol{\mu}}}
\newcommand{\bfeta}{{\boldsymbol{\eta}}}
\newcommand{\bftau}{{\boldsymbol{\tau}}}
\newcommand{\bfomega}{{\boldsymbol{\omega}}}
\newcommand{\bfsigma}{{\boldsymbol{\sigma}}}

\newcommand{\wh}{\widehat}
\newcommand{\wt}{\widetilde}

\newcommand{\ps}{\mathrm{ps}}  
\newcommand{\CPT}{\mathrm{CPT}}  
\newcommand{\NCT}{\mathrm{NCT}}  
\newcommand{\RBD}{\mathrm{RBD}}
\newcommand{\ETD}{\mathrm{ETD}}
\newcommand{\OCT}{\mathrm{OCT}}
\newcommand{\SRT}{\mathrm{SRT}}
\newcommand{\RBT}{\mathrm{RBT}}

\newcommand{\tdpm}[1][{\Gamma}]{\mathfrak{W}_{\operatorname{pm}}(#1)}
\newcommand{\tdps}[1][{\Gamma}]{\mathfrak{W}_{\operatorname{ps}}(#1)}

\newlength{\halfbls}\setlength{\halfbls}{.8\baselineskip}

\newcommand*{\Teichmuller}{Teich\-m\"uller\xspace}

\DeclareDocumentCommand{\MSfun}{ O{\mu} }{\mathbf{MS}({#1})}
\DeclareDocumentCommand{\MSgrp}{ O{\mu} }{\mathcal{MS}({#1})}
\DeclareDocumentCommand{\MScoarse}{ O{\mu} }{\mathrm{MS}({#1})}
\DeclareDocumentCommand{\tMScoarse}{ O{\mu} }{\widetilde{\mathbb{P}\mathrm{MS}}({#1})}

\newcommand{\kmin}{\kappa_{(2g-2)}}
\newcommand{\ktop}{\kappa_{\mu_\Gamma^{\top}}}
\newcommand{\kbot}{\kappa_{\mu_\Gamma^{\bot}}}


\title[Kodaira dimension of moduli of Abelian differentials]
      {On the Kodaira dimension of moduli spaces of Abelian differentials}

\author{Dawei Chen}
      \thanks{Research of D.C. is supported by the National Science Foundation Grant DMS-2001040
        and the Simons Foundation Collaboration Grant 635235.}
\address{Department of Mathematics, Boston College, Chestnut Hill, MA 02467, USA}
\email{dawei.chen@bc.edu}

\author{Matteo Costantini}
\address{Universität Duisburg-Essen
Fakultät für Mathematik
45117 Essen, Germany}
\email{matteo.costantiuni-due.de}
  \thanks{Research of M.C.  is supported by the DFG Research Training Group 2553.}

\author{Martin M\"oller}
\address{Institut f\"ur Mathematik, Goethe-Universit\"at Frankfurt,
Robert-Mayer-Str. 6-8,
60325 Frankfurt am Main, Germany}
\email{moeller@math.uni-frankfurt.de}
\thanks{Research of M.M.\ is supported
by the DFG-project MO 1884/2-1, by the LOEWE-Schwerpunkt
``Uniformisierte Strukturen in Arithmetik und Geometrie''
and the Collaborative Research Centre
TRR 326 ``Geometry and Arithmetic of Uniformized Structures''.}

\begin{abstract}
This paper lays the foundation for determining the Kodaira dimension of
the projectivized strata of Abelian differentials with prescribed zero and
pole orders in large genus. We work with the moduli space of multi-scale
differentials constructed in~\cite{LMS} which provides an orbifold
compactification of these strata. We establish the projectivity of the moduli
space of multi-scale differentials, describe the locus of canonical
singularities, and compute a series of effective divisor classes. 
Moreover, we exhibit a perturbation of the canonical class which allows
the corresponding pluri-canonical forms to extend over the locus of
non-canonical singularities.
\par
As applications, we certify general type for strata with few zeros
as well as for strata with equidistributed zero orders when $g$ is sufficiently
large. In particular, we show general type for the odd spin components of the minimal strata for $g\geq 13$.
\end{abstract}
\maketitle
\tableofcontents

\section{Introduction}

Let $\bP\omoduli[g,n](\mu)$ be the moduli space of holomorphic (or meromorphic) Abelian differentials (up to scale) with labeled singularities of orders prescribed by a partition~$\mu = (m_1,\ldots,m_n)$ of $2g-2$. This space is called the (projectivized) stratum of Abelian differentials of type $\mu$. 
\par 
The study of $\bP\omoduli[g,n](\mu)$ is important for at least two reasons. On one hand, an Abelian differential induces a flat metric with conical singularities at its zeros such that the underlying Riemann surface can be realized as a polygon with edges pairwise identified by translations. Varying the shape of such polygons by affine transformations induces an action on the strata of differentials (called Teichm\"uller dynamics), whose orbit closures (called affine invariant subvarieties) govern intrinsic properties of surface dynamics.  On the other hand, an Abelian differential (up to scale) corresponds to a canonical divisor in the underlying complex curve. Hence the union of $\bP\omoduli[g,n](\mu)$ stratifies the (projectivized) Hodge bundle over the moduli space of curves, thus producing a number of remarkable questions to investigate from the viewpoint of algebraic geometry, such as compactification, enumerative geometry, and cycle class calculation. The interplay of these aspects has brought the study of differentials to an exciting new era (see e.g., \cite{zorichsurvey}, \cite{wrightsurvey}, \cite{chensurvey} as well as the references therein for an introduction to this fascinating subject).  
\par
Despite the aforementioned advances, not much is known about the birational geometry of $\bP\omoduli[g,n](\mu)$. This is the focus of the current paper. 
A fundamental birational invariant for a variety is the Kodaira dimension, which measures the growth rate of pluri-canonical forms and controls the size of the canonical model of the variety. When the variety has a modular interpretation, determining the Kodaira dimension is closely related to the boundary behavior, singularity analysis, and decomposition of the cone of effective divisors. The study of the Kodaira dimension and related structures has covered many classical moduli spaces and their variants (see e.g., \cite{HarrisMumford}, \cite{TaiKodAg}, \cite{HarrisKodII}, \cite{CR15}, \cite{eh23}, \cite{OG}, \cite{CRsmall}, \cite{kondo1}, \cite{LiKod}, \cite{kondo2}, \cite{logan}, \cite{GHSK3}, \cite{farkaseven}, \cite{FLPrym},  \cite{BFV}, \cite{FVJacobian}, \cite{FVdifference}, \cite{FVoddspin}, \cite{CMKV}, \cite{TVA}, \cite{petok}, \cite{FJP}, \cite{schwarzhyp}, \cite{schwarznodal}, \cite{BMCY}, \cite{AB}, \cite{FJP13}, \cite{FMK3}, \cite{FMHurwitz}). 
It is a general expectation that for sequences of moduli spaces the Kodaira dimension should be 
negative for small complexity (in terms of genus or level covering), but
it should become maximal for large complexity. This is known, e.g., for moduli spaces of
curves (\cite{HarrisMumford}), for moduli spaces of abelian varieties (\cite{TaiKodAg}), and for moduli spaces of K3 surfaces
(\cite{GHSK3}).
\par
The (projectivized) strata of holomorphic Abelian differentials are uniruled for low genus $g \leq 9$ and all zero types $\mu$ as well as for $g \leq 11$ if moreover the number of zeros~$n$ is large (\cite{Barros}, \cite{BudGonality}). On
the other hand, when $n \geq g-1$ these strata can be viewed as generically finite covers of the
moduli space of pointed curves (\cite{GendronThesis}) and thus of general
type for large genus.\footnote{Strictly speaking in \cite{GendronThesis} only the case $n = g-1$ was considered so as to obtain a generically finite map to $\moduli[g]$. But the same argument works for $n > g-1$ by projecting to $\moduli[g, n-g+1]$ and checking finiteness of the fiber over a boundary point parameterizing a general chain of elliptic curves with $n-g+1$ marked points in a tail.} However, in the case of few zeros the large genus behavior of the Kodaira dimension of the strata is wide open, which is one of the main motivations for the techniques developed in this paper.
\par
In order to study the Kodaira dimension of a (non-compact) moduli space, one often needs a good compactification.  The notion of multi-scale differentials from \cite{LMS} gives rise to two compactifications of the moduli stack $\bP\omoduli[g,n](\mu)$.
First, the stack $\bP\MSgrp$ of multi-scale differentials admits a local blowup description compared to the naive
'incidence variety compactification' (\cite{BCGGM1}). Second, there is the
smooth Deligne--Mumford stack $\bP\LMS$ with normal crossing boundary
divisors. They both have the same underlying coarse moduli space $\bP\MScoarse$.
We will recall aspects of the relevant constructions and quotient maps in
\autoref{sec:coarse}. 
\par
A standard method of showing general type is to write the canonical divisor class $K$ as the sum of an ample divisor class and an effective divisor class (i.e., to prove that $K$ is a big divisor class), where the existence of an ample divisor class already requires the underlying space to be projective. Note that  $\bP\LMS$ 
was constructed by a complex-analytic gluing approach, and
that the blowup construction for the stack $\bP\MSgrp$ is also local (i.e., a global ideal sheaf to be blown up is unknown in general).  
Such local operations might destroy the projectivity of the resulting complex-analytic
varieties (see e.g.,\ Hironaka's examples in~\cite[Appendix~B]{Hartshorne}).
Nevertheless, our first result below verifies the projectivity of $\bP\MScoarse$. 
\par
\begin{theorem} \label{thm:LMSisproj}
The coarse moduli space $\bP\MScoarse$ associated with the stack
of multi-scale differentials of type $\mu$ is a projective variety for any~$\mu$.
\end{theorem}
\par
In order to prove the above result, in \autoref{sec:projectivity} we explicitly give a linear combination of boundary divisors
that is relatively ample for the forgetful map from the multi-scale compactification to the incidence variety
compactification (where the incidence variety compactification is projective since it is the closure of the strata in the projective Hodge bundle). 
\par
The canonical class of the stack $\bP\LMS$ was computed in 
\cite[Theorem~1.1]{EC}. We then need to analyze the ramification divisor
of the map from the stack to the coarse moduli space $\bP\MScoarse$, both in the interior and at the boundary.
This is carried out in \autoref{sec:tocoarse}. We remark that for strata
of type $\mu = (m,2g-2-m)$ with $m$ even, the map to the coarse moduli space can actually have a  
ramification divisor in the interior.
\par
We are now in a position to run the aforementioned strategy of expressing the
canonical class as ample plus effective divisor classes. However, there are a number of 
new obstacles comparing to the work of Harris--Mumford for $\barmoduli[g]$ and subsequent
works. The first one of them is about canonical singularities. 
\par
\begin{theorem} \label{thm:introsingofMScoarse}
The interior $\bP\Omega \mathrm{M}_{g,n}(\mu)$ of the coarse moduli space of Abelian differentials 
 with labeled zeros and poles has canonical singularities for all signatures $\mu$ 
except for $\mu = (m, 2-m)$ in $g=2$ with $1\neq m \equiv 1 \mod 3$.
\par
In contrast, the coarse moduli space of multi-scale differentials 
$\bP\MScoarse$ has non-canonical singularities in the boundary for
all but finitely many~$g$.
\end{theorem}
\par
A significant part of this paper deals with these non-canonical singularities
and how to overcome their presence. By the Reid--Tai criterion, the absence
of non-canonical singularities is certified by bounding from below the age of
automorphisms acting on the tangent space of the moduli stack. At the boundary
of $\bP\LMS$  the tangent space decomposes into two parts, as recalled in
\autoref{sec:newage}. One part is the tangent space of the strata determined
by the vertices of the level graph at the corresponding boundary stratum.
In \autoref{nprop:can_interior} and~\autoref{nprop:age-vertex} we compile
tables listing the cases where the  age is small enough to allow for non-canonical
singularities. The other part in the tangent space describes the
opening of nodes in terms of level passages (see \autoref{sec:DMtaut} for
the background on level graphs).
\par
In order to control non-canonical singularities, 
as a preliminary step we recast the action of the stacky structure related groups   
on the level passages in terms of toric geometry. This was only implicitly described
in \cite{LMS}. In \autoref{sec:canonical} we explicitly determine the
cone and fan structure to encode this information, which allows to measure the
failure of a singularity from
being canonical.
\par
Next we define a non-canonical compensation divisor $D_{\NC}$ which is a linear
combination of boundary divisors given explicitly
in~\autoref{eq:non-canonical-term}. It allows us to prove the following
criterion (see \autoref{subsec:singgeneralities} for relevant definitions
and \autoref{sec:refofProp13} for a refinement for strata with
$\mu = (m,2g-2-m)$ that takes the ramification divisor into account).
\par
\begin{prop} \label{prop:GenTypeCrit}
For $g\geq 2$ and all $\mu$ except for $\mu = (m, 2-m)$ in $g=2$ with $1\neq m \equiv 1 \mod 3$, there exists 
an explicit effective divisor class $D_{\NC}$ such that pluri-canonical forms associated to the perturbed canonical class 
$K_{\bP\MScoarse} - D_{\NC}$ in the smooth locus of $\bP\MScoarse$ extends completely in a desingularization. 
\par
In particular if one can write
\be \label{eq:KDAE}
K_{\bP\MScoarse} - D_{\NC} \= A + E
\ee
with~$A$ an ample divisor class and $E$ an effective divisor class,
then~$\bP\MScoarse$ is a variety of general type.
\end{prop}
\par
The description of $D_{\NC}$ builds on fairly delicate statements about
automorphisms of small age acting on tangent spaces (see
\autoref{nprop:ageort} and~\autoref{prop:compensation}). Moreover, the
description of $D_{\NC}$ reflects a strong tension: its divisor class has to
be large enough to compensate non-canonical singularities, while it cannot
be too large to make the desired expression~\autoref{eq:KDAE} unrealizable
(see \autoref{rem:nc-optimal}). Therefore, the ideas and techniques
involved in the study of $D_{\NC}$ can be of independent interests for
applications to other moduli spaces with non-canonical singularities.  
\par
\medskip
We now turn to the class computation of effective divisors suitable
to fulfill~\autoref{eq:KDAE}. The first type of divisors are pullbacks
of a series of effective divisors from $\barmoduli[g,n]$ (such as divisors
of Brill--Noether type). The formulas for such pullbacks and the conversion
formulas between standard divisor classes $\lambda_1$,
$\kappa_1$ and $\xi = c_1(\cO(-1))$ on $\bP\LMS$ are provided in
\autoref{sec:taut}. We remark that the verification of some pullback divisors
not containing the entire stratum requires non-trivial degeneration techniques (e.g., using curves of non-compact type 
in the proof of~\autoref{lem:Hurclass}). The second type of divisors, which are not induced via pullback, will be
called \emph{generalized Weierstrass divisors} and they can be defined for
  all connected
components of the strata except hyperelliptic and even spin components.
For a partition $\alpha = (\alpha_1, \ldots, \alpha_n)$ of $g-1$ such that
$0\leq \alpha_i \leq m_i$ for all~$i$ we define the divisor in the interior as
\be
W_{\mu}(\alpha) \= \bigl\{ (X, \bfz, \omega) \in \bP\omoduli[g,n](\mu)\,:\,
h^0\bigl(X, \alpha_1 z_1 + \cdots + \alpha_n z_n \bigr) \geq 2\bigr \}\, 
\ee
in~$\bP\omoduli[g,n](\mu)$. We will define the generalized Weierstrass divisor in $\bP\LMS$ via a Porteous type setting and study its class, as well
as its boundary behavior, in \autoref{sec:genWP}. Indeed we will use a novel twisted version of Porteous' formula in order to reduce 
extraneous contributions from the boundary divisors.    
\par 
Using these divisor classes we can prove our main results on Kodaira dimensions for strata of various types. We start with the case when there is a unique zero.\footnote{A point $z$ is called subcanonical if $(2g-2)z$ is a canonical divisor. Subcanonical points are among the most special points in algebraic curves. For $g\geq 4$ the locus of subcanonical points in $\mathcal M_{g,1}$ consists of three components, hyperelliptic, odd spin and even spin by~\cite{kozo1}. The hyperelliptic component of subcanonical points parameterizes Weierstrass points in hyperelliptic curves, hence it is a rational variety. In contrast, \autoref{intro:minimal} reveals that the spin components of subcanonical points can behave very differently.}  
\par
\begin{theorem} \label{intro:minimal}
The odd spin components of the coarse moduli spaces
of Abelian differentials with a unique zero $\bP\Omega \mathrm{M}_{g,1}(2g-2)$ are of general type for $g \geq 13$.
\end{theorem}
\par
Recall that the odd spin components of the minimal strata are known to be uniruled for $g\leq 9$.
Hence in this case our result is nearly optimal. For the remaining cases
of $g=10$, $11$ and~$12$, one needs either a more extremal effective divisor
class or a certificate of non-general type. We remark that the proof of this theorem
for $13 \leq g \leq 44$ relies on a computer verification of the
constraints given by \autoref{prop:GenTypeCrit} (see the end of
\autoref{sec:proofandalgo} for the description of the algorithm).
\par
Next we treat the case of strata with 'few zeros'. 
\par
\begin{theorem} \label{intro:fewzero}
Given a constant $M$, consider all holomorphic signatures $\mu = (m_1, \ldots, m_n)$
with even entries such that $m_i \geq M$ for all $i$ and $n \leq 5(M+1)$. Then
the odd spin components of the strata $\bP\Omega \mathrm{M}_{g,n}(\mu)$ are of general type for
all but finitely many such~$\mu$. 
\par
In particular (for $M = 1$), the odd spin
components of the strata $\bP\Omega \mathrm{M}_{g,n}(\mu)$ with at most $10$ zeros are of general
type for all but finitely many $\mu$.  
\par 
Moreover for holomorphic signatures $\mu = (m_1, m_2)$ with two odd entries,
the (non-hyperelliptic) strata $\bP\Omega \mathrm{M}_{g,n}(\mu)$ are of general type for all
but finitely many~$\mu$.
\end{theorem}
\par
The above statement is aimed at simplicity. The precise condition on
'few zeros' for which we prove the theorem is given
in \autoref{subsec:fewzero}. Note that the meaning of 'few zeros' is
relative, e.g., an integer tuple close to 
$(\sqrt{g}^{\sqrt{g}}, g-2)$ with approximately $\sqrt{g} + 1$ zeros is indeed
a signature of 'few zeros'.  On the other hand, the fact that for strata with
zeros of odd order the range of our result is more limited is due to the constraint
of the parameters $\alpha_i$ being integers in the definition of 
generalized Weierstrass divisors (since a natural choice for $\alpha_i$ is $m_i / 2$).  
\par
Finally for strata with many zeros, our method can also be applied to the
following zero type. We say that  a stratum is equidistributed if the zero orders
are all the same, i.e., $\mu = (s^n)$ with $n$ entries of equal value $s$.  
\par
\begin{theorem} \label{intro:stothen}
All but finitely many of connected equidistributed strata $\bP\Omega \mathrm{M}_{g,n}(s^n)$ (and the 
odd spin components in the disconnected case) are of general type. 
\end{theorem}
\par
The above result is again stated for simplicity rather than completeness.
For instance, the result also holds for nearly equidistributed strata when $\mu$ is close to $(s^n)$, e.g.,\
for strata of type $(s-1,s+1,s^{n-2})$ with $n$ large enough. All these
results are proven by using divisors of Brill--Noether type and the
generalized Weierstrass divisor for $\alpha = \mu/2$ (or a rounding of $\mu/2$ if there are zeros of odd order).
\par
Contrary to the case of constructing effective divisors with low slope in $\barmoduli$, a new phenomenon we have discovered is that it does not suffice to control
the usual slope involving the boundary divisor $\Delta_{\irr}$ (whose analogue in $\bP\LMS$ is the 
horizontal boundary divisor $D_h$). Instead, even for $\mu = (s^n)$ 
more boundary divisors are critical, e.g., boundary divisors
consisting of 'vine curves' with two vertices and various numbers of edges.  
These boundary divisors impose tight bounds on the convex combination of divisors of Brill--Noether type and the generalized Weierstrass
divisor 
for constructing the desired effective divisor class~$E$ in~\autoref{eq:KDAE}.
\par 
\medskip
This paper opens the gate for exploring comprehensively the birational geometry of moduli spaces of differentials. In what follows we elaborate on further directions.  
\par 
First, note that strata of holomorphic differentials with very unbalanced zero orders (such as $\mu=(g-2,2^{g/2})$) are not covered by the current method of 
using a single generalized Weierstrass divisor (combined with a divisor of Brill--Noether type), 
which we will explain in \autoref{sec:gentype}. Nevertheless, we have obtained evidences that using a mixed version of generalized Weierstrass divisors might work
(by varying the parameters $\alpha$ and taking a weighted average of all generalized Weierstrass divisors). The remaining challenges are the choice of weights, rounding $\alpha$ to be an integer tuple, and the combinatorial complexity of estimating.  
\par
Next, we have excluded the even spin components of the strata.  This is due to the construction of the generalized Weierstrass divisor, e.g., for the minimal strata when $(2g-2)z$ is a canonical divisor of even spin, i.e., $h^0(X, (g-1)z) = 2$ (or a higher even number), the locus of $h^0(X, (g-1)z) > 1$ used for defining the generalized Weierstrass divisor would contain entirely the even spin component.  A revised approach is to quotient out the (generically) two-dimensional subspace $H^0(X, (g-1)z)$ in the setting of Porteous' formula for the Hodge bundle. The remaining issue is caused by the locus where this subspace jumps dimension and hence the quotient can fail to be a vector bundle.  Nevertheless, we expect that using certain blowup of this locus can help extend and complete the desired calculation.     
\par
Moreover, this paper deals exclusively with the strata where the zeros are marked. When there are zeros of the same order, 
an unmarked stratum is a finite quotient of the corresponding marked stratum induced by permuting the marked points of the same
order, which can thus have distinct birational geometry. For instance, the unmarked principal stratum with $2g-2$ simple zeros is uniruled for all~$g$ because it
is an open dense subset of the (projectivized) Hodge bundle, while the marked principal strata are of general type for large $g$.  
Many ideas and techniques in this paper can readily be adapted
to treat the unmarked strata, e.g.,~in the singularity analysis in
\autoref{nprop:can_interior} leading to
\autoref{thm:singofMScoarse} we also consider the case of unmarked zeros and poles (since those with the same order can appear as indistinguishable edges in a level graph). 
\par
Another generalization is for strata of meromorphic differentials.  
This paper paves the way to treat them as well, e.g., the ramification
divisor in \autoref{sec:coarse}, the projectivity of
$\bP\MScoarse[\mu]$ in \autoref{sec:projectivity}, and the singularity
analysis in \autoref{sec:newage} and \autoref{sec:canonical} cover the
meromorphic case, too. It is interesting to note that the behavior of general type for
meromorphic strata starts as early as for $g=1$ from the corresponding geometry of modular curves.  
\par
Finally, our calculation for the classes of pullback divisors and generalized Weierstrass divisors provides the first step towards understanding the effective cone of $\bP\MScoarse[\mu]$, which together with the ample divisor class we constructed can shed light on the chamber decomposition of the effective cone and other birational models of $\bP\MScoarse[\mu]$. We plan to treat these questions in future work.  
\par 
\subsection*{Acknowledgments} We thank Jonathan Zachhuber
for his support in using the {\tt diffstrata}-package that led to the proof
of projectivity. We also thank Ignacio Barros, Andrei Bud, and Gavril Farkas for helpful discussions.  



\section{From the stack to the coarse moduli space} \label{sec:coarse}

In this section we first recall some background on the geometry of
the moduli stack of multi-scale differentials and its
coarse moduli space. Our main
goal here is to determine in \autoref{prop:ramifInt} and
in \autoref{prop:ramifBound} the ramification loci of the
map from the stack of multi-scale differentials to the coarse moduli space. 
Throughout we use $\CH^\bullet(\cdot)$ to denote the Chow groups with {\em rational} coefficients. 

\subsection{The smooth Deligne--Mumford stack $\bP\LMS$ and its boundary structure}
\label{sec:DMtaut}

We recall some notation from \cite{LMS} and \cite{EC}, summarizing
notions of multi-scale differentials and enhanced level graphs. For simplicity 
we often abbreviate $B=\bP\omoduli[g,n](\mu)$ and $\ol{B}= \bP\LMS$. We will also denote by $\varphi\colon \bP\LMS \to \bP\MScoarse$ the map from the smooth Deligne--Mumford
stack to the coarse moduli space.
\par
Boundary strata in~$\bP\LMS$ are encoded by \emph{enhanced level graphs}, which by definition are dual graphs of stable curves together with additional data. They are provided with a level structure, i.e.\ a total order on the vertices
with equality permitted. Edges between vertices on the same level are
called \emph{horizontal}, and they are called \emph{vertical} otherwise.
Usually the top level is labeled by zero, and the levels below are labeled by
consecutive negative integers. We also refer to the edges starting above level~$-i$
and ending at or below level~$-i$ as the edges crossing the \emph{i-th level passage}. 
An \emph{enhancement} is an assignment of an integer~$p_e \geq 0$ to each edge, with $p_e = 0$ if and only if the
edge is horizontal.\footnote{These enhancements were denoted by 
$\kappa_e$ in \cite{LMS}. We avoid that notation in view of the clash with constants
derived from $\kappa$-classes. Our symbol reflects that these are the prongs of the differentials, comparing to \cite{CMSZ} where the same
notation was used, but called 'twist'.}
The enhancement encodes the number of prongs (real
positive rays) emanating from the zeros and poles the multi-scale
differentials have at the branches of the nodes corresponding
to~$e$. We usually omit 'enhanced' for level graphs. 
\par
Throughout the paper we rely on the nice boundary combinatorics of
the moduli stack of multi-scale differentials. Many computations
happen on boundary strata without horizontal nodes. We denote by
$\LG_L(\ol{B})$ the set of enhanced level graphs with $L$ levels
below the top level and no horizontal nodes. We will also use the notation $\LG_L(\mu)$ in order to emphasize the signature, or simply $\LG_L$ when it is clear in the context. These level graphs correspond to subvarieties in the boundary of~$\ol{B}$ with codimension $L$ in~$\ol{B}$. We denote by~$D_\Gamma$ the
closed subvariety corresponding the boundary stratum with level graph~$\Gamma$
together with its degenerations. Each level of an enhanced level
graph defines a generalized stratum (\cite[Section~4]{EC}).
We denote by $d^{[i]}_\Gamma $ the projectivized dimension of level~$i$.
The normal crossing boundary structure implies that $\dim(\ol{B})
= L + \sum_{i=-L}^0 d^{[i]}_\Gamma$ for every level graph~$\Gamma$ with $L$ levels
below zero.
\par
Adjacency of boundary strata is encoded by the undegeneration map
\be \label{eq:bdundegeneration}
\delta_{i_1,\dots,i_n} \colon \LG_L(\ol{B})\to \LG_{n}(\ol{B})\,,
\ee
which contracts all the passage levels of a non-horizontal level graph
$\Gamma$ except for the passages between levels $-i_{k}+1$ and $-i_{k}$ for $k=1,\dots,n$. With this notation, $D_\Gamma \in \LG_L(\ol{B})$ is
a union (due to prong-matchings) of connected components of the intersection
of $\delta_{j}(D_\Gamma)$ for $j=1,\ldots,L$. For $I = \{i_1,\dots,i_n\}$ we also define $\delta_I^{\complement} = \delta_{I^{\complement}}$ for notation convenience. 
\par
Next we recall the notion of a \emph{multi-scale differential}. This
is a tuple $(X,\bfz,\bfomega,\bfsigma,\Gamma)$ consisting of a pointed
stable curve $(X,\bfz)$, a level graph~$\Gamma$, a twisted differential $\bfomega$ 
compatible with~$\Gamma$ and a collection of prong-matchings~$\bfsigma$.
Here a twisted differential is a collection of differentials for each vertex 
of~$\Gamma$ with zeros and poles as prescribed by the marked points and
enhancements, subject to the residue conditions, as given more precisely
in \cite{LMS}. The prong-matching~$\bfsigma$ is a collection of prong-matchings
for each vertical edge, i.e.\ an orientation-reversing bijection of the
prongs of the differential at the branches of the corresponding node. For simplicity 
we often omit certain part of the tuple  $(X,\bfz,\bfomega,\bfsigma,\Gamma)$ when it is
clear from the context.
\par
The stack $\bP\LMS$ parameterizes equivalence classes of multi-scale differentials,
where two equivalent multi-scale differentials differ by the action of the \emph{level rotation
torus}. This is a (multiplicative) torus that acts simultaneously by
rotating the differential and turning the prong-matching. The level rotation torus
should be considered as the quotient of its universal covering $\bC^L$ by the
subgroup that fixes the differential on each level and brings all prongs
back to themselves, where the subgroup is called the \emph{twist group} and denoted by $\text{Tw}_\Gamma$. Not all
elements in this group can be written as a product of twists that act on one
level passage only, and those that can form the \emph{simple twist group}
$\text{Tw}^s_\Gamma$ which is important since the normal crossing boundary structure
of $\bP\LMS$ stems from compactifying each level passage. The quotient group
$K_\Gamma = \text{Tw}_\Gamma/\text{Tw}^s_\Gamma$ is thus part of the stack structure
as we see in the sequel. We call the elements of $K_\Gamma$ the \emph{ghost
automorphisms}. We will frequently use that if $\Gamma \in \LG_1$
or if $\Gamma$ has only horizontal nodes, then $K_\Gamma = \{e\}$ is
trivial by definition.

\subsection{Coordinates at the boundary}
\label{sec:Coordbd}

Recall from \cite{LMS} that a coordinate system near the boundary is given by
perturbed period coordinates. Consider a boundary stratum $D_\Gamma$ with~$L$
levels below zero, possibly also with horizontal nodes. Then the perturbed
period coordinates around a multi-scale differential $(X,\bfomega,\bfz,\bfsigma)$
compatible with~$\Gamma$ can be described as a product of three groups of coordinates:
\begin{itemize}
	\item A parameter $t_i$ parameterizing the opening-up of the level passage
	above level~$-i$. We group them to a point $\bft = (t_i) \in \bC_{\lev} \cong
	\bC^{L}$.
	\item The level-wise projectivized period coordinates $\bP
	H_\rel^1(X_{(-i)})^{\frakR_i}$ of the subsurfaces $X_{(-i)}$ on each level,
	where $\frakR_i$ is the constraint imposed by the global residue condition
	to level~$-i$. We define $\bA_\rel(X_{(-i)})\subseteq \bP
	H_\rel^1(X_{(-i)})^{\frakR_i}$ to be an affine chart containing the image of
	the level-wise flat surfaces $(X_{(-i)},[\omega_{(-i)}])$ under the  level-wise
	period coordinates, and denote by 
	\begin{equation}\label{eq:affinechart}
		\bA_\rel(X)= \prod_{i=0}^L \bA_\rel(X_{(-i)})
	\end{equation}
	the product of the level-wise affine charts.
	\item A parameter $x_i$ for each horizontal node. We group them to a point
	$\bfx = (x_j) \in \bC_{\hor} \cong \bC^{H(\Gamma)}$ where $H(\Gamma)$ is the number
	of horizontal edges of~$\Gamma$.
\end{itemize}
\par
We revisit the map $\varphi\colon \bP\LMS \to \bP\MScoarse$ near~$D_\Gamma$.
It can be factored first as a quotient by the group of ghost
automorphisms~$K_\Gamma$ and then by the group $\Aut(X,\bfomega)$.
An important conclusion from the construction of $\bP\LMS$ is that the action
of~$K_\Gamma$ on period coordinates is on $\bC_{\lev}$ only and that $\Aut(X,\bfomega)$
maps prongs to prongs and thus acts on $\bC_{\lev} / K_\Gamma$. As a result, 
a neighborhood~$U$ of $(X,\bfomega,\bfz,\bfsigma)$ can be described as
\be \label{eq:bdlocalstructure}
(\bA_\rel(X) \times  \bC_{\hor} \, \times \, \bC_{\lev} / K_\Gamma )
/ \Aut(X,\bfomega)\,\,\cong \,\,U \quad  \subset \quad \bP\MScoarse\,.
\ee
\par
\begin{rem} \label{rem:notsemidirect}
	{\rm The exact sequence
		\bes
		0 \to K_\Gamma \to \mathrm{Iso}(X,\bfomega) \to \Aut(X,\bfomega) \to 0
		\ees
		describing the isotropy group in $\bP\LMS$ of a multi-scale differential $(X,\bfomega)$
		does not split in general. In particular, it is not a semidirect product in general.
		Consider for example a triangle graph (with three levels and one vertex on each level) with the prong $p_e=2$ on the
		long edge~$e_2$ and $p_e=1$ on the
		short edges~$e_1$ and $e_3$. Standard coordinates $x_i$ and $y_i$ that put the differentials 
		at the upper and lower end of the edges in normal form are related to the level
		parameters~$t_i$ by (\cite[Equation~12.5]{LMS})
		\bes
		x_1y_1 = t_1^2, \qquad x_2y_2 = t_1t_2, \qquad x_3y_3 = t_3^2\,.
		\ees
		In this case $K_\Gamma \cong \bZ/2$, comparing \cite[Example~3.3]{EC}, as we will also retrieve
		in the sequel. Consider an automorphism of order two on the middle level,
		that acts by $y_1 \mapsto -y_1$ and $x_3 \mapsto -x_3$ while fixing the vertices on
		the other levels. This can be easily realized by a hyperelliptic involution, which
		moreover acts trivially on $\bA_\rel(X)$ thanks to level-wise projectivization and $\bC_{\hor}$
		is void here. The lifts of this action to an action on $\bA_\rel(X) \times  \bC_{\hor}
		\times \bC_{\lev}$, i.e.\ to elements of $\mathrm{Iso}(X,\bfomega)$ are given
		by the action on $\bC_{\lev}$, namely by
		\bes
		t_1 \mapsto \zeta_4^a t_1, \qquad t_2 \mapsto \zeta_4^{4-a} t_2, \qquad a \in \{1,3\}\,. 
		\ees
		The cases of $a=1$ and $a=3$ both have order four and differ
		by the action of the non-trivial element in $K_\Gamma$, thus ruling out the possibility of a splitting. 
}\end{rem}
\par
Here we analyze the action of $\Aut(X,\bfomega)$ on $\bC_{\lev} / K_\Gamma$, in the
simple case that $\Gamma$ has two levels, which implies
that $K_\Gamma$ is trivial. We recall the essential step of the plumbing
construction from \cite[Section~12]{LMS}. At the upper and lower ends of each
edge~$e$ of~$\Gamma$ we choose one pair (of the $p_e$ possible choices) of 
coordinates~$x_e$ and~$y_e$ that puts the level-wise components~$\omega_{(0)}$
and~$\omega_{(-1)}$ of~$\bfomega$ in standard form and such that the collection
of local prong-matchings $(dx_e \otimes dy_e)_{e \in E(\Gamma)}$ represents the
global prong-matching~$\bfsigma$. Then the surfaces in a neighborhood are given
by gluing in the plumbing fixture $x_ey_e =t^{m_e}$ where $m_e = \ell_\Gamma/p_e$ with $\ell_\Gamma = \lcm (p_e)_{e \in E(\Gamma)}$. 
\par
Suppose $\tau$ is an automorphism of $(X,\bfomega,\bfz,\bfsigma)$, say
mapping the edge~$e'$ to~$e$ (with $p_e = p_{e'}$). Then $\tau^*x_e$ is a coordinate near the upper end
of~$e'$, which puts $\tau^* \omega_{(0)} = \zeta_{(0)}\omega_{(0)}$ in standard
form. Consequently $\tau^*x_e/x_{e'} = \zeta_{e^+}$ for some root of unity
$\zeta_{e^+}$ with $\zeta_{e^+}^{p_e} = \zeta_{(0)}$. Similarly,
$\tau^* y_e/y_{e'} = \zeta_{e^-}$ for some root of unity~$\zeta_{e^-}$
with $\zeta_{e^-}^{-p_e} = \zeta_{(-1)}$. The hypothesis that $\tau$ fixes the
equivalence class of the prong-matching~$\bfsigma$ implies that there is some
$c \in \bC$ such that $\tau^*x_e \cdot \tau^* y_e = (c t)^{m_e}$ for all~$e \in E(\Gamma)$.
This~$c$ is in fact a root of unity and describes the action of~$\tau$ on the
coordinate~$t$ transverse to the boundary. If it exists, $c$ is uniquely
determined by
\be \label{eq:zetacond2lev}
\zeta_{e^+} \cdot \zeta_{e^-} \= c^{m_e} \quad \text{for all} \quad e \in E(\Gamma)\,.
\ee
\par

\subsection{Ramification from the stack to the coarse moduli
  space} \label{sec:tocoarse}

In this section we are mainly interested in the map 
$\varphi\colon \bP\LMS \to \bP\MScoarse$ from the smooth Deligne--Mumford
stack to the coarse moduli space. We want to import intersection
theory computations from \cite{EC} on $\bP\LMS$ and then pass to 
study the birational geometry of~$\bP\MScoarse$. We thus need to study
the ramification divisor of this map.
\par
We will also consider the factorization $\varphi = \varphi_2 \circ \varphi_1$
where $\varphi_1\colon \bP\LMS \to \bP\MSgrp$ is the map to the orderly blowup
constructed in \cite{LMS} and where $\varphi_2\colon \bP\MSgrp \to \bP\MScoarse$ is the map
to its coarse moduli space. Since $\varphi_1$ is locally given by the map
$[U/K_\Gamma] \to U/K_\Gamma$, where $U$ is a neighborhood of a generic point in $D_\Gamma$ and $K_\Gamma$ was defined above as the group of ghost automorphisms, it implies
that $\bP\LMS$ and $\bP\MSgrp$ have the same coarse moduli space
and we thus get a factorization as claimed. Since $K_\Gamma$
is trivial if $\Gamma \in \LG_1$ or if $\Gamma$ has horizontal edges only, the map $\varphi_1$ has
no ramification divisor.
\par
Next we focus on $\varphi_2$. Recall that $\bP\MSgrp$ is locally obtained
as the normalization of the blowup of an ideal sheaf (the 'orderly blowup',
see \cite[Section~7]{LMS}) in the normalization of the incidence variety compactification, which by definition is the closure of strata in the projective Hodge bundle.  Hence the isomorphism
groupoids of $\bP\MSgrp$ are contained in the isomorphism groupoids of $\barmoduli[g,n]$.
They are given by the automorphism groups of pointed stable curves
that respect the additional data encoded in the enhanced level graph,
i.e., the enhancements need to be taken into account additionally.
\par
From now on we will often encounter the notion of hyperelliptic 
differentials $(X, \omega)$ where $X$ is hyperelliptic and~$\omega$ is
anti-invariant under the hyperelliptic involution. We remark that this notion is
stronger than only requiring~$X$ to be a hyperelliptic curve. Moreover,
a hyperelliptic component of a stratum means that the locus of hyperelliptic
differentials forms a connected component of the stratum. Since in our setup the zeros and poles are labeled, among all hyperelliptic components only the ones for 
$\mu = (2g-2)$ (holomorphic) or $(2g-2+2m, -2m)$ with $m > 0$ (meromorphic) have the hyperelliptic involution as a non-trivial automorphism for a generic differential (up to sign) contained in them (in contrast for e.g. $\mu = (g-1, g-1)$ the hyperelliptic involution of a generic element in the hyperelliptic component swaps the two zeros). For this reason when analyzing the ramification of the map $\varphi$, we will exclude these special hyperelliptic components (whose birational geometry is much better known anyway, e.g. being unirational). 
\par
The result below describes the ramification divisor of the map $\varphi\colon \bP\LMS \to \bP\MScoarse$ in the interior of
the strata. 
\par
\begin{prop} \label{prop:ramifInt}
Suppose $\mu$ is of holomorphic type (and the hyperelliptic component is excluded if $\mu = (2g-2)$). Then the ramification
divisor of the map $\varphi$ in the interior of the stratum is empty unless $\mu = (m,2g-2-m)$ consists of two zeros of even order (i.e. $m$ is even).
In this case the ramification divisor in the interior arises from the locus of
canonical double covers of quadratic differentials in the stratum 
$\cQ_{0,2g+2}(m-1, 2g-3-m, -1^{2g})$.
\par
Suppose $\mu$ is of (stable) meromorphic type (and the hyperelliptic component is excluded if $\mu = (2g-2+2m, -2m)$). Then the ramification divisor of the map $\varphi$ in the interior of the stratum is empty unless $\mu = (m_1, m_2, 2g-2-m_1-m_2)$ consists of three zeros and poles of even order (i.e. $m_1$ and $m_2$ are both even). In this case the ramification divisor arises from the locus of
canonical double covers of quadratic differentials in the stratum 
$\cQ_{0,2g+2}(m_1-1, m_2-1, 2g-3-m_1-m_2, -1^{2g-1})$.
\end{prop}
\par
\begin{proof} To determine the ramification divisor we only need to consider
automorphisms stabilizing pointwise a divisorial locus in $\bP\LMS$.
We can moreover restrict to automorphism groups of prime order, since any non-trivial group has
such a subgroup.
\par
First consider $(X, \omega, z_1, \ldots, z_n) \in \bP\omoduli[g,n](\mu)$
in the stratum interior for $\mu = (m_1, \ldots, m_n)$.  Let $\tau$ be an
automorphism of~$X$ of prime order~$k$, so that $\tau$ induces a cyclic cover $\pi\colon X\to Y$ of degree $k$ with the quotient curve $Y$ of genus $h$. Then $\omega^k$
is $\tau$-invariant, hence there exists a $k$-differential~$\eta$ in $Y$
such that $\pi^{*}\eta = \omega^k$. The marked zeros and poles $z_1, \ldots, z_n$ of~$\omega$ are fixed by~$\tau$, hence they are totally ramified under $\pi$. Suppose $\pi$ has additional ramification points at $x_1, \ldots, x_\ell \in X$, each of which must also be totally ramified since $k$ is prime. 
Then the signature of $\eta$ is $(a_1, \ldots, a_n, 1-k, \ldots, 1-k)$ where $a_i = m_i + 1 - k$.
We have the Riemann--Hurwitz relation
\be\label{eq:RHrel}
2g-2 \= k (2h-2) + (n+\ell)(k-1).
\ee
\par
First suppose $\ell = 0$.  For $\mu$ of holomorphic type (and hence $g\geq 1$), the projective dimension of the stratum of
such $(Y, \xi)$ is (at most) $2h - 2 + n$ (where the maximum dimension is attained if all $m_i + 1 - k \geq 0$).
If $2h- 2 + n \geq 2g - 3 + n$, then  $h \geq g$.  But this is impossible for the 
branched cover $\pi$ with at least one totally ramified point. This argument also works for $\mu$ of meromorphic type, using the inequality $2h- 3 + n \geq 2g - 4 + n$ instead. 
\par
Next suppose $\ell > 0$.  Consider first the case that $\mu$ is of holomorphic type. 
Since $1 - k < 0$, the projective dimension of
the stratum of such $(Y, \eta)$ is $2h - 3 + n + \ell$.
Suppose $2h- 3 + n + \ell \geq 2g - 3 + n$, i.e., if such locus has at most
codimension one in $\bP\omoduli[g,n](\mu)$. Then $2h + \ell \geq 2g$ and it follows from \autoref{eq:RHrel} that   
$$(2k-2)h + (n-2)(k-1) + \ell(k-2) \leq 0\,. $$ 
\begin{itemize}
\item Suppose $k \geq 3$. Then $n = 1$, $h = 0$ and $\ell \leq 2$. Hence $2g \leq 2h + \ell \leq 2$ and $g\leq 1$. The only possibility is $g = 1$ and $\mu = (0)$, which gives the hyperelliptic component excluded in the assumption. 
\item Suppose $k = 2$.  Then $2h + n -2 \leq 0$, hence $h = 0$ and $n \leq 2$.
If $n = 1$, then it gives the hyperelliptic component of the
minimal stratum in genus~$g$ which is excluded in the assumption.  
If $n = 2$, then we obtain the locus of
hyperelliptic differentials in the stratum $\bP\omoduli[g,2](m, 2g-2-m)$ that arises via canonical double covers from quadratic differentials in the stratum 
$\cQ_{0,2g+2}(m-1, 2g-3-m, -1^{2g})$ of projectiveized dimension $2g-1$, i.e., of codimension
one in the stratum $\bP\omoduli[g,2](m, 2g-2-m)$. In this case the two zeros are ramified, hence they are
Weierstrass points. Consequently the zero orders $m$ and $2g-2-m$ have to be
even. 
\end{itemize}
\par 
For $\mu$ of meromorphic type, the above inequalities become $2h + \ell \geq 2g-1$ and $(2k-2)h + (n-2)(k-1) + \ell(k-2) \leq 1$. Moreover in this case $n \geq 2$, as a (stable) meromorphic differential has at least one zero and one pole. Then a similar analysis as above leads to the locus of hyperelliptic differentials in the meromorphic stratum $\bP\omoduli[g,3](m_1, m_2, 2g-2-m_1-m_2)$ with $m_1$ and $m_2$ both even. 
\par 
Finally we have to make sure that the branching order of $\varphi$ along the
locus of hyperelliptic differentials in the stratum
$\bP\omoduli[g,2](m, 2g-2-m)$ and in $\bP\omoduli[g,3](m_1, m_2, 2g-2-m_1-m_2)$ 
is just $k=2$, not a higher power of two. This follows from the fact that a general
hyperelliptic curve has no non-trivial automorphisms except the hyperelliptic involution. 
\end{proof}
\par
\begin{rem}{\rm 
Using the definition of theta characteristics we can determine the spin parity of the differentials in the ramification divisors in \autoref{prop:ramifInt}.  For a holomorphic differential of type $(2m, 2g-2-2m)$ with both zeros as Weierstrass points in the underlying hyperelliptic curve, the parity of the spin structure is given by $\lfloor m/2\rfloor + \lfloor(g+1-m)/2 \rfloor  \mod 2$.  For a meromorphic differential of type $(2m_1, 2m_2, -(2m_1+2m_2+2 - 2g))$ with both zeros and the pole as Weierstrass points in the underlying hyperelliptic curve, the parity is given by 
$\lfloor m_1/2\rfloor + \lfloor m_2/2\rfloor - \lfloor(m_1 + m_2 - g)/2 \rfloor  \mod 2$.  
}
\end{rem}
\par
To control the ramification at the boundary of $\bP\LMS$ we need a variant of the above 
proposition, allowing marked points to be permuted, but with automorphisms on the full stratum
rather than just on a divisor.
\par
\begin{lemma} \label{le:ramifbd1}
Let $\mu$ be an arbitrary signature of a (stable) stratum of differentials, possibly
of meromorphic type and with unlabeled singularities. If each $\omega$ in the stratum (component) 
$\bP\omoduli[g,\{n\}](\mu)$ admits a non-trivial automorphism~$\tau$ of order~$k$
fixing $\omega$ up to a $k$-th root of unity, then $k =2$, the automorphism~$\tau$
is the hyperelliptic involution, and the stratum (component) is hyperelliptic for 
$\mu = (2g-2)$, $\{ g-1, g-1\}$, $(2g-2+2m, -2m)$ with $m > 0$, $(2g-2 + 2m, \{-m, -m \})$ with $m > 0$ or $m < 1-g$, and $(\{m_1, m_1\}, \{-m_2, -m_2\})$ with $m_i > 0$ and $m_1 - m_2 = g-1$. 
\end{lemma}
\par
\begin{proof}
As before $\tau$ induces a cyclic cover $\pi\colon X\to Y$ of degree~$k$ with~$Y$
of genus~$h$, and there exists a $k$-differential~$\eta$ in~$Y$ of signature $\mu'$
such that $\pi^{*}\eta = \omega^k$. By the Riemann--Hurwitz relation we have 
$$ 2g-2 \= k(2h-2) + \sum_{i=1}^b (d_i-1)r_i + \sum_{j=1}^{b'} (d'_j - 1) r'_j $$
where the singularities of $\omega$ are distributed into $b$ orbits under $\tau$,
each having cardinality $r_i$ with $d_i = k /r_i$, and  in addition there are~$b'$
special (unmarked) orbits, each having cardinality $r'_j < k$ with $d'_j = k / r'_j$.
(Note that we do {\em not} require $r_i < k$.) With these notations we have
$n = \sum_{i=1}^b r_i$, hence the above relation can be rewritten as 
\be \label{eq:RHbb'}
2g - 2 + n  \= k (2h - 2 + b + b') - \sum_{j=1}^{b'} r'_j\,.
\ee
\par
By assumption, the dimension of $\bP\omoduli[g,\{n\}](\mu)$ agrees with the
dimension of the corresponding projectivized stratum of $k$-differentials
$\bP\omoduli[h,\{b+b'\}]^k(\mu')$, i.e. $2g - 2 + n = 2 h - 3 + b + b'$ (if $\mu$
is of holomorphic type) or $2g - 3 + n = 2 h - 3 + b + b'$ (if $\mu$ is of meromorphic type).
In particular, it implies that $2h -2 + b + b' \geq 2g - 2 + n$. Hence combining
with~\autoref{eq:RHbb'} it gives 
$$ (k-1)(2h - 2 + b + b') \,\leq\, \sum_{j=1}^{b'} r'_j\,. $$
Since $k\geq 2$ and $r'_j\leq k/2$, we deduce that $(k-1)b' \geq kb' / 2 \geq \sum_{j=1}^{b'} r'_j$. 
Then the preceding inequality implies that $(k-1)(2h - 2 + b) \leq 0$.
Since $b \geq 1$, it follows that $h = 0$ and $b = 1$ or $2$. 
\par 
Suppose $b = 2$. Then $(k-1)b' \leq \sum_{j=1}^{b'} r'_j \leq kb' / 2$. If $k > 2$,
then $k-1 > k/2$, and the only possibility is $b' = 0$.  In this case
$2g-2 = - r_1 - r_2$, hence $g = 0$ and $r_1 = r_2 = 1$, which leads to the unstable
signature $\mu = (-1, -1)$. If $k = 2$, then together with $h = 0$ and $b=2$ we
obtain those meromorphic hyperelliptic components as claimed. 
\par
Suppose $b = 1$. Then all $n = r_1$ singularities are in one fiber,
i.e.\ $\mu = (m,\ldots,m)$. Moreover $(k-1)(b' - 1) \leq \sum_{j=1}^{b'} r'_j
\leq kb' / 2$, i.e. $(k-2)(b'-2) \leq 2$. If $k > 2$, then $b' \leq 4$. It
follows that 
$$2g-2 +n \= k(b'-1) - \sum_{j=1}^{b'} r'_j \,\leq\, b' - 1 \leq 3\,.$$ 
If $g = 0$, then $n = r_1 > 2$ by stability, but $m n = 2g-2 = -2$, leading to a
divisibility contradiction. If $g = 1$, it is the case of elliptic curves with some ordinary markings, but the only non-trivial automorphism of a generic elliptic curve is the involution of order $k = 2$. The remaining case is $g = 2$ and $n = 1$, and consequently $k = 3$ and $b' = 4$, which occurs for $\mu = (2)$. But the unique zero $z$ is a Weierstrass point and $3z$ cannot be a fiber of a cyclic triple cover of $\bP^1$ since $z$ is a base point of the line bundle $\cO(3z)\cong K(z)$. Finally if $ k = 2$, since $h = 0$ and $b =1$, it leads to those holomorphic hyperelliptic components as claimed. 
\end{proof}
\par
We can now describe the ramification at the boundary. For holomorphic signatures $\mu$, the following boundary strata will be used for this purpose. We say that a two-level
graph~$\Gamma \in \LG_1(\mu)$ is a \emph{hyperelliptic bottom tree (HBT)},
if it is a tree and all the vertices on bottom level belong to
a stratum with signature $(2m_0, -2m_1,\ldots,-2m_{\ell})$ for some integers $m_i > 0$. In particular, every bottom vertex 
has exactly one labeled zero.  For an HBT graph~$\Gamma$ we denote by $D_\Gamma^{\rm H} \subset D_\Gamma$ the union
of irreducible components where the differentials on each vertex on bottom level admit a hyperelliptic involution which fixes the labeled zeros 
and poles (i.e.~the edges of the graph, and hence in particular the residues at the edges are all zero, thus satisfying the GRC as required in~\cite{LMS}). In particular, these bottom differentials in $D_\Gamma^{\rm H}$ are hyperelliptic, i.e. they are anti-invariant under the hyperelliptic involution.
\par
We say that a two-level graph~$\Gamma \in \LG_1(\mu)$ is a \emph{hyperelliptic
top backbone (HTB)}, if it is a tree with a unique bottom vertex (i.e.\ a
backbone graph) such that every top vertex is of type $ (2g_i-2)$
where $g_i$ is the genus of the vertex. In particular, all labeled zeros
are on bottom level. For an HTB graph~$\Gamma$ we denote by $D_\Gamma^{\rm H} \subset
D_\Gamma$ the union of irreducible components where the differentials on
each vertex on top level belong to the hyperelliptic component of the
stratum with signature $(2g_i-2)$. 
\par
Note that a graph~$\Gamma\in \LG_1(\mu)$ for holomorphic signature $\mu$ can be of type HBT and HTB at the same time only
if $\mu = (2g-2)$, and in that case $\Gamma$ is a backbone graph with a
unique bottom vertex carrying the unique labeled zero.
\par
We will also encounter graphs~$\Gamma \in \LG_1(2g-2)$ that we call
\emph{hyperelliptic banana backbones (HBB)}, where $\Gamma$ has a unique bottom vertex (carrying the unique zero), there exists an involution~$\tau$
 fixing the vertices of $\Gamma$, the signature
at each vertex admits a hyperelliptic component (taken the GRC into account for the bottom vertex) with $\tau$ as the hyperelliptic involution, and the quotient graph by~$\tau$ is a backbone graph. In particular,
the edges of~$\Gamma$ are either fixed or pairwise permuted by~$\tau$ (where
each permuted pair of edges looks like a banana in the drawing). We further require an HBB graph to contain at least one banana (otherwise it is of type 
HBT and HTB). 
\par 
For an HBB graph $\Gamma$ we denote by
$D_\Gamma^{\rm H} \subset D_\Gamma$ the union of irreducible components where
the differentials on each vertex belong to the hyperelliptic component
and where moreover the prong-matchings are chosen so that the hyperelliptic
involution does not extend to a neighborhood (see~\autoref{eq:zetacond2lev} and the surrounding paragraphs in 
\autoref{sec:Coordbd} below for more details). We remark that without considering prongs the hyperelliptic component and the spin component of a reducible stratum can actually intersect along the boundary (see~\cite[Corollary 7.10]{GendronThesis} and~\cite[Theorem 5.3]{chendiff} for an example). 
\par
	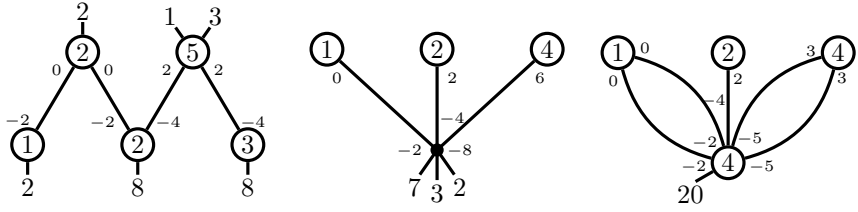
\begin{figure}[ht]
		\begin{tikzpicture}[
		baseline={([yshift=-2.4ex]current bounding box.center)},
		scale=2,very thick,
		bend angle=30]
		\node[circled number] (T) [] {$1$};
		\node[circled number] (T-1) [right=of T] {$2$};
		\node[circled number] (T-2) [right=of T-1] {$3$};
		\path (T) -- (T-1) node[circled number, midway, above=of T] (B) {$2$};
		\path (T-1) -- (T-2) node[circled number, midway, above=of T] (B-1) {$5$};
		\path (B) edge 
		node [order bottom left] {$0$} 
		node [order top left,xshift=0pt,yshift=1pt] {$-2$} (T); 
		\path (B) edge 
		node [order bottom right] {$0$} 
		node [order top left,xshift=0,yshift=0pt] {$-2$} (T-1);
		\draw (B) -- +(90: .2) node [xshift=0,yshift=4] {$2$};
		\draw (B-1) -- +(125: .2) node [xshift=-2,yshift=4] {$1$}; 
		\draw (B-1) -- +(55: .2) node [xshift=2,yshift=4] {$3$}; 				
		\draw (T) -- +(270: .2) node [xshift=0,yshift=-5] {$2$};
		\draw (T-1) -- +(270: .2) node [xshift=0,yshift=-5] {$8$};	
		\draw (T-2) -- +(270: .2) node [xshift=0,yshift=-5] {$8$};				
		\path (B-1) edge 
	    node [order bottom left] {$2$} 
	    node [order top right,xshift=-2pt,yshift=0pt] {$-4$} (T-1);	    
	    \path (B-1) edge 
	    node [order bottom right] {$2$}
	    node [order top right,xshift=-2pt,yshift=0pt] {$-4$} (T-2);
		\end{tikzpicture}
		\quad
		\begin{tikzpicture}[
		baseline={([yshift=-.5ex]current bounding box.center)},
		scale=2,very thick,
		bend angle=30]
		\node[circled number] (T) [] {$1$};
		\node[circled number] (T-1) [right=of T] {$2$};
		\node[circled number] (T-2) [right=of T-1] {$4$};
		\node[comp, below=of T-1] (B) {};
		\path (T-1) -- (B) node[] (B-1) {};
		\path (T-2) -- (B) node[] (B-2) {};
		\path (B.center) edge 
		node [order top left,xshift=1pt,yshift=-3pt] {$0$} 
		node [order bottom left,xshift=0pt,yshift=-2pt] {$-2$} (T);
		\path (B-1.center) edge 
		node [order top right,xshift=0pt,yshift=-1pt] {$2$} 
		node [order bottom right,xshift=-3pt,yshift=10pt] {$-4$} (T-1);
		\path (B-2.center) edge 
		node [order top right,xshift=0pt,yshift=-3pt] {$6$} 
		node [order bottom right,xshift=-2pt,yshift=-2pt] {$-8$} (T-2);
		\draw (B-1.center) -- +(270: .2) node [xshift=0,yshift=-5] {$3$}; 	
		\draw (B-1.center) -- +(235: .2) node [xshift=-2,yshift=-5] {$7$};
		\draw (B-1.center) -- +(305: .2) node [xshift=2,yshift=-5] {$2$};	
		\end{tikzpicture}
		\quad
		\begin{tikzpicture}[
		baseline={([yshift=0.5ex]current bounding box.center)},
		scale=2,very thick,
		bend angle=30]
		\node[circled number] (T) [] {$1$};
		\node[circled number] (T-1) [right=of T] {$2$};
		\node[circled number] (T-2) [right=of T-1] {$4$};
		\node[circled number, below=of T-1] (B) {$4$};
		\path (T) -- (T-1) node[ midway, below=of T,xshift=-10pt,yshift=-32pt] (BL) {};
		\path (T-1) -- (T-2) node[ midway, below=of T,xshift=10pt,yshift=-32pt] (BR) {};
		\path (T-1) -- (B) node[] (B-1) {};
		\path (T-2) -- (B) node[] (B-2) {};
		\path (T-1) edge 
		node [order top left,xshift=3pt,yshift=14pt] {$-4$} 
		node [order bottom right,xshift=-2pt,yshift=-2pt] {$2$} (B);
		\path (T) edge [bend left] 
		node [order bottom left, xshift=-5,yshift=-8] {$0$} 
		node [order top left, xshift=-1,yshift=-13] {$-2$} (B)
		edge [bend right] 
		node [order bottom right, xshift=2,yshift=11] {$0$} 
		node [order top right, xshift=-6,yshift=4] {$-2$} (B);
		\path (B) edge [bend left] 
		node [order bottom right, xshift=-3,yshift=0] {$-5$} 
		node [order top left, xshift=7,yshift=4] {$3$} (T-2)
		edge [bend right] 
		node [order bottom right, xshift=-5,yshift=-4] {$-5$} 
		node [order top right, xshift=-1,yshift=2] {$3$} (T-2);
		\draw (B) -- +(210: .25) node [xshift=-2,yshift=-5] {$20$};
		\end{tikzpicture}
		\caption{A hyperelliptic bottom tree graph (HBT), a hyperelliptic top backbone graph (HTB), and a hyperelliptic banana backbone graph (HBB).}
		\label{figure:ram}
	\end{figure}	

\par  
\begin{prop} \label{prop:ramifBound}
Suppose~$\mu$ is of holomorphic type (and the hyperelliptic component is excluded if $\mu = (2g-2)$).
At boundary divisors the map $\varphi\colon \bP\LMS \to \bP\MScoarse$ is ramified at most of order two. 
\par
More precisely, $\varphi$ is ramified at the components $D_\Gamma^{\rm H}$ for
$\Gamma$ of type HTB or HBT (except for $\mu = (2g-2)$ and $\Gamma$ of both HTB and HBT types), 
$\varphi$ is ramified at the components $D_\Gamma^{\rm H}$ for
$\Gamma$ of type HBB (where the hyperelliptic involution does not extend to a neighborhood), and $\varphi$ is not ramified at any other components of the boundary. 
\par
In particular, $\varphi$ is unramified at the horizontal boundary divisor $D_h$.
\end{prop}
\par
For meromorphic signatures $\mu$, the ramification situation can be similarly described at the boundary. We add a prime, e.g. HTB', to denote the corresponding types of graphs in the meromorphic case, with some extra allowances or requirements as follows. HTB' graphs allow additional meromorphic signatures of type $\mu_i = (-2m_i, 2g_i - 2 + 2m_i)$ and $(- 2m_i, \{g_i - 1 + m_i, g_i-1+m_i \})$ for top level vertices.  HBB' graphs occur for meromorphic signatures $\mu = (-2m, 2g-2+2m)$ with the marked zero on the unique bottom vertex and the marked pole on one of the top vertices. HBT' graphs require all marked poles to concentrate on one top vertex (and the other top vertices are of holomorphic type). 
\par 
\begin{prop} \label{prop:ramifBoundmero}
Suppose~$\mu$ is of meromorphic type (and the hyperelliptic component is excluded if $\mu = (-2m, 2g-2+2m)$).
At boundary divisors the map $\varphi\colon \bP\LMS \to \bP\MScoarse$ is ramified at most of order two. 
\par
More precisely, $\varphi$ is ramified at the components $D_\Gamma^{\rm H}$ for
$\Gamma$ of type HTB' or HBT' (except for $\mu = (-2m, 2g-2+2m)$ and $\Gamma$ of both HTB' and HBT' types), 
$\varphi$ is ramified at the components $D_\Gamma^{\rm H}$ for
$\Gamma$ of type HBB' (where the hyperelliptic involution does not extend to a neighborhood), and $\varphi$ is not ramified at any other components of the boundary. 
\par
In particular, $\varphi$ is unramified at any horizontal boundary divisor.
\end{prop}
\par
In the proof below we will use the following observation implicitly. Suppose $\omega$ is anti-invariant under the hyperelliptic involution $\tau$, i.e. $\tau^{*}\omega = -\omega$. Then the residue of $\omega$ is zero at any hyperelliptic Weierstrass point, and the sum of the residues of $\omega$ is zero at any pair of hyperelliptic conjugate points. In particular, if the GRC imposes the residue-zero condition to an edge of $\Gamma$ fixed by $\tau$ or to a pair of edges swapped by $\tau$ (i.e. a banana), then this residue condition is automatically satisfied by $\omega$. 
\par
\begin{proof}[Proof of \autoref{prop:ramifBound} and~\autoref{prop:ramifBoundmero}]
For holomorphic signatures $\mu$, consider first the horizontal divisor $D_h$
given by $\bP\LMS[\mu, -1, -1]$ identifying the simple poles~$p_1$ and~$p_2$ to a
node~$q$. If an automorphism~$\tau$ fixes the generic point in this boundary divisor
up to a rescaling of~$\omega$, then \autoref{le:ramifbd1} implies that $\mu=(2g-2)$,
a case which has been excluded.
\par
Second, for meromorphic signatures $\mu$, there can be another type of horizontal
divisors where the horizontal edge corresponds to a separating node.  But in that case
one component separated by the node must admit a non-trivial automorphism that fixes
the simple pole at the node, which is impossible by \autoref{le:ramifbd1}.
\par
Next we treat boundary divisors $D_\Gamma$ for $\Gamma \in \LG_1(\mu)$. If $D_\Gamma$
is in the ramification locus of $\varphi$, then every multi-scale differential
$(X, \bfomega = (\omega_{(0)}, \omega_{(-1)}), \bfsigma)$
in $D_\Gamma$ admits a non-trivial automorphism~$\tau$. Since every top level
vertex has either positive genus (with generically distinct moduli) or contains a labeled pole, $\tau$
cannot permute top level vertices. Similarly every bottom level vertex contains
a labeled zero, hence $\tau$ cannot permute them either.  Therefore, $\tau$
acts as an automorphism on each vertex of $\Gamma$. 
\par
Recall from \cite{LMS} that projectivized multi-scale differentials
are represented by $(X, \bfomega, \bfsigma)$ up to projectivization (rescaling
all levels simultaneously). Moreover, the action of the level rotation torus
rescales $\omega_{(-1)}$ and acts on the prong-matching~$\bfsigma$
simultaneously. Suppose that~$\tau$ has order~$k_{(i)}$ when restricted
to level~$i$. We conclude that $\tau^{*}\omega_{(i)} = \zeta_{(i)} \omega_{(i)}$
for $i=0,-1$, where $\zeta_{(i)}$ is a $k_{(i)}$-th root of unity
(not necessarily primitive). We remark that in the sequel a non-trivial action on a
vertex~$v$ means~$\tau$ restricted to the underlying marked surface~$X_v$
is non-trivial, and it does {\em not} necessarily imply that $\zeta \neq 1$
(but if $\zeta \neq 1$ then clearly $\tau$ must act non-trivially on $X_v$).
We also make a useful observation that if $\tau$ does not fix every edge
of $\Gamma$, then $\tau$ acts non-trivially on some vertices in both levels.   
\par 
Consider first the case that $\tau$ acts non-trivially on the top level
and we enumerate the vertices on that level by~$X_i$. 
By \autoref{le:ramifbd1}, the action of $\tau$ restricted to each top level vertex, if non-trivial, can only be the hyperelliptic involution which maps $\omega$ to $-\omega$ on that vertex. Since the top level differentials are projectivized simultaneously, if $\tau$ induces a hyperelliptic involution on one top level vertex, then it must act in the same way for every top vertex. Therefore, if $\mu$ is a holomorphic signature, 
each top level vertex carries a hyperelliptic differential of type
$\mu_i = (2g_i - 2)$ or $\{g_i -1, g_i-1 \}$, where $g_i$ is the genus of the vertex,
i.e.\ it has no labeled zero and admits either one edge fixed by $\tau$ or a
pair of edges (i.e. a banana) interchanged by $\tau$. Similarly if $\mu$ is a meromorphic signature, we allow in addition meromorphic signatures of type $(-2m_i, 2g_i - 2 + 2m_i)$ and $(-2m_i, \{ g-1 + m_i, g-1 + m_i \})$ for top level vertices.  In both cases $\Gamma$ is a
banana tree with a unique bottom vertex that contains all labeled zeros.  If $\tau$ acts trivially on the bottom level, then there is no banana, and we get HTB and HTB' (hyperelliptic top backbone) in the holomorphic and meromorphic cases, respectively.  If $\tau$ acts non-trivially on the unique bottom vertex (e.g. if there exists at least one banana), one can verify (by a similar but simpler argument as in the next paragraph) that the bottom (generalized) stratum must be hyperelliptic with a unique labeled zero, and all labeled poles belong to one top level vertex. This case corresponds to HBB and HBB' (hyperelliptic banana backbone) with $\mu = (2g-2)$ in the holomorphic case and $\mu = (-2m, 2g-2+2m)$ in the meromorphic case, respectively. 
\par 
Next suppose $\tau$ acts trivially on the top level and non-trivially on the
bottom level, and we enumerate now the vertices on bottom level
by~$X_i$. In particular, all edges are fixed by $\tau$. 
Suppose the bottom level has $v^{\bot}$ many vertices $X_i$, each of genus~$g_i$, 
with~$n_i$ labeled zeros and poles, and $e_i$ edges. We separate the discussion in two cases. Consider first the case $\tau^{*}\omega_i = \zeta_{(-1)} \omega_i$ for
some $\zeta_{(-1)} \neq 1$ for all differentials~$\omega_i$ on~$X_i$ in the
bottom level. This assumption implies that $\tau$ restricted to the marked
surface $X_i$ is a non-trivial action of order $k_i \geq 2$. Let~$Y_i$ be the
quotient of~$X_i$ by~$\tau$ and denote its genus by~$h_i$. Suppose that
$\pi_i \colon X_i \to Y_i$ has in addition $\ell_i$ branch points (not from the images of labeled zeros and poles and edges), and that over each such branch point
the fiber cardinality is $c_{i,j}$ with multiplicity $d_{i,j}$ for each fiber point, i.e. $c_{i,j} d_{i,j} = k_i$ with $c_{i,j} \leq k_i/2$.  Then we have the Riemann--Hurwitz relation 
$$ 2g_i - 2 + n_i + e_i \= k_i(2h_i - 2 + n_i + e_i  + \ell_i) - \sum_{j=1}^{\ell_i} c_{i,j}\,.  $$
For a holomorphic signature $\mu$, the bottom generalized stratum has (unprojectivized) dimension equal to 
$$N^{\bot} \= \Big(\sum_{i=1}^{v^{\bot}} (2g_i - 2 + n_i + e_i ) \Big) - (v^{\top} - 1)\,, $$
where $v^{\top}$ is the number of top vertices, and we subtract $v^{\top} - 1$
because the GRC imposes this many independent conditions (besides the Residue Theorem condition on each vertex). For a meromorphic signature $\mu$, 
the bottom generalized stratum has (unprojectivized) dimension bigger than or equal to the above formula, since a top vertex with marked poles does not impose a GRC, and the equality is attained if and only if all top level marked poles belong to the same vertex (so the GRC is imposed by the other $v^{\top}-1$ holomorphic top vertices independently).  The dimension of the (unprojectivized) locus of those $Y_i$ is 
$$N' \= \sum_{i=1}^{v^{\bot}} (2h_i - 2 + n_i + e_i + \ell_i)\,. $$
By assumption we have $N^{\bot} \leq N'$, which implies that 
$$ \sum_{i=1}^{v^{\bot}} (k_i-1)(2h_i - 2 + n_i + e_i  + \ell_i) \leq v^{\top} - 1 + \sum_{i=1}^{v^{\bot}}\sum_{j=1}^{\ell_i} c_{i,j}\,. $$ 
Since $k_i \geq 2$ and $c_{i,j}\leq k_i/2$, we have 
$$\sum_{i=1}^{v^{\bot}}\sum_{j=1}^{\ell_i} c_{i,j} \leq  \sum_{i=1}^{v^{\bot}} (k_i/2)\ell_i   \leq  \sum_{i=1}^{v^{\bot}}(k_i-1)\ell_i\,.$$ 
Moreover since $\Gamma$ is connected, we have $\sum_{i=1}^{v^{\bot}} e_i \geq v^{\top} -1 + v^{\bot}$, where the equality holds if and only if $\Gamma$ is a tree graph. It follows that $\sum_{i=1}^{v^{\bot}} (2h_i - 1 + n_i) \leq 0$. Since every $n_i > 0$, the only possibility for $2h_i  - 1 + n_i \leq 0$ is $(h_i, n_i) = (0, 1)$ for which $2h_i  - 1 + n_i = 0$, and hence all inequalities involved above must be equalities. Therefore, we conclude that all $k_i = 2$, $h_i = 0$, $n_i = 1$, $c_{i,j} =1$, and that $\Gamma$ is a tree graph with each bottom vertex carrying a hyperelliptic differential and having exactly one labeled zero, and any two adjacent top and bottom vertices are joined by a single edge fixed under the hyperelliptic involution of the bottom vertex. We thus conclude that this case corresponds to HBT and HBT' (hyperelliptic bottom tree) in the holomorphic and meromorphic cases, respectively. 
\par 
Now consider the other case when $\tau$ acts trivially on the top level, non-trivially on the bottom level, but $\tau^{*}\omega_i = \omega_i$ for all
differentials~$\omega_i$ on the lower level vertices. Using the above
notation for the quotient map, $\omega_i$ being $\tau$-invariant implies
that $\omega_i = \pi_i^{*} \eta_i$ for an Abelian differential $\eta_i$
in each $Y_i$. Moreover, since all edges and labeled zeros and poles are fixed by $\tau$,
the quotient differential $(Y_i, \eta_i)$ has the same number of zeros and
poles as $(X_i, \omega_i)$ (despite possibly a different genus
and different orders of zeros and poles). We also infer that the residue
of~$\omega_i$ at any polar edge is equal to~$k_i$ times the residue
of~$\eta_i$ at the image pole under $\pi_i$. For the purpose of dimension count 
we can thus replace each $(X_i, \omega_i)$ by $(Y_i, k_i \eta_i)$ in the bottom
level of~$\Gamma$, as residue constraints imposed by the GRC are linear and depend
only on the graph topology of~$\Gamma$, not its decoration by genera and
enhancements. It follows
that the GRC imposes the same number of conditions to the lower level of
the graph before and after the replacement. Hence by assumption we conclude that 
$$ \sum_{i=1}^{v^{\bot}} (2g_i - 2 + n_i + e_i)  \=  \sum_{i=1}^{v^{\bot}} (2h_i - 2 + n_i + e_i)\,. $$ 
By the Riemann--Hurwitz relation 
$$2g_i - 2 + n_i + e_i = k_i (2h_i - 2 + n_i + e_i) \geq 2h_i - 2 + n_i + e_i$$ 
and the inequality is strictly if $k_i > 1$. But there exists at least one $k_i > 1$ by the assumption that $\tau$ acts non-trivially on the bottom level, thus leading to a contradiction to the preceding identity. 
\par
The ramification orders at the components~$D_\Gamma^{\rm H}$ for those special graphs $\Gamma$ 
will be justified in \autoref{sec:Coordbd} below.
\end{proof}
\par

\begin{proof}[Proof of~\autoref{prop:ramifBound} and~\autoref{prop:ramifBoundmero}, ramification
orders] For simplicity we use the graph notations for holomorphic signatures $\mu$.  The argument works identically for the case of meromorphic signatures.  
We start with the case of $(X,\bfomega,\bfz,\bfsigma)$ in a component
$D_\Gamma^{\rm H}$ for~$\Gamma$ of type HBT or HTB, i.e., the involution~$\tau$ fixes all the
edges. In particular $\zeta_{e^\pm}^2 = 1$ for every edge $e$. Moreover, $p_e$ is odd for all~$e$,
hence~$\zeta_{e^\pm} = -1$ if and only if the action on the corresponding level
is non-trivial. This implies that the ramification order is two at $D_\Gamma^{\rm H}$ for these graphs, except for the simultaneous
involution of the intersection of HBT and HTB, where the action on the $t$ parameter is given by multiplication by $c = (-1)(-1) = 1$ (this is because the system of equations \autoref{eq:zetacond2lev} has a unique solution and $c=1$ is a valid one), and hence the map is not ramified at $D_\Gamma^{\rm H}$ when $\Gamma$ is of both HBT and HTB types. 
\par
Next consider an HBB graph~$\Gamma$. If an edge~$e$ is
fixed by the involution $\tau$, then~$p_e$ is odd and hence as above $\zeta_{e^+} \zeta_{e^-} = (-1)(-1) = 1$.
Suppose two edges~$e_1$ and $e_2$ are swapped. Being an involution implies
that $\zeta_{e_1^+}  \zeta_{e_2^+} = 1 = \zeta_{e_1^-}  \zeta_{e_2^-}$. Since
these two edges have the same~$p_e$ and hence the same $m_e$, the system~\autoref{eq:zetacond2lev} is
solvable only if $\zeta_{e_1^+} \zeta_{e_1^-}  = \zeta_{e_2^+} \zeta_{e_2^-}$. Since $\tau$ is an involution,
$c^2=1$, and hence $\zeta_{e_i^+} \zeta_{e_i^-} = \pm 1$. Moreover, any component of~$D_\Gamma$ is either unramified (if $c = 1$) 
or has ramification order two (if $c = -1$), depending on whether or not the hyperelliptic involution on the boundary component can extend to a neighborhood in the interior.  
\par
If~$p_{e_i}$ is even, then indeed both possibilities can occur. Suppose for simplicity that there are no other edges besides the swapped pair.
Given a solution of the above system with $\zeta_{e_i^+} \zeta_{e_i^-} = -1$,
we can replace the coordinate~$x_1$ by~$-x_1$. This still puts~$\omega_{(0)}$
in standard form there (as $d(x_1^{p_{e_1}})$ is unchanged), but swaps the sign of both $\zeta_{e_i^+}$ and
thus the new roots of unity satisfy that $\zeta_{e_i^+} \zeta_{e_i^-} = +1$. 
\par
On the other hand, if all the twists $p_{e}$ are odd, then $\ell_\Gamma$ is odd. 
Since $\zeta_{(0)} = \zeta_{(-1)} = -1$, raising~\autoref{eq:zetacond2lev} to the $p_e$-th power implies that $c^{\ell_\Gamma} = (-1)(-1) = 1$ and hence $c \neq -1$.  In this case the hyperelliptic involution always extends to the interior, i.e., the corresponding $D_\Gamma$ lies in the boundary of the hyperelliptic component.  
\end{proof}


\section{Projectivity of the coarse moduli space}
\label{sec:projectivity}

In this section we recall the background on the geometry of
the moduli stack of multi-scale differentials and its
coarse moduli space. The main goal is to prove the projectivity
announced in \autoref{thm:LMSisproj}.
\par
Denote by $\barmoduli[g,n](\mu)$ the closure of the (projectivized) stratum of Abelian differentials of type $\mu$ in the Deligne--Mumford compactification $\barmoduli[g,n]$, and by $\ol{M}_{g,n}(\mu)$ its coarse moduli space. This is the image of a projection of the
incidence variety compactification originally defined in \cite{BCGGM1}, where the projection contracts boundary strata whose level graphs have at least two vertices on top
level. There is a forgetful morphism of stacks $f\colon \bP\LMS \to \barmoduli[g,n](\mu)$, 
and we denote by $\ol{f}\colon \bP\MScoarse \to \ol{M}_{g,n}(\mu)$ the 
corresponding map of coarse moduli spaces. The space~$\ol{M}_{g,n}(\mu)$,
as a subvariety of $\ol{M}_{g,n}$, is projective and thus has an ample line bundle $\cA$. Recall that a line bundle $\cB$ is called
{\em $\ol{f}$-ample} or {\em relatively ample},  if $\cB$ is ample on
every fiber of~$\ol{f}$.  If $\cA$ is ample and $\cB$ is $\ol{f}$-ample, then
$\ol{f}^*\cA \otimes \ep \cB$ is ample on $\bP\MScoarse$ for small enough~$\ep$ (see e.g.\ \cite[Section~1.7]{LazBook1} for these facts).
It thus suffices to show the existence of such an $\ol{f}$-ample bundle.
\par
Our strategy relies on three observations: First, the fibers of~$\ol{f}$ are finitely
covered by toric varieties. Intuitively, the toric structure stems from
rescaling the differentials on subsets of the components of the stable
curve (such that the global residue condition is preserved).
There we can use the toric Nakai-Kleimann
criterion for ampleness.
\par
\begin{prop}[{\cite[Theorem~6.3.13]{CLS}}] \label{prop:CLScrit}
A Cartier divisor~$D$ on a
proper toric variety~$X$ is ample if and only if $D \cdot C >0$ for every
torus-invariant irreducible curve~$C$ on~$X$.
\end{prop}
\par
Second, those torus-invariant curves map to a class of curves in~$\bP\LMS$
that we call relevant curves and that are easy to describe for a given
level graph. Third, we can verify on~$\bP\LMS$ the required positivity
by showing the following, using the notation that will be introduced
in~\autoref{eq:calL} below. For notation simplicity, we use $\ol{B}$ to denote the moduli space of multi-scale differentials as we did in \autoref{sec:DMtaut}. 
\par
\begin{prop}\label{prop:fampleclass}
There is an effective divisor class $D$ such that $\cL_{\ol{B}}\otimes
\cO_{\ol{B}}(-D)$ has positive intersection numbers with all relevant curves
in all fibers of~$\ol{f}$.
\end{prop}
\par

\subsection{Some line bundles}

Let $\Gamma \in \LG_1(\ol{B})$ be a graph corresponding to a divisor~$D_\Gamma$
in $\ol{B}$. The least common multiple
\bes
\ell_\Gamma \= \lcm\Bigl(p_e \,:\, e \in E(\Gamma)\Bigr)
\ees
appears frequently in formulas, which is the size of the orbit of the
twist group acting on all prong-matchings for~$\Gamma$. One prominent
line bundle is defined by the following sum of boundary divisors
\be \label{eq:calL}
\cL_{\ol{B}} \= \cO_{\ol{B}}\Bigl(\sum_{\Gamma\in \LG_1(\ol{B})} \ell_\Gamma
D_\Gamma \Bigr)\,.
\ee
The compactification of the stratum~$\ol{B}$, being constructed as the $\bC^*$-quotient of the
unprojectivized space $\LMS$ comes with a tautological bundle~$\cO_{\ol{B}}(-1)$,
whose first Chern class is denoted by~$\xi$.
\par
For inductive arguments we need the following generalization. For $\Gamma\in \LG_L(\ol{B})$, 
we denote by $\fraki_\Gamma\colon D_\Gamma \to \ol{B}$ the inclusion that maps the
boundary strata into the total space and by $\frakj_{\Delta, \Gamma}\colon 
D_\Delta \to D_\Gamma$ the inclusion into an undegeneration. As in
\cite{EC} we denote by $\ell_{\Gamma,i}$ the lcm of the
enhancements~$p_e$ of the edges~$e$ of the two-level undegeneration
$\delta_{i}(\Gamma)$ and let $\ell_\Gamma = \prod_{i=1}^{L} \ell_{\Gamma,i}$.
We now define 
\be \label{eq:defLithlev}
\cL_{\Gamma}^{[i]} \=  \cO_{D_\Gamma}
\Bigl(\sum_{\Gamma \overset{[i]}{\rightsquigarrow}
  \wh{\Delta} } \ell_{\wh{\Delta},-i+1}D_{\wh{\Delta}} \Bigr)
\quad \text{for any}  \quad i\in \{0,-1,\dots,-L\}\,,
\ee
where the sum is over all graphs $\wh{\Delta} \in \LG_{L+1}(\ol{B})$
that yield divisors in~$D_\Gamma$ by splitting the $i$-th level. 
\par
Generalizing the definition of~$\xi$ to the level strata of $D_\Gamma$, we define~$\xi_\Gamma^{[i]} \in \CH^1(D_\Gamma)$ to be the first
Chern class of the tautological bundle at level~$i$ on $D_\Gamma$. 
\par

\subsection{Toric covers}

We start with a description of the fibers of~$\ol{f}$ and recall some
more details about the construction in \cite{LMS}.
Let $(X,\bfz,\bfeta)$ be a twisted differential consisting of a pointed
stable curve~$(X,\bfz)$ and a collection $\bfeta =
\{\eta_v\}_{v \in V(\Gamma)}$ of differentials indexed by the vertices of
the dual graph~$\Gamma$ of~$X$, satisfying the conditions of a twisted
differential compatible with some level structure on~$\Gamma$, as given
in~\cite{BCGGM1}. The level structure is not unique, but we can assume that
the level structure has a minimal number of levels. The fiber~$\ol{F}$
of~$\ol{f}$ over (the image in $\ol{M}_{g,n}(\mu)$ of) this twisted differential consists of all multi-scale
differentials~$(X,\bfz,\bfomega,\bfsigma, \Delta)$ where the
collection of differentials~$\bfomega$ is compatible with some enhanced
level graph structure~$\Delta$ on the given
dual graph~$\Gamma$, where~$\bfsigma$ is some prong-matching
and where each of the components~$\omega_v$ of~$\bfomega$ is a multiple
of~$\eta_v$. Recall moreover that two such multi-scale differentials
are equivalent if they differ by the action of the level rotation torus
$T_{\Delta}$ rescaling the differential level-wise and simultaneously
rotating the prong-matchings, see~\cite{LMS}.
\par
The set of all enhanced level graphs~$\Delta$ compatible
with~$\bfeta$ in this way, with arrows given by undegeneration, forms 
a directed graph. The terminal elements in this graph, i.e.\ those
with the minimal number of levels, correspond to the
irreducible components of~$\ol{F}$. By slight abuse of notation, we will  
denote by $\Gamma$ such a terminal element.
\par
Consider the action of the 'big' torus $T^{V(\Gamma)}$ rescaling the differentials
on each vertex individually. Since the multi-scale differentials are
constrained by the global residue condition, which are always of the
form that a sum of residues is zero, the orbit of a subtorus $T_\Gamma^P
\subset T^{V(\Gamma)}$ preserves the differentials~$\bfomega$ in the
fiber~$\ol{F}$. The torus $T_\Gamma^P$ contains the subtorus $T_\Gamma^{\mathrm{np}}
\cong (\bC^*)^{L(\Gamma)}$ that rescales the differentials level by level
(this torus is isogeneous to the level rotation torus $T_\Gamma$, but
so far we have no prong-matchings taken into account, thus explaining the upper index).
\par
Let $F \subset \bP\LMS$ be the fiber of~$f$ over $(X,\bfz,\bfeta)$.
There is a natural map $F \to \ol{F}$, given by passing from the quotient
stack (by the automorphism group of the pointed curve and the local factor group
$K_\Gamma$) to the coarse quotient. Our first goal is to show the following result.  
\par
\begin{prop} \label{prop:Fibaretoric}
For each irreducible component $F_\Gamma$ of~$F$ there exists a proper toric
variety~$\widetilde{F}_\Gamma$ for a torus isogeneous to~$T_\Gamma^P/
T_\Gamma^{\mathrm{np}}$ that admits a cover $\widetilde{F}_\Gamma \to F_\Gamma$
which is unramified over the open torus orbit.
\end{prop}
\par
The irreducible component $F_\Gamma$ of the fiber might be called a toric stack. However there are
various definitions of that notion and we prefer not entering that
discussion. Note that $\ol{F}_\Gamma$ might not be a toric variety due to
graph automorphisms, e.g.\ the quotient of~$\bP^{46}$ by $\bZ/47$ is
not even a rational variety (\cite{Swan}) and such examples can occur
in fibers of~$f$ for sufficiently large genus.
\par
In order to prove \autoref{prop:Fibaretoric} we prove local
versions and piece them together as in \cite{MuEnum}. The idea was also
used in \cite[Section~4.2]{EC} and we use the same notation as there,
except that all objects are restricted to a fiber of~$f$ and that
we need to additionally exhibit a torus action. Let $\Delta$ be
a degeneration of $\Gamma$ in the fiber $F$ and let $U(\Delta) \subset F_\Gamma$
be the open substack of multi-scale differentials
compatible with undegenerations of~$\Delta$.
\par
\begin{lemma} \label{le:Fiblocallytoric}
There exists a toric variety~${U}^s_\Delta$ for the torus~$T_\Gamma^P/T_\Gamma^{\text{np}}$
that admits an unramified cover ${U}^s_\Delta \to U_\Delta$ of stacks.
\end{lemma}
\par
\begin{proof} We use the cover constructed in \cite[Section~14]{LMS}
that provides the smooth DM-stack structure of $\bP\LMS$ by pieces
of the simple Dehn space. This cover first provides sufficiently small open
neighborhoods in $U_\Delta$ with a $\Tw[\Delta]$-marking, i.e. a marking up to
the monodromy in the group $\Tw[\Delta]$. The smooth charts of the stack
are then given by coverings that have a $\sTw[\Delta]$-marking (rather then
just a $\Tw[\Delta]$-marking) and are thus pieces of simple Dehn space.
These pieces glue together to the unramified cover ${U}^s_\Delta \to U_\Delta$.
This cover is indeed a smooth complex variety since the simple Dehn space is.
\par
Since the topology in the whole fiber of~$f$ is constant we can unwind
the definition of the simple Dehn space in \cite[Section~12]{LMS} intersected
with the fiber of~$f$ and obtain the following explicit description.
Let $\sTw[\Delta] = \oplus_{i \in L(\Delta)} \sTw[i]$ be the level-wise
decomposition of the simple twist group and let $T_i = \bC/\sTw[i]$ be the
level-wise constituents of the level rotation tori. Then, as a set
\bes
U_\Delta^s \= \coprod_{\Gamma {\rightsquigarrow} \Pi {\rightsquigarrow}
\Delta} \Biggl(\mathfrak{W}_{\operatorname{pm}}(\Pi) / \bigl(
\bigoplus_{i \in L(\Pi)} T_i \oplus \bigoplus_{i \in L(\Delta) \setminus L(\Pi)} \sTw[i]
\bigr)\Biggr)\,,
\ees
where, with the same deviation of notation compared to the cited sources,
$\mathfrak{W}_{\operatorname{pm}}(\Pi)$ denotes the set of all
prong-matched differentials compatible with~$\Pi$ in the given fiber of~$f$
and with $\sTw[\Pi]$-marking. Since all points in $\mathfrak{W}_{\operatorname{pm}}(\Pi)$
are obtained by rescaling~$\bfeta$ and choosing a prong-matching, we find 
\be \label{def:UDeltas}
U_\Delta^s \= \coprod_{\Gamma {\rightsquigarrow} \Pi {\rightsquigarrow}
\Delta} \Bigl( \bigoplus_{i \in L(\Delta) \setminus L(\Pi)} T_i \Bigr)
\cdot \Bigl(\coprod_\bfsigma (X,\bfz,\bfeta,\bfsigma,\Pi) \Bigr)\,,
\ee
where $\bfsigma$ runs through representatives of the prong-matchings
under the action of $\oplus_{i \in L(\Pi)} T_i$. In particular the torus
\bes
T_{\Delta,\Gamma}^{s} \= \bigoplus_{i \in L(\Delta) \setminus L(\Gamma)} T_i  
\ees
acts on $U_\Delta^s$ by acting on the stratum for~$\Pi$ via the components in
$L(\Delta) \setminus L(\Pi)$. The continuity of this action is clear
from the definition of the topology on $U_\Delta^s$. The torus
$T_{\Delta,\Gamma}^{s}$ is a cover of a factor of $T^P_\Gamma/T_\Gamma^{\text{np}}$
(in fact isogeneous if $\Delta$ is maximally degenerate), since the levels
of~$\Delta$ are obtained by pulling apart the levels of~$\Gamma$ according to
the rescalable pieces. The transitivity of the action of 
$ T_{\Delta,\Gamma}^{s}$ on each irreducible component of $U^s_\Delta$ is
obvious.
\end{proof}
\par
\begin{proof}[Proof of \autoref{prop:Fibaretoric}]
Since the claim is for each irreducible component~$F_\Gamma$
of~$F$, we can thus focus on the degenerations of a fixed~$\Gamma$ with
minimal number of levels. To ease notation, we keep calling this component~$F$.
We take $\widetilde{F}$ to be the normalization of~$F$ in the function
field of the smallest extension of $K(F)$ that contains field extensions
corresponding to $U_\Delta^s \to F$ in \autoref{le:Fiblocallytoric}. Note that there is a finite number of
extensions and that these extensions are unramified. Consequently 
$\widetilde{F} \to F$ is unramified, too. Since each of the $U_\Delta^s$
admits an action of a torus $T_{\Delta,\Gamma}^{s}$ isogeneous to
$T^P_\Gamma/T_\Gamma^{\text{np}}$, so does $\widetilde{F}$. In fact the fiber product
of the $T_{\Delta,\Gamma}^{s}$  over $T^P_\Gamma/T_\Gamma^{\text{np}}$ acts, by
the minimality of the field extension. Consequently $\widetilde{F}$ is toric and
its properness follows from the properness of~$F$.
\end{proof}
\par
\begin{example} \label{ex:cherry}
  {\rm Consider~$\Gamma$ the 'cherry' graph giving the boundary
divisor 
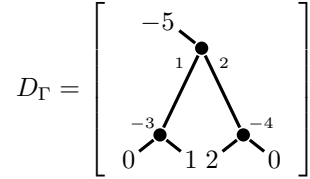
\bes
D_\Gamma=\left[
\begin{tikzpicture}[
baseline={([yshift=-.5ex]current bounding box.center)},
scale=2,very thick,
bend angle=30]
\node[comp,fill] (T) [] {};
\node [order node dis,above left] (T-1) at (T.north east) {$-5$};
\path[] (T) edge [shorten >=5pt] (T-1.center);
\node[comp,fill] (B) [below left= 1cm and 0.4cm of T] {}
edge 
node [order bottom left] {$-3$}
node [order top left] {$1$} (T); 
\node[comp,fill] (C) [below right= 1cm and 0.4cm of T] {}
edge 
node [order bottom right] {$-4$}
node [order top right] {$2$} (T);
\node [minimum width=18pt,below right] (B-2) at (B.south east) {$1$};
\path (B) edge [shorten >=5pt] (B-2.center);
\node [minimum width=18pt,below left] (B-2) at (B.south west) {$0$};
\path (B) edge [shorten >=5pt] (B-2.center);
\node [minimum width=18pt,below right] (C-2) at (C.south east) {$0$};
\path (C) edge [shorten >=5pt] (C-2.center);
\node [minimum width=18pt,below left] (C-2) at (C.south west) {$2$};
\path (C) edge [shorten >=5pt] (C-2.center);
\end{tikzpicture}\right]
\ees
in the stratum with $\mu = (2,1,0,0,-5)$. Since each of the vertices parameterizes a rational curve with three marked points which has unique moduli, the cherry represents a single point in $\ol{M}_{0,5}(\mu)$ (see also \cite[Example~14.5]
{LMS}). We will describe the fiber~$F$ and the toric variety $\widetilde{F}$
in this case.
\par
The residues at all the poles are zero by the Residue Theorem, 
so that $T_\Gamma^P \cong (\bC^*)^2$ rescales independently the vertices on lower
level and $T_\Gamma^{\textrm{np}} \cong \bC^*$ sits diagonally in~$T_\Gamma^P$.
Let~$\Delta_\ell$ (resp.~$\Delta_r$) be the slanted cherry graph with
the left (resp.\ right) edge being shorter. We focus on the case $\Delta := \Delta_\ell$.
Then, as subgroups of the group $\bZ \oplus
\bZ$ generating the Dehn twists around the left and the right nodes 
\bes \sTw[1] \= \langle (6,0) \rangle, \quad \sTw[2] \= \langle (0,3) \rangle,
\quad \sTw[\Delta] \= \sTw[1]  \oplus \sTw[2] \ees
and
\bes \Tw[\Delta] \= \langle \sTw[\Delta], (2,1) \rangle
\quad \text{hence} \quad K_{\Delta} \cong \bZ/3\bZ\,.
\ees
In this case $\mathfrak{W}_{\operatorname{pm}}(\Gamma) = (\bC \times \bC)/
\langle(6,6)\rangle$. The torus $T_1 = \bC/6\bZ$ acts diagonally on this space
and the discrete group~$\sTw[2]$ acts effectively on the second factors.
In this decomposition
\bes U_{\Delta}^s \= U_{\Delta}^s(\Gamma) \, \amalg \, U_{\Delta}^s(\Delta)\,,
\ees
and we have thus identified the first subset
$U_{\Delta}^s(\Gamma) = \mathfrak{W}_{\operatorname{pm}}(\Gamma)/(T_1 \oplus \sTw[2])$.
On the other hand,
\bes \mathfrak{W}_{\operatorname{pm}}(\Delta) \=
(\bC \times \bC)/\langle (6,6), (0,3) \rangle \quad \text{and} \quad 
U_{\Delta}^s(\Delta) = \mathfrak{W}_{\operatorname{pm}}(\Delta)/T_1 \oplus T_2\,,
\ees
where $T_1 = \bC/6\bZ$ acts diagonally as in the preceding case and
$T_2 = \bC/3\bZ$ acts on the second factor. Obviously $T_1 \oplus T_2$
acts faithfully and transitively and $U_{\Delta}^s(\Delta)$ is a single point.
The group $K_\Delta$ acts faithfully on the first subset while fixing the
second and thus produces the non-trivial quotient stack structure of
$U_\Delta$ at the image point of the left slanted cherry.  
\par
The more useful description of this decomposition of~$U_{\Delta}^s$
is~\autoref{def:UDeltas}. Since $T_1$ acts transitively on the prong-matchings
for both subsets, the decomposition into prong-matchings representatives is reduced to a single factor,
a single prong-matching equivalence class. The set~$U_{\Delta}^s$ is
thus a toric variety for $T_2 = (0 \times \bC)/(0 \times 3\bZ)$, which
is a triple cover of~$T_\Gamma^P/T_\Gamma^{\text{np}}$, where
$T_\Gamma^P = (\bC \times \bC)/(\bZ \times \bZ)$ and $T_\Gamma^{\text{np}}$ is
the diagonal in~$T_\Gamma^P$.
\par
The same description holds for $\Delta_r$ exchanging the role of the prong numbers $p_1 = 2$ and $p_2=3$ everywhere.
\par
Consequently, the full fiber~$F$ consists a complex plane $U_{\Delta_\ell}$
with orbifold order three at the origin, glued via $z \mapsto 1/z$ to a
complex plane $U_{\Delta_r}$ with orbifold order two at the origin.
The cover~$\widetilde{F} \to F$ is a cyclic cover of degree six, fully
ramified over the origin and $\infty$, which is the smallest cover that dominates
the cyclic cover of order three ramified at~$0$ and that of order
two ramified at~$\infty$. The fiber product of~$T_2$ and the
corresponding torus for the right slanted cherry
$T_1 = (\bC \times 0)/(2\bZ \times 0)$ over $T_\Gamma^P/T_\Gamma^{\text{np}}$
admits an isogeny of degree six to $T_\Gamma^P/T_\Gamma^{\text{np}}$ and acts
on~$\widetilde{F}$ as requested.
}\end{example}
\par

\subsection{Relevant curves}

We introduce the notion of {\em relevant curves} in~$\bP\LMS$ constructed as follows.
In the first step take a boundary stratum~$D_\Delta$, say with $\Delta \in
\LG_L(\ol{B})$ that has the following two features.
\par
First, one level~$i$ whose vertices~$V^{[i]}$ can
be partitioned into two sets~$V^{[i]}_A$ and~$V^{[i]}_B$ with the following property.
There exists a codimension-one degeneration $\Delta^A \in \LG_{L+1}(\ol{B})$
of $\Delta$ such that
the $i$-th level is split and the vertices in~$A$ go down while those in~$B$ stay up,
and vice versa a codimension-one degeneration $\Delta^B \in \LG_{L+1}(\ol{B})$
of $\Delta$
where those in~$B$ go down and those in $A$ stay up. (It may happen that the
two degenerations produce abstractly isomorphic graphs, see the rhombus graph
in \autoref{ex:rhombus} below. Note that just being able to put~$A$ down while
keeping $B$ up (without the converse) is not a sufficient criterion, see e.g.
the zig-zag graph in~\cite[Figure~2]{EC}.)
\par
Second, there is a unique level~$i$ with such a splitting and  there is no finer
partition of the vertices at level~$i$  that can be moved up and down independently (this also justifies that the vertices involved in the splitting into~$A$ and~$B$ are 
suppressed in the notation of relevant curves). 
\par
In the second step we define the relevant curve $C_{\Delta,i}$ inside $D_\Delta$ given
by specifying a point class on the generalized stratum of each level different
from~$i$, and a point class on the generalized strata corresponding to the
partitions~$A$ and~$B$ (these are the generalized strata at level~$i$
of $\Delta^A$ and $\Delta^B$ respectively). Note that a boundary
stratum can be disconnected due to prong-matching equivalence classes (see e.g.
\cite[Section~3]{EC}), but this second step pins down a component the
relevant curve lies in. We do no record this in the notation, since
the intersection numbers below do not depend on the component of the boundary
stratum.
\par
Obviously relevant curves are contained in the fibers of~$f$.
\par
\begin{lemma} \label{le:torusinvrelevant}
Let $C$ be a torus-invariant irreducible curve in the cover~$\widetilde{F}$ of
any fiber~$F$ of~$f$. Then the image of~$C$ in $\bP\LMS$ is a relevant curve.
\end{lemma}
\par
The converse statement also holds by the same argument, i.e. relevant curves are torus-invariant, but it is not
needed in the sequel. 
\par
\begin{proof} Recall that we covered~${F}$ by the images of the
sets $U^s_\Delta$. For a torus-invariant curve $C$, the generic point of~$C$ thus lies in the preimage
of some open set $U^s_\Delta$, and hence $C$ is given by the closure of certain one-dimensional
subtorus~$T_1$ of~$T:=T_\Gamma^P/T_\Gamma^{\text{np}}$. 
\par
First suppose that 
$T_1$ acts on at least two levels~$i$ and~$j$ non-trivially. Decompose the vertices on level $i$ into three subsets 
$A_i$, $B_i$ and $C_i$, where $T_1$ scales the vertices in $A_i$ down at its parameter $t = 0$, scales the vertices in $B_i$ down 
at $t = \infty$, and does not scale the vertices in $C_i$. Note that $A_i$ and $B_i$ are non-empty and $C_i$ can possibly be empty.  
In the same way we decompose the vertices on level $j$ into three subsets $A_j$, $B_j$ and $C_j$. 
Then the subtorus of~$T$ that rescales $A_i$ with $t$ and $A_j$ with $t^{-1}$ (and does nothing to the other subsets) 
exhibits~$C$ as not $T$-invariant, which contradicts the assumption. 
\par
Next suppose that the action
of $T_1$ is non-trivial on level~$i$ only and partitions the vertices of
that level into three non-empty subsets, the set~$A$ that goes down
at $0 \in \overline{T}_1$, the set~$B$ that goes down at $\infty \in \overline{T}_1$,
and the rest~$S$. Then the subtorus of~$T$ that fixes~$A \cup B$ and
rescales the vertices in~$S$ diagonally exhibits~$C$ as not
$T$-invariant, leading to the same contradiction.  
\par 
Therefore, $T_1$ acts non-trivially only on a single level~$i$ and decomposes the vertices on level~$i$ into two subsets~$A$
and~$B$ which are obtained by considering the limits to~$0$ and to~$\infty \in \ol{T}_1$. Note that this argument justifies the uniqueness and minimality in the second part of the first step defining relevant curves, while the second step in the definition simply cuts down the dimension to one (i.e. to a curve).
\end{proof}
\par
Using \autoref{prop:fampleclass} we can now complete:
\par
\begin{proof}[Proof of \autoref{thm:LMSisproj}]
The bundle $\cB = \cL_{\ol{B}}\otimes \cO_{\ol{B}}(-D)$ from
\autoref{prop:fampleclass} descends to a bundle $\overline{\cB}$
on~$\bP\MScoarse$ since the boundary divisors (and thus $\cL_{\ol{B}}$) are
invariant under the local isomorphism groups. We claim that $\overline{\cB}$
is $\ol{f}$-ample. By definition we need to show that the restriction of $\overline{\cB}$ to
any fiber of~$\ol{f}$ is ample. For this it suffices to prove the ampleness
of the pullback via the finite covering $\widetilde{F} \to F \to \ol{F}$
given by \autoref{prop:Fibaretoric}, which implies that we
can check ampleness via the toric criterion in \autoref{prop:CLScrit}.
Namely,  we need to check the positivity of the pullback of $\overline{\cB}$
to~$\widetilde{F}$ on any torus-invariant curve. Then by
\autoref{le:torusinvrelevant} and push-pull it suffices 
to check the positivity of~$\cB$ on any relevant curve.  
We have thus reduced the claim to that of
\autoref{prop:fampleclass}.
\end{proof}
\par

\subsection{The proof of positivity}

It remains to show \autoref{prop:fampleclass}.
Recall from the description of the relevant curves $C_{\Delta,i}$
above that among the boundary divisors that do not contain $C_{\Delta,i}$,
there are precisely two divisors $D_{\Gamma^A}$ and $D_{\Gamma^B}$ with
non-zero (hence positive) intersection numbers with $C_{\Delta,i}$, namely
the ones containing $D_{\Delta^A}$ and $D_{\Delta^B}$, respectively.  
\par
\begin{lemma}\label{lem:contribatlevel}
Let $C_{\Delta,i}$ be a relevant curve. Then
\bes
\big(\xi_{\Delta}^{[i]}+ \c_1(\cL_{\Delta}^{[i]})\big)\cdot [C_{\Delta,i}]>0\quad \text{and}
\quad (\xi_{\Delta}^{[i]})\cdot [C_{\Delta,i}]<0\,.
\ees
\end{lemma}
\par
\begin{proof}
For the first statement we express $\xi_{\Delta}^{[i]}$ using
\cite[Proposition~8.2]{EC}
as a positive $\psi$-class contribution and a negative boundary contribution,
consisting of certain summands that also appear in $\cL_{\Delta}^{[i]}$. Since
$\psi$-classes are pullbacks from $\barmoduli[g,n]$, they have zero intersection numbers with contracted
curves. The boundary divisors in $\cL_{\Delta}^{[i]}$ do not contain $C_{\Delta,i}$,
so their intersection numbers with  $C_{\Delta,i}$ are non-negative. More precisely, the boundary
terms in \cite[Proposition~8.2]{EC} are (for a chosen leg~$x$ specifying
the $\psi$-class) those where~$x$ goes down in the splitting. This set
contains exactly one of $D_{\Delta^A}$ and $D_{\Delta^B}$ since
in the definition of relevant curves the vertices in the split level are partitioned into~$A$ and~$B$. So the contribution
is positive, as claimed.
\par
The second statement follows from the same relation, since the boundary terms
in this relation do not contain~$C_{\Delta_i}$, and since one of 
$D_{\Delta^A}$ and $D_{\Delta^B}$ appears in the relation (with negative sign).
\end{proof}
\par
\begin{lemma} \label{le:cLwithrelevant}
Let $C_{\Delta,i}$ be a relevant curve with $\Delta \in \LG_L(\ol{B})$. 
Then
\bes \c_1(\cL_{\ol{B}})\cdot [C_{\Delta,i}]>0 \quad \text{if} \quad i \in \{0,-L\}
\ees
and $\c_1(\cL_{\ol{B}})\cdot [C_{\Delta,i}]=0$ otherwise.
\end{lemma}
\par
\begin{proof}
We denote the successive undegenerations of $\Delta$ keeping the
top~$j$ levels by $\Delta_j:=\delta_{1,\dots,j}(\Delta)$ and $\Delta_0:=\ol{B}$.
In particular $\Delta_L=\Delta$. Define moreover the successive pullbacks
of $\cL_{\ol{B}}$ to be $\cE_j:=\frakj_{\Delta_j,\Delta_{j-1}}^*(\cE_{j-1})$,
where $\cE_0:=\cL_{\ol{B}}$. In particular $\c_1(\cE_L)=\c_1(\cL_{\ol{B}})\cdot [D_{\Delta}]$. Applying \cite[Lemma 7.6]{EC} successively we find that
$\c_1(\cE_L) = \c_1(\cL_{\Delta}^{[-L]}) +\xi_{\Delta}^{[-L]}-\xi_{\Delta}^{[0]}$.
The desired claim thus follows from \autoref{lem:contribatlevel}. 
\end{proof}
\par
Finally recall from \cite[Theorem~7.1]{EC} that the normal bundle
of a boundary divisor is given by
\be \label{eq:nb}
\c_1(\cN_\Gamma) \= \frac{1}{\ell_\Gamma} \bigl(-\xi_\Gamma^\top
- \c_1(\cL_\Gamma^\top) + \xi_\Gamma^\bot \bigr)\quad \text{in} \quad
\CH^1(D_{\Gamma})\,.
\ee
More generally, the normal bundle of a codimension-one degeneration
of graphs, say with $\delta_{-i+1}^\complement(\Gamma)=\Pi$ is given
in loc.\ cit.\ by
\be \label{eq:nbinPi}
\c_1(\cN_{\Gamma,\Pi}) \= \frac{1}{\ell_{\Gamma, -i+1}} \bigl(-\xi_\Gamma^{[i]}
- \c_1(\cL_\Gamma^{[i]}) + \xi_\Gamma^{[i-1]} \bigr)\quad \text{in} \quad
\CH^1(D_{\Gamma})\,.
\ee
\begin{proof}[Proof of \autoref{prop:fampleclass}]
We start by numbering the (non-horizontal) boundary divisors~$D_1, D_2, \ldots$
of $\ol{B}$ in such a way that whenever two boundary divisors intersect, the one
with smaller index will be obtained as undegeneration keeping the top level
passage. In symbols, if $\Gamma\in \LG_2(\ol{B})$
is a graph corresponding to a boundary stratum in $D_i\cap D_j$ with $i<j$,
then $D_{\delta_1(\Gamma)}=D_i$ and $D_{\delta_2(\Gamma)}=D_j$. This is possible
because of  \cite[Proposition 5.1]{EC}, or equivalently it is a total order
refining the partial order defined in \cite[Proposition 3.5]{CoMoZadiffstrata}. 
\par
Consider now the set $\frakC$ of relevant curves for~$\ol{B}$, which we
write as the disjoint union of $\frakC_E$ and $\frakC_H$. The first set
$\frakC_E$, 'easy' to deal with, consists of relevant curves $C_{\Delta,i}$
such that $i=0$ or $\Delta \in \LG_i(\ol{B})$, i.e.\ the split level
is the top or bottom level of $\Delta$. 'Hard' to deal with is $\frakC_H = \frakC
\setminus \frakC_E$, i.e.\ the split level is strictly in between.  
\par
\autoref{le:cLwithrelevant} shows that $\cL_{\ol{B}}$ intersects
all relevant curves non-negatively, and intersects the easy ones positively. 
The strategy is to create positive intersections with curves in~$\frakC_H$
without destroying the positivity we already have. For this purpose we next decompose
$\frakC_H$ further.
\par
To a relevant curve $C_{\Delta,i} \in \frakC_{H}$ with $\Delta \in \LG_L(\ol{B})$ and $0 > i > -L$, 
we can associate
a three-level graph~$\wt{\Delta}$ by keeping only the levels right above
and below the critical level~$i$. In symbols $\wt{\Delta} = 
\delta_{-i,-i+1}(\Delta) \in \LG_2(\ol{B})$. We define $\frakC_H^{a,b}$,
for $a<b$, to be the set of relevant curves $C_{\Delta,i}$ whose associated
graph~$\wt{\Delta}$ labels a component of the intersection $D_a\cap D_b$. In other words, if $D_\Delta$ is a component of $D_{k_1}\cap \cdots \cap D_{k_L}$ with $k_1 < \cdots < k_L$, then  
$k_{-i} = a$ and $k_{-i+1} = b$. 
\par
In what follows we only need to consider the boundary divisors $D_a$ where $\frakC_H^{a,b}$ is non-empty for some $b > a$. For notation simplicity we still label them in increasing order as $D_1, D_2, \ldots$. We start by considering~$\frakC_H^{1,b}$ for $b>1$ and the line bundle
$\cL_{\ep_1} = \cL_{\ol{B}} \otimes \cO(-\ep_1 D_1)$ for $\ep_1 >0$. We claim
that for $\ep_1$ small enough
\be \c_1(\cL_{\ep_1}) \cdot [C] > 0 \quad \text{for} \quad C \in \frakC_E \cup
\bigcup_{b > 1}\frakC_H^{1,b}
\ee
and $\c_1(\cL_{\ep_1}) \cdot [C] \geq 0$ for the remaining relevant curves. First, 
to justify the last part it suffices to show that $[D_1] \cdot [C] = 0$
for $C \in \frakC_H^{a,b}$ with $a > 1$. If $D_1$ contains $C = C_{\Delta, i}$, then
$1<a$ implies that $D_1$ can be obtained as an undegeneration keeping the level passage
one or more above level~$i$. Namely, in the above notation $D_\Delta$ is a component of 
$D_{k_1}\cap \cdots \cap D_{k_L}$ with $k_1 = 1$, $k_{-i} =a$, $k_{-i+1} = b$ and $ -i > 1$.  
Pull back the normal bundle $\cN_{D_1}$ successively 
to $D_{1} \cap D_{k_2} \cap \cdots \cap D_{k_j}$ with $j$ varying from $2$ to $L$, and apply~\cite[Corollary 7.7]{EC}. 
We obtain that $[D_1]\cdot [C]$ is equal to the degree of $\cN_{\Delta, \delta_{1}^\complement(\Delta)}$ on $C$, 
which is zero by using~\autoref{eq:nbinPi} (with $i = 0$ therein and noting that 
a relevant curve obtained by splitting level~$i$ has zero intersection with $\c_1(\cL^{[j]}_\Delta)$ and $\xi_\Delta^{[j]}$ for $j \neq i$). 
So the only way the claim can fail is
that~$D_1$ and~$C=C_{\Delta,i}$ intersect in one of the two graphs
$\Delta^A$ or $\Delta^B$. Undegenerating the levels away from those adjacent to 
level~$i$ of $\Delta^A$ or $\Delta^B$, we obtain a four-level graph
$\Delta'$ such that
$\delta_1(\Delta')=a$, $\delta_2(\Delta')=1$ and $\delta_3(\Delta)=b$.
This graph would correspond to an intersection of $D_a$ and $D_1$ 
with $D_a$ on top. By the initial choice of the ordering of the divisors,
this is not possible. Second, for $C = C_{\Delta,i} \in \frakC_H^{1,b}$, in the above notation it means $i = -1$. 
Pull back the normal bundle $\cN_{D_1} $ successively 
to $D_{k_1} \cap D_{k_2} \cap \cdots \cap D_{k_j}$ with $k_1 = 1$, $k_2 = b$ and $j$ varying from $2$ to $L$, and apply~\cite[Corollary 7.7]{EC}. We thus obtain that 
\be
-\ep_1 [D_1] \cdot [C_{\Delta, -1}] \= - \ep_1\,\deg\big(\cN_{\Delta, \delta_{1}^\complement(\Delta)}|_{C}\big)
\= -\frac{\ep_1}{\ell_{\Gamma, 1}} \, \xi_\Delta^{[-1]} \cdot [C_{\Delta, -1}] > 0
\ee
where the last equality and inequality follow from~\autoref{eq:nbinPi} 
and~\autoref{lem:contribatlevel} respectively.  
Third, the positivity for $C\in \frakC_E$ already established is not
destroyed for~$\ep_1$ small enough.
\par
We next consider $\frakC_H^{2,b}$  for $b>2$ and the line bundle
$\cL_{\ep_1,\ep_2}= \cL_{\ep_1} \otimes \cO(-\ep_2 D_2)$. Again, we claim that
for $\ep_2$ small enough
\be \c_1(\cL_{\ep_1,\ep_2}) \cdot [C] > 0 \quad \text{for} \quad C \in \frakC_E \cup
\bigcup_{a \in \{1,2\} \atop b > a}\frakC_H^{a,b}
\ee
and $\c_1(\cL_{\ep_1,\ep_2}) \cdot [C] \geq 0$ for the remaining relevant curves.
As in the previous step, the claim about the remaining curves follows from the
ordering of the divisors $D_i$. The claim about $\frakC_H^{2,b}$ follows
from the description of the normal bundle, and the positivity of the intersection 
pairing with curves in $\frakC_E$ and $\frakC_H^{1,b}$ is not destroyed
for $\ep_2$ small enough.
\par
Iteratively we can define $\cL_{\ep_1,\ldots,\ep_{j}}= \cL_{\ep_1,\ldots,\ep_{j-1}} \otimes
\cO(-\ep_j D_j)$ with $\ep_j$ small enough until at the last step the resulting line bundle class intersects all relevant curves positively.
\end{proof}
\par
\begin{example}\label{ex:rhombus} {\rm
The multi-scale space $\ol{B}=\bP\Xi\cM_{1,4}(2,0,0,-2)$ has six relevant curves.
Consider first the set $\frakC_E$ of 'easy' curves. There are four curves
in this set, which are described in \autoref{figure:rhombuseasy}. The first
two are supported on divisors, so we need to impose an extra codimension-one condition, for example a $\psi$-decoration, in order to make them into curve
classes, while the last two are honest one-dimensional boundary strata.
In all of these cases, the split level is the bottom one. Hence by 
\autoref{le:cLwithrelevant} they intersect positively with $\c_1(\cL_{\ol{B}})$.
One can check that the intersection numbers of $\c_1(\cL_{\ol{B}})$ with these
four curves are given by
	\[\c_1(\cL_{\ol{B}})\cdot [C_{\Delta_1,-1}]=\frac{1}{8}\,,\quad  \c_1(\cL_{\ol{B}})\cdot [C_{\Delta_2,-1}]= \frac{1}{2}\,,\] 
	\[ \c_1(\cL_{\ol{B}})\cdot [C_{\Delta_3,-2}]=\frac{1}{2}\,, \quad \c_1(\cL_{\ol{B}})\cdot [C_{\Delta_4,-2}]=1\,.  \]
\par
\begin{figure}[ht]
	\bes
	C_{\Delta_1,-1}=\left[
	\begin{tikzpicture}[
	baseline={([yshift=-.5ex]current bounding box.center)},
	scale=2,very thick,
	bend angle=30]
	\node[comp] (T) [] {};
	\node[circled number] (T-1) [right=of T] {$1$};
	\path (T) -- (T-1) node[comp, midway, above=of T] (B) {};
	\path (B) edge 
	node [order top left,xshift=2pt,yshift=-4pt] {} (T);
	\path (B) edge 
	node [order bottom right] {} 
	node [order top right,xshift=-2pt,yshift=-4pt] {} (T-1);
	\draw (B) -- +(90: .2) node [xshift=0,yshift=2] {$-2$};
	\draw (T) -- +(245: .2) node [xshift=-2,yshift=-5] {$0$};
	\draw (T) -- +(290: .2) node [xshift=2,yshift=-5] {$0$};
	\draw (T-1) -- +(270: .2) node [xshift=2,yshift=-5] {$2$};
	\node [below left] at (T-1) {\small$\psi$};
	\end{tikzpicture}\right],
	\quad C_{\Delta_2,-1} =\left[
	\begin{tikzpicture}[
	baseline={([yshift=-.5ex]current bounding box.center)},
	scale=2,very thick,
	bend angle=30]
	\node[comp] (T) [] {};
	\node[comp] (T-1) [right=of T] {};
	\path (T) -- (T-1) node[comp, midway, above=of T] (B) {};
	\path (B) edge 
	node [order top left,xshift=2pt,yshift=-4pt] {} (T);
	\node[comp] (C) [right=of T] {}
	edge [bend left] 
	node [order bottom left] {} 
	node [order top left] {} (B)
	edge [bend right] 
	node [order bottom right] {} 
	node [order top right] {} (B);
	\draw (B) -- +(90: .2) node [xshift=0,yshift=2] {$-2$};
	\draw (T) -- +(245: .2) node [xshift=-2,yshift=-5] {$0$};
	\draw (T) -- +(290: .2) node [xshift=2,yshift=-5] {$0$};
	\draw (T-1) -- +(270: .2) node [xshift=2,yshift=-5] {$2$};
	\node [above left] at (B) {\small$\psi$};
	\end{tikzpicture}\right],
	\ees
	\bes
	\quad C_{\Delta_3,-2} =\left[
	\begin{tikzpicture}[
	baseline={([yshift=-.5ex]current bounding box.center)},
	scale=2,very thick,
	bend angle=30]
	\node[comp] (T) [] {};
	\node[comp] (T-1) [right=of T] {};
	\path (T) -- (T-1) node[comp, midway, above=of T] (B) {};
	\node[comp] (D) [above=of B] {};
	\path (D) edge 
	node [order top left,xshift=2pt,yshift=-4pt] {} (T);
	\path (D) edge 
	node [order top left,xshift=2pt,yshift=-4pt] {} (B);
	\node[comp] (C) [right=of T] {}
	edge [bend left] 
	node [order bottom left] {} 
	node [order top left] {} (B)
	edge [bend right] 
	node [order bottom right] {} 
	node [order top right] {} (B);
	\draw (D) -- +(90: .2) node [xshift=0,yshift=2] {$-2$};
	\draw (T) -- +(245: .2) node [xshift=-2,yshift=-5] {$0$};
	\draw (T) -- +(290: .2) node [xshift=2,yshift=-5] {$0$};
	\draw (T-1) -- +(270: .2) node [xshift=2,yshift=-5] {$2$};
	\end{tikzpicture}\right],
	\quad C_{\Delta_4,-2} =\left[
	\begin{tikzpicture}[
	baseline={([yshift=-.5ex]current bounding box.center)},
	scale=2,very thick,
	bend angle=30]
	\node[comp] (T) [] {};
	\node[comp] (T-1) [right=of T] {};
	\path (T) -- (T-1) node[comp, midway, above=of T] (B) {};
	\node[comp] (D) [above=of B] {};
	\path (B) edge 
	node [order top left,xshift=2pt,yshift=-4pt] {} (T);
	\path (D) edge 
	node [order top left,xshift=2pt,yshift=-4pt] {} (B);
	\path (D) edge [bend left] 
	node [order top left,xshift=2pt,yshift=-4pt] {} (T-1);
	\path (T-1) edge 
	node [order top left,xshift=2pt,yshift=-4pt] {} (B);
	\draw (D) -- +(90: .2) node [xshift=0,yshift=2] {$-2$};
	\draw (T) -- +(245: .2) node [xshift=-2,yshift=-5] {$0$};
	\draw (T) -- +(290: .2) node [xshift=2,yshift=-5] {$0$};
	\draw (T-1) -- +(270: .2) node [xshift=2,yshift=-5] {$2$};
	\end{tikzpicture}\right]
	\ees
	\caption{The set $\frakC_E$ of 'easy' relevant curves in the space $\bP\Xi\cM_{1,4}(2,0,0,-2)$ which intersect positively with $\cL_{\ol{B}}$.}\label{figure:rhombuseasy}
\end{figure}
\par 
There are two 'hard' relevant curves, given by boundary strata defined
by three-level graphs where the split level is at level $-1$
(see \autoref{figure:rhombushard}). The
first curve $C_{\Delta_5,-1}$ is of cherry-banana type. The other curve $C_{\Delta_6,-1}$ is of rhombus type.
By \autoref{le:cLwithrelevant}, the
intersection numbers of $\c_1(\cL_{\ol{B}})$ with these two curves are zero. 
Let $D_1$ be the divisor of cherry type given by the undegeneration of 
$C_{\Delta_5,-1}$ keeping the first level passage. Let $D_2$ be the divisor of banana type given by the undegeneration of 
$C_{\Delta_6,-1}$ keeping the first level passage. One can check that 
$$[D_1]\cdot [C_{\Delta_5,-1}] = - \frac{1}{2}\,, \quad [D_2] \cdot [C_{\Delta_6,-1}] = -1\,.$$
\par
\begin{figure}
\bes C_{\Delta_5,-1}=\left[
\begin{tikzpicture}[
baseline={([yshift=-.5ex]current bounding box.center)},
scale=2,very thick,
bend angle=30]
\node[comp] (T) [] {};
\node[comp] (T-1) [right=of T] {};
\path (T) -- (T-1) node[comp, midway, above=of T] (B) {};
\path (B) edge 
node [order top left,xshift=2pt,yshift=-4pt] {} (T);
\path (B) edge 
node [order bottom right] {} 
node [order top right,xshift=-2pt,yshift=-4pt] {} (T-1);
\node[comp] (C) [below=of T-1] {}
edge [bend left] 
node [order bottom left] {} 
node [order top left] {} (T-1)
edge [bend right] 
node [order bottom right] {} 
node [order top right] {} (T-1);
\draw (B) -- +(90: .2) node [xshift=0,yshift=2] {$-2$};
\draw (C) -- +(270: .2) node [xshift=0,yshift=-5] {$2$};
\draw (T) -- +(245: .2) node [xshift=-2,yshift=-5] {$0$};
\draw (T) -- +(290: .2) node [xshift=2,yshift=-5] {$0$};
\end{tikzpicture}\right],
\quad C_{\Delta_6,-1} =\left[
\begin{tikzpicture}[
baseline={([yshift=-.5ex]current bounding box.center)},
scale=2,very thick,
bend angle=30]
\node[comp] (T) [] {};
\node[comp] (T-1) [right=of T] {};
\path (T) -- (T-1) node[comp, midway, above=of T] (B) {};
\path (T) -- (T-1) node[comp, midway, below=of T] (C) {};
\path (B) edge 
node [order top left,xshift=2pt,yshift=-4pt] {} (T);
\path (B) edge 
node [order bottom right] {} 
node [order top right,xshift=-2pt,yshift=-4pt] {} (T-1);
\path (C) edge 
node [order top left,xshift=2pt,yshift=-4pt] {} (T);
\path (C) edge 
node [order top left,xshift=2pt,yshift=-4pt] {} (T-1);
\draw (B) -- +(90: .2) node [xshift=0,yshift=2] {$-2$};
\draw (C) -- +(270: .2) node [xshift=0,yshift=-5] {$2$};
\draw (T) -- +(245: .2) node [xshift=-2,yshift=-5] {$0$};
\draw (T-1) -- +(290: .2) node [xshift=2,yshift=-5] {$0$};
\end{tikzpicture}\right]\ees
\caption{The set $\frakC_H$ of 'hard' relevant curves in the space $\bP\Xi\cM_{1,4}(2,0,0,-2)$ which have zero intersection with $\cL_{\ol{B}}$.}
\label{figure:rhombushard}
\end{figure}
\par
One can also see that besides $C_{\Delta_5,-1}$ among all the relevant curves $D_1$
intersects non-trivially only $C_{\Delta_1,-1}$ 
which evaluates to $-1/8$, while $D_2$
intersects trivially all the relevant curves apart from $C_{\Delta_{6},-1}$.
 Hence in this case 
the line bundle $\c_1(\cL_{\ol{B}})-\ep_1 [D_1]-\ep_2[D_2]$ with $\ep_1>0$ and $\ep_2>0$ is $f$-ample.  
\par
In this example, one can also compute (with the help of
{\tt diffstrata}) the cone of $f$-ample divisors. It is a $13$-dimensional
polyhedron defined as the convex hull of $1$~vertex, $9$~rays, and $5$~lines. 
}\end{example}


\section{The age of automorphisms }
\label{sec:newage}

This section prepares for determining non-canonical singularites
in $\bP\MScoarse$ in the subsequent \autoref{sec:canonical}.
The first step is to determine ages of automorphisms, since we want
to apply a variant of the Reid--Tai criterion for canonical
singularities (\cite{ReidYoung}, \cite{TaiKodAg}). We give age estimates
in two situations, first for interior points of strata allowing
permuted marked points and second at the boundary, restricting there
to fixed marked points since this is our main goal.
\par
From now on let $B = \bP \omoduli[{g,n}](\mu)$ be the moduli space of Abelian differentials of type $\mu\in \bZ^n$, so of possibly meromorphic type, and
$\ol{B}=\bP\LMS$ be its compactification by multi-scale differentials. For
the subsequent propositions we split the set of marked points~$\bfz$
of a multi-scale differential into the set~$Z$ of zeros ($m_i \geq 0$)
and into the set~$P$ of poles ($m_i < 0$).  
\par 
First of all we consider automorphisms of Abelian differentials with possibly unlabeled
marked points. Recall that the local deformation
space of an Abelian differential $(X, \omega, Z, P)$ can be identified
with $H^1(X \setminus P, Z; \bbC)$. The tangent space to the associated
projectivized stratum is thus naturally identified with the affine
space $\bA_\rel(X)$ introduced in \autoref{eq:affinechart}, where in this
case we only deal with a one-level graph consisting of a single vertex. When computing the age of
an automorphism, it will be on this affine space throughout. We write
$\omoduli[{g}](\mu)$ for the spaces with unlabeled zeros and poles.  
\par
\begin{prop} \label{nprop:can_interior}
Let $(X,Z,P)$ be a (stable) pointed smooth curve with an automorphism $\tau$ of order $k \geq 2$,
fixing zeros and poles setwise (but not necessarily pointwise) and fixing
projectively an Abelian differential~$\omega$ of type $\mu$ with the zeros~$Z$
and poles~$P$. Let $\zeta$ be any
primitive $k$-th root of unity, and let $\zeta^{a'_1}, \ldots, \zeta^{a'_d}$
be the eigenvalues of the induced action on $\bA_\rel(X)$,
where $0 \leq a'_i < k$ and $d = \dim(B)$. Then  
$$\age (\tau|_{\bA_\rel(X)}) \,:=\, \sum_{i=1}^{d}
\frac{a'_i}{k} \,\geq\, 1\,, $$
except for the cases listed in \autoref{ncap:noncan}.
\end{prop}
\par
\begin{figure}
    \centering      
$$
\begin{array}{|c|l|c|c|c|}
\hline &&&& \\[-\halfbls]
\multicolumn{1}{|c|}{\text{Case}}&
\multicolumn{1}{|c|}{\text{Stratum}}&
\multicolumn{1}{|c|}{\text{order}}&
\multicolumn{1}{|c|}{\text{eigenvalues}} &
\multicolumn{1}{|c|}{\text{age}\geq }\\
[-\halfbls] & & & &\\
\hline  
&&&& \\ [-\halfbls]
(H) & \text{hyperelliptic differentials}  & k=2 & &0\\
[-\halfbls] & &&&\\ \hline &&&& \\[-\halfbls]
(1) & \bP\omoduli[1](\{0,0\})  & k=2 & 1, 1,-1 &1/2\\
[-\halfbls] & &&&\\ \hline  &&&& \\ [-\halfbls]
(2) & \bP\omoduli[0](\{m,m,m\},-3m-2)  & k=3 & \zeta_3, \zeta_3^2&1/3\\
[-\halfbls] & &&&\\ \hline &&&& \\ [-\halfbls]
(3) &\bP\omoduli[2](m, 2-m), \,\, 1\neq m \equiv 1 \mod 3   & k=3
& \zeta_3, \zeta_3, \zeta_3^2, \zeta_3^2 &2/3\\
[-\halfbls] & &&&\\ \hline  &&&& \\ [-\halfbls]
(4) & \bP\omoduli[1](0)   & k=3 & \zeta_3, \zeta_3^2&1/3\\
[-\halfbls] & &&&\\ \hline  &&&& \\ [-\halfbls]
(5) & \bP\omoduli[1](m,-m),\,\, 0 \neq  m \equiv 0 \mod 3
 & k=3 & \zeta_3, \zeta_3^2&1/3\\
[-\halfbls] & &&&\\ \hline  &&&& \\ [-\halfbls]
(6) & \bP\omoduli[1](\{0,0,0\})  & k=3 & \zeta_3, \zeta_3,
\zeta_3^2, \zeta_3^2 &2/3\\
[-\halfbls] & &&&\\ \hline  &&&& \\ [-\halfbls]
(7) & \bP\omoduli[1](\{m,m,m\},-3m), \,\, m \neq 0  & k=3 &
\zeta_3, \zeta_3, \zeta_3^2, \zeta_3^2 &2/3\\
[-\halfbls] &&&&\\ \hline  &&&& \\ [-\halfbls]
(8) & \bP\omoduli[0](\{m, m, m, m\}, -4m-2)  & k=4 & i,-1,i^3 &3/4\\
[-\halfbls] & &&&\\ \hline  &&&& \\ [-\halfbls]
(9) & \bP\omoduli[1](0) & k=4 & i, i^3&1/2\\
[-\halfbls] & &&&\\ \hline  &&&& \\ [-\halfbls]
(10) & \bP\omoduli[1](m,-m),\,\, 0 \neq m \,\, 
\text{even}  & k=4 & i, i^3&1/2 \\
[-\halfbls] &&&&\\ \hline  &&&& \\ [-\halfbls]
(11) & \bP\omoduli[1](\{0,0\})  & k=4 & i,-1,i^3 &3/4\\
[-\halfbls] &&&&\\ \hline  &&&& \\ [-\halfbls]
(12) & \bP\omoduli[1](\{m,m\},-2m),\,\, 0 \neq m\,\,
\text{even} & k=4 & i,-1,i^3 &3/4\\
[-\halfbls] &&&&\\ \hline  &&&& \\ [-\halfbls]
(13) & \bP\omoduli[1](0) & k=6 & \zeta_6 ,\zeta_6^5 &1/3\\
[-\halfbls] && &&\\ \hline 
\end{array}
$$
\caption{Automorphisms with $\age < 1$. The column of eigenvalues corresponds to the
induced action on the unprojectivized chart in $H^1(X \setminus P, Z; \bbC)$. Here $m$ can be possibly negative.}
\label{ncap:noncan}
\end{figure}
\par
\medskip
Consider now the space $\ol{B}$. Recall from \autoref{sec:Coordbd} the local
coordinate system near
a multi-scale differential $(X,\bfomega,\bfz,\bfsigma)\in \ol{B}$
compatible with the enhanced level graph~$\Gamma$, in particular the
decomposition of the coordinates in~\autoref{eq:bdlocalstructure}. Moreover
there we explained that automorphisms of multi-scale differentials
are only well defined on the quotient of the affine space
$\bA = \bA_\rel(X) \times  \bC_{\hor} \, \times \, \bC_{\lev}$  by the group~$K_\Gamma$.
On the other hand the age of an automorphism
is defined only for a linear action. We abuse this definition
for $\bftau \in \Aut(X,\bfomega,\bfz,\bfsigma)$ and say that $\age(\bftau) \geq C$
if each lift of~$\bftau$ to an automorphism of~$\bA$ has $\age \geq C$. In particular, if the induced action of $\bftau$ 
on $\bA_\rel(X)$ or on $\bA_\rel(X)\times \bC_{\hor}$ has age $\geq C$, then $\age (\bftau) \geq C$. 
\par
\begin{prop} \label{nprop:age-vertex}
 Let $\Gamma$ be a level graph representing a boundary stratum in $\ol{B}$ and let
$\bftau=(\tau_{(-i)})$ be an automorphism of a multi-scale differential
compatible with~$\Gamma$, fixing the labeled points. Suppose moreover
that~$\bftau$ fixes a vertex~$v$. Then $\age(\bftau)\geq 1$, if~$v$ does not
belong to the lists in \autoref{ncap:noncan} and \autoref{ncap:noncanRC}.
\end{prop}
\par
\begin{figure}
    \centering      
$$
\begin{array}{|c|l|c|c|c|}
\hline &&&& \\[-\halfbls]
\multicolumn{1}{|c|}{\text{Case}}&
\multicolumn{1}{|c|}{\text{Stratum}}&
\multicolumn{1}{|c|}{\text{order}}&
\multicolumn{1}{|c|}{\text{eigenvalues}} &
\multicolumn{1}{|c|}{\text{age}\geq } \\
[-\halfbls] & &&&\\ \hline &&&& \\ [-\halfbls]
(RH) &  \text{hyperelliptic differentials} & k=2 &  &0\\
[-\halfbls] & & &&  \\ \hline &&&& \\ [-\halfbls]
(R1) &  \bP\omoduli[0]^{\frakR}(3m_1 + m_2-2, -m_2,\{-m_1\}^3),  & k=3 &
\zeta_3, \zeta_3^2 &1/3\\
[-\halfbls] & &&& \\  &&&& \\ [-\halfbls]
&m_2\not\equiv 1 \mod 3, & \ & \\
[-\halfbls] & &&& \\  &&&& \\ [-\halfbls]
& \frakR = \{ r_2 = 0\}\,\, \text{or}\,\, \frakR = \{r_2=0, \, r_3+r_4+r_5 =0\} &&\\
[-\halfbls] & &&&\\ \hline &&&& \\ [-\halfbls]
(R2) &  \bP\omoduli[0]^{\frakR}(\{m_1+m_2-1\}^3, -3m_1, 1-3m_2), & k=3
& \zeta_3, \zeta_3^2&1/3 \\
[-\halfbls] & &&& \\  &&&& \\ [-\halfbls]
& \frakR = \{r_4=0\} \,\, \text{or}\,\, \frakR = \{r_5=0\} \,\, \text{or}\,\, \frakR = \{r_4=r_5=0\} &&\\
[-\halfbls] && &&\\ \hline &&&& \\[-\halfbls]
(R3) &  \bP\omoduli[2,4]^{\frakR}(\sum_{i=1}^3 m_i +2,-m_1,-m_2,-m_3),  
&  k=3 & \zeta_3, \zeta_3, \zeta_3^2, \zeta_3^2&2/3 \\
[-\halfbls] && &&\\ &&&& \\[-\halfbls]
& m_i\not\equiv 1 \mod 3, \,\, \sum_{i=1}^3 m_i \not\equiv 0 \mod 3, &&  \\
[-\halfbls] && &&\\ &&&& \\[-\halfbls]
& \frakR = \{r_2= r_3 =0\} \,\, \text{or}\,\,\frakR =\{r_2= r_3=r_4=0\}&&\\
[-\halfbls] & &&&\\ \hline &&&& \\ [-\halfbls]
(R4) & \bP\omoduli[2,4]^{\frakR}(m_1+m_2+2,-m_1,-m_2), &  k=3 
&  \zeta_3, \zeta_3, \zeta_3^2, \zeta_3^2  &2/3\\
[-\halfbls] & &&&\\  &&&& \\ [-\halfbls]
& m_i \not\equiv 1 \mod 3, \,\,m_1+m_2\not\equiv 0 \mod 3, &&\\
[-\halfbls] & &&&\\  &&&& \\ [-\halfbls]
& \frakR = \{r_2 =0\} \,\, \text{or}\,\,\frakR = \{r_2=r_3 = 0\} &&\\
[-\halfbls] & &&&\\ \hline &&&& \\ [-\halfbls]
(R5) & \bP\omoduli[1,3]^\frakR( m_1+m_2, -m_1, -m_2),   &  k=3
& \zeta_3,  \zeta_3^2  &1/3\\
[-\halfbls] & &&&\\  &&&& \\ [-\halfbls]
& m_i\not\equiv 1 \mod 3, \,\, m_1+m_2\not\equiv 2 \mod 3,&& \\
[-\halfbls] & &&&\\  &&&& \\ [-\halfbls]
& \frakR = \{r_2 =0\} \,\, \text{or}\,\,\frakR = \{r_2=r_3 = 0\} &&\\
[-\halfbls] & &&&\\ \hline &&&& \\ [-\halfbls]
(R6) &  \bP\omoduli[1]^{\frakR}(\{m_1+m_2\}^2,
-2m_1,-2m_2), &  k=4 &  i, -1, i^3 &3/4\\
[-\halfbls] & &&&\\  &&&& \\ [-\halfbls]
&  m_1+m_2 \,\, \text{even}, && \\
[-\halfbls] & &&&\\  &&&& \\ [-\halfbls]
& \frakR = \{r_3 =0\} \, \text{or}\,\,\frakR = \{r_3=r_4 = 0\} &&\\
[-\halfbls] & &&&\\ \hline &&&& \\ [-\halfbls]
(R7) &  \bP\omoduli[1]^{\frakR}(4m_1 + 2m_2, \{-2m_1\}^2, -2m_2),  & k=4 & i, -1, i^3 &3/4\\
[-\halfbls] & &&& \\  &&&& \\ [-\halfbls]
&  \frakR = \{r_4 =0\} \, \text{or}\,\,\frakR = \{r_2 + r_3 = 0\} &&\\
[-\halfbls] &&& &\\ \hline 
\end{array}
$$
\caption{Vertices with residue conditions that can yield $\age(\bftau) < 1$. Here the $m_i$ are always positive. Each $r_i$ denotes the residue of the $i$-th entry in the signature. Alternative versions of the GRC are equivalent by the Residue Theorem and relabeling the poles.}
\label{ncap:noncanRC}
\end{figure}
\par
\begin{rem} \label{rem:hypeigen}
{\rm 
Let~$X$ be a hyperelliptic curve of genus~$g$ with $f_z$ fixed zeros,
$f_p$ fixed poles, $c_z$ conjugate pairs of zeros and $c_p$ conjugate
pairs of poles, under the hyperelliptic involution $\tau$. It is easy to check that 
the eigenvalue decomposition of the $\tau$-action on the (unprojectivized) relative periods
(if unconstrained by the GRC) is
\bes
(-1)^{2g}\,, \ (+1)^{c_z+ f_z-1}\,, \ (-1)^{c_z}\,,\   (+1)^{c_p+ f_p-1}\,,\  (-1)^{c_p}\, 
\ees
where the fourth term is empty if $c_p+f_p = 0$ (i.e., when there is no pole).  
In particular, if (the projectivized) $\age(\tau) < 1$ then $c_z+ f_z \leq 2$
and $c_p+ f_p \leq 2$, and not both of these are equal to~$2$. For later use we also need to consider hyperelliptic involutions with GRC constraints. In particular if $\age(\tau) = 0$, then every fixed pole of $X$ is constrained by the GRC to have zero residue and every conjugate pair of poles of $X$ is constrained by the GRC to have the sum of the residues equal to zero. 
}
\end{rem}
\par
\begin{rem} \label{rem:congruence}
{\rm 
Conversely, all the cases listed in the tables can be realized via canonical covers of $k$-differentials
on $X/\langle \tau \rangle$, and all the congruence conditions of the signatures are
necessary as well. We go over one case in each table and leave the others for the
reader to verify.  For example for $\bP\omoduli[0](\{m, m, m, m\}, -4m-2)$ in
Case (8), take a quartic differential of signature $(4m, -4m-5, -3)$ in $\bP^1$
and pull it back via the canonical quartic cover totally ramified at the last
two singularities, where the second ramification point over the pole of order three becomes an unmarked ordinary
point. Next for the cases in \autoref{ncap:noncanRC}, first note that if
$\omega$ is not a $\bftau$-invariant form, then any $\bftau$-fixed pole of
$\omega$ must have zero residue. In addition, if there is a totally ramified point
of multiplicity $k$ under the $\bftau$-action, then $\omega$ is an invariant form
if and only if the singularity order at the ramification point is
$\equiv k-1 \mod k$. Using these, consider Case (R1) as an example.  One can
pull back a cubic differential of signature $(3m_1+m_2-4, -m_2-2, -3m_1)$ on~$\bP^1$
via the canonical triple cover totally ramified at the first two singularities,
whose cubic root thus gives the desired $\omega$ satisfying the residue condition
$\frakR$ as well. 
}
\end{rem}
\par
\medskip
In order to prove the above propositions, we make some preparation first.
Let $(X,\omega,Z,P)$ be an Abelian differential and $\tau$ be an automorphism of order $k$,
either in the context of \autoref{nprop:can_interior} or as the restriction
of $\bftau$ from \autoref{nprop:age-vertex} to the surface at the vertex~$v$.
We fix a primitive $k$-th root of unity~$\zeta = \zeta_k$
and $a \in \{0,\dots,k-1\}$ such that
\be \label{neq:adef}
\tau^{*}\omega \= \zeta^{a} \omega\,.
\ee
Let~$\gamma$ be a homology class such that $\tau_{*}\gamma = \zeta^{a_i} \gamma$. Then 
$$\zeta^{a_i}\int_{\gamma} \omega \= \int_{\tau_{*}\gamma} \omega \=
\int_{\gamma}\tau^{*}\omega \= \zeta^a \int_{\gamma} \omega\,. $$
Therefore, the eigen-period $\int_{\gamma} \omega$ must be zero if $a_i \neq a$.
Since the periods of $\omega$ cannot be all equal to zero, there must exist
some $a_i$ equal to~$a$ such that the corresponding eigen-period of~$\omega$ is
nonzero, and hence we can use it to projectivize the domain of periods as well as
the induced action. Then each of the exponents of the projectivized
action is $a'_i = a_i - a \mod k$, where
we use representatives with $0 \leq a'_i < k$ throughout.
\par
\medskip
The first lemma below gives a lower bound for the age contribution of a
$\tau$-orbit of zeros or poles in the context of \autoref{nprop:can_interior}.
\par
\begin{lemma}\label{lem:ageinterior}
Let $\tau$ be a non-trivial automorphism of order $k$ of a (stable) pointed smooth 
curve~$(X,Z,P)$, fixing zeros and poles setwise (but not necessarily pointwise)
and fixing projectively an Abelian differential~$\omega$ with the zeros~$Z$ and
poles~$P$. Let $\{x_1,\dots,x_{k'}\}$ be a $\tau$-orbit of unlabeled zeros or poles,
where $k=k' \ell$. Consider the subspace $U\subseteq \bA_\rel(X)$ generated by
the $k'-1$ relative periods joining the $x_i$ if they are zeros, or by the $k'-1$
loops at each of the $x_i$ if they are poles. Then for $a\in \{0,\dots,k-1\}$ as
in~\autoref{neq:adef} we have 
\begin{equation}\tag{I}\label{eq:ageinterior}
\age(\tau|_{U})\, \geq \,\frac{1}{k}  \sum_{i=1}^{k'-1} (\ell i - a) \mod k
\,\geq \, \frac{1}{2k'}(k'-2)(k'-1)\,.
\end{equation}
In particular, in this case $\age(\tau) \geq 1$ if $k' \geq 5$. 
\end{lemma}	
\begin{proof}
The unprojectivized eigenvalues of the $\tau$-action restricted to the subspace~$U$
are $\zeta_{k'},\dots,\zeta_{k'}^{k'-1}$, where $\zeta_{k'}$ is a primitive $k'$-th root
of unity. In other words, these eigenvalues are  $\zeta_{k}^\ell,\dots,
\zeta_{k}^{\ell(k'-1)}$, which thus implies the first inequality. Note that the sum
in the bound consists of $k'-1$ distinct numbers in $[0, k-1]$ that belong to
the same congruence class mod~$\ell$. Thus its minimum is attained at
$\frac{1}{k'}\sum_{j=0}^{k'-2} j$.
\end{proof}
\par
In the next lemma we are in the context of an automorphism $\bftau$ of a multi-scale differential compatible with $\Gamma$, where $\Gamma$ is a level graph corresponding to a boundary component of $\ol{B}$. Consider a vertex $v$ of the graph and denote by $(X, \omega, Z, P)$ the differential associated to it. We define the \emph{multi-vertex RC-independent subspace} $V\subset H_1(X\setminus P, Z;\bC)\oplus \bbC_{\hor}(X)$ to be the largest subspace such that periods of $(X, \omega)$ in $V$ (including residues of the poles and plumbing parameters for horizontal edges of $X$) can vary independently without being constrained by the other vertices of $\Gamma$ due to global and matching residue conditions. In other words, $V$ is the largest subspace of parameters associated to $X$ such that the dimension of $\bigoplus_{i=0}^{d-1} \tau^i (V)$ does not drop after imposing the residue conditions ${\frak R}$ and matching horizontal plumbing parameters (in case a horizontal edge joins two permuted vertices). Note that $V$ always contains the absolute periods of $X$ and relative periods that join between zeros of $X$. We denote by $M$ the dimension of $V$.  
\par
\begin{lemma} \label{lem:permute}
Suppose~$\bftau$ permutes~$d$ vertices of~$\Gamma$ and let $\tau$  be its
restriction to the homology of these vertices. Let $M$ be the dimension of the
multi-vertex RC-independent subspace of each vertex. Then 
\bes
\age(\tau) \,\geq\, \frac{M}{2}(d-1)\,. 
\ees
In particular, $\age (\tau)\geq 1$ if $d\geq 3$, $M > 0$ and if $d = 2$, $M > 1$.  
\end{lemma}
\par
If a permuted vertex has genus $g$ with $n$ zero edges, then $M \geq 2g + n -1$ because $V$ contains the subspace of absolute and relative periods. In particular, the case $M \leq 1$ can only occur for $g = 0$ with $n\leq 2$. 
\par
\begin{proof}[Proof of \autoref{lem:permute}]
The automorphism $\bftau$ permutes these $d$ vertices and for each one we have
an $M$-dimensional subspace $V$ of homology which is cyclically permuted among the
vertices.  Then the restricted action of $\tau$ to the sum of these $M$-dimensional
subspaces can be described in a suitable basis as
\be \label{eq:permmat}
\begin{pmatrix}
		0_{M\times M(d-1)} &A_{M \times M}\\
		I_{M(d-1)\times M(d-1) } &0_{M(d-1)\times M}
\end{pmatrix}
\ee
where~$A$ is the matrix representing the automorphism $\tau^d$ acting on $V$. Suppose the action of $\tau^d$ on $V$ has order $k$. If the eigenvalues of $A$ are given by $\zeta_{k}^{a_i}$ for $i=1,\dots,M$, then the eigenvalues of the full matrix are given by $\zeta_{d k}^{a_i+jk}$ for all $i=1,\dots,M$ and $j=0,\dots,d-1$. In the case of $M>0$, there exists a nonzero
eigen-period corresponding to certain eigenvalue $\zeta_{d k}^{a_{i_0}+jk}$, which can be
used to projectivize the action. Therefore, $\age(\tau)$ is bounded below by $M$ sums of type 
$\frac{1}{dk}\sum_{j=0}^{d-1} (jk + a') \mod dk$ for some $a' = a_i - a_{i_0}$.  Since each sum consists of the representatives of the same congruence class ($a' \mod k$) in consecutive subintervals of length $k$ in $[0, dk-1]$, its minimum is attained at $\frac{1}{dk} \sum_{j=0}^{d-1} jk = (d-1)/2$. We thus conclude that $\age (\tau) \geq M(d-1)/2$. 
\end{proof}
The following third lemma gives a lower bound for the age contribution of
a $\bftau$-orbit of zeros or poles from the permutation representation on the space of residues
and relative cycles on a $\bftau$-fixed vertex and also on the
rest of a multi-scale differential whose marked points are fixed.
We define $ \bA_\rel(X_{>-J}) = \prod_{-i>-J}
\bA_\rel(X_{-i})$ and similarly for~$X_{<-J}$ to parameterize periods of the vertices above or below level $-J$ (up to level-wise projectivization).   
\par
\begin{lemma}\label{nlem:ageGRC}
Let  $\Gamma$ be a level graph representing a boundary stratum in $\ol{B}$ and
let $\bftau=(\tau_{(-i)})$ be an automorphism
of a projectivized multi-scale differential compatible with~$\Gamma$.
Suppose $\bftau$ fixes a vertex~$v$ at level~$-J$ of~$\Gamma$. Denote the
restriction of the multi-scale differential to $v$ by~$(X,\omega,Z,P)$ and
suppose that the order of $\tau = \bftau|_v$ is $k = k'\ell$. 
Then the following estimates hold:
\begin{itemize}
\item[(i)] Let $\{x_1,\dots,x_{k'}\}$  be a $\bftau$-orbit of unlabeled
non-simple poles
of $(X,\omega)$ corresponding to the lower ends of edges ending at~$v$.
Consider the subspace $U\subseteq \bA_\rel(X_{(-J)})$ generated by the loops around these poles constrained by
the residue conditions~$\frakR$ imposed by the higher levels of~$\Gamma$.
Suppose that these edges are adjacent to $d$ connected components of the graph
$\Gamma_{>-J}$ at higher level. Then 
\begin{equation}\tag{P1}\label{neq:tauup}
	\age(\bftau|_{U})\, \geq \,\frac{1}{k}
	\sum_{\substack{i = 1 \\ (k'/d)\nmid i}}^{k'-1} (\ell i - a) \mod k\,, 	 
\end{equation}
\begin{equation}\tag{P2} \label{neq:tauotherp1} 
\age(\bftau|_{\bA_\rel(X_{(>-J)})}) \,\geq\,
d-1\,. 
\end{equation}
\item[(ii)] Let $\{x_1,\dots,x_{k'}\}$  be a $\bftau$-orbit of unlabeled zeros
of $(X,\omega)$ corresponding to higher ends of edges adjacent to~$v$.
Let $U\subseteq  \bA_\rel(X_{(-J)})$ be the subspace generated by the relative
periods between these zeros. Then 
\begin{equation}\tag{Z}\label{neq:tauuz}
\age(\bftau|_{U})\,\geq\, \frac{1}{k}  \sum_{i=1}^{k'-1} (\ell i - a)
\mod k \,\geq \,  \frac{1}{2k'} (k'-2)(k'-1)\,. 
\end{equation} 
\item[(iii)] Let $\{x_1,\dots,x_{k'}\}$  be a $\bftau$-orbit of unlabeled
simple poles of $(X,\omega)$ corresponding to horizontal edges, each of which has exactly one end adjacent to~$v$ at the $x_i$. Let $U\subseteq  \bA_\rel(X_{(-J)})$ be the subspace generated by the residue cycles at these poles. Then
\begin{flalign}
\age(\bftau|_{\bC_{\hor}})& \,\geq\,
\frac{1}{k'}\,\,\sum_{j=0}^{k'-1} j \mod k' \= \frac{k'-1}{2}\,, 
\tag{SPH1}\label{neq:tauhor1}
\end{flalign}
\begin{flalign}
\age(\bftau|_{U})& \,\geq\,
\frac{1}{k'}\,\,\sum_{j=1}^{k'-1} (j - a) \mod k' \,\geq\, \frac{1}{2k'}(k'-2)(k'-1)\,.
\tag{SPH2}\label{neq:tauhor2}
\end{flalign}
\item[(iv)] Let $\{x_1,\dots,x_m,y_1,\ldots,y_m\}$  be a
$\bftau$-orbit (i.e.\ $k' = 2m$) of unlabeled simple
poles of $(X,\omega)$ corresponding to horizontal edges
both of whose ends $x_i$ and $y_i$ are adjacent to~$v$ (i.e. they form self-nodes as $m$ loops in the dual graph at $v$). 
Let $U\subseteq  \bA_\rel(X_{(-J)})$ be the subspace generated by the residue cycles at these poles. Then 
\begin{flalign}
\age(\bftau|_{\bC_{\hor}})& \,\geq\,
\frac{1}{k'}\,\,\sum_{j=0}^{m-1} 2j \mod k' = \frac{m-1}{2}\,, 
\tag{SPS1}\label{neq:tauself1}
\end{flalign}
\begin{flalign}
\age(\bftau|_{U})& \,\geq\,
\frac{1}{k'}\,\,\sum_{j=1}^{m-1} (2j - a) \mod k' \geq \frac{1}{2m}(m-1)^2 \,.
\tag{SPS2}\label{neq:tauself2}
\end{flalign}
In all of the above~$a$ is analogously defined as in~\autoref{neq:adef} and in the last inequality it is an odd number. 
\end{itemize}
\end{lemma}
\par 
\begin{proof}
In the setting of~(i), the poles~$x_i$ can be grouped accordingly
into $d$ sets each of which has $k'/d$ elements adjacent to the same connected component of~$\Gamma_{>-J}$. 
Let $\gamma_i$ be a loop around each $x_i$ and consider the space generated by all the $\gamma_i$. 
Applying the GRC, the subspace $U$ is cut out by the equations
$\sum_{j=1}^{k'/d} \gamma_{i+jd} = 0$, for $i=1,\dots,d$. Then $(\tau_{(-J)})|_{U}$
is the regular representation of the cyclic group of order~$k'$ with all
representations induced from the subgroup of order~$k'/d$ removed.
This means that the eigenvalues of $(\tau_{(-J)})|_{U}$ are $\zeta_{k'}^{i}
= \zeta_{k}^{\ell i}$, for all $i\neq j k'/d$ with $j=1,\dots,d$. Hence
$\age (\tau_{(-J)})|_{U} \geq \frac{1}{k}\sum_{i\in \{1,\dots,k'\}: (k'/d)\nmid i}
(\ell i - a) \mod k$ for
some $a\in \{0, \dots,k-1\}$ where $\zeta_{k}^{a}$ is the eigenvalue for the
eigenform $\omega$ under the action of $\tau_{(-J)}$. This shows~\autoref{neq:tauup}.
\par
For the age estimate on higher level, the automorphism $\bftau$ permutes the top vertices
of the $d$~connected components of~$\Gamma_{>-J}$. Hence there exist $d'$ disjoint (locally) top vertices on some level $-J' > J$ that are permuted by $\tau_{(-J')}$ where $d \mid d'$. Note that (locally) top vertices do not admit any vertical polar edges, hence they are not constrained by the GRC. If the multi-vertex RC-independent subspace of these $d'$ vertices has dimension $M \geq 2$, then~\autoref{neq:tauotherp1} follows from \autoref{lem:permute}. The only possibility for $M \leq 1$ is when these $d'$ vertices have genus zero and each admits some horizontal edges joining between them (so that we need to match horizontal residues and plumbing parameters between them). Since these $d'$ vertices belong to $d$ connected components, in this case $d' \geq 2d$ and we obtain at least $d$ independent residue cycles and $d$ independent plumbing parameters (from at least $d'/2 \geq d$ horizontal edges) that are permuted by $\tau_{(-J')}$, which yields the desired bound by applying the argument of  \autoref{lem:permute} for the case $M = 2$. 
\par
In case~(ii) of zeros the argument for~\autoref{neq:tauuz} is the same as for
\autoref{neq:tauup}, using relative periods joining the zeros instead of loops
around the poles. Since there are no residue conditions, only the trivial
representation ($i=0$) has to be omitted. 
\par
In case~(iii) of horizontal nodes, suppose each horizontal edge has local equation $x_i y_i = h_i$ for $i = 1, \dots, k'$, where $h_i$ is the plumbing parameter, such that $x_1 \mapsto x_2, x_2\mapsto x_3, \dots, x_{k'}\mapsto \zeta_\ell^x x_1$ under $\bftau$, where the last one is due to that the action of $\bftau^{k'}$ multiplies $x_1$ by an $\ell$-th root of unity. Similarly suppose 
$y_1 \mapsto y_2, y_2\mapsto y_3, \dots, y_{k'}\mapsto \zeta_\ell^y y_1$ under $\bftau$. It follows that $h_1 \mapsto h_2, h_2 \mapsto h_3, \dots, h_{k'}\mapsto \zeta_\ell^{b} h_1$ under $\bftau$ for some $b = x+y \in \{ 0, \dots, \ell -1\}$. Then the associated eigenvalues as $k$-th roots of unity have exponents $b, b + \ell, \dots, b + (k'-1)\ell \mod k$ for $k = k' \ell$. Since in this case $\bftau^{k'}$ fixes a simple pole branch, it implies that $(\bftau^{k'})^{*}\omega = \omega$ restricted to $v$, and hence $a$ is divisible by $\ell$. Altogether it gives 
$$ \age(\bftau|_{\bC_{\hor}}) \,\geq\, \frac{1}{k} \sum_{j=0}^{k'-1} (b + j \ell) \mod k\,. $$
Clearly $b=0$ minimizes the bound, which gives~\autoref{neq:tauhor1}. For~\autoref{neq:tauhor2}, the proof is similar. The only possible difference is that the sum of the residue cycles might be trivial due to the Residue Theorem (combined with GRC from higher level if $v$ has other poles). Hence we can only obtain a bound by summing up $k'-1$ terms (instead of $k'$), which minimizes at $\frac{1}{k'}\sum_{j=0}^{k'-2} j$ as seen before.  
\par
In case~(iv) we can assume that $\bftau(x_i)=x_{i+1}$, $\bftau(y_i)= y_{i+1}$,
$\bftau(x_m) = y_1$ and $\bftau(y_m) = x_1$ by labeling the points appropriately.
In particular $\bftau^m$ swaps $x_i$ and $y_i$ for all~$i$. The opposite
residue condition implies that $(\bftau^m)^* \omega = -\omega$. Consequently
$a$ is odd. Using similar arguments as in (iii), we thus obtain~\autoref{neq:tauself1}
and~\autoref{neq:tauself2} respectively from the cyclic group action on the~$m$
horizontal plumbing parameters and from the residue cycles at the self-nodes
formed by gluing each $x_i$ with $y_i$.  
 \end{proof}
\par
Finally we need the following estimate for sums analogous to those appearing in
the preceding lemmas.  
\par
\begin{lemma} \label{nle:eulerphi}
Let $\phi_k(a) = \sum_{n} (n - a) \mod k$ where the sum ranges over all
the $\varphi(k)$ integers $1\leq n < k$ that are relatively prime to $k$.
Then $\phi_k(a) \geq k$ for any $a$ as long as $k \neq 2,3,4,6$.
\end{lemma}
\par
\begin{proof}
The set of integers relatively prime to $k$ equidistributes
in the intervals $(0,a]$ and $(a,k]$, in fact in any interval, with
an effective error rate of the number~$d(k)$ of divisors of~$k$ (see \cite[Lemma~1.4]{BIR} and simplify the argument therein by removing the extra congruence condition). Consequently, $\phi_k(a) \to k\varphi(k)/2$
for any~$a$ as $k \to \infty$ with controlled error terms. It thus suffices to
check small values of~$k$, which gives the list above.  
\end{proof}
\par
We can now prove \autoref{nprop:can_interior} and \autoref{nprop:age-vertex} at the same time. The main difference to keep in mind is that we can use the inequality \autoref{eq:ageinterior} of \autoref{lem:ageinterior} in the context of  \autoref{nprop:can_interior}, while we have to use the inequalities of \autoref{nlem:ageGRC} in the context of \autoref{nprop:age-vertex} due to possible residue conditions. We will also skip the verification of 
realizability and congruence condition of each case in the following already lengthy proof (see \autoref{rem:congruence} if the reader is interested). 
\par 
Throughout the proof we stick to the following notations. Let $X'$ be the quotient of $X$ by the group action generated by an automorphism $\tau$ of order $k$, and let $g'$ be the genus of $X'$.  Then the associated map $\pi\colon X \to X'$ is
a cyclic cover of degree~$k$ with the deck transformation group generated
by~$\tau$. Let~$Z'$ and~$P'$ be the $\pi$-images of $Z$ and $P$ respectively.
Let $b$ be the number of branch points and $s_i$, for $i=1,\dots,b$, be the
cardinality of each ramified fiber. Note that every $s_i$ divides $k$ since $\pi$ is a cyclic cover. Moreover in this case the
Riemann--Hurwitz formula gives 
\begin{equation}\label{neq:RHcyclic}
	2g-2 \= k(2g'-2)+bk-\sum_{i=1}^b  s_i\,.
\end{equation} 
We say that a fiber of $\pi$ is \emph{special} if it consists of zeros or 
poles of~$\omega$. Moreover, a special fiber is called a zero (resp. pole)
fiber if it consists of zeros (resp. poles) of $\omega$.
Note that a special fiber does not have to consist of ramification points,
and conversely, a ramified fiber does not have to be special.  
\par
\begin{proof}[Proof of \autoref{nprop:can_interior} and
\autoref{nprop:age-vertex}] Let  $\Gamma$ be a level graph representing
a boundary stratum in $\ol{B}$ and let $\bftau=(\tau_{(-i)})$ be an automorphism
of a projectivized multi-scale differential compatible with~$\Gamma$ fixing a
vertex~$v$. We denote by $\tau$ the restriction of $\bftau$ to $v$, and by
$(X,\omega,Z,P)$ the restriction of the multi-scale differential to $v$.
Let $\zeta^{a_1},\ldots, \zeta^{a_N}$ be
the eigenvalues for the induced action of $\tau$ on the homology
$H_1(X \setminus P, Z; \bbC)^\frakR$ (and on thus also the cohomology
$H^1(X \setminus P, Z; \bbC)^\frakR$), where~$\frakR$ are the residue
conditions in the general situation of \autoref{nprop:age-vertex}. 
\par 
Recall that if $\tau^{*}\omega = \zeta^a \omega$, then we can use a nonzero
eigen-period corresponding to the eigenvalue $\zeta^a$ to projectivize the induced
action (restricted to the level of $v$).  Then each of the exponents of the
projectivized action is $a'_i = a_i - a \mod k$, where $0 \leq a'_i < k$. 
\par
The $\tau$-invariant subspace of $H^1(X \setminus P, Z; \bbC)^\frakR$ can be
identified with a subspace of the cohomology of the quotient surface 
$H^1(X' \setminus P', Z'; \bbC)$, cut out by some residue constraints~$\frakR'$.
Independently of~$\frakR'$, the action of $\tau$ preserves the symplectic
paring of the absolute homology $H_1(X; \bbC)$. Hence the $2g$ eigenvalues
from the absolute part split into $g$ conjugate pairs of type $(\zeta^{a_i},
\zeta^{k-a_i})$ for $0< a_i \leq [k/2]$ and pairs of type $(1,1)$ if $a_i = 0$.
In particular, the sum of the exponents from the absolute part is divisible
by $k$, and is at least $k$ unless all absolute periods are $\tau$-invariant.
\par
\subsubsection*{\textbf{Case} $\bf{g'\not =0}$, $\bf{a=0}$} First consider
the case $a = 0$, i.e. $\omega$ is a $\tau$-invariant form. If
$g > g'$, then $H^1(X'; \bbC) \to H^1(X; \bbC)$ is not onto for the dimension
reason, hence the absolute periods are not all invariant. By the preceding
paragraph, we conclude that $\age(\tau) \geq k/k = 1$.  
\par
The opposite case $g \leq g'$ is only possible by Riemann--Hurwitz if
$g = g' = 1$ or $g = g' = 0$. Suppose $g = g' = 1$.  Then $\pi$ is an elliptic
isogeny with no ramification, and consequently there is at least one unramified zero fiber.
Applying to this fiber \autoref{eq:ageinterior} or \autoref{neq:tauuz} (for $a = 0$, $\ell = 1$ and $k'= k$ therein), we obtain that $\age(\tau)\geq1$ if $k\geq 3$.  For $k=2$, we also get $\age(\tau)\geq1$ if there are at least two special unramified zero fibers by the same reason, or if the second special unramified fiber consists of poles (with possible residue constraints) by \autoref{neq:tauup}, \autoref{neq:tauotherp1} or \autoref{neq:tauhor1}. For $k = 2$ and only one special fiber, the map $\pi$ is a bielliptic cover and the corresponding stratum is
$\bP\omoduli[1](\{0,0\})$, where the two zeros $z_1$ and $z_2$ (of order zero) are exchanged by $\tau$. 
In this case $\tau$ induces a quasi-reflection, listed
as \textbf{Case~(1)}. 
\par 
The other case is $g = g' = 0$, and we leave it to be discussed later.
\par
\subsubsection*{\textbf{Case} $\mathbf{g'\not =0}$, $\mathbf{0<a<k}$}
Next consider the case $0 < a < k$.  If $a \leq k/2$, from the subspace of
invariant periods we obtain 
$\age (\tau) \geq \frac{1}{2} \dim H^1(X' \setminus P', Z'; \bbC)^{\frakR'}$ 
which is at least one except for $g' = 0$, and we again leave this exceptional
case to be discussed later. 
\par
If $a > k/2$, consider any conjugate pair of eigenvalues with exponents~$a_i$
and $k - a_i$ from the absolute
part, where $a_i \leq k/2$ (or the pair $(1,1)$ with $\zeta$-exponent~$0$).
Then after subtracting $a$ and normalizing to the range $[0, k-1]$, this pair
contributes at least $2/k$ to $\age(\tau)$.  Hence the $g$ pairs contribute
at least $2g/k$.  By Riemann--Hurwitz, $2g-2 \geq (2g'-2)k$, hence
$\age(\tau) \geq 2g/k \geq  1$ if $g' > 1$.
\par
It remains to discuss $g' = 1$. In this case $2g-2 = (k - s_1) + \cdots + (k - s_b)$
by Riemann--Hurwitz \autoref{neq:RHcyclic}, where recall that~$b$ is the number
of branch points of $\pi$, $s_i$ is the cardinality of each ramified
fiber, and every $s_i$ divides~$k$. If $b \geq 2$, then $2g - 2 \geq k/2 + k/2 = k$ and consequently $\age (\tau) \geq 2g/k > 1$.
\par
Suppose $b = 1$ and we set $s = s_1$. Then $2g-2 = k - s \geq k/2$
and consequently $g \geq 2$. If the unique ramified fiber is not a special zero fiber, then there must exist a special unramified zero fiber. The contribution~\autoref{eq:ageinterior} or \autoref{neq:tauuz} (with $k = k'$) from this unramified zero fiber gives $\age(\tau)
\geq 1$ for $k \geq 4$ in the case of $a > k / 2$. For $k \leq 3$, the
estimate $2g/k\geq 4/3 > 1$ justifies $\age(\tau) > 1$ in this case. 
\par
Now suppose that the ramified fiber is a special zero fiber. Applying~\autoref{eq:ageinterior} or~\autoref{neq:tauuz} to this fiber (with~$s = k'$), if $s \geq 3$, then it contributes at least $(s-2)/s \geq (s-2)/k$, and hence $\age (\tau) \geq (2g + s-2) / k = 1$ in this case. If $s \leq 2$, then $2g \geq k$ and from the absolute periods we obtain $\age (\tau) \geq 2g/k \geq 1$ in this case. 
\par
\subsubsection*{\textbf{Case} $\bf{g' = 0}$} Finally consider the
case $g' = 0$, i.e. $\pi$ is a cyclic cover of $\bP^1$ of degree~$k$
with~$b$ branch points. We discuss various cases according to the number of the branch points~$b$. For convenience we denote by $S_{\text{ram}}$ and $S_{\text{un}}$ the sets of ramified and unramified special fibers respectively.
\subsubsection*{\textbf{Case} $\bf{g' = 0,\ k\geq 3,\ b = 2}$}
Suppose $b = 2$. Then $X \cong \bP^1$ and $\pi$ is totally ramified at
two points. Since the stability of~$X$ in a stratum of genus zero implies
$|Z \cup P| \geq 3$, there is at least one special unramified fiber,
i.e.\  $|S_{\text{un}}|\geq 1$. If an unramified special fiber consists of zeros, the age contribution from~\autoref{eq:ageinterior} or~\autoref{neq:tauuz} implies that $\age(\tau) \geq 1$ for $k\geq 5$. 
Moreover, if an unramified special fiber consists of simple poles, then using~\autoref{neq:tauhor1} or~\autoref{neq:tauself1} we
obtain $\age(\bftau) \geq 1$ for $k \geq 3$. 
If an unramified special fiber consists of higher order poles, recall that for such
a fiber $d$ is the number of connected components of $\Gamma_{> -J}$ adjacent
to these $k$ poles as in \autoref{nlem:ageGRC}. If $d\geq 3$, then
by~\autoref{neq:tauotherp1} we obtain $\age(\bftau)\geq 1$. If $d\leq2$, by combining~\autoref{neq:tauup} and~\autoref{neq:tauotherp1} we obtain
$\age(\bftau)\geq 1$ for $k\geq 5$. The remaining cases are thus $k \leq 4$.
\par
\subparagraph{\textbf{Subcase} $\mathbf{k = 3}$} By the estimates in \autoref{lem:ageinterior} and \autoref{nlem:ageGRC}, any unramified special fiber in all situations (with or without residue constraints) contributes at least~$1/3$ to the age. Therefore, if $|S_{\text{un}}|
\geq 3$, we thus obtain that $\age (\bftau) \geq 1$. Moreover, the case $|S_{\text{ram}}|=0$ is impossible since the sum of entries of $\mu$ (which is $-2$) is not divisible by $3$. Hence we only need to consider the cases $|S_{\text{un}}| = 1$ or $2$ and $|S_{\text{ram}}| = 1$ or $2$. 
\par 
If $|S_{\text{un}}| = 2$ and
$|S_{\text{ram}}| = 1$ or $2$, then by the dimension reason there is  
an extra subspace of periods (besides the two unramified fiber contributions) on which~$\tau$ acts with eigenvalues~$1$ or $(1,1)$ or
$(\zeta, \zeta^2)$ since the total determinant of the $\tau$-action is one.
One checks that $\age(\tau) \geq 1$ in all these cases.
\par
Consider now the case $|S_{\text{un}}| =1$ and $|S_{\text{ram}}| =2$. Assume
there is a totally ramified zero. If the second totally ramified fiber is
a (GRC free) pole or a zero, then besides the eigenvalues $\zeta, \zeta^2$
from the unramified fiber we get an eigenvalue~$1$ from the residue or
a relative period, since the determinant of the unprojectivized action
of~$\tau$ is one. This gives $\age(\tau) \geq 1$. If the second ramified fiber
is a (GRC constrained) pole with $d=1$, we get $\omoduli[0]^\frakR(3m_1+m_2-2, -m_2,
\{-m_1,-m_1,-m_1\})$, with $m_i>0$ and $\frakR=\{r_3+r_4+r_5=0,\ r_2=0\}$,
which is \textbf{Case~(R1)}.  
\par
Suppose both ramified fibers are poles. Then the unramified special fiber
must consist of zeros, hence it gives signatures of type $(\{m, m, m\}, -m_1, -m_2)$
with $3m - m_1 - m_2 = -2$ and $m_1, m_2 > 0$, where the unramified zero fiber
contributes eigenvalues $\zeta$ and $\zeta^2$. If any of the two poles can have
a nonzero residue, then the residue gives an extra eigenvalue~$1$, which together
with $\zeta$ and $\zeta^2$ makes $\age(\tau) \geq 1$. The remaining case leads to
the GRC-constrained  stratum with $\frakR=\{r_4 = 0\}$ or $\frakR=\{r_4 = r_5 = 0\}$
listed as~\textbf{Case~(R2)}. 
\par
Finally if $|S_{\text{un}}| =1$ and $|S_{\text{ram}}| =1$, it gives
$\omoduli[0](\{m,m,m\},-3m-2)$ with~$m \in \bZ$, which is~\textbf{Case~(2)}.
\par
\subparagraph{\textbf{Subcase} $\mathbf{k = 4}$} Suppose first an unramified fiber
consists of poles that are subject to residue conditions. As noticed above, these
are vertical edges joined to $d$ higher connected components with $d > 1$. 
Then combining inequalities~\autoref{neq:tauup} and~\autoref{neq:tauotherp1}
gives $\age(\bftau) \geq 1$ in this case. 
\par 
Now we can assume that any unramified special fiber consists of either zeros or poles without extra residue constraints (i.e. $d=1$). Using \autoref{eq:ageinterior}, \autoref{neq:tauuz}, or \autoref{neq:tauup}, the relative periods or residue cycles from this fiber contribute at least $3/4$ to $\age(\tau)$. Hence we can further assume that $|S_{\text{un}}|= 1$. If $|S_{\text{ram}}|=2$, then $X$ has genus $g = 0$ with six zeros and poles, and hence the (unprojectivized) stratum of $X$ has dimension equal to four. Besides the eigenvalues $\zeta, \zeta^2, \zeta^3$
from the special unramified fiber, the remaining eigenvalue must be~$1$
or $\zeta^2 = -1$, since the determinant of the $\tau$-action
is $\pm 1$ (from an invertible integral matrix). In both cases one checks that
$\age(\tau) \geq 1$. Therefore, the remaining possibility is
$|S_{\text{un}}| = |S_{\text{ram}}|=1$, i.e. $\bP\omoduli[0](\{m, m, m, m\}, -4m-2)$
where the four zeros (or poles) are permuted by $\tau$ and the pole (or zero)
is fixed, with (unprojectivized) eigenvalues $\zeta, \zeta^2, \zeta^3$
and (projectivized) $\age(\tau) = 3/4 < 1$ if $a = 1$. This is \textbf{Case~(8)}. 
\subsubsection*{\textbf{Case} $\mathbf{g' = 0,\ k\geq 3,\ b \geq 3}$}
It is well-known (e.g.\ \cite[Proposition~2.3.1]{DeligneMostow}) that any
primitive $k$-th root of unity appears as an eigenvalue with multiplicity
$b-2 \geq 1$ for the action of $\tau$ on the $2g$-dimensional absolute
part $H^1(X; \bbC)$, and hence $g\geq 1$.
\par
\subparagraph{\textbf{Subcase} $\mathbf{k \not \in \{2,3,4,6\}}$}
For  $k \not \in \{2,3,4,6\}$, using \autoref{nle:eulerphi} and its notation
we find that $\age (\tau) \geq (b-2)\phi_k(a) / k \geq 1$.
Hence the remaining cases are cyclic covers of $\bP^1$ with
degree $k = 2, 3, 4, 6$, with $b \geq 3$ branch points, and $g \geq 1$. 
\par
\subparagraph{\textbf{Subcase} $\bf{k=3}$} In this case the cover is totally
ramified at $b = g+2$ points. If moreover $g \geq 3$, since $\phi_3(a) \geq 1$
for any $a$, we find $\age (\tau) \geq (b-2) / 3 \geq 1$. Hence we only need to
consider $g=2$ and $g=1$. 
\par
Consider first $g = 2$. Then the absolute periods already give a contribution
of at least~$2/3$ to $\age(\tau)$. If $|S_{\text{un}}|\geq 1$, then an unramified special fiber
contributes to the age at least $1/3$ by the estimates in \autoref{lem:ageinterior} and \autoref{nlem:ageGRC}, and hence $\age (\bftau) \geq 1$ in this case. We can thus assume that $|S_{\text{un}}|=0$. Moreover, if 
there are two ramified zero fibers, then the relative paths joining them give an eigenvalue $1$. If there is any 
ramified pole fiber whose residue is not constrained to be zero, then its residue gives an eigenvalue $1$. In 
both cases $\age (\tau) \geq 1$ since $1, \zeta, \zeta^2$ all appear as eigenvalues. The remaining cases 
are $|S_{\text{ram}}|=2, 3, 4$ with exactly one ramified zero fiber and all residues constrained to be zero, giving 
\textbf{Cases (3)}, {\bf (R3)} and {\bf (R4)}. 
\par
Now we deal with $g = 1$. The absolute periods give eigenvalues $\zeta, \zeta^2$,
and hence contribute at least $1/3$ to $\age (\tau)$. Note that any two special zero fibers contribute an eigenvalue~$1$
using their relative periods, which makes the total age at least one together with
the absolute periods. Similarly if there is a totally ramified pole fiber whose
residue is not constrained to be zero, then its residue gives an eigenvalue~$1$,
which again makes the total age at least one. We can thus assume that there is
exactly one special zero fiber and all ramified poles are constrained to have zero
residue. Moreover if $|S_{\text{un}}|\geq 2$, since each unramified fiber contributes
at least $1/3$ to the age  by the estimates in \autoref{lem:ageinterior} and
\autoref{nlem:ageGRC}, then altogether we obtain $\age (\bftau)\geq 1$. Hence we
only need to consider $|S_{\text{un}}| = 1$ or $0$.  
\par
Suppose $|S_{\text{un}}| = 1$, which contributes at least $1/3$ to the age as explained above. If this is an unramified pole fiber giving three edges that permute three connected components in upper level, then \autoref{neq:tauotherp1} contributes $1$ to the age, hence we can assume that the unramified fiber is either a zero fiber or a pole fiber not constrained by the GRC. The remaining cases give signatures  $(\{0,0,0\})$ or $(\{m, m, m \}, -3m)$
without residue conditions and those with residue conditions as follows: 
\begin{itemize}
\item $(\{m, m, m \}, -m_1, -m_2)^{\frakR}$ with $\frakR = \{r_4 =0\}$ or 
$\frakR = \{r_4 = r_5 = 0\}$,
\item $(\{m, m, m \}, -m_1, -m_2, -m_3)^{\frakR}$ with residue condition
$\frakR = \{r_4 = r_5 = 0\}$ or  $\frakR  = \{r_4 = r_5 = r_6 = 0\}$,
\item $(3m_1+m_2, \{-m_1, -m_1, -m_1 \}, -m_2)^{\frakR}$ with $\frakR = \{ r_5 = 0\}$, 
\item $(3m_1+m_2+m_3, \{-m_1, -m_1, -m_1 \}, -m_2, -m_3)^{\frakR}$ with
$\frakR = \{ r_5 = r_6 = 0\}$.
\end{itemize}
The first two cases without residue conditions correspond to {\bf Cases (6)} and {\bf (7)}. 
However in the cases with non-trivial~$\frakR$, the vertex $v$ already contributes
at least $2/3$ to $\age (\bftau)$. If the three edges in the unramified special fiber
join three vertices in $\Gamma$ that are permuted by $\bftau$, then
$\age (\bftau) \geq 1$ by the argument in the proof of \autoref{nprop:ageort} below.
If they join the same vertex $v'$, then their relative periods or residues
(without GRC) for $v'$ contribute at least $1/3$ to $\age (\bftau)$, thus making
$\age (\bftau) \geq 1$. Hence these cases do not appear in the tables.
\par
Suppose $|S_{\text{un}}| = 0$, i.e. there is no unramified special fiber. If there is
no pole, then the special fibers consist of a unique zero fiber, thus giving the
signature $\mu = (0)$ listed as \textbf{Case (4)}. 
If there are (ramified) poles, then each of them is constrained to have residue zero
as explained above, hence we obtain signatures  $(m, -m)$  and
$(m_1+m_2, -m_1, -m_2)^{\frakR}$ with $\frakR=\{r_2=0\}$ or $\{r_2=r_3=0\}$, which
correspond to \textbf{Cases (5)} and {\bf (R5)}.  
\par
\subparagraph{\textbf{Subcase} $\bf{k=4}$} We have $\phi_4(a) \geq 2$
for any $a$, and hence $\age (\tau) \geq 2(b-2) / 4 \geq 1$ for $b\geq 4$.
We thus need to consider the case $b = 3$. Since $2g+ 6 = \sum_{i=1}^3(4-s_i) \leq 9$ with~$s_i$ a proper divisor of $4$, we conclude that $g\leq 1$. Since there exist
eigenvalues $\zeta$ and $\zeta^3$  from the absolute homology, we only need to
consider the case $g = 1$.   
\par
In this case $\pi$ has two totally ramified fibers and the third ramified fiber consists of two simply ramified points. Since the absolute periods
already contribute eigenvalues $\zeta$ and $\zeta^3$ (thus at least $1/2$ to the age), if there exists an unramified special
fiber, it contributes by \autoref{nlem:ageGRC} enough to make $\age(\bftau)
\geq 1$. If there are two special zero fibers, then their relative periods give an eigenvalue $1$, which contributes enough to make $\age(\bftau)
\geq 1$. Similarly if there is a totally ramified pole without residue constraint, then its residue gives an eigenvalue $1$, which also makes 
$\age(\bftau) \geq 1$. Hence we can assume in the sequel that $|S_{\text{un}}|= 0$, there is a unique special zero fiber, and any totally ramified pole is constrained to have residue zero.   
\par
If the fiber with two simply ramified points is not special, we get the strata $\bP\omoduli[1](0)$ and $\bP\omoduli[1](m,-m)$ 
which are listed as \textbf{Cases~(9)} and \textbf{(10)}.
\par
If the fiber consisting of two simply ramified points is a special zero fiber, then it contributes at least an eigenvalue $\zeta^2$. As said, any remaining special fiber must be a ramified pole fiber with residue constraint. We thus get $(\{0,0\})$, $(\{m, m \}, -2m)$ and $(\{m, m \}, -m_1, -m_2)^{\frakR}$ with $\frakR = \{ r_3 = 0\}$ or $\{ r_3 = r_4 = 0\}$, listed as \textbf{Cases~(11)}, \textbf{(12)} (with two permuted zeros) and \textbf{(R6)}. 
\par
Now suppose the fiber with two simply ramified points is a special pole fiber. 
A GRC with the case $d=2$ for this fiber 
(as defined in \autoref{nlem:ageGRC}~(i)) gives enough to make $\age(\bftau) \geq 1$
by~\autoref{neq:tauotherp1}. The case $d=1$ leads to the strata 
$\bP\omoduli[1] (\{-m, -m\}, 2m)$ and $\bP\omoduli[1]^\frakR(2m_1+m_2, \{-m_1,-m_1\},-m_2)$ with $\frakR=\{r_4=0\}$ or $\frakR = \{r_2+r_3=0\}$,
listed as \textbf{Cases (12)} (with two permuted poles) and \textbf{(R7)}. 
\par
If the fiber with two simply ramified points consists of simple poles adjacent
to other vertices, then the residues at the two nodes give
eigenvalues~$\pm 1$. Together with the eigenvalues $\zeta, \zeta^3$ from the absolute periods
we obtain $\age(\tau) \geq 1$ for any~$\zeta^a$ used for projectivization. If the fiber with two simply ramified points has two simple poles that form a self-node, then $\omega$ is $\tau$-anti-invariant, i.e. $a=2$.
Projectivization of the eigenvalues of the absolute periods already
gives $\age(\tau) \geq 1$.
\par
\subparagraph{\textbf{Subcase} $\bf{k=6}$} In this case $\phi_6(a) \geq 4$
for any $a\neq 5$ and $\phi_6(5) = 2$, hence $\age (\tau) \geq 2(b-2) / 6
\geq 1$ for $b\geq 5$ and any $a$. Moreover,
$\age(\tau) \geq 1$ for $b = 4$ and any $a \neq 5$.
\par
Consider first $b = 4$ and $a = 5$. Then the eigenvalues of~$\tau$ from
the absolute homology of $X$ contain $\zeta, \zeta, \zeta^5, \zeta^5$,
contributing $2/3$ to $\age(\tau)$, and $g\geq 2$ in this case. If $g > 2$,
an additional conjugate pair of eigenvalues from the absolute part of
type $(\zeta^{a_i}, \zeta^{6-a_i})$ for $1 \leq a_i \leq 3$ or $(1,1)$,
after dividing by $\zeta^5$, can contribute at least $1/3$ to $\age(\tau)$,
which is enough. For $g = 2$, if there is a special zero or pole fiber with cardinality $2$ or $3$ or $6$,
using the respective estimates in \autoref{lem:ageinterior} and \autoref{nlem:ageGRC}, we can still obtain 
that $\age(\bftau) \geq 1$.
If all zeros and poles are totally ramified, then by Riemann--Hurwitz the only case is $\bP\omoduli[2](2)$ for $s_1 = 1$ and $s_2 = s_3 = s_4 = 3$ where the unique zero $z$ is totally ramified. But in that case $\tau^2$ induces a triple cover of $\bP^1$ totally ramified at the Weierstrass point~$z$, which is impossible because the linear system $|3z|$ has a base point at $z$. 
\par
Next consider $b = 3$. By Riemann--Hurwitz $g=1$ or $2$ in this case.
For $ b = 3$ and $g = 1$, we get $s_1 + s_2 + s_3 = 6$ with $1\leq s_i \leq 3$
dividing $6$. The only possibilities are $(2,2,2)$ or $(1,2,3)$. The former
is impossible for a connected cyclic cover of $\bP^1$ as the gcd of $s_1, s_2, s_3$
is not relatively prime to $6$. For the latter, note that the
eigenvalues of $\tau$ from the absolute homology of $X$ contain
$\zeta$ and $\zeta^5$. If the concerned stratum for $X$ has any extra dimension (from non-absolute periods or residues), since the  determinant of the $\tau$-action is $\pm 1$, the extra eigenvalue(s) together with $(\zeta, \zeta^5)$ will make $\age(\tau) \geq 1$. If there is no extra dimension, then there is a unique zero and all poles are constrained to have zero residue. Moreover if $d$ polar edges in a special fiber joining $d$ higher connected components get permuted for $d > 2$, then we obtain $\age{\bftau} \geq 1$ by~\autoref{neq:tauotherp1}, and for $d =2$ we obtain an extra eigenvalue $-1 = \zeta^3$ which combined with the absolute eigenvalues still makes $\age{\bftau} \geq 1$. Therefore, the only remaining case is $\bP\omoduli[1](0)$ where the marked point is totally ramified, which gives \textbf{Case~(13)}. 
\par
For $ b = 3$ and $g=2$, we get $s_1 + s_2 + s_3 = 4$ with $1\leq s_i \leq 3$
dividing~$6$. The only possibility is $(1,1,2)$. The eigenvalues of $\tau$
from the absolute homology of $X$ are $(\zeta, \zeta^5)$ together with
another conjugate pair.  Since the target is $\bP^1$, there is no invariant
absolute period, hence the other pair is either $(\zeta^2, \zeta^4)$
or $(\zeta^3, \zeta^3)$. Both cases give $\age (\tau) \geq 1$
after taking projectivization by $\zeta^a$ for any $a$.  
\end{proof}
\par
In the sequel we need to bound the number of vertices in a level graph that can be permuted by an automorphism of small age.  
\par
\begin{lemma} \label{lem:boundpermuted}
Suppose $\bftau$ is an automorphism of a multi-scale differential
of type~$\mu$ (possibly meromorphic) with $\age(\bftau)<1$. Then any $\bftau$-orbit of cyclicly permuted vertices in the associated level graph has cardinality at most two. Moreover if two vertices are swapped, then $\bftau^2$ acts trivially on their underlying stable curves. 
\end{lemma}
\par
\begin{proof}
Suppose there is an orbit of~$d>1$ cyclicly permuted
vertices (which have to be on the same level). We need to show that $d = 2$. By \autoref{lem:permute} $\age(\bftau)$ is at least one 
if $d \geq 3, M>0$ or $d = 2, M > 1$, where recall that $M$ is the dimension of the multi-vertex RC-independent subspace. Therefore, any permuted vertex must have genus zero with at most two zero edges (otherwise $M \geq 2$ from the relative periods of the zeros). Moreover if $d\geq 3$, then each of the permuted vertices can admit only one zero edge (otherwise $M\geq 1$).  
\par 
Suppose there exists an orbit of $d \geq 3$ permuted vertices. Among all of such orbits we choose the largest $d$, and further choose the highest level containing the orbit if there are multiple orbits of $d$ permuted vertices. We denote the chosen $d$ permuted vertices by $v_1, \ldots, v_d$. If they admit horizontal edges (in all situations of self-loops, or joining between them, or joining some other vertices on the same level), then there exist at least $d$ independent horizontal plumbing parameters that are cyclicly permuted, making $\age (\bftau) \geq 1$ 
by \autoref{lem:permute} (using $d\geq 3$ and $M\geq 1$), which contradicts the assumption. We can thus assume that all polar edges of the $v_i$ are vertical only. By the preceding paragraph, each $v_i$ admits only one zero edge and hence at least two vertical polar edges (by stability and since it does not have simple poles).  Take a polar edge $e_1$ of $v_1$ joining to a higher level and consider part of its $\bftau$-orbit $e_1, e_2, \ldots, e_d$ (where the next one $e_{d+1}$ does not have to be $e_1$ as it might be another polar edge of $v_1$).  Suppose the upper ends of these $d$ edges are adjacent to $d'$ vertices $v'_1, \ldots, v'_{d'}$ that are permuted by $\bftau$, where $d = d' k'$.  By assumption $d' < d$ (as we chose the $d$ permuted vertices to be both largest and highest) and hence each $v'_j$ admits at least $k' = d/d'$ zero edges joining to some of the $v_i$.  Using $M' \geq k'-1$ from the relative periods of these zero edges of $v'_j$, we obtain at least $(d'-1)(k'-1)/2$ for the age by \autoref{lem:permute}.  Since $k' = d/d' > 1$, if $d' \geq 3$, then $\age (\bftau)\geq 1$ which contradicts the assumption.  Hence the only possibilities are $d' = 1$, or $d' = 2$ with $k' = 2$.    
\par 
Consider first $d' = 1$, i.e. $k' = d$. Then $e_1, \ldots, e_d$ are the zero edges of a single higher-level vertex and they are permuted by $\bftau$. Then we can apply~\autoref{neq:tauuz} to conclude that $d \leq 4$.  The case $d=4$ is only possible when $e_1, \ldots, e_4$ are the zero edges of a single higher-level vertex in Case (8) of \autoref{ncap:noncan}, which contributes at least $3/4$ to the age.  Since each of the $v_i$ admits at least another polar edge, consider its orbit with $d' = 1$ or $d'=k'=2$ in the same notation.  The former contributes at least another $3/4$ and the latter contributes at least an extra $1/2$, both of which make the total age bigger than one, contradicting the assumption.  For the case $d = 3$, the estimate in~\autoref{neq:tauuz} gives at least $1/3$.  Since each of $v_1, v_2, v_3$ admits at least another polar edge, consider its orbit with $d' = 1$ in the same notation and we gain another $1/3$ (here $d=3$ is not divisible by $2$, so $d' = 2$ cannot occur). 
Moreover, $v_i$ cannot admit more than two polar edges (as we have already got $2/3$ for the age).  Then the residue cycles $r_i$ of the two polar edges (up to sign) on the $v_i$ can either freely vary or are constrained by $r_1 + r_2 + r_3 = 0$ only, hence the action of $\bftau$ on the subspace generated by $r_1, r_2, r_3$ contains eigenvalues $\zeta_3, \zeta_3^2$, contributing an extra $1/3$ and making the total age $\geq 1$, leading to a contradiction.  Alternatively, the three zero edges of $v_1, v_2, v_3$ go down to another three permuted vertices (otherwise the adjacent single vertex in lower level contributes at least $1/3$ by checking the cases of $k=3$ in \autoref{ncap:noncan} and \autoref{ncap:noncanRC} which makes the total age $\geq 1$). These three new permuted vertices can only admit one zero edge each. Hence we can go down along them again, until we find a single vertex with three permuted edges, which leads to the same contradiction as before. 
\par 
Consider the remaining possibility that $d' = k' = 2$, and hence $d = 4$.  In this case the age contribution from $v'_1, v'_2$ is already at least $(d'-1)(k'-1)/2 = 1/2$.  Since each $v_i$ admits at least another polar edge, consider its orbit and the associated upper adjacent vertices (again with $d' = k' = 2$ as the only possibility).  Then we obtain another contribution $1/2$, which altogether makes the total age $\geq 1$, leading to a contradiction.
\par 
We have thus proved that $d=2$, i.e., any $\bftau$-permuted vertices must appear in pairs. Take two mutually permuted vertices $v_1$ and $v_2$. Next we will show that $\bftau^2$ acts trivially on this pair, i.e., $\bftau^2$ fixes every edge of $v_1$ and $v_2$. Prove by contradiction. 
First suppose $\bftau^2$ does not fix all zero edges in the pair. Then we can find an orbit of zero edges 
$z_1 \mapsto z_2 \mapsto z'_1 \mapsto z'_2$ ($\mapsto z_1$, since we have seen that each permuted vertex can admit at most two zero edges), where $z_i$ and $z'_i$ are the zero edges of $v_i$ for $i = 1,2$. If these four zero edges belong to a single higher vertex, then it can only be Case (8) in \autoref{ncap:noncan} which contributes at least $3/4$ to the age. On the other hand, the permuted pair contributes at least $1/2$ from the relative periods of the zero edges, which altogether makes the total age $\geq 1$ and contradicts the assumption.  The remaining possibility is that $v_i$ joins a lower vertex $v'_i$ via $z_i$ and $z'_i$, and hence the two banana pairs $(v_1, v'_1)$ and $(v_2, v'_2)$ are swapped by $\bftau$.  But the upper relative periods and the lower residues each contribute at least $1/2$, in total making the age $\geq 1$ and thus leading to a contradiction. The same argument can be used to show that $\bftau^2$ fixes all vertical polar edges of $v_1$ and $v_2$. 
For any horizontal edge of $v_i$, it cannot be a self-loop at $v_i$ (otherwise its residue cycle and plumbing parameter can make the multi-vertex RC-independent dimension $M \geq 2$). Similarly we can rule out the case that two or more horizontal edges of $v_i$ are permuted by $\bftau$. Therefore, every horizontal edge of $v_1$ joins $v_2$, which is fixed by $\bftau$ (with the two ends swapped). Consequently $\bftau^2$ fixes the ends of every horizontal edge at both $v_1$ and $v_2$. In summary, we have shown that $\bftau^2$ fixes all zero and polar edges on each $v_i$ (whose number is at least three by stability since $v_i$ has genus zero). Hence $\bftau^2$ acts trivially on $v_1$ and $v_2$. 
\end{proof}
\par
\begin{rem}\label{rem:trivial-pair}
{\rm Suppose $v_1$ and $v_2$ are two vertices of genus zero such that ${\bftau}$ swaps them and ${\bftau}^2$ is the identity restricted to each of them. 
If their relative periods and residues contribute zero to $\age (\bftau)$, then the above proof implies that each $v_i$ admits a unique zero edge, and moreover, any two permuted (vertical) polar edges must be constrained by the GRC to have the sum of the residues equal to zero. We call such permuted $(v_1, v_2)$ of age zero a {\em trivial pair}. In particular if we view a trivial pair as a 'single' vertex, then it behaves the same as hyperelliptic vertices of age zero as described in \autoref{rem:hypeigen}. 
}
\end{rem}
\par
Finally we can show the following result about automorphisms of multi-scale differentials with small age.  Recall the coordinates 
$\bA_{\rel}(X)$, $\bC_{\hor}$ and $\bC_{\lev}$ introduced in \autoref{sec:Coordbd}, and $\bA = \bA_{\rel}(X) \times \bC_{\hor}\times \bC_{\lev}$. 
\par
\begin{prop} \label{nprop:ageort}
Let $\bftau$ be a lift to~$\bA$ of an automorphism of a multi-scale differential
of type~$\mu$ (possibly meromorphic). 
If $\age(\bftau)<1$ and $\bftau$
does not induce a trivial action or a quasi-reflection on $\bA$, then the action of $\bftau$ on the subspace $\bC_{\lev}$ of level parameters 
is non-trivial for all $\mu$ in $g \geq 2$ except for Case (3) in $g=2$ in \autoref{ncap:noncan}.  
\end{prop}
\par
\begin{proof} 
Suppose that $\bftau$ acts trivially
on $\bC_{\lev}$. We remark that in this case any two $\bftau$-fixed
vertices connected by vertical edges must have the same $\bftau$-restricted order.
To see it, let $X$ and $Y$ be two $\bftau$-fixed vertices joined by vertical
edges $e_1, \ldots, e_d$ that are permuted in one orbit (there could be more edges
between them in other orbits).  Suppose the $\bftau$-orders restricted to~$X$ and~$Y$
are $k_1$ and $k_2$ respectively, and $k_i = d \ell_i$ where $\ell_i$ is the
ramification multiplicity of the quotient map of $\bftau$ restricted to each end
of the edges. Each $e_j$ gives a relation $x_j  y_j$ equal to some product of the
level parameters (as in~\autoref{eq:zetacond2lev}), where $x_j$ and $y_j$ are local
standard coordinates at $e_j$. Note that $\bftau^d$ maps $x_1$ to $\zeta_{\ell_1} x_{1}$
and maps $y_1$ to $\zeta_{\ell_2}  y_{1}$ for some primitive $\ell_i$-th roots of unity.
Since by assumption $\bftau$ acts trivially on the level parameters, it implies that
$\zeta_{\ell_1}  \zeta_{\ell_2} = 1$, and hence $\ell_1 = \ell_2$. Since $k_i = d \ell_i$, 
we thus conclude that $k_1 = k_2$. 
\par
If $\bftau$ permutes any vertices, by \autoref{lem:boundpermuted} $\bftau^2$ fixes every vertex, 
acts trivially on all permuted pairs, and acts trivially on the level parameters (since $\bftau$ does). Applying the remark in the preceding paragraph to $\bftau^2$, it implies that $\bftau^2$ is a trivial action. Therefore, $\bftau$ has order two and it induces either a trivial action or a quasi-reflection (i.e., with age $0$ or $1/2$) as the only possibilities of having age smaller than one, which is ruled out by the assumption.  
\par 
From now on suppose that $\bftau$ has order at least three and that~$\bftau$ fixes every vertex. If there are at least three permuted horizontal
edges, then by~\autoref{neq:tauhor1} we get $\age (\bftau) \geq 1$. Moreover if a
horizontal edge~$e$ is fixed at a vertex~$v$, then the local standard form $dx/x$
is invariant at~$e$, hence $\omega$ restricted to~$v$ is an invariant form, which
cannot occur for any cases in \autoref{ncap:noncan} and \autoref{ncap:noncanRC},
except for possibly hyperelliptic vertices or vertices on which $\bftau$ acts
trivially. If two horizontal edges (or the two ends of a horizontal self-loop) are
permuted, looking at the cases, only hyperelliptic vertices are possible. In summary,
horizontal edges are only adjacent to vertices on which the restricted
$\bftau$-orders are one or two. Combining with the remark that vertical edges join
vertices that have the same restricted $\bftau$-order, we thus conclude
that $\bftau$ restricted to every vertex has order exactly $k$ for some $k\geq 3$.
\par 
From the Cases of $k\geq 3$ in \autoref{ncap:noncan}
and \autoref{ncap:noncanRC} we see that the age contribution is at least $1/3$ from
any $\bftau$-fixed vertex of order $k\geq 3$. Since $\age(\bftau) < 1$, there can be at most two vertices in the graph.  If there is a single vertex, the only case with $g\geq 2$ in \autoref{ncap:noncan} is Case (3) in genus two. Suppose there are exactly two vertices $v_1$ and $v_2$. Then the only possibilities are 
Cases (2), (4), (5), (R1), (R2), (R5) for $k=3$, and Case (13) for $k = 6$. The latter cannot appear twice as its signature does not have a pole.  Hence we can assume that $k=3$. Since all zeros and poles of the ambient stratum are labeled, any permuted zeros and poles in these Cases must be edges. Then Cases (4), (5), (R5) can appear in pairs and Cases (2), (R1), (R2) can appear in pairs (including possibly self-pairing). Due to the congruence conditions and residue conditions in these Cases, the only possibilities are (4) paired with (R5) and (2) paired with (2). In both cases the signatures of the ambient strata belong to Case (3).   
\end{proof}


\section{Singularities of the coarse moduli space}
\label{sec:canonical}

The purpose of this section is to control the singularities of the coarse
moduli space $\bP\MScoarse$ in order to show that the usual strategy for
proving general type---the canonical bundle is ample plus effective---can
be used after an appropriate modification due to non-canonical singularities
at the boundary as stated in \autoref{prop:GenTypeCrit}, which
we will prove at the end of this section.
\par
We start with the interior of the moduli space and the digression on the
logarithmic viewpoint:
\par
\begin{theorem} \label{thm:singofMScoarse}
For any signature~$\mu$, except 
of type $\mu = (m, 2-m)$
in genus $g=2$ for $1\neq m \equiv 1 \mod 3$ (Case (3) in \autoref{ncap:noncan}),
the interior of the coarse moduli space $\bP\MScoarse$ has canonical singularities.
\par
The pair $(\bP\MScoarse, \mathbf{D})$ consisting of the coarse moduli space
and the total boundary $\mathbf{D} = \partial \bP\MScoarse$ is a log canonical pair.
\par
Both statements hold as well for the strata with unlabeled zeros and poles,
after further ruling out Cases (6), (7), (8), (11), (12) in \autoref{ncap:noncan}. 
\end{theorem} 
\par
This situation is quite parallel to the moduli space of curves $\overline{M}_{g}$.
The singularities in the interior $M_{g}$ are also canonical, as shown by
\cite{HarrisMumford}. The pair $(\overline{M}_{g}, \partial M_{g})$ being
log canonical is a general fact for the coarse moduli space of a Deligne--Mumford stack
with a normal crossings boundary divisor (see e.g.,~\cite[Appendix A]{HHlogcan}). 
\par
The proof of \autoref{thm:singofMScoarse} is given in \autoref{sec:SingCurveAut}.
In contrast, there can be non-canonical singularities at the boundary, see
\autoref{rem:cherry} for an easy example of such a singularity induced by
ghost automorphisms. These are discussed in general in \autoref{sec:singMS}.
The effect of curve automorphisms and ghost automorphisms are combined in
\autoref{sec:Singbd}. There we define the compensation divisor~$D_{\NC}$ and
prove \autoref{prop:GenTypeCrit}.
\par

\subsection{Canonical sheaf and singularities on quotient stacks}
\label{subsec:singgeneralities}

We recall several well-known facts on singularities and the canonical
sheaf of an irreducible normal variety~$W$, with a focus on the case of coarse
moduli spaces of smooth Deligne--Mumford stacks. In particular, the spaces
we consider are $\bQ$-factorial, i.e., every Weil divisor is $\bQ$-Cartier. 
\par
On a singular variety there are three competing definitions of sheaves
of differential forms. First, $\Omega^1_W$ denotes the sheaf of 
K\"ahler differentials and $\Omega^p_W = \wedge^p \Omega^1_W$ 
its tensor power. It is badly behaved near the singular points and plays
hardly any role in the sequel. Second, let 
$i\colon U:= W_{\mathrm{reg}} \to W$ be the inclusion of the regular part
and let
\be
\Omega^{[p]}_W   \,:=\, i_*(\Omega^p_U) = (\Omega^{p}_{W})^{**}
\ee
be the sheaf of reflexive differential forms. Its top power
$\omega_W = \Omega^{[p]}_W$ for $p = \dim(W)$ is a line bundle, called 
Grothendieck's dualizing sheaf (even though~$W$ is not 
Cohen--Macaulay in general, so Serre duality does not hold with~$\omega_W$).
Third, let $ \pi\colon \widetilde{W} \to W$ be a resolution of singularities
and we define
\be
\widetilde{\Omega}^p_W \= \pi_* \Omega^p_{\widetilde{W}}\,
\ee
which is useful for computing global sections on~$\widetilde{W}$. 
Finally we define the canonical divisor $K_W = \pi_* K_{\widetilde{W}}$ using the
pushforward of cycles. Note that $\omega_W = \cO(K_W)$.
\par
We now discuss the logarithmic situation, or more generally the case of pairs.
Let $W$ still be a normal variety and $\Delta = \sum a_i D_i$ a sum of prime
divisors with $a_i \in \bQ$. Choose the smooth resolution~$\pi$
so that moreover the preimage $\widetilde{{D}} = \pi^{-1} {D}$
is a normal crossings divisor. Now let $i\colon U \to W$ be the inclusion of
the open subset where~$W$ is smooth and ${D}$ is normal crossings.
As above we define the sheaf of reflexive logarithmic differentials and
the pushforward of logarithmic differentials on the resolution to be
\be
\Omega^{[p]}_W(\log {D})  \,:=\, i_*(\Omega^p_U(\log {D}|_U))
\quad \text{and} \quad \widetilde{\Omega}^p_W(\log {D})
\= \pi_* \Omega^p_{\widetilde{W}}(\log \widetilde{{D}})\,.
\ee
The logarithmic canonical divisor is defined to be $K_W + {D}$.
This is consistent with the above notation since $\cO(K_W + {D})
= \Omega^{[p]}_W(\log {D})$ for $p = \dim(W)$.
\par
Next we review some types of singularities of pairs. Recall that the
discrepancy $\discrep(W,\Delta)$ is the infimum over all exceptional
divisors~$E$ in all birational morphisms $\widetilde{W} \to W$ of the
coefficient $a(E,\Delta,\widetilde{W})$ of~$E$ in the pullback of $K_W
+ \Delta$. The pair $(W,\Delta)$ has \emph{canonical singularities} if
$\discrep(W,\Delta) \geq 0$ and \emph{logarithmic canonical singularities}
if $\discrep(W,\Delta) \geq -1$. More details can be found in
\cite[Section~2.3]{KMbook}. 
\par
In particular if $W = (W,\emptyset)$ has canonical singularities, then
sections of the canonical bundle restricted to the regular locus of~$W$ can
be extended across the singularities (\cite{ReidYoung}).
\par

\subsection{Singularities from curve automorphisms} \label{sec:SingCurveAut}

We apply this discussion to $W = \bP\MScoarse$, keep the notation $D_\Gamma$ for
the boundary divisors of the stack and write $\mathbf{D}_\Gamma$ etc\ for the
boundary divisors of the coarse moduli space. Since  $\bP\MScoarse$ has only
finite quotient singularities, every subvariety of codimension~one is $\bQ$-Cartier.
Thus the rational Picard groups of $\bP\MScoarse$ and $\bP\LMS$ are identical and
computations are performed mostly in terms of the classes of $D_\Gamma$ in the sequel.
\par
\begin{proof}[Proof of \autoref{thm:singofMScoarse}]
For the first statement we examine the cases in \autoref{ncap:noncan}.
Except Case (3), all the strata in that table either involve
permuted marked points (not allowed for labeled strata), or are one-dimensional
with smooth quotients.
\par
To prove the second statement we use that $(\bP\MScoarse, \mathbf{D})$ is the coarse moduli space
associated with a smooth Deligne--Mumford stack with normal crossings boundary divisor
(\cite{LMS}) and that $\mathbf{D}$ has each boundary term with coefficient one. 
In fact, \cite[Proposition~A.13]{HHlogcan} shows that in this situation
there is some boundary divisor~$\Delta$ such that the pullback of
$m(K_{\bP\MScoarse} + \Delta)$ equals $m(K_{\bP\LMS} + D)$ and furthermore
that $(\bP\MScoarse, \Delta)$ is log canonical. Since all boundary divisors in~$D$
appear with coefficient one, the same holds for $\Delta$, i.e., it implies
that $\Delta = \mathbf{D}$. This reflects the fact that log canonical divisors are 
insensitive to branching (see \cite[Proposition~20.2]{FA}  for a general
comparison formula). 
\par 
Finally for unlabeled strata we only need to rule out those of (projectivized) dimension at least two in \autoref{ncap:noncan}. 
\end{proof}
\par
\begin{rem} \label{rem:cherry}
  {\rm Consider the left slanted cherry $\Delta = \Delta_\ell$ from
\autoref{ex:cherry}, i.e., with the twist $p = 2$ on the short edge. The stack
structure is given by the group $K_\Delta = \bZ/3\bZ$. As we verified there, the
generator $(2,1)$ of this group acts by $(\zeta_6^2, \zeta_3^1)
= (\zeta_3, \zeta_3)$ on the coordinates corresponding to opening up the
levels. (See also Equations~(6.7) and~(12.6) in \cite{LMS} for the construction.
The local parameter called $t_i$ in loc.\ cit.\  raised to the lcm of $p_e$
for all edges~$e$ crossing the level passage gives the coordinate called $s_i$
that rescales the differential.) This group action does not satisfy the
Reid--Tai criterion and \cite[Example~1.8 (2)]{ReidYoung} shows explicitly why
the corresponding singularity is not canonical. We will elaborate on
this in the next subsection.
}\end{rem}
\par
\begin{proof}[Proof of \autoref{thm:introsingofMScoarse}]
Besides \autoref{thm:singofMScoarse} we need to show that non-canonical
singularities occur in the boundary of all strata with possible exceptions in
low genus only. A more elaborate version
of \autoref{rem:cherry} can be used to show that as long as
$g \geq 5$. Consider a triangle graph with one vertex at each level,
the top vertex in the stratum $\omoduli[2,2](0,2)$, the middle
vertex in the stratum $\omoduli[2,2](-2,4)$ and the bottom level vertex
with the remaining genus and all the marked points. The prongs are~$3$
on the long edge from the top to bottom level, and $1$ and $5$ on the short
edges. The element $(\zeta_{3}^1, \zeta_{15}^5)$ fixes all prongs, and it 
defines a ghost automorphism with $\age < 1$.
\end{proof}

\subsection{Singularities induced by ghost automorphisms}
\label{sec:singMS}

The singularities induced at the boundary of $\bP\LMS$, say at
an enhanced level graph~$\Gamma \in \LG_L$, stem from the action
of the ghost automorphisms $K_\Gamma = \Tw/\sTw$. These are toric
singularities. We explain
here how to fit the data of the graph and the enhancements into
the standard framework of toric geometry. The goal is to give a
formula for the discrepancies and, more generally, a formula for
the pullback of torus-invariant divisors in terms of these graph
data. 
\par
\medskip
We start by recalling some well-known toric terminologies. 
An affine toric variety is given by a $\bZ$-module $N$ that we
view as a lattice inside $N_{\bR} = N \otimes_{\bZ} \bR$ and a
convex rational polyhedral cone~$\sigma \subset N_{\bR}$. We let $M = N^\vee
= \Hom(N,\bZ)$ and view the dual cone~$\sigma^\vee$ as a subset of~$M_{\bR}$.
Then the group algebra $\bC[M \cap \sigma^\vee]$ is a finitely generated
algebra and the associated (affine) toric variety is defined as 
\bes X_{N,\sigma} \= \mathrm{Spec}(\bC[M \cap \sigma^\vee])\,. \ees
The spanning rays of~$\sigma$ generated by the primitive elements
$r_1,\ldots,r_L$  are in bijection to
the torus-invariant divisors $D_1,\ldots,D_L$ of~$X_{N,\sigma}$.
We omit~$\sigma$ from the notation, if it is the positive cone
for some implicitly chosen basis of~$N_{\bR}$.
\par
The affine toric variety~$X_{N}$ is non-singular if $\sigma\cap N$ is generated
by a subset of a basis of~$N$. If this is not the case, we can resolve
the singularities by subdividing~$\sigma$ as a union of subcones such that 
each of the subcones satisfies the above condition. Let~$F$ be the fan obtained
from the cone subdivision.
The additional rays of the subcones are given by the primitive
interior points in $\sigma^{\circ} \cap N$, which we list as $v_1,\ldots,v_s$.
Each of the rays~$v_i$ corresponds to a torus-invariant exceptional
divisor~$E_i$ in the resolution $\pi\colon \widetilde{X}_{F} \to X_N$.
\par
We state the next proposition in the case of interest to us, namely
that~$X_{N}$ has only abelian finite quotient singularities, which is equivalent
to $\sigma$ being a simplicial cone by \cite[Theorem~11.4.8]{CLS}, i.e.,
$L = \dim N_{\bR}$. Consequently there are elements $m_{\sigma,i} \in M_\bR$
such that $\langle m_{\sigma,i}, r_j \rangle = \delta_{ij}$.
Let $m_\sigma = \sum_{i=1}^L m_{\sigma,i}$. We denote the non-exceptional
torus-invariant divisors, the strict transforms
of the~$D_i$ by $\widetilde{D}_i$. 
\par
\begin{prop} \label{prop:toricpullback}
The canonical divisor $K_{X_N} = - \sum_{i=1}^L D_i$
is the negative sum of the torus-invariant divisors. Moreover if $X_{N}$ has
only abelian quotient singularities, then  the discrepancy
of~$E_i$ is given by
\bes
K_{\widetilde{X}_F} - \pi^* K_{X_N} \= \sum_{j=1}^s (\langle m_\sigma, v_j
\rangle -1) E_j\,
\ees
and more generally 
\bes
\pi^* D_i - \widetilde{D}_i \= \sum_{j=1}^s \langle m_{\sigma,i}, v_j
\rangle  E_j\,
\quad \text{for all}\quad  i=1,\ldots,L\,.
\ees
\end{prop}
\par
\begin{proof}
The claims follow from combining~\cite{CLS} Theorem~8.2.3 on the description
of the canonical bundle, Lemma 11.4.10 for its pullback to resolutions,
Theorem~4.2.8 for the conversion of divisors into support functions
and Proposition~6.2.7 for the pullback of divisors written in these
terms.
\end{proof}
\par
In general, not all pluricanonical forms extend from $X_N$ to its
resolution $\widetilde{X}_{F}$, since they can acquire poles along the exceptional
divisors $E_i$. However, we can consider only a subset of pluricanonical forms
having high enough order of vanishing along the divisors $D_i$. We now show
a criterion about how high the order of vanishing of pluricanonical forms
along~$D_i$ has to be in order to ensure that they extend to the resolution. 
\par
\begin{prop} \label{prop:extendtocan}
	Let $(b_1,\ldots,b_L) \in \bN^L$ be a tuple such that 
	\begin{equation}\label{eq:conditionextension}
	 \sum_{i=1}^L (b_i + 1)
	\langle m_{\sigma,i}, v_j \rangle\geq 1,\quad \text{for } j=1,\dots, s.
	\end{equation}
Then for all $a \in \bN$ we have the inequality
\[	h^0(\widetilde{X}_F, aK_{\widetilde{X}_F})
	\,\geq \, h^0\Bigl(X_N, a\Big(K_{X_N} - \sum_{i=1}^L b_i D_i\Big)\Bigr)\,.
	\]
\end{prop}
\par
\begin{proof} It suffices to show that the exceptional
divisors~$E_j$ occur with non-negative coefficients in the difference
\bas &\phantom{\=}
K_{\widetilde{X}_F} \,-\, \pi^*\Bigl(K_{X_N} - \sum_{i=1}^L b_i Di\Bigr) \\
&\= \sum_{i=1}^L b_i \widetilde{D}_i
\+\sum_{j=1}^s (\langle m_\sigma, v_j \rangle -1) E_j
\+ \sum_{j=1}^s E_j \cdot \sum_{i=1}^L b_i \langle m_{\sigma,i}, v_j \rangle \\
&\= \sum_{i=1}^L b_i \widetilde{D}_i
\+\sum_{j=1}^s E_j \cdot \Bigl(-1  + \sum_{i=1}^L (b_i + 1)
\langle m_{\sigma,i}, v_j \rangle \Bigr)\,,
\eas
which is ensured by the standing assumption.
\end{proof}

Recall that an inclusion of lattices $N' \hookrightarrow N$ with quotient
group~$G$, and the same cone~$\sigma$ in both lattices, gives rise to
a quotient map $X_{N'} \mapsto X_{N'}/G = X_{N}$, see e.g.\
\cite[Sections~1.5 and~3.3]{CLS}. 
\par
We focus now on toric varieties obtained by the the quotient of affine space
via a cyclic group of order~$n$. We say that a singularity is of
type $\tfrac1n(a_1,\ldots,a_L)$ if it is the quotient of~$\bC^L$ by a
cyclic group $G=\langle\tau\rangle$  of order~$n$ acting by $\zeta_n^{a_i}$
on the $L$~coordinates.  Consider then the case $X_{N'}=\bC^L$, so~$N'$ is a
lattice in $\bR^L$  generated by vectors~$e_i$ and $\sigma$ is the standard
cone generated by the basis $(e_i)$ of $N'$. If we define $N$ to be the
lattice generated by the basis of $N'$ and by $v_\tau=\sum_{i=1}^L a_i/n\cdot e_i$,
then $N/N'=G$ and $X_{N}=\bC^L/G$. 
\par
We specialize further to the toric varieties~$X_{N}$. Since the cone~$\sigma$
is generated by the coordinate vectors $e_i=v_i$, using the notation previously
introduced, the $m_{\sigma,i}$ are simply the dual vector $e_i^{\vee}$. The only
primitive interior points in $\sigma^{\circ} \cap N$ are the vectors $v_{\tau^j}$,
for $j=1,\dots,n$. Hence, in this setting, \autoref{prop:extendtocan}
specializes to the following statement.
\par
\begin{cor}\label{cor:cyclicsing}	
  Let $X=\bC^d/\langle\tau\rangle$, where $\tau$ acts by multiplication
  of $\zeta_n^{a_i}$ on the $i$-th coordinate. If 
\begin{equation}\label{eq:conditionextensioncyclic}
	\sum_{i=1}^L \frac{a_i}{n}(b_i+1)\geq 1
\end{equation}
then the inequality 
\[	h^0(\widetilde{X}, aK_{\widetilde{X}})
\,\geq \, h^0\Bigl(X, a\Big(K_{X} - \sum_{i=1}^d b_i D_i\Big)\Bigr)\,
\]
holds for all $a \in \bN$, where $D_i$ are the image in $X$ of the
divisors $\{x_i=0\}\subset \bC^d$ and where $\widetilde{X}$ is a smooth
resolution of $X$.
\end{cor} 
\par
Note that if the age of $\tau$ is indeed greater or equal to one, then we
take $b_i=0$, and so the singularities of $X_{N}$ are canonical.
\par
\begin{example} \label{ex:cyclictoric}{\rm
The resolution $\pi\colon \widetilde{X} \to X$ of a singularity of
type $\tfrac13(1,1)$ has a single exceptional divisor~$E$. In this case we have
\be
K_{\widetilde{X}} - \pi^* K_{X} \= -\frac13 E
\quad \text{and} \quad \pi^*D_i - \widetilde{D}_i \= \frac13 E
\ee
where $D_i$ for $i=1,2$ are the two coordinate axes. If we consider for example $b_1=1$ and $b_2=0$, then by \autoref{eq:conditionextensioncyclic} we have that all sections of $K_{X}-D_1$ extend to the resolution $\widetilde{X}$.
\par
The resolution $\pi\colon \widetilde{X} \to X$ of a singularity of
type $\tfrac14(1,2)$ also has a single exceptional divisor~$E$. In this case we have
\be
K_{\widetilde{X}} - \pi^* K_{X} \= -\frac14 E
\quad \text{and} \quad \pi^*D_1 - \widetilde{D}_1 \= \frac14 E,
\quad \pi^*D_2 - \widetilde{D}_2 \= \frac12 E\,.
\ee
In this case, if we set for example $b_1=1$ and $b_2=0$, by \autoref{eq:conditionextensioncyclic} we have that all sections of $K_{X}-D_1$ extend to the resolution $\widetilde{X}$.
		This is a special case of the resolutions in \autoref{ex:cherryresolution}.
}\end{example}
\par
\medskip
Even though not needed in the sequel, it is instructive to describe in standard
toric geometry language the structure near a boundary component of the orderly
blow-up  $\bP\MSgrp$, where only ghost automorphisms are quotiented out. Recall
indeed that the map $\varphi_1\colon \bP\LMS \to \bP\MSgrp$ to the orderly blowup
is locally given by the map $[U/K_\Gamma] \to U/K_\Gamma$, where $U$ is a
neighborhood of a generic point of a boundary component $D_\Gamma$
and $K_\Gamma= \text{Tw}_\Gamma/\text{Tw}^s_\Gamma$ is the group of ghost automorphisms.
\par
Recall from \autoref{sec:Coordbd} the coordinate system
near the boundary. We analyze the toric geometry of the part
$\bC^{\lev}$ of this coordinate system. This is the affine toric
variety with the (dual) lattice
\bas
M' &\= \, \Bigl\langle \frac1{p_e} \cdot w_i, \quad
e \in E(\Gamma),\, \text{$i$-th level passage crossed by~$e$}
\Bigr\rangle_{\bZ}  \\ 
& \= \, \Bigl\langle \frac1{\ell_i} \cdot w_i, \quad i = 1, \ldots, L \Bigr\rangle_{\bZ}
\eas
and $\sigma^{\vee}$ the positive (dual) cone generated by the~$w_i$ in
$M'_{\bR} \cong \bR^L$, where $w_i$ is the $i$-th unit vector and $\ell_i$ is the $\lcm$ of
all enhancements~$p_e$ crossing the
$i$-th level passage.
Moreover,   
\be \label{eq:NwithsTw}
N' \= (M')^\vee \= \sTw[\Gamma] \= \langle \ell_i \cdot e_i,\, i=1,\ldots, L \rangle_{\bZ}\,
\ee
is the simple twist group by definition, where $e_i$ is the dual vector of $w_i$. 
\par
In order to define the twist group similarly, recall from
\cite[Section~5]{LMS} that it depends only on the level passages
crossed by the edges and the enhancements, not on the vertices the
edges are attached to. For $0 \leq i<j \leq L$ we let
\[w_{-j}^{-i} \= (0,\ldots,0,1,\ldots,1,0,\ldots,0) \in \bR^L\]
be the vector where the string of ones goes from~$i+1$ to~$j$. Then
\[
N \,:= \Tw[\Gamma] \= (M')^\vee \quad \text{and} \quad
M \= \Bigl \langle
\frac{1}{p_e} \cdot w^{e^+}_{e^-}, \quad e \in E(\Gamma)
\Bigr \rangle_{\bZ}\,
\]
where $e^\pm$ are the upper and lower ends of~$e$. Note that the
explicit computation of a basis of $M$ is in general not possible
in closed form, which requires working with gcds, i.e., computing
a Smith normal form. 

\begin{example}\label{ex:cherryresolution} {\rm
We continue with the running example of the cherry graph, now generalizing
to enhancements $a$ on the short edge and~$b$ on the long edge. We
let $\ell_1 = \lcm(a,b)$ and $\ell_2=b$. Then
\bas
N' &\= &\sTw[\Gamma]  &\=& \Bigl\langle \bigl(\frac1a, 0\bigr),
\bigl(\frac1b, 0\bigr), \bigl(0, \frac1b \bigr)   \Bigr\rangle^\vee
&\=& \Bigl\langle \bigl(\ell_1, 0\bigr),
\bigl(0,\ell_2 \bigr)  \Bigr\rangle \phantom{\,.} \\
N &\= &\Tw[\Gamma]  &\=& \Bigl\langle \bigl(\frac1a, 0\bigr),
\bigl(\frac1b, \frac1b \bigr)  \Bigr\rangle^\vee
&\=&
\Bigl\langle \bigl(a, -a \bigr), \bigl(0,\ell_2 \bigr)  \Bigr\rangle\,
\eas
so that $n:= [N:N'] = \ell_1/a = b/\gcd(a, b)$.
\par
We restrict moreover to $b \geq a$. We see
that this is a cyclic quotient singularity of order $n =\ell_1/a$ and
of type $\frac1n (1, q)$ where $q = \tfrac{b-a}{\gcd(a, b)}$. This generalizes
\autoref{rem:cherry}. To resolve this singularity minimally we have
to insert the rays generated by boundary points of the lower convex hull
of $N$ in the positive quadrant~$\sigma = \sigma'$. These are the
rays generated by 
\[v_j=j\cdot (\tfrac1n, \tfrac{q}n)\in \sigma^{\circ}\cap N\] for $j=1,\ldots,q'$ with $0< q' < n$ and $qq' \equiv 1 \mod n$, see
\cite[Section~II.6]{vanderGeer} or \cite[Section~10]{CLS} for another version
of this resolution ('Hirzebruch-Jung continued fraction').
}\end{example}
\par

\subsection{Singularities at the boundary}
\label{sec:Singbd}

Here we combine the previous two subsections to analyze the singularities
at the boundary with the goal of proving \autoref{prop:GenTypeCrit}.
We start with the definition of the non-canonical compensation divisor $D_{\NC}$.
\par
First we distinguish several special edge types in a two-level graph (for vertical edges only).  
If an edge corresponds to a separating node, we say that it is of {\em compact type (CPT)}. Otherwise we say that it is of {\em non-compact type (NCT)}. 
If the lower part of the graph separated by a CPT edge consists of a single rational vertex, the edge type is called a {\em rational bottom tail (RBT)}. 
Recall that a {\em (vertical) dumbbell (VDB)} graph 
is defined to be a graph of compact type with a unique (separating)
edge which is vertical. If the graph contains a unique (vertical) edge (i.e., a VDB graph) and if one end of the edge is of genus one, we say that it is an {\em elliptic dumbbell  (EDB)}. 
An edge of compact type which is neither RBT nor EDB is called {\em other compact type (OCT)}.
\par
Let $E_\Gamma$ be the number of (vertical) edges of a two-level graph $\Gamma$. We then define (for $g\geq 2$) that 
\ba \label{eq:non-canonical-term}
D_{\mathrm{NC}} & \,:= \,\sum_{\Gamma \in \LG_1} b_{\NC}^\Gamma [D_\Gamma] \,:=\, \sum_{\Gamma \in \LG_1}
(\ell_\Gamma R_{\NC}^\Gamma -1) [D_\Gamma]\, \quad \text{where}  \\
R_{\NC}^\Gamma &=\, \sum_{\rm NCT} \frac{1}{2}\frac{1}{p_e} + 
\sum_{\rm RBT} \frac{1}{p_{e}}  + \sum_{\rm OCT}
\frac{2}{ p_{e}} +\sum_{\rm EDB} \frac{4}{p_{e}} \,.
\ea
In the above each sum runs over the edges of the corresponding type. Note that the edge type EDB is exclusive, i.e., if it appears in a two-level graph, then the graph has a unique edge and all other edge types do not appear. In particular, an edge cannot be both RBT and EDB due to the assumption that $g\geq 2$. 
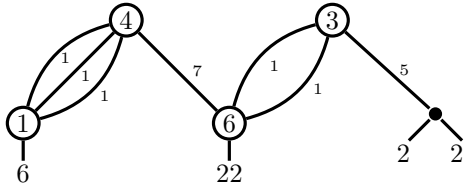
\begin{figure}[ht]
	\begin{tikzpicture}[
		baseline={([yshift=0.5ex]current bounding box.center)},
		scale=2,very thick,
		bend angle=30]
		\node[circled number] (T) [] {$4$};
		\node[] (T-1) [right=of T] {};
		\node[circled number] (T-2) [right=of T-1] {$3$};
		\node[] (T-3) [right=of T-2] {};
		\node[] (T-4) [left=of T] {};
		\node[circled number, below=of T-1] (B) {$6$};
		\node[comp, below=of T-3] (C) {};
		\node[circled number, below=of T-4] (D) {$1$};
		\path (T) -- (T-1) node[ midway, below=of T,xshift=-10pt,yshift=-32pt] (BL) {};
		\path (T-1) -- (T-2) node[ midway, below=of T,xshift=10pt,yshift=-32pt] (BR) {};
		\path (T-1) -- (B) node[] (B-1) {};
		\path (T-2) -- (B) node[] (B-2) {};
		\path (T-2) edge 
		node [order top left,xshift=3pt,yshift=14pt] {} 
		node [order bottom right,xshift=15pt,yshift=-13pt] {$5$} (C);
	\path (T) edge 
	node [order top left,xshift=3pt,yshift=14pt] {} 
	node [order bottom right,xshift=15pt,yshift=-13pt] {$7$} (B);
		\path (B) edge [bend left] 
		node [order bottom right, xshift=-3,yshift=0] {} 
		node [order top left, xshift=-5,yshift=-14] {$1$} (T-2)
		edge [bend right] 
		node [order bottom right, xshift=-5,yshift=-4] {} 
		node [order top right, xshift=-8,yshift=-15] {$1$} (T-2);
		\path (T) edge [bend left] 
		node [order bottom left, xshift=0,yshift=-20] {$1$} 
		node [order top left, xshift=-1,yshift=-13] {} (D)
		edge [bend right] 
		node [order bottom right, xshift=-20,yshift=-12] {$1$} 
		node [order top right, xshift=-6,yshift=4] {} (D);
				\path (T) edge 
		node [order top left,xshift=3pt,yshift=14pt] {} 
		node [order bottom right,xshift=-15pt,yshift=-15pt] {$1$} (D);
		\draw (B) -- +(270: .25) node [xshift=0,yshift=-5] {$22$};
		\draw (C) -- +(225: .25) node [xshift=-2,yshift=-5] {$2$};
		\draw (C) -- +(-45: .25) node [xshift=-2,yshift=-5] {$2$};
		\draw (D) -- +(270: .25) node [xshift=0,yshift=-5] {$6$};
	\end{tikzpicture}
	\caption{In this graph there are five NCT edges with prong $1$, one OCT edge with prong $7$, and one RBT edge with prong $5$. }
	\label{figure:Rgamma}
\end{figure}	
\begin{rem}{\rm 
For certain range of genera and signatures, one can alter the definition of $R_{\NC}^\Gamma$ in \autoref{eq:non-canonical-term} slightly. Indeed for certifying general type for the minimal strata in low genera, we need another version of $R_{\NC}^\Gamma$ (see \autoref{prop:compensationv2} and \autoref{prop:nc-refinement}). 
}
\end{rem}
\par 
We are now ready to present the main proof of this section.
\par
\begin{proof}[Proof of \autoref{prop:GenTypeCrit}]
The content of the proposition is that global sections of
$a(K_{\bP\MScoarse} - D_{\mathrm{NC}})$ extend to $a$-canonical sections on
a smooth resolution~$\tMScoarse$ of~$\bP\MScoarse$, for any $a\in \bN$. We revisit the
argument~\cite[Proposition~3.1]{TaiKodAg} for this purpose.
\par
Suppose such a section~$\eta$ does not extend to~$\tMScoarse$ and
suppose this happens near some multi-scale differential
$(X,\bfomega,\bfz, \bfsigma)$ compatible with some level graph~$\Pi$, in fact
necessarily in the boundary by \autoref{thm:singofMScoarse}, i.e.\ $\Pi$ is
non-trivial The
section~$\eta$ thus acquires poles near a divisor~$E$ of $\tMScoarse$. Using the
same notation as in \autoref{sec:Coordbd}, we consider the local covering~$\bA$
of $\bP\MScoarse$, where  $\bA =  \bC_{\lev} \times \bA_\rel(X) \times  \bC_{\hor} $,
the tangent space to the orbifold chart near $(X,\bfomega,\bfz,\bfsigma)$.  
Consider finally the normalization of~$\tMScoarse$ in the function field
of $\bA$, i.e., the normalization of the corresponding
fiber product. For each component~$E'$ of the preimage
of~$E$ in this normalization, the stabilizer (in the full Deck group of the
cover, the extension $ \mathrm{Iso}(X,\bfomega)$ of~$\Aut(X,\bfomega)$
by~$K_\Pi$ as in \autoref{rem:notsemidirect}) is cyclic,
say generated by an element~$\bftau$, a lift of an automorphism in
$\Aut(X,\bfomega)$ to an automorphism acting on~$\bA$ as considered in
\autoref{sec:newage}. Consequently the pullback of~$\eta$ to
$\bA/\langle \bftau \rangle$ does not extend to its smooth resolution.
We will show that this does not happen for sections under consideration, i.e.,
with enough vanishing along some boundary divisors, using \autoref{cor:cyclicsing}.
\par
If $\bftau$ is a quasi-reflection on $\bA$, the quotient is
smooth and the extension of $a$-canonical sections is automatic. In the case
of $\age(\bftau) \geq 1$, the singularities
of $\bA/ \langle \bftau \rangle$  are canonical and all $a$-canonical sections
of $K_{\bP\MScoarse}$ extend to $\bA/ \langle \bftau \rangle$ by the original argument
of Tai's criterion. 
\par
In the remaining cases, by \autoref{nprop:ageort} and as our genus and signature 
hypothesis, $\bftau$ acts non-trivially on at least one of the level
coordinates~$t_i$. We may thus assume that $\bftau$ acts by $\exp(2\pi i a_j/k)$
in an appropriate basis of $\bA_\rel(X) \times  \bC_{\hor}$ and by
$\exp(2\pi i c_j/k_j)$ on the coordinates~$t_j$, where at least one of the
entries~$c_i$ of~$n_\bftau$ is non-zero by our assumption.  Using
\autoref{cor:cyclicsing} with $b_i=b_{\NC}^{\delta_i(\Gamma)}$ for $i=1,\dots,L$
(where $\delta_j(\Gamma) \in \LG_1$ is the undegeneration of a level graph $\Gamma$
obtained by compressing all level passages but the $j$-th one) and with $b_i=0$
for $i=L+1,\dots, d$, we need to show that 
\begin{equation}\label{eq:condition}
c(\bftau)\,:=\,\sum_{i=1}^L \frac{c_j}{k_j}\left(b_{\NC}^{\delta_i(\Gamma)}+1\right)
+\sum_{i=L+1}^d \frac{a_j}{k}\geq 1\,.
\end{equation}	
This is the statement of \autoref{prop:compensation} below, which is long
and technical so we separate it. 
\end{proof}
\par 
It remains to introduce and justify \autoref{prop:compensation} used in the above
proof. This requires some additional preparation. Denote by 
\[0\leq s_j \= c_j / k_j < 1\]
the (rational) argument of the action of $\bftau$ on $t_j$ (mod $2\pi i$) as used
in \autoref{eq:condition}. For a (vertical) edge $e$, denote by $[e]$ the interval
of level passages crossed by $e$. We say that a level passage is {\em non-trivial}
if the corresponding $s_j > 0$ and that an edge is {\em non-trivial} if it crosses
a non-trivial level passage. We also say that a vertex has {\em order~$k$} if the
order of~$\bftau$ restricted to that vertex has order $k$. Finally for an edge~$e$
we define its {\em contribution}  
$$ c_{e} \= \sum_{j\in [e]} (\ell_j/p_e) s_j $$
which depends on $\bftau$ but we skip it in the notation when the context is clear.  
\par 
\begin{lemma}
\label{le:nc-edge-fix}
Let $e$ be an edge fixed by $\bftau$ and joining two vertices $v_1$ and $v_2$ where each $v_i$ has order $k_i$. Suppose either $k_1 \neq k_2$, or $k_1 = k_2$ and $e$ is non-trivial. Then $ \lcm (k_1, k_2) c_e$ is a positive integer.  In particular, 
\begin{eqnarray}
\label{eq:nc-edge-fix}
c_e \,\geq\, \frac{1}{\lcm (k_1, k_2)}\,.
\end{eqnarray}  
\end{lemma}
\par
\begin{proof}
Recall the local equation $x_1 x_2 = \prod_{j \in [e]} t_j^{\ell_j / p_e}$ at the node represented by~$e$ in the universal family over the moduli space of multi-scale differentials, where $x_i$ is a standard coordinate at $e$ in $v_i$. Consider first the case $k_1 \neq k_2$. Since $\bftau(x_i) = \zeta_{k_i}x_i$ where $\zeta_{k_i}$ is a primitive $k_i$-th root of unity, $\bftau(x_1x_2)$ differs from $x_1x_2$ by a non-trivial root of unity of order at most $\lcm (k_1, k_2)$. Hence in this case $c_e \geq 1/\lcm (k_1, k_2)$. If $k_1 = k_2 = k$, then $\bftau(x_1x_2)$ differs from $x_1x_2$ by a root of unity of order at most $k$, which implies that $c_e$ is either zero or at least $1/k$. However the former is impossible since by assumption $e$ crosses some non-trivial level passage, i.e., at least some $s_j > 0$ in the sum. This thus verifies the inequality~\autoref{eq:nc-edge-fix}. 
\end{proof}
\par 
\begin{lemma}
\label{le:nc-edge-permute}
Let $e_1, \ldots, e_h$ be edges which are cyclically permuted by $\bftau$ and which  join two vertices $v_1$ and $v_2$ fixed by $\bftau^h$, where each $v_i$ has order $k_i$. Suppose either $k_1 \neq k_2$, or $k_1 = k_2$ and the $e_i$ are non-trivial. Then $\lcm (k_1, k_2) c_{e_i}$ is a positive integer. In particular
\begin{eqnarray}
\label{eq:nc-edge-permute}
c_{e_i} \geq \frac{1}{\lcm (k_1, k_2)}\, 
\end{eqnarray} 
for all $i$.  
\end{lemma}
\par
\begin{proof}
Note that $\bftau^h$ fixes each $e_i$ and the two vertices, and that it has order $k_i / h$ on each $v_i$. Then 
the proof of the previous lemma implies that $h c_e$
is a multiple of $1/\lcm (k_1/h, k_2/h)$. It is in fact a non-zero multiple
by the previous proof in the case of $k_1 \neq k_2$, and directly by the existence of
a non-trivial level passage in the case of $k_1 = k_2$. 
\end{proof}
\par 
Finally we let 
$$r_{j} = b_{\rm NC}^{\delta_j(\Gamma)} + 1$$ 
and rewrite the contribution $c(\bftau)$ in \autoref{eq:condition}
 in terms of the above notations as
$$c(\bftau) \= \sum_{j=1}^L r_j s_j + \age (\bftau)|_{\bA_{\rel} \times \bC_{\hor}}\,$$ 
where by definition $\age (\bftau)|_{\bA_{\rel} \times \bC_{\hor}}=\sum_{i=L+1}^d \frac{a_j}{k}$. 
\par
\begin{prop}
	\label{prop:compensation}
	For $g\geq 2$, suppose that $\bftau$ does not induce a quasi-reflection on~$\bA$
	and that not all level passages are trivial under $\bftau$. Then~$c(\bftau)\geq 1$. 
\end{prop}

Before showing the proof of the above proposition, we present an alternative version of $R_{\NC}^\Gamma$ which we will need to use to prove  \autoref{intro:minimal} in low genera. (More precisely, we will need the following version of $R_{\NC}^\Gamma$, together with the improvement in \autoref{prop:nc-refinement}, to show that the minimal strata with odd spin parity are of general type for $13\leq g\leq 43$.)

\begin{prop}
	\label{prop:compensationv2} 
	For $g\geq 2$,  suppose that $\bftau$ does not induce a quasi-reflection on~$\bA$
	and that not all level passages are trivial under $\bftau$. Let $v^\top$ be the number of top level vertices in $\Gamma$. 
	Then substituting $R_{\NC}^\Gamma$ in the definition \autoref{eq:non-canonical-term} with
	\[R_{\NC}^\Gamma =\, \sum_{\rm NCT} \frac{1}{E_\Gamma}\frac{1}{p_e} + 
	\sum_{\rm RBT} \frac{1}{p_{e}}  + \sum_{\rm OCT} 
\frac{2}{ p_{e}} +\sum_{\rm EDB} \frac{4}{p_{e}}+(v^\top-1) \,\]
still satisfies that $c(\bftau)\geq 1$.
\end{prop}

Our strategy is to prove first \autoref{prop:compensationv2}. Since $1/E_\Gamma\leq 1/2$ in the presence of NCT edges, the same proof will work for \autoref{prop:compensation} if we can show that the additional $v^\top$-term, which is not present in \autoref{eq:non-canonical-term}, is not needed if the coefficient for NCT edges is $1/2$ instead of  $1/E_\Gamma$. 

\begin{proof}[Proof of \autoref{prop:compensationv2}]
We denote by $\Pi$ the level graph on which we perform the analysis. If $\bftau$ permutes some vertices of $\Pi$, then $\age (\bftau)|_{\bA_{\rel} \times \bC_{\hor}}\geq 1$ unless the permuted vertices consist of two vertices of genus zero swapped by $\bftau$  as described in \autoref{lem:boundpermuted}. For such a permuted pair we can combine the two vertices as one 'hyperelliptic' vertex, which does not affect the analysis of edge contributions when we apply \autoref{le:nc-edge-permute} (for $h = 2$). Moreover, a trivial pair of age zero described in \autoref{rem:trivial-pair} behaves the same way as a single hyperelliptic vertex of age zero. In this sense from now on we can assume that \emph{$\bftau$ fixes every vertex of $\Pi$}. Moreover, we can also assume that \emph{$\Pi$ has no horizontal edges} (since higher order vertices with $\age < 1$ do not
admit such edges according to the tables and as a consequence of the
trivial pair discussion in the proof of \autoref{prop:compensation}).
\par
Let $\mathrm{H}$ be a non-trivial level passage of $\Pi$ such that it is crossed by the maximum number of edges among all non-trivial level passages, where we denote by $E_H$ the number of edges crossing $\mathrm{H}$. By~\autoref{eq:nc-edge-fix} or~\autoref{eq:nc-edge-permute}, any non-trivial edge joining two vertices of order $k_1$ and $k_2$ contributes at least $c_e \geq 1/\lcm (k_1, k_2)$ times the corresponding $(\ell/p_e)$-coefficient in~\autoref{eq:non-canonical-term}, which is then at least $1/(E_{\mathrm{H}}\lcm (k_1, k_2))$ for all edge types. Hence we can sum up the contributions of the $E_{\rm H}$ edges and obtain that 
\begin{equation}
	\label{eq:levelHcontrib}
c(\bftau) \geq \frac{E_{\mathrm{H}}}{E_{\mathrm{H}}\lcm (k_1, k_2)}=\frac{1}{\lcm (k_1, k_2)}\,.
\end{equation}
 If $k_1 = k_2 = 1$ for these edges, then we obtain enough contribution to $\age(\bftau)$.
\par
Another preliminary remark is that, if there is a level passage $j_0$ whose
corresponding two-level graph has $v^\top=v_{j_0}^\top>1$ vertices on top level, then
there is a \emph{special edge}~$e$ crossing this level passage such that all the level
passages $j\in [e]$ satisfy $v_j^\top>1$, where $v_j^\top$ denotes the
number of top level vertices of $\delta_j(\Pi)$. To see the existence of such a special edge $e$, consider all
non-backtracking paths that start and end with an edge crossing~$j_0$,
with exactly these two edges crossing $j_0$, and that connect two disjoint
connected components above~$j_0$.
By hypothesis this set is nonempty, and we can orient the paths such that the
starting level of the paths is not above the ending level.
This means that the first edge~$e_1$ of the path has the property that
all level passages in $[e_1]$ above $j_0$ satisfy $v_j^\top>1$. 
Consider now a path where the lowest level touched by the path is as high as
possible among all paths. Then the first edge~$e_1$ of this path also has the property
that all level passages in $j\in [e_1] $ satisfy $v_j^\top>1$. Indeed if there
is a level passage $j'\in [e_1]$ below $j_0$ with  $v_{j'}=1$, then we would find
another path with~$j'$ as the lowest touched level. In summary, one can use $e = e_1$ as a special edge. Note that if the starting level passage $j_0$ is non-trivial, e.g., $j_0=H$, then the special edge $e$ is non-trivial, since it crosses $j_0$. Hence by \autoref{le:nc-edge-fix} or \autoref{le:nc-edge-permute},
we obtain a contribution of at least 
\begin{equation}
	\label{eq:vtopcontrib}
	(c_e/E_{\mathrm{H}}+1)/\lcm(k_1,k_2)\,
\end{equation} where~$k_i$ are the order of the
vertices joined by~$e$. We call this  the \emph{$v^\top$-contribution} of the
special edge.
\par
From now on we can assume that there is at least one vertex of order $k_i > 1$. In this case the edge contribution can become smaller, but the vertex age can make an extra contribution to $c(\bftau)$ by using \autoref{ncap:noncan} and \autoref{ncap:noncanRC}. We thus need to analyze a number of cases depending on the orders of vertices. 

 \subsection*{\em Case that all vertices are of order one or two}
 Assume that all vertices are of order one or two, and there is at least one vertex of order two. For any such vertex, we can assume that it is hyperelliptic under $\bftau$ (otherwise the vertex age is at least one). In this case an edge is either fixed or permuted with another edge under the hyperelliptic involution. Since the total edge contribution is at least $1/2$ from the contribution \autoref{eq:levelHcontrib} (for $k_1, k_2 \in \{1,2\}$), we can assume that the age of every hyperelliptic vertex is zero, since $\age(\bftau) \geq 1$ and we are done otherwise. By \autoref{rem:hypeigen} and \autoref{nlem:ageGRC}, a hyperelliptic vertex of age zero (including the age from $\bC_{\rm hor}$) does not admit horizontal edges, has at most one fixed zero or one pair of permuted zeros, has every fixed pole with residue zero constrained by GRC, and has every permuted pair of poles with the zero sum of the residues constrained by GRC. Therefore, we can assume that every hyperelliptic vertex belongs to a hyperelliptic (banana) tree of age zero, which is a connected subgraph consisting of hyperelliptic vertices of age zero whose edge configurations are described above. If a non-hyperelliptic vertex is adjacent to one end of a hyperelliptic tree, we say that the connecting edge is a {\em handle}. 
\par 
Among the $E_{\mathrm{H}}$ edges crossing ${\mathrm{H}}$, suppose $E_1$ edges are fixed
and have both endpoints trivial (\emph{type one}), $E_2$ edges are fixed and have at
least one endpoint hyperelliptic (\emph{type two}), and $2E_3$ edges are permuted in
pairs with both endpoints hyperelliptic (\emph{type three}), where
$E_1 + E_2 + 2E_3 = E_{\mathrm{H}}$. 
By~\autoref{eq:nc-edge-fix} and~\autoref{eq:nc-edge-permute},
edges $e_1$, $e_2$, $(e_3, e'_3)$ of the three types satisfy that $c_{e_1} \geq 1$,
$c_{e_2} \geq 1/2$, and $c_{e_3} + c_{e'_3} \geq 1/2 + 1/2 = 1$, respectively. If there are edges of type two or three constrained by GRC, then the undegeneration of the level passage $\mathrm{H}$ gives rise to a two-level graph with $v_{\mathrm{H}}^\top>1$, and the $v^\top$-contribution \autoref{eq:vtopcontrib} is at least $1/2$ which is enough. Hence we can assume that for every edge of type two which is not a handle, or for every pair of edges of type three, there is at least one associated CPT handle or at least two NCT handles.
 Note that a handle~$e$ of a hyperelliptic tree has endpoints of order one and
two respectively, hence $c_e \geq 1/2$. In addition, any permuted edge of type three
is not a handle.    
\par 
Consider first the case that every handle of a hyperelliptic tree is not CPT. Then every hyperelliptic tree has at least two handles, one at each end. 
It follows that $c(\bftau) \geq E_1 / E_{\mathrm{H}}  + (1/2 + 1/2) E_2 / E_{\mathrm{H}}  + (1+1/2+1/2) E_3 / E_{\mathrm{H}}  = 1$. Next consider the case that there are at least two CPT handles. Then $c(\bftau) \geq 1/2 + 1/2 = 1$. Finally suppose there is exactly one CPT handle $e$. Then $c(\bftau) \geq 1/2 + E_1 / E_{\mathrm{H}} + E_2  / E_{\mathrm{H}} + 2E_3 / E_{\mathrm{H}} - 1/E_{\mathrm{H}} \geq 1$ for $E_{\mathrm{H}}\geq 2$, where we subtract $1/E_{\mathrm{H}}$ for not counting redundantly the contribution from the hyperelliptic tree with the handle $e$. If $E_{\mathrm{H}} = 1$, then the graph consists of a single hyperelliptic tree (modulo trivial level passages), hence either $c(\bftau) \geq 1$ or it induces a quasi-reflection. 

\smallskip
From now on we assume that besides vertices of order one and two, all other vertices are of order three, four, or six from \autoref{ncap:noncan} and \autoref{ncap:noncanRC}.
\par
	\subsection*{\em Case that a vertex is of order four}
We start by assuming that there is at least one vertex  of order four. Since all vertices of order four have age at least $1/2$, we can assume that there is a unique vertex $v$ of order four, which rules out Case (8) as it has four cyclicly permuted edges and thus requires another vertex of the same type to pair with it, giving enough vertex age.  If $v$ has a pair of permuted edges, then they can only join a hyperelliptic vertex of age zero. Since $\age(v) = 3/4$ in the relevant Cases (11), (12), (R6) and (R7), we can assume that all other vertices are of order one or two. As before, since the $E_{\mathrm{H}}$ edges contribute at least $1/4$, we conclude that $c(\bftau) \geq \age(v) + 1/4 = 1$. 
\par 
Now suppose all edges of $v$ are fixed. The relevant cases are~(9) and~(10), for which $\age(v)=1/2$. Case~(9) has genus one with a unique edge $e$ going down. Suppose the lower end~$v'$ of~$e$ has order $k'$ which can be one, two, or three (as in Case~(13) of order $6$ cannot be a lower vertex). Since $e$ is CPT (and not RBT as $g\geq 2$), we obtain that $c(\bftau) \geq \age(v) + \age(v') + 2/\lcm(4, k') \geq 1$ for any $k' = 1,2,3$. 
\par 
For Case (10), if it admits one edge and one leg, then the same argument works as in Case (9). Suppose it admits two edges $e_1$ and $e_2$, one going up to $v_1$ and the other going down to $v_2$ where $v_i$ is of order $k_i\neq 4$. As before we can reduce to the case that all non-trivial edges are not CPT, any hyperelliptic vertex has age zero, and there is at most one other vertex with non-zero age. Note that $e_1$ and $e_2$ cross different level passages. If there is no vertex of order three or six, then $c(\bftau) \geq \age(v)  + (1/ 4 + 1/4 + (E_{\mathrm{H}}-1)/2) / E_{\mathrm{H}} = 1$. If there is exactly one vertex of order three or six, then $c(\bftau) \geq \age(v) + 1/3 + (1/12 + 1/4 + (E_{\mathrm{H}}-1)/6) / E_{\mathrm{H}} \geq 1$. 
\par 
\subsection*{\em Case that a vertex is of order six} 
Now we assume that all vertices are of order one, two, three, or six, and there is at least one of order six. Consider Case (13) for a vertex $v$ of order six where $\age(v) = 1/3$ and $v$ admits a unique edge $e$ going down. Let $v'$ be the lower end of $e$ with order $k'$. Since $e$ is CPT (and not RBT as $g\geq 2$), its contribution to $c(\bftau)$ is $2c_e\geq 2/\lcm (6, k') \geq 1/3$, hence combining with $\age(v)$ we can assume that all vertices (except $v$) are of order one or two and that any hyperelliptic vertex has age zero. If there are other non-trivial edges besides $e$, then either they are CPT and we obtain enough age, or similarly as before we obtain $c(\bftau) \geq 2/3 + (E_{\mathrm{H}}-1)/(2
E_{\mathrm{H}}) \geq 1$ for $E_{\mathrm{H}}\geq 3$. Otherwise it reduces to EDB, which implies that $c(\bftau)\geq \age (v) + 4/\lcm(1,6) = 1$.
\par
\subsection*{\em Case that a vertex is of order three}
Finally we assume that all vertices are of order one, two, or three, and there exists at least one vertex $v$ of order three.  If $v$ admits three cyclicly permuted edges, then it must be paired with another vertex $v'$ of order three, which includes Cases (2), (6), (7), (R1) and (R2).  The only combination for $\age(v) + \age (v') < 1$ is where $v$ and $v'$ are of type (2), (R1) or (R2), and then $\age(v) + \age (v') = 2/3$. We can assume that any hyperelliptic vertex has age zero and that there are no CPT edges.  Either the three edges do not cross the level passage ${\rm H}$, and then
$c(\bftau)\geq 2/3+(1/6+1/6+(E_{\mathrm{H}}-1)/2)/E_{\mathrm{H}}\geq 1$, or we obtain that 
$c(\bftau) \geq 2/3 + (1/3+1/3+1/3 + (E_{\mathrm{H}}-3)/2)
/E_{\mathrm{H}} \geq 1$ since $E_{\mathrm{H}}\geq 3$ in this case.
\par 
Now suppose all vertices of order three have fixed edges only, which includes
Cases (3), (4), (5), (R3), (R4) and (R5).  For Cases (3), (R3) and (R4), $\age (v) = 2/3$, and they can be treated similarly as in the preceding case. For Case (4), the vertex~$v$ admits a unique edge and we can use the same argument as in Case (13). 
\par 
For Case (5), if $v$ admits one leg and one edge, then the argument is still the same as Case (13). Next suppose $v$ admits two edges $e_1$ and $e_2$, going up and down respectively to $v_1$ and $v_2$.  Again we only need to consider the case that all other vertices (except $v$) are of order one or two and any hyperelliptic vertex has age zero. It follows that $c(\bftau) \geq 1/3 + (1/6 + 1/6 + (E_{\mathrm{H}}-1)/2 ) /
E_{\mathrm{H}}$. If there is an edge crossing ${\mathrm{H}}$ constrained by GRC, then the corresponding two-level graph has $v_{\mathrm{H}}^\top>1$, and so the $v^\top$-contribution \autoref{eq:vtopcontrib} gives an additional $1/2$, which  is enough. Hence, as in the case of vertices of order two, the edges crossing the level ${\mathrm{H}}$ and the handles of the hyperelliptic trees crossing ${\mathrm{H}}$, give a contribution of at least $(E_{\mathrm{H}}-1)/E_{\mathrm{H}}+1/6$, where we subtract one since the contributing edge might be $e_1$ or $e_2$. If $E_{\mathrm{H}} >1$, then the preceding estimate together with the vertex age is enough.
Consider finally the case that  $E_{\mathrm{H}} = 1$. Then the graph (modulo trivial level passages) reduces to a tree where both $e_1$ and $e_2$ are CPT. Moreover if $v_2$ is hyperelliptic, then $e_2$ cannot be RBT by stability of $v_2$. Therefore, the contribution of $e_i$ to $c(\bftau)$ is at least $1/3$, which is from either $2c_{e_i} \geq 2/\lcm (2,3) = 1/3$ or $c_{e_i} \geq 1/\lcm (1,3) = 1/3$. It follows that  $c(\bftau) \geq \age(v) + 1/3 + 1/3 = 1$.  
\par
For Case (R5) a similar argument as the previous one works. 
\end{proof}
A slight modification of the previous proof can be used to show \autoref{prop:compensation}.
\begin{proof}[Proof of  \autoref{prop:compensation}]
Since $1/E_\Gamma\leq 1/2$ in the presence of NCT edges, the same proof will work for \autoref{prop:compensation} if we can show that the $v^\top$-contribution \autoref{eq:vtopcontrib} is not needed if we change the NCT coefficient to $1/2$ as in \autoref{prop:compensation}. Indeed in the proof of \autoref{prop:compensationv2}, the $v^\top$-contribution is only used in two instances, one in the case of vertices of order one or two and the other in the case of vertices of order three, to ensure that if there are edges crossing $\mathrm{H}$ constrained by GRC, then we obtain a contribution of at least $1/2$. In both cases since each edge crossing $\mathrm{H}$ constrained by GRC gives a contribution of at least $(1/2)\cdot (1/2)= 1/4$ and there are at least two such edges, these edges contribute at least $1/2$, hence giving enough without using the $v^\top$-contribution.    
\end{proof}
\par
In order to prove \autoref{intro:minimal} for low genera (the general type result for the minimal odd strata with $g\leq 43$), we need to further reduce the $(\ell/p_e)$-coefficients in $R_{\rm NC}^{\Gamma}$ for certain graphs as follows.  
\par
We say that a two-level graph for the minimal stratum in genus $g$ is
a {\em rational multi-banana (RMB)} if it has two vertices only, one
on top and one on the bottom,
joined by~$E$ edges for $E\geq 2$, where the bottom vertex is of genus zero.
\par
\begin{prop}   \label{prop:nc-refinement}
For the minimal strata $\mu=(2g-2)$, we can refine the $R_{\NC}^\Gamma$
coefficients of \autoref{prop:compensationv2} by setting
$R_{\NC}^\Gamma= 1/\ell_\Gamma$ for  RMB graphs with prongs of type $(1^{E_\Gamma-1}, p)$
where $p\geq 2E_\Gamma-3$ or with prongs of type $(1^{E_\Gamma-2}, 2, p)$ where $p\geq 2E_\Gamma-2$, and by
setting  $R_{\NC}^\Gamma= P_{-1}/(E_\Gamma+1)$ for RMB graphs with at least one
prong of order one and at least one prong of order greater than seven, where 
$P_{-1}=\sum_{e\in E(\Gamma)} 1/p_e$.
\end{prop}
\par
\begin{proof}
We will use the same notation as in the proof of \autoref{prop:compensationv2}.
\par
\subsection*{\em Case that there exists an RMB with prongs
$(1^{E_\mathrm{R}-1}, p)$ or $(1^{E_\mathrm{R}-2}, 2, p)$}
Consider first the case in a minimal stratum that $\Pi$ has a level passage $\mathrm{R}$ whose undegeneration gives an RMB graph of prong type $(1^{E_\mathrm{R}-1}, p)$ with $p \geq 2$. Let $e_p$ be the edge of prong $p$. Let $v_i^{\pm}$ be the upper and lower ends of the edges of prong order one and $u^{\pm}$ the upper and lower ends of $e_p$, where some of the vertices can coincide. We make some observations first. We can assume that
$\mathrm{R}$ is non-trivial, since otherwise we can apply the initial argument. Since in the minimal strata $\Pi$ has no local minima other than the unique bottom level vertex, any vertex $v_i^+$ lower than~$u^+$ admits a unique polar edge of prong one. It follows that $p$ and $1$ are the only prong values for the level passages between $\mathrm{R}$ and $u^+$. On the other hand, the subgraph below $\mathrm{R}$ is a tree with rational vertices where each vertex goes down to the bottom vertex via a unique path. In particular, other than~$p$, the maximum prong value between $\mathrm{R}$ and $u^-$ can be at most $2 E_{\mathrm{R}} - 3$, which is less than or equal to $p$ by assumption.  Moreover the prong values between $\mathrm{R}$ and $u^-$ are
all greater than two. It follows that $P_{-1,j}/E_j \geq 1/p$ for $j \in [e_p]$ and any
level passage in $[e_p]$ cannot be RMB with at least one prong of order one (except for
the prong type $(1^{E-1}, p)$ where $1/\ell_{\mathrm{R}}=1/p$ is the desired no-compensation
coefficient). Since $e_p$ crosses the non-trivial level passage $\mathrm{R}$, we thus conclude from \autoref{le:nc-edge-fix} and \autoref{le:nc-edge-permute} that 
\[c(\bftau) \geq \sum_{j\in [e_p]} \frac{ \ell_jP_{-1,j}}{E_j} s_j \geq c_{e_p} \geq 1/\lcm (k^+, k^-)\] 
where $k^{\pm}$ are the vertex orders of $u^{\pm}$ under $\bftau$. 
\par
If $u^{\pm}$ are both of order one, then $c_{e_p}\geq 1$ and we are done. If one of them is hyperelliptic, then $c_{e_p} \geq 1/2$, and we can assume that any hyperelliptic vertex is of age zero and that $e_p$ is not constrained by GRC, since otherwise the $v^\top$-contribution \autoref{eq:vtopcontrib} would give another $1/2$. First suppose $u^-$ is hyperelliptic. Since it has age zero and $e_p$ is not constrained by GRC, the vertex $u^-$ admits a unique polar edge (which is $e_p$) and a unique zero edge (or leg), which contradicts stability as $u^-$ is of genus zero. Next suppose that~$u^+$ is hyperelliptic and let $e_q$ be a polar edge of maximal prong $q$ (among all the polar edges of $u^+$). Using as before that there is no local minima other than the bottom level vertex, we conclude that all prong values in $[e_q]$ are $q$ and $1$. Therefore, if $e_q$ is not trivial, we obtain from these level passages that $c_{e_q} \geq 1/2$ and hence $c(\bftau) \geq c_{e_p} + c_{e_q} \geq 1$. If $e_q$ is trivial, then its upper vertex~$u_1^+$ is hyperelliptic (of age zero), and all other vertices of order one connected to~$u^+$ are on higher level than $u_1^+$. Hence we can iterate the same argument we used for~$u_+$ to~$u_1^+$, and continue this procedure until we reach a non-trivial edge. The procedure has to stop since the top level of $\Pi$ cannot have hyperelliptic vertices of age zero (since those have to be of genus zero).
\par
The argument for an RMB graph of prong type $(1^{E_\mathrm{R}-2}, 2, p)$ with $p \geq 2$ is similar. Indeed similarly as before,  $p$, $1$ and $2$ are the only prong values for the level passages between~$\mathrm{R}$ and $u^+$. Moreover, arguing as before, the maximum prong value between~$\mathrm{R}$ and $u^-$ can be at most $2 E_{\mathrm{R}} -2$, which is less than or equal to $p$ by assumption. 
One additional observation is that in this case $-2 (E_\mathrm{R} - 2) - 3 - (p+1) + 2g-2 = -2$, hence $p$ is even. In particular $1/\ell_{\mathrm{R}}=1/p$ as before. So the previous case of having only order one vertices and the previous case where $u^-$ is hyperelliptic are clear.  If  $u^+$ is hyperelliptic of age zero, then under the above notation assume that $e_q$ is non-trivial. Then $e_q$ is the only polar edge and $q$ is even, since  $-q - 1 + p-1 = - 2$ and $p$ is even. Indeed, if there are polar edges constrained by GRC, there is a non-trivial level passage with $v^\top>1$, and hence the $v^\top$-contribution \autoref{eq:vtopcontrib} gives another $1/2$. If there are a pair of permuted edges with prong $q$, we have $2(-q-1)+p-1= -2$, which is impossible since $p$ is even. It follows that $q\geq 2$ and $q$ is still the largest prong for the level passages in $[e_q]$ (since there might be edges with prongs of order~$2$ crossing a level passage in $[e_q]$, we had to rule out the case of $q=1$ in order for the previous sentence to be true).  If $e_q$ is trivial, as before we can iterate the procedure until we find a non-trivial edge giving a contribution of $1/2$.
 \par
 Consider now RMB with prongs of type $(1^{E_\mathrm{R}-1}, p)$ or $(1^{E_\mathrm{R}-2}, 2, p)$ in the situation of having higher order vertices from \autoref{ncap:noncan} and \autoref{ncap:noncanRC}. 
 \par
 First we treat the cases in \autoref{ncap:noncanRC} with GRC of order four. If there is such a vertex, then $\Pi$ has a non-trivial level passage whose corresponding divisor has $v^\top>1$. Then in this case the $v^\top$-contribution \autoref{eq:vtopcontrib} is at least $1/\lcm(k_1,k_2)$, where~$k_i$ are the orders of the vertices joined by the special edge with the property that all the level passages crossed by it have $v^\top>1$. It is easy to check that this contribution, together with the age contribution from the vertices, is enough.
 \par
Next consider the cases in  \autoref{ncap:noncan}. Since the Cases (4), (9) and (13) give CPT edges, which cannot cross $R$, we do not need to consider these cases. Indeed we have seen that, using the above notation, it is enough to consider the contribution from the edges $e_p$ and $e_q$. Since $e_p$ cannot be a CPT edge, the only possibility would be that $e_q$ is the CPT edge, but in this case $c_{e_q}\geq 2/6$, which, together with $c_{e_p}$ and the age contribution from the vertices, is enough. Moreover, Case (8) needs to join another vertex of order four, thus giving enough vertex age. The remaining cases of order four are (10) and (11). Since these cases correspond to vertices of positive genus, the vertex $u^{-}$ cannot be of this type. If $u^{+}$ is not of type (10) or (11), then we obtain enough age from $c_{e_p}\geq 1/2$ and the vertex age. If $u^{+}$ is of type (10) or (11), we obtain that $c_{e_p}+c_{e_q}+\age(u^{+})\geq 1/4+1/4+1/2=1$.
 \par
 We are left with considering the cases of order three. The only cases of genus zero are (R1) and (R2), both having three permuted edges which need to join another vertex of order three. Then the age contribution of at least $2/3$ from these two vertices, together with the additional $v^\top$-contribution \autoref{eq:vtopcontrib} of at least $1/6$ and the edge contribution of at least $1/6$, is enough. The remaining cases are all of positive genus, hence they cannot correspond to the vertex $u^{-}$. If $u^{-}$ is hyperelliptic, by the argument above we have $\age(u^{-})\geq 1/2$. Then the vertex age is at least $1/2+1/3$ and $c_{e_p}\geq 1/6$, hence we obtain  enough age. Consider finally the case where $u^{-}$ is of order one, and then $c_{e_p}\geq 1/3$. If  the vertex age is at least $2/3$,  we obtain enough age. Hence we are left to consider the case where only one  vertex of order three is present and it has age $1/3$, which means Cases (5) and (R5) (Case (4) admits a non-trivial CPT edge which gives enough contribution). If the upper end of $e_q$ is of order one, then we obtain $c_{e_p}+c_{e_q}\geq 2/3$, which together with the vertex age, is enough. If the upper end of $e_q$ is hyperelliptic (necessarily of age zero), then by the same analysis for $u^+$ we obtain another edge giving a contribution of at least $1/2$, hence we obtain enough together with the contribution of $c_{e_p}+c_{e_q}\geq 1/3+1/6$ and the vertex age.
 
\subsection*{\em Case that there exists an RMB with a prong of order one and a prong of large order}
Consider a graph $\Pi$ such that one of its undegenerations is an RMB graph with at least one prong of order one and one prong of order greater than seven. We can also assume that there is no special RMB level of prong type $(1^{E_\mathrm{R}-1}, p)$
or $(1^{E_\mathrm{R}-2}, 2, p)$. As before, denote by~$\mathrm{H}$ the level passage with the largest number of edges crossing and by~$\mathrm{R}$ the RMB level passage, which we can assume to be non-trivial.
\par
Consider first the case of having only vertices of order one. Since any non-trivial edge in $\Pi$ yields a contribution of at least $1/(E_\mathrm{H}+1)$, the $E_\mathrm{H}$ edges crossing~$\mathrm{H}$ give a contribution of at least $E_\mathrm{H}/(E_\mathrm{H}+1)$. If there is a non-trivial level passage different from $\mathrm{H}$,  then we gain the additional contribution to reach age one. If $\mathrm{H}$ is the only non-trivial level passage, then by \autoref{le:nc-edge-fix} and \autoref{le:nc-edge-permute} we obtain $s_\mathrm{H}\ell_{\mathrm{H}}/p_e\in \bN$ for all edges $e$ crossing $\mathrm{H}$, which is impossible for $0<s_{\mathrm{H}}<1$.
\par
Consider now the case where the vertices have order one or two. Since the edge contribution is at least $E_\mathrm{H}/2(E_\mathrm{H}+1)$, if there is a contribution of $1/2$ from a vertex of order two or from a non-trivial level passage with $v^\top>1$, we can argue similarly as before. Indeed any edge yields a contribution of at least $1/2(E_\mathrm{H}+1)$. Hence if there is a non-trivial level passage different from $\mathrm{H}$, we obtain enough. As before we can also rule out the case that  $\mathrm{H}$ is the only non-trivial level passage. Therefore, it suffices to consider the case where all vertices of order two are hyperelliptic of age zero and all non-trivial level passages have $v^\top=1$.
\par
Since we can assume that $v_\mathrm{H}^\top=1$, using the same analysis as in the general case, i.e., the case with vertices of order two in the proof of \autoref{prop:compensation}, we obtain a contribution of at least $c(\bftau)\geq E_\mathrm{H}/(E_\mathrm{H}+1)$ by considering the contributions of the $E_\mathrm{H}$ edges together with the contributions of the handles of hyperelliptic trees crossing $\mathrm{H}$. We thus need to find an additional contribution of at least $1/(E_\mathrm{H}+1)$. 
\par 
Let $e$ be the special edge with  prong one and let $v^\pm$ be the vertices at its upper and lower
ends.  If $e$ is a handle of type two, i.e., it is a fixed edge and at least one of its endpoints is hyperelliptic, then $v^+$ cannot be hyperelliptic of age zero and hence $v^-$, being hyperelliptic of age zero,  has a pair of permuted zero edges of prong one. This is impossible since the graph below $v^-$ is a tree with rational vertices by the RMB assumption for the level $\mathrm{R}$.  Hence if $e$ does not cross $\mathrm{H}$, we obtain the additional contribution we need. If $e$ crosses $\mathrm{H}$ and is of type three, i.e., it belongs to a pair of permuted edges with prong one, then $v^+$  admits a single polar edge $e_1$, which is a handle of prong one. Moreover, $v^-$ also admits a single zero edge  $e_2$ of prong three, which (modulo trivial level passages) is a handle since the graph below $v^-$ is a tree.  It follows that $v^+$ (resp. $v^-$) is (modulo trivial level passages) on the top (resp. bottom) level of $\Pi$ (since otherwise we would gain an additional contribution from edges joining trivial vertices or hyperelliptic trees not crossing $\mathrm{H}$). Moreover, they are the only vertices on the top and bottom levels, since $v_{\mathrm{H}}^\top=1$ and $\Pi$ has a unique bottom level. Let $e_p$ be the special edge crossing $\mathrm{R}$ with prong $p\geq 8$. If $e_p$ crosses only the level passages crossed by $e$ or the upper handle $e_1$ (which have both prong order one), then $c_e+c_{e_1}\geq pc_{e_p}\geq p/2\geq 2$. Since the original calculation only used a contribution of $1$ for $c_e+c_{e_1}$, we have found the desired missing contribution. Therefore, we can assume that $e_p$ joins the top and bottom vertices (since otherwise we would get enough contribution from additional edges joining trivial vertices or hyperelliptic trees). In this case $e_p$ joins two trivial vertices and since $p\geq 8$, then $c_e+c_{e_1}+c_{e_2}\geq pc_{e_p}/3\geq p/3\geq 3/2+1$. This is enough since in the original computation the contribution given by $c_e+c_{e_1}+c_{e_2}$ was $3/2$.
\par
Now we have reduced to the situation where $e$ crosses $\mathrm{H}$ and it is of type one, i.e., the vertices $v^\pm$  are trivial. In the following arguments we can assume that edges joining trivial vertices are non-trivial (e.g., by collapsing those trivial edges and merging the corresponding trivial vertices). If $v^+$ admits a polar edge going up  to a hyperelliptic vertex of age zero, then the hyperelliptic tree associated to this vertex does not cross $\mathrm{H}$ and hence yields an additional contribution of at least $1/(E_\mathrm{H}+1)$ given by its (at least) two handles.
 If $v^+$ admits a polar edge going up to another trivial vertex, then we gain the desired extra $1/(E_\mathrm{H}+1)$ from this extra edge, which is neither a handle nor crossing~$\mathrm{H}$. 
From now on we can assume that $v^+$ admits zero edges only.  If there is an edge crossing $\mathrm{H}$ not joining $v^+$,  since $v^\top_\mathrm{H}=1$ and since~$v^+$ is above $\mathrm{H}$ (as $e$ crosses $\mathrm{H}$), then there is a path completely above~$\mathrm{H}$ joining the upper vertex of this edge and~$v^+$.  Then this path gives enough contribution either because there is one edge joining two trivial vertices or because it is a hyperelliptic tree. Hence we can  assume that $v^+$ is the common upper vertex for all edges crossing~$\mathrm{H}$. Suppose $v^-$ is the lowest vertex reached by edges crossing~$R$. In this case, since by assumption there is an edge $e_p$ with prong $p \geq 8$ crossing~$\mathrm{R}$, using that $e$ has prong one we obtain $c_e \geq p  c_{e_p} \geq 8  (1/2) \geq 2$. Since in the original estimate we only used $c_e\geq 1$, we gain the extra contribution we needed. We are left to show that if there is an edge crossing $\mathrm{R}$, then its ending vertex is above $v^-$ or we get enough compensation. If $v^-$ is on the bottom level of $\Pi$, then we are done.  If not, $v^-$ joins a lower level via an edge $e'$.   If the lower end of $e’$ is trivial, we obtain the desired extra contribution of at least $1/(E_\mathrm{H}+1)$.  If the lower end of $e’$ is hyperelliptic (which cannot be the lowest level of $\Pi$), then the hyperelliptic tree starting with handle $e’$ gives enough contribution.
\par
Consider finally the situation where higher order vertices appear. If $v_\mathrm{H}^\top>1$, we obtain a contribution of at least $(E_{\mathrm{H}}/(E_{\mathrm{H}}+1)+1)/\max(\lcm(k_1,k_2))$, where $k_i$ run among the order of the vertices. First, the case of $\mathrm{H}$ being the only non-trivial level passage is impossible since $\mathrm{R}$ is non-trivial with $v_\mathrm{R}^\top=1$. If now there is at least another non-trivial level passage, we have an additional edge contribution of $1/(E_\mathrm{H}+1)/\max(\lcm(k_1,k_2))$. By inspecting \autoref{ncap:noncanRC} and considering also the vertex age contribution,  the only possibility is to have at most one higher order vertex with age $1/3$, which means Cases (5) and (R5) (Cases (4) and (13) can be excluded since a non-trivial CPT edge yields enough contribution). In this case the $v^\top$-contribution is at least $1/6$. Let $v_1$ and $v_2$ be the two vertices joined to the special vertex of order three. If $v_1, v_2$ are both trivial and there are in total at least $E_\mathrm{H} + 2$ non-trivial edges, then we have the age estimate $c(\bftau)\geq 1/3 + 1/6 + (1/3 + 1/3 + E_\mathrm{H}/2) / (E_\mathrm{H}+1) > 1$.  If there are exactly $E_\mathrm{H} + 1$ non-trivial non-CPT edges, then $v_\mathrm{H}^\top=1$.  
Suppose one of the $v_i$  is hyperelliptic and the other is trivial.  Then the hyperelliptic vertex needs to go up or down further (to avoid CPT and since the bottom minimal rational vertex cannot be hyperelliptic).  Hence there are in total at least $E_\mathrm{H} + 2$ non-trivial edges.  Then we have $c(\bftau)\geq 1/3 + 1/6 + (1/3 + 1/6 + E_\mathrm{H}/2) / (E_\mathrm{H}+1) = 1$.  
Suppose finally that $v_1$ and $v_2$ are both hyperelliptic.  Then they need to go up and down respectively for the same reason. Hence there are in total at least $E_\mathrm{H} + 3$ non-trivial edges, and we have $c(\bftau)\geq 1/3 + 1/6 + (1/6 + 1/6 + (E_\mathrm{H}+1)/2) / (E_\mathrm{H}+1) > 1$.  
\par
At this point we have reduced to the situation of $v_\mathrm{H}^\top=1$ for any level passages, and we can argue similarly as in the non-special RMB situation, and use the same notation as in the case of vertices of order two in the proof of \autoref{prop:compensation}, where we considered three types  of edges. Besides these  three types, there can be now $E_4$ edges of a fourth type joining at least one vertex of order greater than two. If $E_4=0$ and the handles of the hyperelliptic trees crossing $\mathrm{H}$ do not join higher order vertices, then as before the level passage $\mathrm{H}$ gives a contribution of at least $ E_\mathrm{H}/(E_\mathrm{H}+1)\geq 2/3$ for $ E_\mathrm{H}\geq 2$, which is enough together with the additional age from the higher order vertices. One can argue similarly if there are some handles that join higher order vertices. If $E_\mathrm{H}=1$, then we have only CPT edges, and we have already discussed this case. If $E_4>0$, then we can only have $E_4=1,2,3$. If $E_4=3$, then there are at least two vertices of order three and the vertex age is at least $2/3$, which together with $ (E_\mathrm{H}-E_4+3/3)/(E_\mathrm{H}+1)\geq 1/3$ for $E_\mathrm{H}\geq 4$ is enough. If $E_\mathrm{H}=E_4=3$, then ${\rm H}$ cannot be ${\rm R}$, since ${\rm R}$ is crossed by at least two edges with different prongs. If $E_4=2$, one can argue similarly. If $E_4=1$, then we get a contribution from ${\rm H}$ of at least $ (E_\mathrm{H}-E_4)/(E_\mathrm{H}+1)$, which is at least $1/3$ for $E_\mathrm{H}\geq 2$ and at least  $2/3$ for $E_\mathrm{H}\geq 5$. Hence if  $E_\mathrm{H}\geq 5$ or $E_\mathrm{H}=1$ we get enough contribution. If $2\leq E_\mathrm{H}\leq 4$, the only two possibilities left are Cases (5) and (10) (as before Cases yielding CPT edges do not need to be considered). For Case (10), which gives a vertex age of $1/2$, we obtain a refined estimate $ (E_\mathrm{H}-1+1/4)/(E_\mathrm{H}+1)+1/E_\mathrm{H}\geq 1/2$ for $E_\mathrm{H}\geq 2$, where the last $1/E_\mathrm{H}$ term comes from the additional edge attached to the special vertex. Similarly for Case (5), which gives a vertex age of $1/3$, we obtain $ (E_\mathrm{H}-1+1/6)/(E_\mathrm{H}+1)+1/E_\mathrm{H}\geq 2/3$ for $E_\mathrm{H}\geq 2$.
\end{proof}
\par
\begin{rem}
\label{rem:nc-optimal}{\rm 
We give some examples to illustrate that the choices of the coefficients in $R_{\rm NC}^{\Gamma}$ are delicate. 
\par
Take three vertices, each on a different level, where the middle vertex joins the top and the bottom each by a single edge of prong $p_1$ and $p_2$, respectively, the top joins the bottom by $E-1$ edges of prong all equal to $p_3 = p_1 + p_2$, and $p_1, p_2$ are relatively prime.  
  Suppose $\bftau$ acts trivially on the vertices and acts on the level parameters $t_1$ and $t_2$ by $t_i\mapsto e^{2\pi i / p}t_i$, i.e., the arguments of the action (mod $2\pi i$) are $s_1 = s_2 = 1/p_3$ as in the proof of \autoref{le:nc-edge-fix}. 
  In this case $\ell_1 = p_1p_3$, $\ell_2 = p_2p_3$, and they satisfy the requirements that $(\ell_1 / p_1)  s_1 = 1 $, $(\ell_1 / p_2) s_2 = 1$, and  $(\ell_1 / p_3)  s_1 + (\ell_2 / p_3)  s_2 = 1$. We have $\ell_1 (1/p_1 + (E-1)/p_3) s_1 + \ell_2 (1/p_2 + (E-1)p_3) s_2 = E+1$, hence we cannot use a coefficient smaller than $1/(E+1)$ (and indeed we use $1/E$ for general NCT edge types in \autoref{prop:compensationv2}). 
\par 
For the RBT coefficient, since $\ell = p$ for the unique RBT edge and $\ell / p - 1 = 0$, i.e., we do not make any compensation, it is clearly sharp.  
\par 
For the general CPT coefficient (i.e., OCT), take Case (5) of order 3 and age $1/3$ as the middle vertex sitting in between two hyperelliptic trees of age zero. Let  $s_1 = s_2 = 1/6 = 1/\lcm (2,3)$ for the action on the level parameters of the two handles.  Then $c(\bftau)  = 1/3 + 2\cdot (1/6+1/6) = 1$, which implies that we cannot reduce the OCT coefficient to be smaller than $2$. Note that in this case the genus of the top or the bottom can be almost arbitrary after a divisorial undegeneration (except that it cannot be zero since Case (5) of genus one is contained in a middle level, which singles out RBT). 
\par 
For the EDB coefficient, take Case (5) as a top or a bottom vertex joining a hyperelliptic tree of age zero, and let $s = 1/6$ for the action on the level parameter of the handle. Then $c(\bftau)  = 1/3 + 4\cdot (1/6) = 1$, which shows that the EDB coefficient~$4$ is necessary.  
\par
For the minimal strata RMB coefficients in \autoref{prop:nc-refinement}, consider a triangle graph with one top vertex $v_1$, one middle vertex $v_2$, and one bottom vertex $v_3$, where $v_2$ is hyperelliptic of age zero, $v_1$ and $v_3$ are trivial, and moreover $v_3$ is rational with the unique marked zero.  Let $p_1$ be the higher short edge, $p_2$ the lower short edge, and $p_3$ the long edge satisfying that $2p_3 = p_1 + p_2$.  Suppose $p_1, p_2, p_3$ are odd and pairwise relatively prime.  Then $\ell_1 = p_1p_3$ and $\ell_2 = p_2p_3$. The top level undegeneration is a non-RMB banana graph (as $v_2$ cannot be rational with only two edges and no legs) and the bottom level undegeneration is an RMB.  Take $s_1 = s_2 = 1/(2p_3)$ which satisfies the half-integer requirement along the short edges and the integer requirement along the long edge.  Suppose we use the coefficient $1/2$ for the top non-RMB banana graph and the coefficient $1/(E_{\rm R}+1) = 1/3$ for the bottom RMB.  Then we obtain  $c(\bftau) = (1/2) \cdot \ell_1 s_1 (1/p_1 + 1/p_3) + (1/3) \cdot \ell_2  s_2  (1/p_2 + 1/p_3) = 3/4 + p_1/(12p_3)$ which can be smaller than $1$.  However in this case $p_2 \neq 1$ since $v_2$ cannot be rational, which also implies that $p_3 \neq 1$ by the relation $2p_3 = p_1 + p_2$.  Therefore, imposing a prong of order one when reducing the coefficients of RMB graphs in the minimal strata makes sense.  
\par
Finally we show that an extra contribution from multiple top vertices (i.e., the $v^\top$-term) is necessary in \autoref{prop:compensationv2}. Consider a three-level graph where $u_1, \ldots, u_n$ are trivial vertices on the top level, $v$ is a hyperelliptic vertex (of age zero) on the middle level, and $w$ is a trivial vertex on the bottom level. Suppose each $u_i$ joins $v$ by a short edge of prong $p_1$, $v$ joins $w$ by a short edge of prong $p_2$, and each $u_i$ joins $w$ by a long edge of prong $p_3$, where $p_1 + p_2 = 2p_3$ and $p_1, p_2, p_3$ are pairwise relatively prime.  The GRC imposes zero residue for each polar edge of $v$, hence the age of $v$ can be zero.  Moreover, every edge in every undegeneration of this graph is NCT.  Take $s_1 = s_2 = 1/(2p_3)$ which satisfies the half-integer requirement along every short edge and the integer requirement along every long edge.  The top level passage has $2n$ edges and the bottom level passage has $n+1$ edges. Then we obtain $c(\bftau) = 1/(2n) \cdot \ell_1 s_1 (n/p_1 + n/p_3) + 1/(n+1) \cdot \ell_2 s_2  (1/p_2 + n/p_3)  \sim  3/4 < 1$ for $n \gg 0$ and $p_3 \gg p_2$. Hence in this case an extra contribution from $v^{\top}$ is needed.  
}\end{rem}
\par 


\section{Pullback classes and the canonical class} \label{sec:taut}

In this section we recall basic properties of the tautological ring of
the moduli space of multi-scale differentials and express some divisor classes
in terms of standard generators. We apply this to the formula for the canonical
bundle and also compute the classes of divisors of Brill--Noether type that are
pulled back from the moduli space of curves. 
\par
The divisorial part of the tautological ring is generated by
the $\psi$-classes and boundary classes. It contains standard tautological divisor
classes $\xi, \lambda_1, \kappa_1$, which are all proportional in the strata interior
but can differ along the boundary (see e.g.,~\cite{ChenTauto}).  
For their conversion we define rational numbers
\be \label{eq:defkappamu}
\kappa_\mu \= \sum_{ m_i\neq -1} \frac{m_i(m_i+2)}{m_i+1} 
\= 2g - 2 + s + \sum_{m_i\neq -1} \frac{m_i}{m_i+1}\, 
\ee
for any signature $\mu = (m_1, \ldots, m_n)$, where $s$ is the number of entries equal
to~$-1$ (i.e. the number of simple poles).
The values  $\kappa^\bot := \kappa_{\mu_{\Gamma}^{\bot}}$ and $\kappa^\top :=
\kappa_{\mu_{\Gamma}^{\top}}$ are similarly defined for
the bottom and top level strata of $\Gamma$, \emph{including the edges
as legs}. In particular 
$$ \kappa^\bot + \kappa^\top \= \kappa_\mu\,.$$
This constant $\kappa_\mu$ previously
appeared (with an additional factor $\tfrac1{12}$) in \cite{EKZ} for relating sums of
Lyapunov exponents and (area) Siegel-Veech constants in the strata of holomorphic
differentials.  Our definition here includes meromorphic signatures as well. We also
define the $\psi$-class over simple poles to be 
$$\psi_{-1} \= \sum_{ m_i = -1} \psi_i\,.$$ 
\par 
The main conversion result of divisor classes we need is the following relation,
whose proof is given in \autoref{sec:kappa}. 
\par
\begin{prop} \label{prop:xiinlambda} In the tautological ring we have
  the relation
 \be
 \label{eq:eta-xi}
\kappa_\mu  \xi \=  \psi_{-1} + 12 \lambda_1 - [D_h]
- \sum_{\Gamma \in \LG_1} \ell_\Gamma \kbot [D_\Gamma]\,.
\ee
\end{prop}
\par
The canonical class of the coarse moduli space follows from the canonical class of the stack computed in \cite{EC}, using the above conversion formula and taking into account the branching behavior at the boundary. Recall
from \autoref{sec:tocoarse} the definition of graphs $\Gamma$ of type HTB, HBT and HBB (and their prime versions for meromorphic strata) as well as the corresponding ramification divisors $D_\Gamma^{\rm H}$. Recall also that $N$ denotes the dimension of the (unprojectivized) strata.  
\par
\begin{prop} \label{prop:canonicaloncoarse}
Let $\mu$ be a holomorphic signature not of type $(2m,2g-2-2m)$ (and the hyperelliptic component is excluded if $\mu = (2g-2)$), 
or a meromorphic signature not of type $(2m_1, 2m_2, 2g-2-2m_1 - 2m_2)$ (and the hyperelliptic component is excluded if $\mu = (2g-2+2m, -2m)$).
Then the class of the canonical bundle of the coarse moduli space
$\bP\MScoarse$  is given by 
\ba \label{eq:canformula}
\frac{\kappa_\mu}{N} &\c_1\bigl(K_{\bP\MScoarse}\bigr) 
\= \psi_{-1} + 12 \lambda_1 - \Bigl(1 + \frac{\kappa_\mu}{N}\Bigr)[D_h] \\
& - \sum_{\Gamma\in \LG_1} \Big(\ell_\Gamma  \kappa_{\mu_{\Gamma}^{\bot}} - \frac{\kappa_\mu}{N} (\ell_\Gamma N_\Gamma^{\bot}-1)
 \Big) [D_\Gamma]
- \frac{\kappa_\mu}{N} \sum_{\Gamma \, \text{is} \, \mathrm{HTB} \, \text{or}\atop
 \mathrm{HBT} \,\text{or}\, \mathrm{HBB}}
 [D^{\rm H}_\Gamma]
\ea
in $\CH^1(\bP\MScoarse) = \CH^1(\bP\LMS)$, where for meromorphic signatures the last
sum is for graphs of type HTB', HBT' and HBB'.  
\end{prop}
\par
We remark that for $\mu$ of holomorphic type $(2m,2g-2-2m)$ or of meromorphic type
$(2m_1, 2m_2, 2g-2-2m_1 - 2m_2)$ the above expressions have to be modified
by the class of the ramification divisor in the interior described by
\autoref{prop:ramifInt}. 
\par 
The proof of \autoref{prop:canonicaloncoarse} together with the variants for
connected components of the strata is given in \autoref{sec:spinhyp}.
\par

\subsection{Relations among tautological classes} \label{sec:kappa}

We denote by 
\be
R^\bullet(\bP\LMS) \,\subset\, \CH^\bullet(\bP\LMS)
\ee
the \emph{tautological ring}, being the smallest subring
that contains the $\psi$-classes, that is closed under the push-forward of level-wise clutching and forgetting a marked point, and moreover that contains the class of the horizontal divisor~$D_h$ (including the class of each component if $D_h$ is reducible).\footnote{A version of (small) tautological ring without $D_h$ was also considered in \cite{EC} for the purpose of running {\tt diffstrata}.}
Note that $D_h$ is irreducible in each holomorphic stratum (component), but in general it can be reducible for the meromorphic strata. 
\par
The tautological ring contains all boundary strata classes and 
standard tautological classes such as the $\kappa$-classes. Let $\pi\colon \cX \to \ol{B}$ be the universal family with $s_i\colon \ol{B}
\to \cX$ the universal sections, $S_i \subset \cX$ their images and $\omega_\pi$ the relative cotangent bundle. 
Define the Miller--Morita--Mumford class $\kappa_1 = \pi_{*} (\c_1 (\omega_{\pi})^2)$.
We give a closed expression for~$\kappa_1$ as follows.  
\par 
\begin{prop} \label{prop:kappexpr}
The class $\kappa_1$ can be expressed in terms of the standard generators
of $R^\bullet(\bP\LMS)$ as
\bes
\kappa_1 \= -\psi_{-1} + \kappa_\mu \xi + \sum_{\Gamma \in \LG_1} \ell_\Gamma
\Big(\kbot - \sum_{e\in E(\Gamma)} \frac{1}{p_e}\Big) [D_\Gamma]\,
\ees
where $\psi_{-1}$ is the sum of $\psi$-classes associated to marked simple poles.  
\end{prop}
\par
We remark that the above expression can be converted to the
Arbarello--Cornalba $\kappa$-class via the relation 
\bes
\kappa_1^{{\rm AC}}  \= \pi_{*} \Big(\c_1 \Big(\omega_{\pi}\Big(\sum_{i=1}^n S_i\Big)\Big)^2\Big)  \= \kappa_1 + \psi\,
\ees
where $\psi = \sum_{i=1}^n \psi_i$ is the total $\psi$-class. 
\par
\begin{proof} The prescribed vanishing of the differentials in the universal family along the sections $S_i$ implies that 
$$ \c_1 (\omega_{\pi}) \= \pi^{*}\xi + \sum_{i=1}^n m_i S_i
+ \sum_{\Gamma\in \LG_1} \ell_\Gamma [\cX_\Gamma^{\bot}]\,$$ 
where $\cX_\Gamma^{\bot}$ is the vertical vanishing divisor over the locus
with level graph $\Gamma$. We compute
\bas
\kappa_1 &\=  \pi_{*}  \Big(\Big( \pi^{*}\xi + \sum_{i=1}^n m_i S_i
+ \sum_{\Gamma\in \LG_1} \ell_\Gamma [\cX_\Gamma^{\bot}]\Big)^2\Big) \\
&\=  -\sum_{i=1}^n m_i^2 \psi_i +  \pi_{*} \Big(\Big(\sum_{\Gamma\in \LG_1} \ell_\Gamma [\cX_\Gamma^{\bot}]\Big)^2\Big) + (4g-4) \xi + 2\sum_{i=1}^n \sum_{\Gamma\in {}_i{\LG_1}} m_i \ell_\Gamma [D_\Gamma] \\
&\=   (4g-4) \xi -\sum_{i=1}^n m_i^2 \psi_i +  2\sum_{i=1}^n \sum_{\Gamma\in {}_i{\LG_1}} m_i \ell_\Gamma [D_\Gamma]  + \sum_{\Gamma\in \LG_1} \ell_\Gamma^2 \pi_{*} ([\cX_\Gamma^{\bot}]^2)\,
\eas
where ${}_i{\LG_1}$ means the $i$-th marked point is in lower level.
\par 
Next we evaluate 
$\pi_{*}([\cX_\Gamma^{\bot}]^2)$. Since $[\cX_\Gamma^{\bot}] + [\cX_\Gamma^{\top}] = \pi^{*} [D_\Gamma]$, it suffices to evaluate $\pi_{*} ([\cX_\Gamma^{\bot}]\cdot [\cX_\Gamma^{\top}])$, which is a class supported on $D_\Gamma$ with suitable multiplicity. Geometrically $\cX_\Gamma^{\bot}$ and $\cX_\Gamma^{\top}$ intersect along the (vertical) edges of $\Gamma$. For $e \in E(\Gamma)$, to figure out its contribution to the multiplicity, it suffices to take a general one-parameter family $C$ crossing through $D_\Gamma$.  The local singularity type of $\cX|_C$ at $e$ is $\ell_\Gamma / p_e$, i.e., locally the corresponding node is defined by 
$xy = t^{\ell_\Gamma / p_e}$ where $t$ is the base parameter. Hence the local contribution is $p_e / \ell_\Gamma$ (as the reciprocal of the exponent of $t$). We conclude that 
 $$\pi_{*}([\cX_\Gamma^{\bot}]^2) \=  - \pi_{*}([\cX_\Gamma^{\bot}]\cdot [\cX_\Gamma^{\top}])
 \=  - \frac{1}{\ell_\Gamma}\Big(\sum_{e\in E(\Gamma)} p_e\Big)[D_\Gamma]\,. $$  
Moreover by the relation at the beginning we have 
 \be \label{eq:xitopsi}
 \xi \= (m_i+1)\psi_i - \sum_{\Gamma\in _i\LG_1} \ell_\Gamma [D_\Gamma]
 \ee
(see e.g.,~\cite[Proposition 8.2]{EC}) and consequently  
$$ (2g-2) \xi \= \sum_{i=1}^n (m_i^2 + m_i) \psi_i -  \sum_{i=1}^n \sum_{\Gamma\in {}_i{\LG_1}} m_i \ell_\Gamma [D_\Gamma] \,.$$
 In particular, $\psi_i$ can be converted to $\xi$ with boundary classes as long as $m_i \neq -1$, and $\xi = - \sum_{\Gamma\in _i\LG_1} \ell_\Gamma [D_\Gamma]$ for $m_i = -1$. It follows that 
\begin{align*}
\kappa_1 + \psi_{-1} 
&\=  (2g-2+s)\xi + \sum_{m_i\neq -1} m_i \psi_i + \sum_{\Gamma\in \LG_1} \ell_\Gamma 
\Big( \sum_{i \in L(\Gamma)\atop m_i\neq -1} m_i  - \sum_{e\in E(\Gamma)} p_e\Big) [D_\Gamma] \\
&\=  \kappa_\mu \xi + \sum_{ m_i\neq -1} \frac{m_i}{m_i+1} \sum_{\Gamma\in _i\LG_1} \ell_\Gamma [D_\Gamma] + \sum_{\Gamma\in \LG_1} \ell_\Gamma 
\Big( \sum_{i \in L(\Gamma)\atop m_i\neq -1} m_i  - \sum_{e\in E(\Gamma)} p_e\Big) [D_\Gamma] \\
&\= \kappa_\mu \xi + \sum_{\Gamma \in \LG_1} \ell_\Gamma \Big( \sum_{i\in L(\Gamma)\atop m_i\neq -1} \frac{m_i^2+2m_i}{m_i+1} -   \sum_{e\in E(\Gamma)} p_e    \Big)[D_\Gamma] \\
&\=  \kappa_\mu \xi + \sum_{\Gamma \in \LG_1} \ell_\Gamma \Big(\kbot - \sum_{e\in E(\Gamma)} \frac{1}{p_e}\Big) [D_\Gamma]\,
\end{align*}
as claimed, where $L(\Gamma)$ denotes the lower level of $\Gamma$ and in the last step we used 
that a vertical edge $e$ with pole order $-p_e - 1$ contributes $-p_e + 1/p_e$ in $\kbot$.   
\end{proof}
\par
Next we compute the first Chern class of the Hodge bundle~$\lambda_1$ as shown in \autoref{prop:xiinlambda}. For that purpose we need to pull back boundary divisor classes from the moduli space of curves to the moduli space of multi-scale differentials.   
Let $\delta_{[n]}$ be the total boundary divisor class in $\barmoduli[{g,n}]$. Denote by $f_{[n]}\colon \bP\LMS \to \barmoduli[g,n]$ the natural map remembering the underlying pointed stable curves only.  
\par
\begin{lemma} \label{le:deltapullback}
The pullback of $\delta_{[n]}$ to the moduli space of multi-scale differentials has divisor 
class 
\bes
f_{[n]}^*(\delta_{[n]}) \= [D_h] + \sum_{\Gamma\in \LG_1} \ell_\Gamma \Big(\sum_{e\in E(\Gamma)}
\frac{1}{p_e}\Big)[D_\Gamma]\,.  
\ees
\end{lemma}
\par
\begin{proof}  The local equation of the universal curve over $\barmoduli[{g,n}]$ is $xy=t$ for the node defining each boundary divisor. In the family over~$\bP\LMS$ the local equation is
$xy=t$ for $D_h$ by the horizontal plumbing \cite[Equation~(12.6)]{LMS}
and $xy = t^{\ell_\Gamma / p_e}$ for each (vertical) edge~$e$ in
a graph~$\Gamma\in \LG_1$ by Equation~(12.8) in loc.\ cit. 
We thus obtain the sum on the right-hand side of
the desired equation and each summand contributes $\ell_\Gamma / p_e$. 
\end{proof}
\par
\begin{proof}[Proof of \autoref{prop:xiinlambda}]
Recall the well-known relation $12\lambda_1 = \kappa_1 + \delta_{[n]}$ on $\barmoduli[{g,n}]$ (see e.g.~\cite[Chapter XIII, Equation~(7.7)]{acgh2}). Combining it with 
\autoref{le:deltapullback} and \autoref{prop:kappexpr} we thus conclude the desired formula. 
\end{proof}
\par
For convenience of later calculations we combine \autoref{prop:xiinlambda}  and~\autoref{eq:xitopsi} to get
\be \label{eq:psitolambda}
\psi_i \=\frac{1}{\kappa_\mu(m_i +1)} \Bigl(
\psi_{-1} + 12 \lambda_1 - [D_h] -\sum_{\Gamma \in \LG_1} \bigl(\kappa_{\mu_\Gamma^\bot} - \delta_{i,\bot} \kappa_\mu
\bigr) \ell_\Gamma [D_\Gamma]
\Bigr)
\ee
if $m_i \neq -1$, where we define the 'Kronecker' symbol $\delta_{i,\bot}$ to be 1 if
the $i$-th leg is on bottom level and zero otherwise.

\subsection{The canonical bundle formula}
\label{sec:spinhyp}

Our goal here is to prove the canonical bundle formula in
\autoref{prop:canonicaloncoarse}. We continue to use $\mathbf{D}_\Gamma$,
$\mathbf{D}_h$, etc, for the reduced boundary divisors of the coarse moduli space $\bP\MScoarse$, to
distinguish from the divisors $D_\Gamma$, $D_h$, etc, in the stack $\bP\LMS$. We use ${\bf D}$ and $D$ to denote respectively the total boundary divisors.  
We will express the canonical classes in standard generators in $\CH^1(\bP\LMS)$, both for the stack and the
coarse moduli space, so that even for the latter no boldface objects appear
in the final formulas. We start with:
\par
\begin{prop} \label{prop:canonicalonstack}
The class of the log canonical bundle of the smooth Deligne--Mumford stack $\bP\LMS$ is
given by 
\bas
\frac{\kappa_\mu}{N} \c_1\Bigl(\Omega^{d}_{\bP\LMS}(\log {D})\Bigr) \= \psi_{-1} + 
12 \lambda_1 - [D_h] + \sum_{\Gamma\in \LG_1} \ell_\Gamma
\Big(\kappa_\mu \frac{N_\Gamma^{\bot}}{N} -
\kappa_{\mu_{\Gamma}^{\bot}}\Big) [D_\Gamma]
\eas
in $\CH^1(\bP\LMS)$, where $d = \dim \bP\LMS$ and $N_\Gamma^{\bot}$ is the dimension of the unprojectivized bottom level stratum in $D_\Gamma$. 
\end{prop}
\par
\begin{proof} This follows from \cite[Theorem~1.1]{EC} and~\autoref{eq:eta-xi}.
\end{proof}
\par
Next we pass to the coarse moduli space, still in the logarithmic context:
\par
\begin{prop} \label{prop:logcanonicaloncoarse}
Let $\mu$ be a holomorphic signature not of type $(2m, 2g-2-2m)$ (and the hyperelliptic component is excluded if $\mu = (2g-2)$), 
or a meromorphic signature not of type $(2m_1, 2m_2, 2g-2-2m_1 - 2m_2)$ (and the hyperelliptic component is excluded if $\mu = (2g-2+2m, -2m)$). Then 
the class of the reflexive log canonical bundle of the coarse moduli space
$\bP\MScoarse$  is given by 
\bas
\frac{\kappa_\mu}{N} \c_1\Bigl(\Omega^{[d]}_{\bP\MScoarse}(\log \mathbf{D})\Bigr) \= \psi_{-1} + 
12 \lambda_1 - [D_h] + \sum_{\Gamma\in \LG_1} \ell_\Gamma
\Big(\kappa_\mu \frac{N_\Gamma^{\bot}}{N} - \kappa_{\mu_{\Gamma}^{\bot}} \Big) [D_\Gamma]
\eas
in $\CH^1(\bP\MScoarse) = \CH^1(\bP\LMS)$, where $d = \dim \bP\MScoarse$. 
\end{prop}
\par
\begin{proof}[Proof of \autoref{prop:logcanonicaloncoarse} and
\autoref{prop:canonicaloncoarse}]
By the assumption~\autoref{prop:ramifInt} says that the map 
$\varphi\colon \bP\LMS \to \bP\MScoarse$ is unramified in the interior.  
For any linear combination of boundary divisors $\sum a_i D_i$ in the stack, suppose $\varphi$ is ramified with order~$e_i$ at~$D_i$ and with image $\mathbf{D}_i$ in the coarse moduli space. Then the ramification
formula passing from the stack to the coarse moduli space (e.g.\ \cite[Proposition~A.13]{HHlogcan}) gives an equality of $\bQ$-Cartier
divisors
\be
\varphi^*\Bigl(K_{\bP\MScoarse} + \sum_i \frac{e_i - 1 + a_i}{e_i} \mathbf{D}_i\Bigr)
\= K_{\bP\LMS} + \sum_i a_i D_i\,.
\ee
Setting all $a_i = 1$, i.e., taking the combination $D_h + \sum_\Gamma D_\Gamma$, \autoref{prop:logcanonicaloncoarse} 
thus follows from the above 
together with \autoref{prop:canonicalonstack}.  Next we can set $a_i=0$
for unramified boundary divisors and $a_i = -1$ for ramified boundary divisors (of order two). Then 
\autoref{prop:canonicaloncoarse} follows from 
\autoref{prop:ramifBound} and \autoref{prop:ramifBoundmero}. 
\end{proof}
\par
Next we explain how to apply the above formulas to each connected component of a disconnected stratum. 
Recall that the connected components of strata of holomorphic differentials have
been classified by Kontsevich--Zorich (\cite{kozo1}). For special signatures these connected components are distinguished by the parity of the spin
structure, odd or even, and by consisting entirely of hyperelliptic differentials. 
We denote these components by an upper index odd,
even or hyp respectively. Note that if a holomorphic stratum has three components, then the hyperelliptic component has a fixed spin parity
depending on~$g$. We emphasize that the superscript odd or even excludes
the hyperelliptic component with that spin parity. The smooth
compactification $\bP\LMS$ still separates these components and we distinguish
them by the same superscripts. In this case we also add the same superscript to
decompose a boundary divisor $D_\Gamma$ or $D_h$, e.g.,
\be
D_\Gamma \= D_\Gamma^{\hyp} \cup  D_\Gamma^{\odd} \cup D_\Gamma^{\even}
\quad \text{with} \quad
D_\Gamma^\bullet \,\subset\, \bP\LMS^\bullet
\ee
as a disjoint union, where $\bullet \in \{\hyp, \odd, \even\}$. Similar
decompositions occur if there are only two components, odd and even, 
or hyperelliptic and non-hyperelliptic.
\par
All the steps of the proof of \cite[Theorem~1.1]{EC} can be performed on
each connected component separately. We thus deduce that the component-wise
version of \autoref{prop:canonicalonstack}
\ba \label{eq:canclassstackspin}
\frac{\kappa_\mu}{N} \c_1\Bigl(\Omega^{d}_{\bP\LMS^\bullet}(\log {D})\Bigr) \= \psi_{-1} + 
12 \lambda_1 - [D_h^{\bullet}] + \sum_{\Gamma\in \LG_1} \ell_\Gamma
\Big(\kappa_\mu \frac{N_\Gamma^{\bot}}{N}
- \kappa_{\mu_{\Gamma}^{\bot}}\Big) [D_\Gamma^\bullet]
\ea
holds with $\bullet \in \{\hyp, \odd, \even\}$ for holomorphic signatures $\mu$. As before this formula
can be converted to 
\ba \label{eq:canclasscoarsespin}
\frac{\kappa_\mu}{N} &\c_1\bigl(K_{\bP\MScoarse^\bullet}\bigr) 
\= \psi_{-1} + 12 \lambda_1 - \Bigl(1 + \frac{\kappa_\mu}{N}\Bigr)[D_h^{\bullet}] \\
&+ \sum_{\Gamma\in \LG_1} \Big(\frac{\kappa_\mu}{N} (\ell_\Gamma N_\Gamma^{\bot}-1)
- \ell_\Gamma  \kappa_{\mu_{\Gamma}^{\bot}} \Big) [D_\Gamma^\bullet]
- \frac{\kappa_\mu}{N} \sum_{\Gamma \, \text{is} \, \mathrm{HTB}\, \text{or} \atop
 \mathrm{HBT} \,\text{or}\, \mathrm{HBB}}
 [D^{H,\bullet}_\Gamma]
\ea
in $\CH^1(\bP\MScoarse) = \CH^1(\bP\LMS)$ for the (non-logarithmic) canonical 
bundle. Here $D^{H,\bullet}_\Gamma$ is the (possibly empty) component of $D^{\rm H}_\Gamma$
in the component indicated by $\bullet \in \{\hyp, \odd, \even\}$. We remark that $D_\Gamma^{\bullet}$ and $D^{H,\bullet}_\Gamma$ can be further reducible due to connected components of the strata in each level of $\Gamma$ and also due to the equivalence classes of prong-matchings of the multi-scale differentials they encode. However for our purpose we do not need to classify their irreducible components.  
\par 
The connected components of strata of meromorphic differentials have been classified by Boissy (\cite{boissymero}). These connected components are similarly distinguished by the spin and hyperelliptic structures, with the only exception in genus one where the distinction is given by divisors of the gcd of the entries in~$\mu$ (also called the rotation or torsion numbers, see~\cite{ChenEED} and~\cite{ChenGendron}). We can analogously add the corresponding superscripts in order to apply \autoref{prop:canonicalonstack} and \autoref{prop:canonicaloncoarse} to each connected component of a meromorphic stratum. 

\subsection{Divisors of Brill--Noether type}
\label{sub:BN}

Recall that two kinds of effective divisors in $\barmoduli[g]$ were used by 
Harris and Mumford in \cite{HarrisMumford} and \cite{HarrisKodII} 
which depend on the parity of~$g$. We proceed similarly. First for $g$ odd,
consider the admissible cover compactification of the locus
\ba \label{eq:defBN}
\wt{BN}_g & \=  \{ X \in \moduli[g]\,: X \, \text{has a} \,\,\frakg^1_k\},
\qquad (k=(g+1)/2)\\
&\= \{ X \in \moduli[g]\,: \text{there is a cover}\,\,
\pi\colon X \to \bP^1, \, \deg(\pi) = k\}\,.
\ea
This is a classical \emph{Brill--Noether divisor} by the description
using linear series. We normalize the class of the Brill--Noether divisor computed
in \cite{HarrisMumford} to be 
$$ [\BN_g] \=  6 \lambda_1 - \frac{g+1}{g+3} \delta_{\irr} -
\sum_{i=1}^{[g/2]} \frac{6 i (g-i)}{g+3} \delta_i \quad \in \CH^1(\barmoduli[g])\,.$$
\par
We define the pullback of the Brill--Noether divisor $\BN_\mu$
(as a $\bQ$-divisor) to be the total transform 
$[\BN_\mu] = f^*[\BN_g]$ where $f\colon \bP\LMS \to\barmoduli[g]$ is the forgetful map. To state the class 
of the pullback Brill--Noether divisor we use the following notation.
For an edge~$e$ in a level graph~$\Gamma$ we write $e \mapsto  \Delta_i$ if
contracting all edges of~$\Gamma$ but~$e$ results in a graph of compact type 
parametrized by the boundary divisor~$\Delta_i$ in $\barmoduli[g]$. We write $e \mapsto
\Delta_{\irr}$ if the edge is non-separating, equivalently if the contraction
results in a graph parameterizing irreducible one-nodal
curves.
\par
\begin{lemma} \label{lem:BNclass} Let $g\geq 3$ be odd. 
If the stratum $\bP\LMS$ is connected then the class of the Brill--Noether divisor
in  $\bP\LMS$ 
$$ [\BN_\mu] \= 6\lambda_1 - \frac{g+1}{g+3} [D_h] - \sum_{\Gamma\in \LG_1}
b_\Gamma [D_\Gamma]\, $$
where 
$$b_\Gamma \= \ell_\Gamma \Big(\sum_{i=1}^{[g/2]}\sum_{e\in E(\Gamma) \atop e \mapsto \Delta_i} \frac{6i(g-i)}{(g+3) p_e} +  
\sum_{e\in E(\Gamma) \atop e \mapsto \Delta_{\irr}} \frac{g+1}{(g+3) p_e} \Big)\, $$
is an effective divisor class.
\par
If the stratum is disconnected and $\bullet \in \{\odd, \even\}$ denotes
a non-hyperelliptic component then similarly
\bes
[\BN_\mu^\bullet] \= 6\lambda_1 - \frac{g+1}{g+3} [D_h^\bullet]
- \sum_{\Gamma\in \LG_1} b_\Gamma [D_\Gamma^\bullet]
\ees
is an effective divisor class.
\end{lemma}
\par
\begin{proof}
It follows from~\cite[Theorem 1.1]{BudGonality} that for odd~$g$ and any
non-hyper\-elliptic connected component of a stratum $\bP\LMS$ not every curve
parameterized therein admits a ${\frakg}^1_{k}$ for $k = (g+1)/2$.\footnote{Strictly speaking \cite{BudGonality} only considered the strata of holomorphic differentials. However the same argument works for the meromorphic case as well by merging all zeros and poles and specializing to the minimal strata, with the exception for $\mu = (2g-2+m, -m)$ with $m > 1$ where the zero and pole cannot be merged due to the GRC. For the exceptional case one can still argue as in loc. cit. by taking a Brill--Noether general curve in $\bP\Omega\cM_{g-1}(2g-4)^{{\rm nonhyp}}$ union an elliptic tail in $\bP\Omega\cM_{1}(2-2g, 2g+2-m, -m)$. } This implies that
the pullback of the Brill--Noether divisor to such strata (components) is an effective divisor.
To compute ($c$ times) its class we can perform the same local computation as in the proof of 
\autoref{le:deltapullback}. 
\end{proof}
We remark that for special $\Gamma$, e.g., if the underlying curves parameterized by a vertex of $\Gamma$ have gonality much smaller than expected, then the total transform $\BN_\mu$ can contain the corresponding boundary divisor $D_\Gamma$, which can be subtracted from $\BN_\mu$ to make the remaining effective divisor class more extremal.
\par
\medskip
Second, to cover $g$ even we need two other types of effective divisors.
For $g$ even Harris used in \cite{HarrisKodII} the closure of 
\ba \label{eq:defHur}
\wt{\Hur}_g &\= \{ X \in \moduli[g]\,: \text{there is a cover}\,\,
\pi\colon X \to \bP^1, \, \deg(\pi) = (g+2)/2, \\
&  \qquad \qquad \qquad \qquad
\text{$\pi$ has a point of multiplicity three}\,\}\,
\ea
where multiplicity being $m$ means locally the cover is given by $x \mapsto x^m$, i.e., being of ramification order $m-1$. We normalize the class of the Hurwitz divisor computed in loc.~cit. to be 
$$  [\Hur_g] \=  6 \lambda_1 - \frac{3g^2 + 12g -6}{(g+8)(3g-1)}
\delta_{\irr} - \sum_{i=1}^{g/2} \frac{6 i (g-i)(3g+4)}
{(g+8)(3g-1)} \delta_i\,. $$
Similarly to the case above we let $[\Hur_\mu] = f^*[\Hur_g]$ be the
pullback {\em Hurwitz divisor class}. 
\par
\begin{lemma} \label{lem:Hurclass}
For even $g \geq 6$, the Hurwitz divisor $\Hur_\mu$ is an effective divisor in every connected $\bP\LMS$ and in every non-hyperelliptic component of disconnected $\bP\LMS$. For $g=4$, $\Hur_\mu$ is an effective divisor in every connected $\bP\LMS$ and in every odd spin  component of disconnected $\bP\LMS$. 
\par 
Moreover, the class of the Hurwitz divisor in a connected  
$\bP\LMS$ is  
$$ [\Hur_\mu] \= 6\lambda_1 -  \frac{3g^2 + 12g -6}{(g+8)(3g-1)} [D_h] - \sum_{\Gamma\in \LG_1}
h_\Gamma [D_\Gamma]\, $$
where 
$$h_\Gamma \= \ell_\Gamma \Big(\sum_{i=1}^{g/2}\sum_{e\in E(\Gamma) \atop e \mapsto \Delta_i} \frac{6 i (g-i)(3g+4)}{(g+8)(3g-1) p_e} +  
\sum_{e\in E(\Gamma) \atop e \mapsto \Delta_{\irr}} \frac{3g^2 + 12g -6}{(g+8)(3g-1) p_e} \Big)\, $$
and the same expression holds for the spin components of $\bP\LMS$ with $D_\Gamma$ decorated by $\bullet \in \{ \odd, \even\}$ respectively. 
\end{lemma}
\par
\begin{proof}
 The divisor class calculation is the same as in the proof of \autoref{le:deltapullback}. In order to prove that no strata (components) are  contained in $\Hur_\mu$ in the claimed range, by merging zeros and poles it suffices to prove it for the minimal strata, i.e., for $(2g-2)^{\odd}$ with even $g \geq 4$ and for $(2g-2)^{\even}$ with even $g\geq 6$ (with the only exception for $\mu = (2g-2+m, -m)$ with $m > 1$ since in this case the zero and pole cannot be merged due to the GRC, which we will treat separately at the end). We will exhibit a (boundary) point in each case that is not contained in $\Hur_\mu$.   
\par
Consider first the odd spin minimal strata. Take a multi-scale differential $(X, \bfomega)$ consisting of an elliptic curve $(E, p_1, \ldots, p_g)\in \bP\Omega\cM_1(0^g)$ union a rational curve $(R, q_1, \ldots, q_g, z) \in \bP\Omega\cM_0(-2^g, 2g-2)$ at $g$ nodes by identifying $p_i \sim q_i$ for all $i$. The GRC is automatically satisfied by the Residue Theorem on $R$. Since all prongs are one at the nodes, the spin parity of $(X, \bfomega)$ (via the Arf invariant) is equal to the parity of the flat torus $E$ which is odd. Therefore, $(X, \bfomega)$ is a boundary point of the odd spin minimal stratum. We claim that for general positions of $p_i$ and $q_j$ in $E$ and $R$, the union $X$ is not contained in the Hurwitz divisor $\Hur_g$ for even $g\geq 4$. To see it, consider an admissible cover of degree $g/2 + 1$ on $X$. Since the number of nodes between $E$ and $R$ is $g > g/2+1$, $E$ and $R$ must map to the same target $\bP^1$-component $C$, such that each pair $p_i$ and $q_i$ is contained in the same fiber.  If all ramification points in $E$ and $R$ are simple, by Riemann--Hurwitz their total number is $g$. Hence the total number of distinct branch points and image points of $p_i, q_i$ in $C$ is (at most) $2g$. However if the admissible cover has a triple point, then this number drops to (at most) $2g-1$. In other words, the parameter space of such admissible covers with a triple point restricted to $E$ and $R$ has dimension bounded by $\dim \cM_{0,2g-1} = 2g-4$. On the other hand, the parameter space of $E$ union $R$ is $\cM_{1,g}\times \cM_{0,g}$ with dimension $2g-3 > 2g-4$. We thus conclude that a general union of $E$ and $R$ does not admit a cover of degree $g/2+1$ with a triple point. 
\par 
Next for the even spin minimal strata, take a multi-scale differential $(X, \bfomega)$ consisting of three components,  
a rational curve $(R, q_1, \ldots, q_g, z) \in \bP\Omega\cM_0(-2^g, 2g-2)$, an elliptic curve $(E, p_1, \ldots, p_{g-1})\in \bP\Omega\cM_1(0^{g-1})$ and another elliptic curve $(E', p_g) \in \bP\Omega\cM_1(0)$ by identifying $p_i \sim q_i$ for all $i$. The GRC requires that the residue of $\omega_R$ at $q_g$ is zero, which can be satisfied by choosing a special position of $z$ with respect to a general choice of $q_1, \ldots, q_g$ in $R$. Since all prongs are one at the nodes, the spin parity of $(X, \bfomega)$ is equal to the sum of the parities of the flat tori $E$ and $E'$ which is even.   Therefore, $(X, \bfomega)$ is a boundary point of the even spin minimal stratum. Since the number of nodes between $E$ and $R$ is $g-1 > g/2+1$ for $g \geq 6$, one can argue similarly as in the preceding paragraph to show that a general such $X$ is not contained in the Hurwitz divisor $\Hur_g$ for even $g\geq 6$. 
\par
For the exceptional case $\mu = (2g-2+m, -m)$ with $m > 1$, we can put the unique zero and pole in the rational component $R$ in the above constructions, and the same arguments still work through.  
\par
We remark that the discrepancy of the genus bounds for the odd and even spin cases is necessary, because every $(X, z)\in \bP\Omega\cM_4(6)^{\even}$
admits a triple cover by the linear system $|3z|$ with a triple point at $z$, and hence the even spin minimal stratum in genus four maps entirely into the Hurwitz divisor.  
\end{proof}
In order to show the statement of \autoref{intro:minimal} for even genera minimal strata, the Hurwitz divisor, whose class was computed in \autoref{lem:Hurclass}, is sufficient for all genera apart from $g=14$. For this special case,  in order to show that  $\bP\MScoarse[26]$ is of general type, we  will need to use a pointed Brill--Noether divisor which was studied in the proof
of \cite[Theorem 4.9]{farkaskoszul}. We let  $d = g/2+1$ and define the {\em divisor of $n$-fold points}
(here for $n=2$) as the closure in $\barmoduli[g,1]$ of 
\bas
\NF_{g,1} &\= {\{(X,p) \in \moduli[g,1]: \exists \cL \in W_d^1(X): h^0(X,\cL(-2p)) \geq 1\}}\,. 
\eas
In loc.~cit.\ the class of this divisor was computed to be\footnote{The factor $6/(g+3)$ compensates the different normalizations of the Brill--Noether divisor class here and in~\cite{farkaskoszul}.}
\begin{align}
 \label{eq:NFtoBNW}
[\NF_{g,1}] & \=   \frac{12}{(g+1)(g-2)}{g-2 \choose g/2} \frac{g+3}6 [\BN_g]
+ \frac{6}{g(g-1)(g+1)} {g \choose g/2 + 1} [{\rm W}] \nonumber \\
&\= \frac{2(g+3)}{(g+1)(g-2)}{g-2
\choose g/2} \cdot  \Bigl([\BN_g]+\frac{12}{(g+3)(g+2)}[{\rm W}] \Bigr)
\end{align}
where $[{\rm BN}_g]$ is the pullback of the (normalized) Brill--Noether divisor class from $\barmoduli$ 
and $\rm W$ is the Weierstrass point divisor. The class of ${\rm W}$ in $\barmoduli[g,1]$ is given by 
$$ [{\rm W}] = -\lambda_1 + {g+1\choose 2} \psi -  \sum_{j=1}^{g-1} {g-j+1\choose 2}
\delta_{j;1}\, $$
where $\delta_{j;1}$ is the class of the boundary divisor parameterizing curves of compact type with a component 
of genus~$j$ that carries the marked point and another unmarked component of genus $g-j$. 
\par
We consider the class given by the rescaling of the pullback  $f_1^*\NF_{g,1}$  having the same $\lambda$-coefficient~$6$ as $\BN_\mu$, where  $f_1\colon \bP\Xi \barmoduli[g,1](2g-2) \to\barmoduli[g,1]$ is the forgetful map. 
\par
\begin{lemma} \label{le:NFBNdiff}
Consider even $g\geq 4$ and let $\NF_{(2g-2)} = e^{-1} \cdot f_1^* \NF_{g,1}$ where
$e=\frac{2(g^2 + 3g - 1)} {(g-1)(g^2-4)}{g-2 \choose g/2}$. 
Then the class of $\NF_{(2g-2)}$ is effective for every non-hyperelliptic component
of the stratum $\bP\Xi \barmoduli[g,1](2g-2)$ for even $g\geq 6$ and for the odd
spin component $\bP\Xi \barmoduli[4,1](6)^{\odd}$ in $g=4$. 
\end{lemma}
\par
\begin{proof}
To prove the claim of effectiveness, we need to show that a general $(X,z) \in \bP\Xi \barmoduli[g,1](2g-2)^{{\rm nonhyp}}$ does not admit a $\frak g^1_{d}$ with multiplicity at least two at $z$ where $d = (g+2)/2$. This was indeed verified in~\cite[Proposition 3.1]{BudGonality}. Again the discrepancy of the genus bounds for the odd and even spin cases is necessary, as every $(X, z)\in \bP\Xi \barmoduli[4,1](6)^{\even}$ admits a triple cover induced by the linear system $|3z|$ which is ramified at $z$. 
\end{proof}
The class of $\NF_{(2g-2)} $ can be explicitly computed using~\autoref{eq:NFtoBNW} and \autoref{eq:psitolambda}. 
\par

\section{Generalized Weierstrass divisors}
\label{sec:genWP}

There are few effective divisors that are directly defined in the strata of
Abelian differentials other than those arising as a pullback from $\barmoduli[g,n]$.
In this section we directly construct a series of such divisors and compute their classes.
These divisors generalize the classical divisor of Weierstrass points
in $\barmoduli[g,1]$. Despite that they can be defined for both holomorphic
and meromorphic signatures, since in this paper we only apply them to certain
holomorphic strata, we limit their definition to the {\em holomorphic} case,
and leave the meromorphic case to future work.  
\par
In the sequel we will mainly work with the 'middle' case of the
generalized Weierstrass divisors, where we use as weights $\mu/2$. Since
this tuple is not always integral, the generalized Weierstrass divisor class
associated to it is not obviously effective. In \autoref{sec:average} we
discuss the quality of approximation by actual effective classes for generalized
Weierstrass divisors given by rounding $\mu/2$.
\par
Finally, with all divisors in place,  we give in \autoref{sec:refofProp13} a
refinement of the general type criterion from \autoref{prop:GenTypeCrit}.

\subsection{The divisor class}
Fix a holomorphic signature $\mu = (m_1, \ldots, m_n)$, i.e., with $m_i \geq 0$.
If the corresponding stratum $\bP \omoduli[g,n](\mu)$ is not connected,
we suppose moreover that in this section $\bP\LMS$ denotes the multi-scale compactification for the connected
component of {\em odd spin}. All the other components (even spin and hyperelliptic)
are disregarded for the construction of generalized Weierstrass divisors.  
\par
Consider a partition $\alpha = (\alpha_1, \ldots, \alpha_n)$ of $g-1$ such that
$0\leq \alpha_i \leq m_i$ for all~$i$. Set-theoretically,  the
\emph{generalized Weierstrass divisor} associated to $\alpha$ in the interior
of a stratum of type~$\mu$ is given by
\be \label{eq:initdefGenW}
W_{\mu}(\alpha) \= \Bigl\{ (X, \bfz, \omega) \in \bP\omoduli[g,n](\mu)\,:\,
h^0\Bigl(X, \sum_{i=1}^n \alpha_i z_i\Bigr) \geq 2\Bigr \}\,.
\ee
By Riemann--Roch and duality we deduce that 
$$h^0\Big(X, \sum_{i=1}^n \alpha_iz_i\Big)
\= h^0\Big(X, \omega_X - \sum_{i=1}^n \alpha_iz_i\Big)
\= h^0\Big(X, \sum_{i=1}^n (m_i - \alpha_i)z_i\Big)$$ 
from the signature of the stratum, where $\omega_X$ denotes the dualizing line
bundle of~$X$. Geometrically speaking, the tautological section $\omega$ gives a
section of the linear system $|\omega_X - \sum_{i=1}^n \alpha_i z_i|$ in the definition
of~$W_{\mu}(\alpha)$, and we thus consider the locus where the linear system has
extra sections. This viewpoint can realize 
the generalized Weierstrass divisor as a degeneracy locus and provide it 
a scheme structure as follows. 
\par
Let  $\omega_{\rel}$  
be the relative dualizing bundle of the universal curve $\pi\colon \cX \to \bP\LMS$ 
and the $S_i\subset \cX$ are the marked sections. Let moreover
$\cH=\pi_*(\omega_{\rel})$ be the Hodge bundle over $\bP\LMS$ and $\cO(-1)$ its
tautological subbundle. Let $\cF_\alpha$ be the bundle whose fiber over
$(X, \bfz, \bfomega)$ is canonically given by
$H^0\left(\omega_X / \omega_X(-\sum_{i=1}^n \alpha_i z_i)\right)$
i.e., 
$$\cF_\alpha \= \pi_*\Big(\omega_{\rel}/\omega_{\rel}\Big(-\sum_{i=1}^n
\alpha_i S_i\Big)\Big).$$ 
Both $\cH/\cO(-1)$ and $\cF_\alpha$ are vector bundles of rank~$g-1$ (for $\cF_\alpha$
this follows from the long exact sequence associated with the inclusion
$\omega_{\rel}(-\sum_{i=1}^n \alpha_i S_i)
\hookrightarrow \omega_{\rel}$ or from the interpretation as the sheaf of
principal parts of length~$\alpha_i$ at the points~$z_i$).
Taking principal parts to order~$\alpha_i$ at each of the
points~$z_i$ defines a bundle map $\cH \to \cF_\alpha$, which fiberwise 
is the map $H^0(X,\omega_X) \to H^0(X,\omega_X /\omega_X(-\sum_{i=1}^n
\alpha_i z_i))$. Since $\omega$ has vanishing order (at least) $m_i \geq \alpha_i$
at each $z_i$, 
the above bundle map factors through the following bundle map
\be \label{eq:phiuntwisted}
\phi\colon \cH/\cO(-1) \to \cF_\alpha\,.
\ee 
We define the substack $\wt{W}_{\mu}(\alpha)$ to be the degeneracy
locus of~$\phi$, i.e.,\ the \emph{generalized Weierstrass divisor} 
\be
\wt{W}_\mu(\alpha) \= \{ \mathrm{rank}(\phi) < g-1\} \subset \bP\LMS\,.
\ee
To express its divisor class we define 
\be
\label{eq:vartheta}
\vartheta \,:=\, \vartheta_{\mu, \alpha}
\= \sum_{i=1}^n \frac{\alpha_i(\alpha_i+1)}{2(m_i+1)}
\ee
for $\mu = (m_1, \ldots, m_n)$ and $\alpha = (\alpha_1, \ldots, \alpha_n)$.
The quantities $\vartheta^\bot := \vartheta_{\mu_\Gamma^{\bot}, \alpha_\Gamma^{\bot}}$ and
$\vartheta^\top := \vartheta_{\mu_\Gamma^{\top},
\alpha_\Gamma^{\top}}$ are similarly defined for the bottom  and top level
strata of $\Gamma$, \emph{but $\alpha_\Gamma^{\bot}$ and $\alpha_\Gamma^{\top}$ assigns
value zero to each leg associated with an edge $e\in E(\Gamma)$.}\footnote{This
is contrary to how edges are treated in the definition of $\kappa_{\mu_\Gamma^\bot}$
in \autoref{sec:taut}.}  With this convention we also have
$$\vartheta_{\mu_\Gamma^{\bot}, \alpha_\Gamma^{\bot}}
+ \vartheta_{\mu_\Gamma^{\top}, \alpha_\Gamma^{\top}} \= \vartheta_{\mu, \alpha}\,.$$
\par
\begin{prop} \label{prop:Wclass}
The substack $\wt{W}_{\mu}(\alpha)$ is an effective divisor in $\bP\LMS$.
The class of this generalized Weierstrass divisor is given by
\bas \ 
[\wt{W}_{\mu}(\alpha)] & \=
\frac{12+12 \vartheta_{\mu, \alpha} - \kappa_\mu}{\kappa_\mu}\lambda_1
- \frac{1+\vartheta_{\mu, \alpha}}{\kappa_\mu} [D_h] \\
& \,-\, \sum_{\Gamma\in\LG_1} \ell_\Gamma \Big( \frac{\kappa_{\mu_\Gamma^{\bot}}}{\kappa_\mu}
(1+\vartheta_{\mu, \alpha}) - \vartheta_{\mu_\Gamma^{\bot}, \alpha_\Gamma^{\bot}}\Big) [D_\Gamma]\,.
\eas 
\end{prop}
\par
The degeneracy locus $\wt{W}_{\mu}(\alpha)$ in general contains extra
boundary divisors and we estimate the boundary contributions in the next subsection.
\par
\begin{proof}[Proof of \autoref{prop:Wclass}]
We first show that $\wt{W}_{\mu}(\alpha)$ is an effective divisor.  The setup as
degeneracy locus implies that it has (local) codimension at most one everywhere.
Hence it suffices to exhibit a (boundary) point not contained in $\wt{W}_{\mu}(\alpha)$.
\par
Consider the boundary divisor~$D_\Gamma$ for~$\Gamma \in \LG_1$ consisting of
curves $X$ of compact type with a differential $(X_0,q,\omega_0) \in
\omoduli[g,1](2g-2)^{\odd}$ on top level, attached via the node~$q$ to a rational tail
$(R,\bfz, q, \omega_1) \in \omoduli[0,n+1](\mu,-2g)$ carrying all the marked
points. This is always possible, including the case of disconnected strata thanks
to our standing odd spin hypothesis (since the rational tail has even spin
and the parity of the spin is additive for compact type divisors).
\par
In a neighborhood of this boundary divisor we consider the (twisted) line bundle
\be
\cK \= \cO_\cX\Bigl(\sum_{i=1}^n \alpha_i S_i + (g-1)\cX^\bot\Bigr)\,
\ee
where $\cX^\bot$ is the lower level component of the universal family~$\cX$.
Since over the interior the bundle $\pi_*(\cK)/\cO(-1)$ is the kernel of the
map $\phi$, it suffices to show that on the fiber~$X$ over a general point
of~$D_\Gamma$, the bundle $\cK$ has a one-dimensional space of sections (i.e., spanned
by the tautological section $\bfomega$ only).
The restriction of $\cK$ to the fiber~$X$ pulled by to its irreducible components
gives the bundle  $K_0 = \cO_{X_0}((g-1)q)$ on~$X_0$ and the degree zero
bundle $K_1 = \cO_R((1-g)q + \sum_i \alpha_i z_i)$ on~$X_1$ (as can be
seen by twisting $\cK$ by $(1-g)$ times a fiber class before restricting).
Note that for general $(X_0, q) \in  \omoduli[g,1](2g-2)^{\odd}$ we have
$h^0(X_0, K_0) = h^0(X_0, K_0(-q)) = 1$ 
(see~\cite{bullocksubcanonical}), which implies that every section of $K_0$ vanishes at $q$.
Moreover $K_1 \cong \cO_{\bP^1}$, hence any section of $K_1$ vanishing at $q$ (in order to
glue with sections of $K_0$) must be identically zero on $R$. We thus conclude that
$h^0(X,\cK|_X) = h^0(X_0, K_0) = 1$. In summary, we have exhibited a (boundary) point not
contained in $\wt{W}_{\mu}(\alpha)$. Consequently $\wt{W}_{\mu}(\alpha)$ is an actual
effective divisor in $\bP\LMS$. 
\par
In view of the next step we compute the first Chern class of $\cF_\alpha$.
Suppose $\alpha_1 > 0$ and define $\alpha'$ by decreasing $\alpha_1$ in~$\alpha$ by~$1$.
Then there is an exact sequence (see \cite[Theorem~11.2~(d)]{EH3264})
\be \label{eq:3264seq}
0 \to  \sigma_1^*(\omega_{\rel})^{\alpha_1} \to
\cF_{\alpha} \to \cF_{\alpha'} \to 0
\ee
where $\sigma_1$ is the map from the base to the section $Z_1$ in the universal
family. It implies that
$$ c_1(\cF_\alpha)  - c_1(\cF_{\alpha'}) \= \alpha_1\psi_{z_1} $$
and thus, proceeding inductively with all marked points, that
\be \label{eq:c1calF}
c_1(\cF_\alpha) \= \sum_{i=1}^n \frac{\alpha_i(\alpha_i+1)}{2}\psi_i\,.
\ee
\par
To compute the class of the generalized Weierstrass divisor we apply 
the Porteous formula and obtain that 
\bas \ 
[\wt{W}_{\mu}(\alpha)] & \=  c_1 (\cF_\alpha) - c_1 (\cH/\cO(-1))  \nonumber \\ 
& \= \sum_{i=1}^n \frac{\alpha_i(\alpha_i+1)}{2}\psi_i - \lambda_1 + \xi  \nonumber \\
& \=  \Big(1 + \sum_{i=1}^n \frac{\alpha_i(\alpha_i+1)}{2(m_i+1)}\Big)\xi - \lambda_1 + 
\sum_{\Gamma\in\LG_1} \ell_\Gamma \sum_{i\in \Gamma^\bot} \frac{\alpha_i(\alpha_i+1)}{2(m_i+1)}
[D_\Gamma]  \nonumber \\
& \=  (1 + \vartheta_{\mu, \alpha})\xi - \lambda_1 + \sum_{\Gamma\in\LG_1}
\ell_\Gamma \vartheta_{\mu_\Gamma^{\bot}, \alpha_\Gamma^{\bot}} [D_\Gamma] \nonumber \\
\eas
using~\autoref{eq:xitopsi} (here $m_i \geq 0$) and the definition
of~$\vartheta_{\mu, \alpha}$. By~\autoref{eq:eta-xi} this agrees with the formula we claimed.
\end{proof}
\par

\subsection{The twisted version}

To reduce the boundary contribution in the generalized Weierstrass divisor we
replace the bundle map~\autoref {eq:phiuntwisted} by a twisted version. For simplicity of notation, in the sequel
we 'pretend' that there is only one boundary divisor~$D_\Gamma$, i.e., we work locally in its neighborhood so that in codimension-one there is no other boundary divisor seen. This will improve the boundary coefficient in \autoref{prop:Wclass} for this particular~$\Gamma$. To obtain the global improvement for all $D_\Gamma$ we can just twist simultaneously by the divisors $V = V_\Gamma$ constructed in the sequel.
\par
We work over a relatively minimal semi-stable model with smooth total space of the
universal family $\cX \to \bP\LMS$ near $D_\Gamma$, i.e., for a node of local type
$xy = t^a$ with $a > 1$, we blow it up by inserting $a-1$ semi-stable rational bridges
to make the resulting new nodes smooth in the universal family. Let~$X_i$ for
$i\in I = I_\Gamma$ be the irreducible components of the vertical divisors (including
the rational bridges) and let $V = \sum_{i\in I} s_i X_i$ be an effective vertical
Cartier divisor supported on some components
of the central fiber, with chosen multiplicities $s_i \geq 0$ on each component. Moreover, we require that $s_i = 0$ for all top level components, $s_i \leq \ell_\Gamma$ for all bottom level components, and $s_i \leq k p_e$ if $X_i$ is the $k$-th rational bridge from the upper end of $e$ to the lower end of $e$ after blowup.  
\par 
We define the twisted relative dualizing line bundle 
\[\cL=\omega_\pi (-V)\]
and the associated bundles on $\bP\LMS$ given by 
\bes
  \cH^\cL=\pi_{*}\cL \quad \text{and} \quad
\cF_{\alpha}^\cL = \pi_{*}\Big(\cL /\cL \Big(-\sum_{i=1}^n \alpha_i Z_i\Big)\Big)\,.
\ees
Consider the evaluation map as before 
\be \label{eq:phitwisted}
 \phi_\cL \colon \cH^\cL / \cO(-1) \to \cF^\cL_{\alpha}\,.
\ee
By construction of the multi-scale space the tautological form $\bfomega$ has vanishing order $\ell_\Gamma$ on $\cX_\bot$ and vanishing order $k p_e$ on the $k$-th rational bridge after blowing up a node with prong $p_e$ (also see~\cite[Section 4]{chendiff} from the twisting viewpoint). The assumption on the range of the twisting coefficients $s_i$ thus ensures that $\cO(-1)$ is a sub-bundle of $ \cH^\cL$, and hence $\phi_\cL$ is well-defined by the assumption on the range of $\alpha_i$. 
\par
We denote by $\wt{W}_\cL$ the degeneracy locus of this map, called the {\em twisted
generalized Weierstrass divisor} associated to the twisted relative dualizing line
bundle $\cL$. Note that $\cH^\cL$ is locally free of rank~$g$ away from the
codimension-two locus where two or more boundary divisors meet. (In fact,
$\cH^\cL$ is torsion-free, and hence locally free over any discrete valuation
ring transverse to the boundary. The claim on the rank follows by considering
the interior. Away from that codimension-two locus we are complex-analytically
locally in a product situation and can apply the DVR-argument.)
Working away from this codimension-two locus is sufficient to compute
divisor classes by using Porteous' formula.
\par
We need a Grothendieck--Riemann--Roch (GRR) computation before we can fully
exploit Porteous' formula. 
\par
\begin{lemma} \label{le:GRR}
Suppose $V = \sum_{i\in I} s_i X_i$ is effective and does not contain an entire fiber, i.e., $s_i \geq 0$ for all $i$ and at least one of $s_i$ is zero. Then $\pi_{*}\cO_{\cX}(V) = \cO_B$
and
$$ c_1(\cH^\cL) \= \lambda_1 + \frac{1}{2}\pi_{*} \big([V]^2 - c_1(\omega_\pi) \cdot [V]\big)\,.$$
\end{lemma}
\par
\begin{proof} 
For the first statement, note that if all $s_i$ are zero then 
$\pi_{*}\cO_\cX = \cO_B$ since the fibers of $\cX$ are connected. If some $s_i$ is positive, consider the exact sequence 
 $$ 0 \to \cO_{\cX}(V - X_i) \to \cO_{\cX}(V) \to \cO_{X_i}(V) \to 0 $$
and its push-forward by $\pi$. Twisting the sequence by $-s_i$-times a
$\pi$-fiber we compute that the degree of $\cO_{X_i}(V)$ is $ \sum_{e\in X_i\cap X_j} (s_j - s_i)$  where the sum runs over each edge $e$ of $X_i$. We can choose $X_i$ among those with the largest twisting coefficient such that this degree is negative and hence the push-forward term is zero.  Then the fist statement follows from applying induction to $V' = V - X_i$.  
\par
For the second statement, denote by $\cN$ the nodal locus in $\cX$ and let
$\gamma = c_1(\omega_\pi)$ for notation simplicity. 
Then we can apply the first statement, duality and GRR
(and the exact sequence $0 \to \Omega_\pi \to \omega_\pi \to \omega_\pi
\otimes \cN \to 0$ to evaluate $\mathrm{td}^\vee(\Omega_\pi)$)
to obtain that 
\bas 
\phantom{x} & \phantom{\=} \ch(\cH^\cL) \= \ch (\pi_{*}\cL) - \ch (\pi_{*}\cO_{\cX}(V))
\= \ch (\pi_{*}\cL) - \ch (R^1\pi_{*}\cL) \\
& \=  \pi_{*} \Big( \ch (\cL) \cdot \Big(1 - \frac{\gamma}{2} + \frac{\gamma^2 + \cN}{12}
+ \cdots \Big) \Big)  \\
& \=  \pi_{*} \Big(\Big(1 + (\gamma - [V]) + \frac{(\gamma - [V])^2}{2} + \cdots\Big)
\cdot  \Big(1 - \frac{\gamma}{2} + \frac{\gamma^2 + \cN}{12} + \cdots \Big) \Big) \\
& \=  \pi_{*}\Big(1 +  \Big(\frac{\gamma}{2} - [V]\Big) + \Big(\frac{\gamma^2 + \cN}{12}
+ \frac{\gamma^2 - 2\gamma[V] + [V]^2 + \gamma [V] - \gamma^2}{2}\Big) + \cdots \Big) \\
& \=  \pi_{*}\Big(1 +  \Big(\frac{\gamma}{2} - [V]\Big) + \Big( \frac{\gamma^2 + \cN}{12}
+ \frac{[V]^2 - \gamma [V]}{2}\Big) + \cdots \Big) \\
& \=  (g-1) + \lambda_1 + \frac{1}{2}\pi_{*} \big([V]^2 - \gamma [V]\big) + \cdots\,   
\eas
using Noether's formula $\pi_*(\gamma^2 + \cN)/12 = \lambda_1$, 
which implies the claimed formula.  
\end{proof}
\par
Combining Porteous' formula with \autoref{le:GRR}
and using that~\autoref{eq:3264seq} turns into
\bes 
0 \to  \sigma_1^*(\omega_{\rel}(-V))^{\otimes \alpha_1} \to
\cF^\cL_{\alpha} \to \cF^\cL_{\alpha'} \to 0
\ees
we find that the degeneracy locus $\wt{W}_\cL$ of the map $\phi_\cL$
in~\autoref{eq:phitwisted} has class 
\begin{eqnarray*}
[\wt{W}_\cL] & = &  c_1 (\cF^\cL_{\alpha}) - c_1(\cH^\cL) + \xi \\
& = & \sum_{i=1}^n \frac{\alpha_i(\alpha_i+1)}{2} \psi_{i}
- \Big(\sum_{i \in I_\Gamma \atop z_j\in X_i} \alpha_j s_i\Big)
[D_\Gamma] - c_1(\cH^\cL) + \xi \\
& = & \sum_{i=1}^n \frac{\alpha_i(\alpha_i+1)}{2} \psi_{i}
- \Big(\sum_{i \in I_\Gamma \atop z_j\in X_i} \alpha_j s_i\Big)
[D_\Gamma] + \xi - \lambda_1 - 
\frac{1}{2} \pi_{*} [V]^2 
+ \frac{1}{2} \pi_{*}\left(c_1(\omega_\pi) [V] \right) \\
& = & [\wt{W}_\mu(\alpha)] - \Big(\sum_{i \in I_\Gamma \atop z_j\in X_i} \alpha_j s_i\Big)
[D_\Gamma]
-  \frac{1}{2} \pi_{*} [V]^2 + \frac{1}{2} \pi_{*}\left(c_1(\omega_\pi) [V] \right)\,,
\end{eqnarray*}
where we recall that we pretend to work with on $\Gamma$, instead of writing
the sum over all $\Gamma \in \LG_1$. Let $\nu_i$ be the number of edges of $X_i$
(e.g., $\nu = 2$ if $X_i$ is a rational bridge) and $\nu_{i,j}$ the number of edges
joining $X_i$ and $X_j$ (note that $\nu_{i,i}=0$ since $\Gamma$ has no horizontal nodes).
Decomposing~$V$ into its components and using that $X_iX_i = (X_i-F)X_j$ for
a full fiber~$F$ of~$\pi$ we find
\begin{eqnarray*}
\pi_{*}[V]^2 \= \Big( - \sum_{i \in I} s_i^2 \nu_i + 2\sum_{i,j \in I} s_i s_j  \nu_{i,j} \Big) [D_{\Gamma}] 
\end{eqnarray*}
and 
\bes
\pi_*\left(c_1(\omega_\pi) [V]\right) \=  \sum_{i \in I} s_i \pi_*\left(c_1(\omega_\pi) [X_i]\right)  
 \=  \Big(\sum_{i\in I} s_i (2g_i - 2 + \nu_i)\Big) [D_\Gamma]\,. 
\ees
The conclusion of this discussion is:
\par
\begin{lemma}
For integer coefficients $s_i$ of the twisting divisor~$V$ we obtain the coefficient difference 
\ba 
\label{eq:twist-diff}
\Delta_\cL \wt{W}_\Gamma &\,:=\,  \big([\wt{W}_\mu(\alpha)] - [\wt{W}_\cL]\big)_{[D_\Gamma]} \\
& \,\= \ \sum_{i \in I \atop z_j\in X_i} \alpha_j s_i +  \sum_{i,j \in I} s_i s_j  \nu_{i,j} 
- \frac{1}{2}  \Big( \sum_{i \in I} s_i^2 \nu_i +\sum_{i \in I} s_i (2g_i - 2 + \nu_i) \Big)\,.
\ea
\end{lemma}
In what follows we want to maximize this difference by choosing suitable twisting coefficients in the allowed ranges. 
Denote by $\Delta \wt{W}_\Gamma$ the {\em maximum} of $\Delta_\cL \wt{W}_\Gamma$ among all possible choices of the twisted relative dualizing line bundle $\cL$.  
\par
\medskip
We now relabel the vertices of~$\Gamma$ according to levels and single out the irreducible
rational components of the central fiber that stem from blowups.  Let $X_1, \ldots, X_{v^{\top}}$ be the
top level vertices and $Y_1, \ldots, Y_{v_\bot}$ the bottom level vertices. Let
$E_{j}$ be the set of edges adjacent to $Y_j$. For any edge~$e$ we denote by
$R_e^{(k)}$ the rational bridges that stem from the
resolution of the node corresponding to the edge~$e$ for $k=1,\ldots, a_{e}-1$ where $a_e= \ell_\Gamma/p_e$. 
Recall that the twisting coefficients~$s_i$ are zero for all top level components. We also rename them as $\sigma_j$ for the bottom level
components~$Y_j$ and as $s_{e,k}$ for the rational bridges $R_e^{(k)}$, with the convention that $s_{e,0}=0$ and $s_{e, a_e} = \sigma_j$ for $e\in E_j$. As before we require $\sigma_j \leq \ell_\Gamma$ and $s_{e,k}\leq k p_e$ for $e\in E_j$. 
\par 
We introduce the notation $e_{Y_j} = |E_j|$ and
\bes
m_{Y_j} \= \sum_{z_i \in Y_j} m_i\,, \qquad \alpha_{Y_j}\= \sum_{z_i \in Y_j} \alpha_i\,,
\qquad  p_{Y_j} \= \sum_{e \in E_{j}} p_e
\ees
for the total sum of $m_i$, the total sum of $\alpha_i$ and the total sum of prongs that are adjacent to~$Y_j$, respectively.  
The rational bridges do not carry any marked points and the top level gets no twist.
It implies that only those $\alpha_j$ on the bottom level contribute to $\Delta_\cL \wt{W}_\Gamma$.
We thus obtain that 
\ba
\label{eq:DWl1} 
\Delta_\cL \wt{W}_\Gamma &\,\=\, \sum_{j=1}^{v_\bot} \sigma_j \alpha_{Y_j} +\sum_{e \in E} \sum_{k=1}^{a_e-1} (s_{e,k}s_{e,k+1} - s_{e,k}^2)  \\
&\phantom{\=}\,-  \frac{1}{2}\sum_{j=1}^{v_\bot} \sigma_j^2 e_{Y_j}  -
\frac{1}{2}\sum_{j=1}^{v_\bot} \sigma_j (2g({Y_j}) - 2 + e_{Y_j}) \\ 
&\= -\frac{1}2 \sum_{j=1}^{v_\bot} \sum_{e \in E_j} \sum_{k=1}^{a_e} (s_{e,k} - s_{e,k-1})^2 +\frac{1}{2}\sum_{j=1}^{v_\bot} \sigma_j (2\alpha_{Y_j} - m_{Y_j} + p_{Y_j})\,. 
\ea
Given $\sigma_j$, to maximize the above expression, we minimize the quadratic terms
in the first summand by choosing the $s_{e,k}$ nearly equidistant,
i.e., roughly $s_{e,k} \sim  k\sigma_j/a_e$. Working with this possibly fractional approximation
we find that 
\ba
\Delta_\cL \wt{W}_\Gamma \,\,\gtrsim \,\,
\frac{1}{2\ell_\Gamma}\sum_{j=1}^{v_\bot} \sigma_j \Bigl(\ell_\Gamma (2\alpha_{Y_j} - m_{Y_j} + p_{Y_j})
- \sigma_jp_{Y_j}\Bigr)\,. 
\ea
This shows that this is a quadratic optimization problem. In particular for the
natural choice $\alpha_i = m_i/2$ the optimal correction term is obtained for
an integer approximation of $\sigma_j = \ell_\Gamma /2$.
\par
We now take care of the fractional parts in detail. For $e\in E_j$, dividing $\sigma_j$ by $a_e$ we write 
\be
\sigma_j \= q_e a_e + r_e \= (a_e-r_e)q_e + r_e(q_e+1) 
\ee
with $0\leq r_e < a_e$ to compute the numbers of $s_{e,k}-s_{e,k-1}$ equal to~$q_e$ and $q_e+1$, respectively.  
Note that 
$$(a_e-r_e)q_e^2 + r_e(q_e+1)^2 = a_e q_e^2 + 2r_eq_e +r_e = q_e\sigma_j + r_e(q_e+1)\,.$$ 
We thus obtain from~\autoref{eq:DWl1} that
\ba
\Delta_\cL \wt{W}_\Gamma\geq  \frac{1}{2}\sum_{j=1}^{v_\bot} \sigma_j (2\alpha_{Y_j} - m_{Y_j} + p_{Y_j})
- \frac12 \sum_{j=1}^{v_\bot} \sum_{e \in E_j} \big(q_e\sigma_j + r_e(q_e+1)\big)\,.
\ea
\par
The following lemma optimizes this lower bound if $\alpha_i \sim m_i/2$ for all~$i$.
\par
\begin{lemma} \label{le:optdifference}
The maximum difference between the twisted and untwisted Weierstrass divisors 
for the coefficient of the boundary divisor~$D_\Gamma$ is at least 
\be
\Delta_\cL \wt{W}_\Gamma \,\geq\, \Bigl\lfloor \frac{\ell_\Gamma} 2 \Bigr \rfloor \sum_{j=1}^{v_\bot}
 \Bigl(\alpha_{Y_j} - \frac12 m_{Y_j}\Bigr)
+ \frac{\ell_\Gamma}8 \Bigl(P - P_{-1}\Bigr)\,
\ee
where $P  = \sum_{e \in E} p_e$ is the sum of all prongs and $P_{-1} = \sum_{e \in E} 1/p_e$ is the sum of their reciprocals.  
\end{lemma}
\par
\begin{proof} The idea is to show that 
$$q_e\sigma_j + r_e(q_e+1) \lesssim \frac{\ell_\Gamma}{4}(p_e+1/p_e)$$
for an integral choice of~$\sigma_j\sim \ell_\Gamma/2$ and apply this to each edge~$e$ individually.
\par
First suppose $\ell_\Gamma$ is even. In this case we take $\sigma_j=\ell_\Gamma / 2$ for all~$j$. If $p_e$ is
even, then $a_e = \ell_\Gamma / p_e$ divides $\sigma_j$, hence $r_e=0$ and the above inequality estimate literally holds (even without
the term $1/p_e$). If $p_e$ is odd, then $q_e = (p_e-1)/2$, $r_e = a_e/2$, and the inequality estimate becomes an equality.
\par 
Next suppose $\ell_\Gamma$ is odd. In this case we take $\sigma_j = (\ell_\Gamma-1)/2$ for all~$j$ so that
$q_e = (p_e-1)/2$ and $r_e = (a_e-1)/2$. Then we obtain that 
$$q_e\sigma_j + r_e(q_e+1) = \frac{\ell_\Gamma}{4}(p_e+1/p_e) - p_e/2\,$$ 
which implies the desired bound since the last term compensates the rounding of $\sigma_j$ from $\ell_\Gamma / 2$ to $(\ell_\Gamma-1) / 2$. 
\end{proof}

\subsection{Improvement for multiple top level components}

Suppose~$\Gamma \in \LG_1$ is a level graph with $v^\top>1$ vertices
on top level. Here we show that the degeneracy locus $\wt{W}_\cL$
of the map~$\phi_\cL$ even in the twisted setup contains extra copies of the boundary
divisor~$D_\Gamma$ and we estimate the multiplicity. For this purpose
it suffices to work over a small disc~$\Delta_t$ with parameter~$t$
transverse to the boundary divisor. Recall that the tautological section $\bfomega$ vanishes at any top
level zero $z_i$ to the zero order~$m_i$ (hence at least to order $\alpha_i$) over~$\Delta$ 
and it vanishes on the bottom level to order~$\ell = \ell_\Gamma$. Using a
plumbing construction we show that besides $\bfomega$ there exist other such sections:
\par
\begin{prop} \label{prop:vanishingsub}
There is a subbundle~$\cT \subset \pi_*(\omega_\pi)$ of rank~$v^\top$
whose sections vanish
\begin{itemize}
\item along the zero sections~$z_i$ at any top level vertex to order~$m_i$, and
\item along the bottom level components of the central fiber to order~$t^\ell$.
\end{itemize}
\end{prop}
\par
\begin{proof}
By the main theorem of \cite{LMS} the universal family of multi-scale differentials in
a neighborhood of the boundary divisor $D_\Gamma$ can be obtained by plumbing, and thus
the family $\pi\colon \cX \to \Delta$  over the fixed disc with parameter~$t$
and tautological differential~$\omega_t$ can also be obtained by plumbing. We review the essential
steps of the construction in order to show that the plumbing can be performed simultaneously
for a $v^\top$-dimensional space of differentials on the central fiber~$X$.
\par
For the plumbing construction each of the nodes of~$X$ (corresponding to an edge~$e$)
is replaced by the plumbing annulus with differential form
\bas
\bV_e &\= \Bigl\{ (u_e,v_e) \in \Delta_\ve^2\,:\, u_ev_e = t^{\ell /p_e} \Bigr\} \\ 
\Omega &\= C\cdot (u_e^{p_e}+ t^\ell r)du_e/u_e \= (-C) \cdot t^\ell (v^{-p_e} + r) dv_e/v_e 
\eas
for $\ve$ small and for $r,C$ to be specified.
In order to glue in this annulus, we need charts $u_e$ and $v_e$ at the level zero end
and level~$-1$ end of the node that put $\omega_0$ into the standard form
\be
\omega_0^{(0)} \= u_e^{p_e} du_e/u_e\,, \qquad  \omega_0^{(-1)} \= -(v^{-p_e} + t^\ell r)du_e/u_e
\ee
and add a modification differential $\xi(t)$ locally given by $\xi(t) = t^\ell r du_e/u_e$
supported on level zero to compensate for the missing residue term. The sum
$\omega_0 + \xi(t)$ glues with $\Omega$ for $C=1$.
\par
Before we proceed, we remark that for an arbitrary differential~$\eta_0$ supported
on the top level $X_{(0)}$  of the special fiber it is not obvious (and sometimes impossible)
to extend it to the plumbed family. In fact in the given chart~$u_e$, in general $\eta_0$ is given by an arbitrary power series, hence possibly by a series with arbitrary negative powers in~$v_e$. In this case the existence of a differential $\eta_0^{(-1)}$
on the lower level~$X_{(-1)}$ having this prescribed polar part in~$v_e$ is unclear. 
However the situation is better for the following class of differentials.
\par
Let $\bfc=(c_1,\ldots,c_{v^\top})$ be a tuple of non-zero complex numbers and define
the differential $\eta_0(\bfc)$ on $X_0^{(0)}$ to be equal to $c_i \omega_0$ on
the $i$-th component $X_{(0),i}$ of this top level curve for some fixed numbering of these
components. For an edge~$e$ whose upper end goes to the $i$-th component, it locally looks
like $c_i\cdot (u_e^{p_e}+ t^\ell r)du_e/u_e$. Consequently, together with the
modification differentials $c_i \xi|_{X_{(0),i}}$ this glues with $t^\ell$ times
a differential $\eta_1(\bfc)$ whose local form at the lower end of the plumbing fixture
is given by $-c_i \cdot  (v^{-p_e} + r) dv_e/v_e$.
\par
It remains to show that a differential $\eta_1(\bfc)$ on $X_{(-1)}$ with this prescribed
principal part near the lower end of each edge exists. By the solution to the
Mittag-Leffler problem (e.g.,~\cite[Theorem 18.11]{forster}) it suffices to check that the sum of the residues is zero. 
Indeed the tautological differential $\bfomega$
satisfies the global residue condition, i.e., for each~$i$ the sum of residues at
all the edges connecting to $X_{(0),i}$ is equal to zero.  Here all these residues are
multiplied by the same constant~$c_i$. Consequently the sum of the 
residues required in the polar parts of $\eta_1(\bfc)$ is zero, and hence $\eta_1(\bfc)$ exists.  
\par
We now take $\cT$ to be the subbundle generated by the plumbings of the differentials
$(\eta_0(\bfc), \eta_1(\bfc))$. Since the top level is just a rescaling of~$\eta_0$ on each top level vertex 
the first condition holds, and the second condition also holds as mentioned in the above 
plumbing process.
\end{proof}
\par
\begin{cor} \label{cor:2ndcorrection}
The degeneracy locus $\wt{W}_\cL$, where $\cL = \omega_\pi(-V)$
and $V$ is an effective vertical divisor containing the bottom level components 
with multiplicity $\sigma$, contains the boundary divisor~$D_\Gamma$
with multiplicity at least $(v^\top-1)(\ell -\sigma)$.
\end{cor}
\par
\begin{proof} Recall that $\wt{W}_\cL$ is defined as the vanishing locus of the
determinant of the Porteous matrix with rows indexed by a basis of sections of~$\cL/\cO(-1)$ and
columns indexed by the local expansions up to order $\alpha_i$ at the
zeros $z_i$. Consider the $v^\top-1$ rows corresponding to a basis of sections 
of ~$\cT/\cO(-1)$ given by \autoref{prop:vanishingsub},
which is a subbundle of~$\cL/\cO(-1)$. The properties listed in the proposition
imply that the entries in the $\alpha_i$ columns for the zero $z_i$ vanish for all
$z_i$ on top level, and moreover, that the entries in the remaining columns for $z_i$ on the bottom level 
are divisible by $t^{\ell - \sigma}$ (since we have already twisted off $t^{\sigma}$ for the bottom level in~$\cL$).
Taking the determinant of the matrix thus implies the claim.
\end{proof}
\par
We apply the previously obtained improvements from twisting and from top level to the 'middle case' of the generalized Weierstrass divisor by taking $\alpha = \mu / 2$. 
\par 
\begin{cor} \label{cor:genWevenM} If all zero orders $m_i$ are even in $\mu$, then the class
\ba \label{eq:Wmid}
[W^{\mid}] := 
w^{\mid}_\lambda(\mu) \lambda - w^{\mid}_{\hor}(\mu)[D_h] -
\sum_{\Gamma \in \LG_1}  w^{\mid}_{\Gamma}(\mu) \ell_\Gamma [D_\Gamma]
\ea
is an effective divisor class where
\bes
w^{\mid}_\lambda(\mu) \= \frac{12+ \kappa_\mu/2}{\kappa_\mu}\,, \quad
w^{\mid}_{\hor}(\mu) \= \frac{1+\kappa_{\mu}/8}{\kappa_\mu}\,,  \quad
w_\Gamma^{\mid}(\mu) \=  \frac12(v^\top-1)+ \frac{\kappa^{\bot}}{\kappa_\mu}\,.
\ees
\end{cor}
\par
We do not claim that $[W^{\mid}]$ is the class of the closure of the interior locus ${W}_\mu(\mu/2)$
defined in~\autoref{eq:initdefGenW} at the beginning of this section, as 
there might exist special boundary divisors~$D_\Gamma$ that can be subtracted further from $[W^{\mid}]$ such that 
 the remaining class is still effective. However, for certain $\mu$ and~$\Gamma$ there are evidences for the 
 sharpness of our bound (which we do not discuss here to avoid making the paper too long).  With these in mind, 
 we call $[W^{\mid}]$ the \emph{class of the mid-range generalized Weierstrass divisor}.
\par
\begin{proof} For $\mu$ even and $\alpha = \mu / 2$ we have $\vartheta_{\mu,\alpha} = \kappa_\mu/8$.
This converts the expressions from \autoref{prop:Wclass} into the desired forms of $w^{\mid}_{\hor}(\mu)$ and $w^{\mid}_\lambda(\mu)$ for the corresponding coefficients. For the $D_\Gamma$-coefficient
we rewrite the $\vartheta$'s and $\kappa$'s in terms of the top level versions. Then we have 
\bas
8\Bigl(\frac{\kappa^{\bot}}{\kappa_\mu}\vartheta_{\mu,\alpha}
- \vartheta_{\mu,\alpha}^{\bot}\Bigr) \,=
8\vartheta_{\mu_\Gamma^{\top}, \alpha_\Gamma^{\top}} -
{\kappa_{\mu_\Gamma^\top}} \= P_{-1} - P\,.
\eas
Using the twisted version of the mid-range generalized Weierstrass divisor from \autoref{le:optdifference} as well as the improvement from \autoref{cor:2ndcorrection}, we conclude that the class with $-\ell_\Gamma[D_\Gamma]$-coefficient  
equal to 
\bas
\, & \frac{\kappa^{\bot}}{\kappa_\mu} + \Bigl(\frac{\kappa^{\bot}}{\kappa_\mu}\vartheta_{\mu,\alpha} - \vartheta_{\mu,\alpha}^{\bot}\Bigr) +\frac{P - P_{-1}}{8}+\frac{\ell - \sigma}{\ell}(v^\top-1) \\
\= & \frac{\kappa^{\bot}}{\kappa_\mu} +  \frac{\ell - \sigma}{\ell}(v^\top-1)
\eas
is effective. Since in our setting $\sigma = \ell/2$ or $(\ell-1)/2$, then $- (\ell - \sigma) / \ell \leq -1/2$, hence 
the class with $-\ell_\Gamma[D_\Gamma]$-coefficient given by $w^{\mid}_{\Gamma}(\mu)$ is (possibly more) effective.
\end{proof}
\par

\subsection{Odd order zeros and rounding approximations}
\label{sec:average}

For general strata we may still \emph{define} the class~$[W^{\mid}]$ by
the formula~\autoref{eq:Wmid} in \autoref{cor:genWevenM}. However if the zero orders $m_i$ are not
all even, then this class is not obviously effective as $\alpha = \mu/2$ is not an integer tuple. In this case we approximate it
by taking the average of rounding up and down. 
\par 
Let $R(\mu/2)$ be the set of \emph{admissible roundings} for $\mu/2$, defined as follows. For $\alpha
\in R(\mu/2)$ we require that $\alpha_i = m_i/2$ if $m_i$ is even, that
$\alpha_i \in \{(m_i \pm 1)/2\}$ if $m_i$ is odd, and that the total sum
of $\alpha_i$ is $g-1$. That is, we round up and down in precisely half
of the cases. Define the effective divisor class
\be
[W^{\app}] \= |R(\mu/2)|^{-1} \cdot \sum_{\alpha \in R(\mu/2)} [{W}_\mu(\alpha)]
\ee
and thus $[W^{\app}] = [W^{\mid}]$ if all entries are even.
To compute the difference of these two classes in general we write
\be
[{W}_\mu(\alpha)] = w^\alpha_\lambda(\mu) \lambda - w^\alpha_{\hor}(\mu)[D_h] -
\sum_{\Gamma \in \LG_1}  w^\alpha_{\Gamma}(\mu) \ell_\Gamma [D_\Gamma]
\ee
and similarly with the upper index by~$\mid$ or~$\app$. 
We summarize that so far we have computed
\ba 
w^{\alpha}_\lambda(\mu) &\= \frac{12+12 \vartheta_{\mu, \alpha} - \kappa_\mu}{\kappa_\mu}\,, \quad
w^{\alpha}_{\hor}(\mu) \= \frac{1+\vartheta_{\mu, \alpha}}{\kappa_\mu}\,,
\label{eq:genWsummary}
\\
w^\alpha_{\Gamma}(\mu) &\, \geq\, \frac{\kappa_{\mu_\Gamma^{\bot}}}{\kappa_\mu}
(1+\vartheta_{\mu, \alpha})
- \vartheta^\bot
+ \Bigl\lfloor \frac{\ell_\Gamma} 2 \Bigr \rfloor \sum_{i=1}^{v_\bot}
 \Bigl(\frac{2\alpha_{Y_i} -  m_{Y_i}}{2\ell_\Gamma}\Bigr)
+ \frac{P - P_{-1}}8 + \frac12(v^\top-1)\,. \nonumber 
\ea
\par
\begin{lemma} \label{le:wcoeffapp}
Let $\mu$ denote holomorphic signatures. For the $\lambda$-coefficients we have $w^{\app}_\lambda(\mu) \geq w^{\mid}_\lambda(\mu)$ and  $\lim w^{\app}_\lambda(\mu) \geq \lim w^{\mid}_\lambda(\mu) = 1/2$ as $g\to \infty$.  
\par
Moreover for the coefficients of $D_{\hor}$ we have
$w^{\app}_{\hor}(\mu) \geq w^{\mid}_{\hor}(\mu)$ and 
$\lim w^{\app}_{\hor}(\mu) \geq \lim w^{\mid}_{\hor}(\mu) = 1/8$ as $g \to \infty$.
\par
Finally for any $\Gamma \in \LG_1$ and any holomorphic stratum 
\bes
w_\Gamma^{\app}(\mu) - {w}^{\mid}_\Gamma (\mu)
\,\geq\,   \frac{\kappa_{\mu_\Gamma^{\bot}}}{\kappa_\mu}\sum_{m_i\,
	{\rm odd}} \frac{1}{8(m_i+1)}\,- \,\sum_{z_i \in \Gamma_\bot 
	\atop m_i\, {\rm odd}} \frac{1}{8(m_i+1)}.
\ees
\end{lemma}
\par
\begin{proof} We define $\theta(a,m) = a(a+1)/2(m+1)$. The key
observation is that
\bes
\frac12 \Bigl(\theta\Bigl(\frac{m-1}2,m\Bigr) + \theta\Bigl(\frac{m+1}2,m\Bigr)
\Bigr) - \theta\Bigl(\frac{m}2, m\Bigr)
\= \frac{1}{8(m+1)}\,.
\ees
We apply this to the summands of $\vartheta_{\mu, \alpha}$ and observe that
any odd order~$m_i$ is rounded up resp.\ down in $R(\mu/2)$ half of the times. This gives the inequalities for the $\lambda$-coefficients
and the $D_{\hor}$-coefficients. The claim on their limits follows from these inequalities together with
the relation $\vartheta_{\mu, \mu/2} = \kappa_\mu/8$ and the fact that $\kappa_\mu \to \infty$ as $g\to \infty$. 
\par
We now consider the $w_\Gamma$-coefficient. In the comparison between
$\app$ and $\mid$ all the terms involving neither~$\vartheta$
nor~$\vartheta^\bot$ cancel. It implies that 
\be \label{eq:wappstrong}
w_\Gamma^{\app} (\mu) - {w}^{\mid}_\Gamma (\mu) 
\,\geq\, \frac{\kappa_{\mu_\Gamma^{\bot}}}{\kappa_\mu}\sum_{m_i\,
{\rm odd}} \frac{1}{8(m_i+1)}\,- \,\sum_{z_i \in \Gamma_\bot 
\atop m_i\, {\rm odd}} \frac{1}{8(m_i+1)}\,
\ee
which yields the desired inequality.  
\end{proof}
\par 
\begin{rem}\label{rem:wlambdaapp}
{\rm 
Denote by $M_{-1}^{{\rm odd}}:=\sum_{m_i {\rm odd}}\frac{1}{m_i+1}$. Then by the proof of \autoref{le:wcoeffapp} we obtain that 
\[w_{\lambda}^{\app}=w_{\lambda}^{\mid}+\frac{12M_{-1}^{{\rm odd}}}{8\kappa_{\mu}}\,.\]
In particular, if $M_{-1}^{{\rm odd}}$ is negligible comparing to the magnitude of $\kappa_\mu$, then we see that the large genus limits of $w_\lambda^{\app}(\mu)$ and $w_{\hor}^{\app}(\mu)$ coincide with the corresponding limits of the mid-range version, i.e., being $1/2$ and $1/8$ respectively.  For instance, this is the case for signatures $\mu$ whose number of entries is a constant independent of $g$. 

In the case of equidistributed strata with $\mu=(s^n)$, when $s$ is odd, the above specializes to the equality 
\[w_{\lambda}^{\app}=w_{\lambda}^{\mid}+\frac{3}{2s^2 + 4s}\,.\]
}
\end{rem}

\subsection{Refining \autoref{prop:GenTypeCrit} for strata with two zeros}
\label{sec:refofProp13}

For strata of type $\mu = (2m,2g-2-2m)$ the canonical class of the coarse
moduli space (rescaled by the factor $\tfrac{\kappa_\mu}{N}$) is not given
by the right-hand side of~\autoref{eq:canformula}, due to the ramification of the
map from the
stack to the coarse moduli space in the interior, as explained in \autoref{prop:ramifInt}. With the help of the
following proposition we can apply the 'ample+effective'-criterion
\autoref{prop:GenTypeCrit} formally without worrying about the presence of the 
ramification divisor, as long as we use effective divisors containing the
ramification divisor with suitably high coefficients.
\par
\begin{prop} \label{prop:RamGenTypeCrit}
Let $K'_\mu$ denote the right-hand side of~\autoref{eq:canformula}.
Consider the coarse moduli space $\bP\MScoarse$ of multi-scale differentials
of type $\mu = (2m,2g-2-2m)$ with two labeled zeros. If we can write
\be \label{eq:KDAEmod}
K'_\mu - \frac{\kappa_\mu}{N} D_{\NC} \= A + x[B] + y\frac{12}{w^{\app}_\lambda}[W^{\app}]
\ee
with~$A$ an ample divisor class and $B = \BN_\mu$ or $B = \Hur_\mu$ depending
on the parity of~$g$,
and if moreover $x \geq 0$ and $y > 1/24$, then~$\bP\MScoarse$ is a variety of general type for sufficiently large $g$. 
\end{prop}
\par
\begin{proof} Let~$R$ be the interior ramification divisor exhibited
in \autoref{prop:ramifInt}. This proposition implies that
\bes
\frac{\kappa_\mu}{N} K_{\bP\MScoarse} \= K'_\mu - \frac{\kappa_\mu}{N} R\,.
\ees
Note that the Brill--Noether divisor $\wt{\BN}_g$
from~\autoref{eq:defBN}, the Hurwitz divisor~$\wt{\Hur}_g$ from~\autoref{eq:defHur}, and
$W_\mu(\alpha)$ for any~$\alpha$ contain the locus of hyperelliptic curves,
and hence contain~$R$. We write $\wt{B} = f^* \wt{\BN}_g$ or
$\wt{B} = f^* \wt{\Hur}_g$ depending on the parity of $g$, where~$f\colon \bP \LMS
\to \barmoduli[g]$. This implies that
\bes
\frac{\kappa_\mu}{N} (K_{\bP\MScoarse} - D_{\NC}) \= A + \frac{x}{c}[\wt{B}-R] +
y\frac{12}{w^{\mid}_\lambda}[W^{\mid}-R]
+ \Bigl(\frac{x}{c} + y\frac{12}{w^{\mid}_\lambda} -\frac{\kappa_\mu}{N}\Bigr)R\,,
\ees
where~$c$ is the (very large) coefficient of rescaling from $\wt{B}$ to~$B$ (given
in \cite{HarrisMumford} resp.\ in \cite{HarrisKodII}) and the mid and app versions of $W$ coincide for $\mu = (2m, 2g-2-2m)$ with even entries only.
Note that $\frac{\kappa_\mu}{N}
\leq 1$ and $ \frac{12}{w^{\mid}_\lambda(\mu)}\to 24$ as $g\to\infty$.  We thus conclude that the above is a sum of an ample class and an effective class. 
\end{proof}


\section{Certifying general type}
\label{sec:gentype}

In this section we prove \autoref{intro:minimal}, \autoref{intro:fewzero}
and \autoref{intro:stothen} about the Kodaira dimension of strata with few zeros
and equidistributed strata. The general strategy is to apply
\autoref{prop:GenTypeCrit} and its variant
\autoref{prop:RamGenTypeCrit}, using the ample divisor constructed
in the proof of \autoref{thm:LMSisproj} and a combination of effective
divisors introduced in \autoref{sec:taut} and \autoref{sec:genWP}.
\par
On one hand we use the generalized Weierstrass divisor minus its extraneous
boundary components. For simplicity we work exclusively with the
average case $\alpha = \mu/2$ or its nearest integer approximation, i.e.,
with the class $[W^{\app}]$ discussed in \autoref{sec:average}.
On the other hand for $g$ odd we use the Brill--Noether divisor class
given in \autoref{lem:BNclass}. For $g$~even the Brill--Noether divisor
is replaced by the Hurwitz divisor in \autoref{lem:Hurclass} (and for
the minimal stratum with $g=14$, it is replaced by the divisor $\NF_{(2g-2)}$
of \autoref{le:NFBNdiff}).
\par
Technically, we work with a convex combination such that the
$\lambda$-coefficient is zero, and all boundary terms will be shown to be
\emph{strictly} positive. Then we can subtract a small multiple of the ample
class while maintaining the boundary terms positive.
More precisely, we consider the sum
\bas
& \frac{\kappa_{{\mu}}}{N} \Bigl( c_1\bigl(K_{\bP\MScoarse}\bigr)
- D_{\mathrm{NC}} \Bigr)
- y\frac{12\kappa_\mu}{12+12 \vartheta^{\app} - \kappa_\mu}
[{W}^{\app}] - (1-y)\cdot 2[{\BN}_{\mu}] \\
& \qquad \= s_{\hor}(y) [D_h] + \sum_{\Gamma \in \LG_1}  \ell_\Gamma
\bigl(s^{\rm H}_\Gamma(y) [D^{\rm H}_\Gamma] + s^{\NH}_\Gamma(y) [D^{\rm NH}_\Gamma]\bigl)
\eas
with symbols defined as follows. The ramified boundary components
$D^{\rm H}_\Gamma$ of the map from the stack to the coarse moduli space have been
singled out for the boundary divisors of type
HTB, HBT and HBB in \autoref{sec:tocoarse} and they are empty otherwise.
The components $D^{\NH}_\Gamma$ denote the corresponding complement in
$D_\Gamma$ in each case. The coefficients are 
\ba
\label{eq:horstar}
s_{\hor}(y) &\= -1 -\frac{\kappa_\mu}{N} +y
\frac{12(1+\vartheta^{\app})}{12+12 \vartheta^{\app} - \kappa_\mu}
+(2-2y)\frac{g+1}{g+3} \\ 
s_\Gamma^\bigstar(y) &\= c^\bigstar_\Gamma + y \frac{12 w^{\app}_{\Gamma}(\mu)}{w^{\app}_{\lambda}(\mu)}
+ (1-y)b_\Gamma
\qquad \text{for $\bigstar \in \{{\rm H, NH}\}$ }\,
\ea
where for $g$ odd the contributions in $s_\Gamma^\bigstar(y)$ of the canonical
bundle and the Brill--Noether divisor are respectively given by
\ba
\label{eq:cwbcontribution}
c^\bigstar_\Gamma &\= \frac{\kappa_\mu}{N} \Big(N_\Gamma^{\bot} -
R^\bigstar_\Gamma \Big) -  \kappa_{\mu_{\Gamma}^{\bot}}\,,  \\
b_\Gamma &\= \sum_{i=1}^{[g/2]}\sum_{e\in E(\Gamma)
	\atop e \mapsto \Delta_i} \frac{12i(g-i)}{(g+3) p_e} +  
\sum_{e\in E(\Gamma) \atop e \mapsto \Delta_{\irr}} \frac{2(g+1)}{(g+3) p_e}\,,  \\
\ea
and where the coefficients of the Weierstrass divisor are summarized
in~\autoref{eq:genWsummary} and \autoref{le:wcoeffapp}. For $g$ even,
we need to replace $b_\Gamma$ by the corresponding coefficient
of $[\Hur_g]$.
Here $R^\bigstar_\Gamma$  is the (renormalized) contribution
of the NC-compensation divisor plus one (coming from the difference between the canonical and the log-canonical class) and plus the contribution of the ramification
divisor if $\bigstar = H$:
\be \label{eq:defR}
R^\bigstar_\Gamma :=
\frac{b_{\mathrm{NC}}^\Gamma  + 1 +\delta_{\Gamma}^{\rm H}} {\ell_\Gamma}\,.
\ee
Recall that $b_{\mathrm{NC}}^\Gamma$ was defined in \autoref{eq:non-canonical-term} and
$\delta_{\Gamma}^{\rm H} =1$ if $\Gamma$ belongs to HBB, HBT or HTB (see \autoref{figure:ram}) and
$\bigstar = H$, i.e., if~$D_\Gamma$ contains a ramification divisor of the map
to the coarse moduli
space, and zero otherwise. Note that $c_\Gamma^{\rm H} \leq c_\Gamma^{\NH}$ and
likewise for $s^{\rm H}_\Gamma(y)$ for all~$y$. For the purpose of estimates
we thus define $s_\Gamma(y) := \min \{s^{\rm H}_\Gamma(y), s^{\NH}_\Gamma(y)\}$
and $c_\Gamma := \min \{c^{\rm H}_\Gamma, c^{\NH}_\Gamma\}$
if~$\Gamma$ belongs to HBB, HBT or HTB, and we need to control this quantity
only.
\par
For the definition of strata with few zeros, we  refer to the condition \autoref{eq:fewzerodef}. As easily seen in its weaker version \autoref{eq:fewzeroconseq}, this in particular implies the condition stated in \autoref{intro:fewzero}. 
\par
The \emph{proof of \autoref{intro:fewzero} and \autoref{intro:stothen}} can be reduced to showing:
\par
\begin{prop} \label{prop:numericsforstothen}
For all but a finite number of strata with few zeros and even signature, if $y=1/4-\ve$ then $s_{\hor}(y)$  and $s_\Gamma(y)$ are strictly
positive for all~$\Gamma$, where $\ve$ is a constant depending on $g$ defined in \autoref{eq:epsdef}. Moreover, for strata with two zeros and odd signature,  the analogous statement is true for $y=1/6$.  

For all but finitely many equidistributed strata $\mu = (s^n)$ with few zeros, there is a choice of a positive $y<1$ such that $s_{\hor}(y)$  and $s_\Gamma(y)$ are strictly positive for all~$\Gamma$.
\end{prop}
\par
Similarly the \emph{proof of \autoref{intro:minimal}} can be reduced to showing:
\par
\begin{prop} \label{prop:numericsforminimal}
For $g \geq 44$ the coefficients $s_{\hor}(0.19)$  and $s_\Gamma(0.19)$ are
strictly positive for all~$\Gamma$.  For  the range $ 13\leq g \leq 43$, the
coefficients $s_{\hor}(y)$  and $s_\Gamma(y)$ are both
strictly positive for all~$\Gamma$ for $y$ given in \autoref{cap:Rangeygt}.
\end{prop}
\par
\begin{figure}[h]
$$ \begin{array}{|c|c|c|c|c|c|c|}
\hline
g & 13 & 14 & 15 & 16 & 17&18 \\
\hline &&&&& \\ [-\halfbls]
y \in  & [ 0.78,0.79] & [0.67,0.68]
& [0.59,0.73]
& [0.63,0.66]
& [0.47,0.68] 
 & [0.54,0.62] \\
[-\halfbls] &&&&&&\\
\hline  &&&&&&\\ [-\halfbls]
g & 19 & 20 & 30 & 40 & 43 & 44\\
\hline &&&&&& \\ [-\halfbls]
y \in  & [0.39,0.53] & [0.47,0.59]
& [0.29,0.42]
& [0.22,0.35]
& [0.21,0.37] 
& [0.21,0.34]   \\
[-\halfbls] &&&&&&\\
\hline
\end{array}
$$
\caption{Range of~$y$ for showing that minimal strata with odd spin are of general type. For $g=14$ we use the $\NF_{(2g-2)}$ divisor instead of the Hurwitz divisor $\Hur_{(2g-2)}$ to substitute the Brill--Noether divisor $\BN_{(2g-2)}$. } 
\label{cap:Rangeygt}
\end{figure}
\begin{proof}[Proof of \autoref{intro:minimal}, \autoref{intro:fewzero}
    and \autoref{intro:stothen}.]
  If we assume the claims of \autoref{prop:numericsforstothen} and
  \autoref{prop:numericsforminimal}, then we can write 
\[ c_1\bigl(K_{\bP\MScoarse}\bigr)
- D_{\mathrm{NC}} =
C_1\cdot
[{W}^{\app}] +C_2 \cdot 2[{\BN}_{\mu}]+E'\]
where $C_i$ are positive constants and $E'$ is a linear combination of
all boundary divisors with strictly positive coefficients. (In the previous expression, 
the Brill--Noether divisor has to replaced by the Hurwitz divisor or the $\NF$ divisor for even genera.)
Let $A = \lambda_1 + \ep \c_1(\cL_{\ol{B}}\otimes \cO_{\ol{B}}(-D))$ be
the ample class constructed in \autoref{sec:projectivity} (see in particular
the introductory paragraphs and \autoref{prop:fampleclass}).
Hence by slightly perturbing  the coefficient of $[{W}^{\app}]$ 
we obtain that for $\delta_1$ small enough 
\[ c_1\bigl(K_{\bP\MScoarse}\bigr)
- D_{\mathrm{NC}} =
C_1'\cdot
[{W}^{\app}] +C_2 \cdot 2[{\BN}_{\mu}]+ \delta_1 A  + \delta_2 \lambda_1 + E''\]
with $\delta_2 >0$ and where $E''$ is still effective.
We can now apply \autoref{prop:GenTypeCrit} for all strata in the given list,
except for strata with two zeros, for which the canonical class is not given by the
formula \autoref{eq:cwbcontribution} because of the ramification of the map from the
stack to the coarse moduli space in the interior, as explained in \autoref{prop:ramifInt}.  For strata with two zeros, we can nevertheless apply \autoref{prop:RamGenTypeCrit} since the value of $y$ used in \autoref{prop:numericsforstothen} is $1/6$, which is greater than $1/24$.
\end{proof}
\par
We have now reduced to prove \autoref{prop:numericsforstothen}
and \autoref{prop:numericsforminimal}. 

\subsection{A summary of notations}

We work throughout in the stratum with signature $\mu = (m_1, \ldots, m_n)$.
If $\mu$ is a holomorphic signature the unprojectivized dimension of
the corresponding stratum is $N=2g+n-1$. We let
\be
M \= \sum_{i=1}^n  m_i = 2g-2\,, \qquad M_{-1} \= \sum_{i=1}^n \frac{1}{m_i+1}\,.
\ee
Recall from~\autoref{eq:defkappamu} and~\autoref{eq:vartheta} that 
\bes
\kappa \,:= \,\kappa_\mu \= \sum_{i=1}^n \frac{m_i(m_i+2)}{m_i+1} 
\= 2g - 2 + n - M_{-1}\,, \quad
\vartheta \,:=\,
\vartheta_{\mu, \alpha} \= \sum_{i=1}^n \frac{\alpha_i(\alpha_i+1)}{2(m_i+1)}
\ees
if $\alpha = (\alpha_1, \ldots, \alpha_n)$ is a partition of $g-1$
such that $0\leq \alpha_i \leq m_i$ for all~$i$.
\par
If $\Gamma$ is a two-level graph, all these notations have the
corresponding meaning for top and bottom levels, giving rise to
$N^\top, N^\bot$, to $M^\top, M_{-1}^\top, M^\bot, M_{-1}^\bot$ and to  
$\kappa^\top := \kappa_{\mu_{\Gamma}^{\top}}$ etc, so that
$$ \kappa^\bot + \kappa^\top \= \kappa_\mu\,,
\qquad
\vartheta^\bot + \vartheta^\top = \vartheta\,.$$
Finally, level graphs come with the prongs associated with the
edges and we define 
\be \label{eq:recallPPm}
P  = \sum_{e \in E} p_e\,, \qquad P_{-1} = \sum_{e \in E} 1/p_e\,.
\ee
\par
\begin{lemma} \label{le:kappamupositive}
For any meromorphic type~$\mu$ without simple poles we
have $\kappa_\mu \geq 0$. Moreover, $\kappa_\mu = 0$ if and only if
$\mu=(0,\dots,0)$ or $\mu=(-m,m-2,0,\dots,0)$ for $m\geq 2$.
In the remaining cases we have $\kappa_\mu \geq 1/3$.
\end{lemma}
\par
\begin{proof}
Note that $m_i/(m_i+1) \geq 1/2$ for $m_i >0$ and  $m_i/(m_i+1) > 1$ for $m_i < -1$. 
If $g\geq 1$, the claim follows since $n- M_{-1} \geq 1/2 $ if
$\mu \neq (0,\dots,0)$. Consider the case of $g=0$, for which $n\geq 3$. We can assume that $\mu$ has at least three non-zero entries $m_1, m_2, m_3$, with at least one positive and one negative, say $m_1 > 0$ and $m_2 < -1$. If $m_3 > 0$, then 
$\kappa_\mu \geq - 2 + 1/2 + 1/2 + 4/3 = 1/3$ (with equality attained for $\mu = (1,1,-4, 0, \ldots, 0)$).   
If $m_3 < -1$, then $\kappa_\mu > -2 + 1/2 + 1 + 1 > 1/3$. 
\end{proof}
\par
We summarize some more parameters that coarsely classify graphs
$\Gamma \in \LG_1$. These are the number of edges~$E = E_\Gamma$, the
number of vertices $v^\top$ on top and $v^\bot$ on the bottom, the number
of marked points $n^\top$ on top and $n^\bot$ on the bottom and the
genera $g_{i}^\top$ of the vertices on top and their
sum $g^\top = \sum g_{i}^\top$ and similarly for bottom level.
\par
Finally we recall from \autoref{sec:Singbd} that the definition of $D_{\NC}$ in \autoref{eq:non-canonical-term},
together with the definition of $R:= R_\Gamma^\bigstar$ above, gives a bound
\be \label{eq:RfromDNC}
R\leq \frac12 P_{-1}^{\NCT} + 4P_{-1}^{\rm EDB} + P_{-1}^{\RBT} +  2P_{-1}^{\OCT}
+ \frac{\delta_{\Gamma}^{\rm H}}{\ell_\Gamma}
\ee
where ${\rm EDB}$ are \emph{elliptic dumbbells}, compact type edges with
one elliptic end which we consider only if $\Gamma$ has exactly one edge, where $\RBT$ are \emph{rational bottom tails},
tails with a rational vertex on bottom level, and $\OCT$ abbreviates
\emph{other compact type} edges. Moreover $\NCT$ denotes {\em non-compact type} edges.
We define $P_{-1}^{\NCT}$ etc, as in~\autoref{eq:recallPPm}, with the sum restricted to
the corresponding subset of edges.
\par
Under these notations the $\Gamma$-coefficients of the Brill--Noether divisor class 
are estimated by 
\be \label{eq:BNcoarse}
b_\Gamma \,\geq\, 2\frac{g+1}{g+3}P_{-1}^{\NCT} + 12 \frac{g-1}{g+3}P_{-1}^{\OCT} + 
12 \frac{g-1}{g+3}P_{-1}^{{\rm EDB}}\,.
\ee

\subsection{The strategy for a general stratum}

The horizontal divisor, the analogue of $\delta_{\irr}$ for $\barmoduli[g]$,
is not the main concern here, as suggested by a coarse estimate as follows, which is not restricted to the special strata we consider but holds true in the general case of any holomorphic signature.
\par
\begin{lemma} \label{le:chorestimate}
For each $y > 0 $ there are at most finitely many holomorphic strata such that $s_{\hor}(y) > 0$
does not hold.  
\par
For the minimal stratum $\mu = (2g-2)$ we have $s_{\hor}(0.19) >0$ for $g \geq 44$
and $s_{\hor}(y) >0$ for $12 \leq g \leq 43$ and for~$y$ satisfying the lower bound ranges given in \autoref{cap:Rangeygt} (and the lower bound $y> 0.91$ for $g=12$).
\end{lemma}
\par
\begin{proof} We start with the case $g$ odd.
The coefficient $s_{\hor}(y) > 0$ if and only if 
  \[y\geq x:=\Big(1+\frac{\kappa}{N}-\frac{2g+2}{g+3}\Big)\Big(\frac{12(1
    +\vartheta^{\app})}{12+12 \vartheta^{\app}- \kappa}-\frac{2g+2}{g+3}\Big)^{-1}\,.\]
If the entries of $\mu$ are even and we use $\alpha=\mu/2$, then
$\vartheta^{\app} =\vartheta^{\mid}=\kappa/8$, and hence 
\[\frac{12(1+\vartheta^{\app})}{12+12 \vartheta^{\app}- \kappa}
\=\frac{24+3\kappa}{24+\kappa}\to 3\]
as $g \to \infty$. Since
\be \label{eq:knuNestimate}
\frac{\kappa}{N}=1-\frac{1+\sum_{i=1}^n \frac{1}{m_i+1}}{2g+n-1}
\, < \, 1\,,
\ee 
the numerator of~$x$ is smaller than any positive bound as $g \to \infty$.
If the entries of $\mu$ are  odd,  we can use the
relation $\vartheta^{\app}=\vartheta^{\mid}+M_{-1}^{{\rm odd}}$ shown in the proof
of \autoref{le:wcoeffapp} to prove that the denominator of $x$ still converges 
to a positive constant.
\par
We next deal with the case $g$ even and $\mu$ not the minimal stratum,
where the expression  $s_{\hor}^{\app}(y)$ involves the Hurwitz divisor.
Its negative $D_h$-coefficient is within $O(1/g)$ of the coeefficient
of $\BN_\mu$. This implies that whenever we claimed above that
$s_{\hor}^{\app}(y) > 0$ for a fixed~$y$ and all but finitely many strata, the same claim holds for the
corresponding sum $s_{\hor}^{\app}(y)$ involving~$\Hur_g$.
\par
For the minimal strata, one can check that the Hurwitz divisor gives a smaller  $s_{\hor}^{\app}(y)$. By using the expression of $s_{\hor}^{\app}(y)$ involving the Hurwitz divisor and the monotonicity of the lower
bound as $g \to \infty$, one can verify that $y=0.19$  works
for $g \geq 44$. An explicit computation for the remaining cases can be done to check that the lower bounds displayed in \autoref{cap:Rangeygt} work for $13\leq g \leq 43$ and $y> 0.91$ works for $g=12$ (for $g=12$ and $g=14$, we use the divisor $\NF_{(2g-2)}$ instead of the Hurwitz divisor).
\end{proof}
\par
We rewrite now the contributions of $c_\Gamma$ and $w^{\mid}_\Gamma$ in terms of the
parameters characterizing boundary divisors. We frequently drop the index~$\Gamma$
in the sequel to lighten the notation.
\par
\begin{lemma}\label{le:CGamma}
The contribution of the canonical bundle in terms of the parameters
classifying boundary divisors is given by 
\ba \label{eq:CGamma}
c_\Gamma&\= \Big(1-\frac{\kappa}{N}\Big)(M^\top+n^\top+P)
-M_{{-1}}^\top-P_{-1} -\frac{\kappa}{N}\Big(v^\top+R\Big)\,.
\ea
\end{lemma}
\par
\begin{proof}
Rewriting all objects in $c_\Gamma$ in terms of the top
level versions gives
\begin{align*}
c_\Gamma&\=\frac{\kappa}{N}  N^{\bot} - \kappa^{\bot}
-\frac{\kappa}{N}R
\= \kappa^{\top} - \frac{\kappa}{N}  N^{\top}
-\frac{\kappa}{N}R\,.
\end{align*}
Substituting in this expression 
\ba \label{eq:topreplacements}
\kappa^\top &\= M^\top + n^\top - M^\top_{{-1}}
+ P- P_{{-1}}\,, \\
N^\top &\=   M^\top + n^\top + P + v^\top
\ea
gives the claim.
\end{proof}
\par
Using the previous lemmas, we write the full main estimate that we want to control as
\begin{flalign} 
&\phantom{\=}\,\,
s_\Gamma(y) -  12y \frac{w^{\app}_{\Gamma}(\mu) - w^{\mid}_{\Gamma}(\mu)}{w^{\app}_{\lambda}(\mu)}
\nonumber\\
&\= c_\Gamma  + y \frac{12 w^{\mid}_{\Gamma}(\mu)}{w^{\app}_{\lambda}(\mu)} + (1-y)b_\Gamma
\nonumber\\
&\,\geq\, \Big(\frac6{w_{\lambda}^{\app}}
	y-\frac{\kappa}{N}\Big)(v^\top-1)
+(1-y)b_\Gamma- P_{{-1}} -\frac{\kappa}{N} R
\label{eq:gencGammaMainEstimate} \\	
	&\,+\, \frac{1+M_{-1}}{N}(M^\top+n^\top +P)
	\,-\, M_{{-1}}^\top\,+\,y\cdot \frac{12}{w_{\lambda}^{\app}}
\frac{\kappa^\bot}{\kappa}-\frac{\kappa}{N}   \qquad =: T_1 + T_2\,
 \nonumber 
\end{flalign}
where the $T_i$ terms on the right-hand side of the above expression are labeled one for each line.
\par
We start by showing an estimate for $T_1$ in the case of a general stratum. From now on we fix 
\begin{equation}\label{eq:epsdef}
	\ve:=\frac{11g - 2}{4g^2 +16g - 8}\,.
\end{equation}

\begin{lemma} \label{le:fewzeroT1est}
	If we have 
	\be \label{eq:yrange}
	\frac{w^{\app}_\lambda}{6} \leq y\leq \frac{1}{4}-\ve\,,
	\ee
	then $T_1 \geq - 3P_{-1}^{\rm RBT}$ for $g$ large enough.
	Unless $\Gamma$ is a dumbbell graph with a rational bottom that carries all the
	marked points, we have the stronger estimate $T_1 \geq - 2P_{-1}^{\rm RBT}$.
\end{lemma}
\par
\begin{proof} 
The lower bound is obtained by imposing the coefficient of $v^{\top}$ to be non-negative. By using $\kappa/N\leq 1$, we see that if the lower bound holds then we have 
\[\Big(\frac6{w_{\lambda}^{\app}}
y-\frac{\kappa}{N}\Big)(v^\top-1)\geq \Big(\frac6{w_{\lambda}^{\app}}
y-1\Big)(v^\top-1)\geq 0\,.\]
\par
In the case of odd genus $g$, we use the Brill--Noether divisor in the expression
of $b_\Gamma$, while for even genus $g$, we need to use the Hurwitz divisor. 
Using the expression of the Hurwitz divisor obtained in \autoref{lem:Hurclass}, we can
check that $b_\Gamma$ is smaller in the Hurwitz case. More specifically, let us define 
\[\delta_{\text{Hur}}^{\text{BN}}=\begin{cases}
\frac{g^2 - 3g + 2}{3g^3 + 32g^2 + 61g - 24}& g \text{ even},\\
0 &  g \text{ odd}.
\end{cases}\]
Using $\kappa/N <1$ and the estimate~\autoref{eq:BNcoarse} for Brill--Noether (or the
analogous one for Hurwitz) we find, thanks to the lower bound for~$y$, that
\ba \label{eq:T1rewrite}
T_1\geq	& \left(\frac{1-4y}{2}-2(1-y)\left(\frac{2}{g+3}+5\delta_{\text{Hur}}^{\text{BN}}\right)\right)P_{-1}^{\NCT}\\
&+\left(9-12y-2(1-y)\left(\frac{24}{g+3}+60\delta_{\text{Hur}}^{\text{BN}}\right)\right) P_{-1}^{\OCT}  \\
& 		+\left(7-12y-2(1-y)\left(\frac{24}{g+3}+60\delta_{\text{Hur}}^{\text{BN}}\right)\right) P_{-1}^{{\rm EDB}} \\
&- 2P_{-1}^{\rm RBT}  - \frac{\delta_{\Gamma}^{\rm H}}{\ell_\Gamma}\,.
\ea
The upper bound for $y$ in \autoref{eq:yrange} is exactly the one that makes the
coefficient of $P_{-1}^{\NCT}$ in the previous expression positive.  Moreover it also implies that all the other terms in the brackets are
positive for $g$ large enough. This shows the first claim. 
\par 
For the strengthening claim we observe that $\frac{\delta_{\Gamma}^{\rm H}}{\ell_\Gamma} \leq 1/p_e$ for every~$e$.
This ramification term is covered by the $P_{-1}^{{\rm EDB}}$-summand or the
$P_{-1}^{\OCT}$-summand, if at least one such edge exists. The only ramified
boundary divisors whose level graphs do not have such an edge are exactly dumbbell
graphs with a rational bottom that carries all the marked points. Since for these
special graphs~$\Gamma$ we have $\ell_\Gamma=P_{-1}^{\rm RBT}$, we have shown the
full claim.
\end{proof}
\par
\begin{rem} {\rm 
There are sequences of connected strata for which it is impossible to
show bigness of the canonical class or general type for all but finitely many cases by just using the
Brill--Noether (and Hurwitz) divisors and the approximation to the
mid-range Weierstrass divisor~$W^{\app}$.
\par
Indeed if we impose no constraints on a sequence of signatures~$\mu_k$, there are
graphs~$\Gamma_k$ giving a boundary divisor of $\bP\MScoarse[\mu_k]$ for which
the term $T_2$ tends to negative infinity for $k$ growing and the term $T_1$ stays
bounded. Consider the example of $\mu_k=(k-2,2^{k/2})$ and the case of
a graph $\Gamma$ where only the zero of high order $k-2$ is on bottom level and
where~$P$ is independent of~$k$. Then  the only linear terms in~$k$ of~$T_2$
are the positive term $(1+M_{{-1}})(M^\top+n^\top)/N=k/10+O(1)$  and
the negative term $-M^\top_{-1}=-k/6$, while all the other terms are bounded.
Note also that if $v^\top$ is independent of $k$, then also $T_1$ is
independent of $k$, so $s^{\mid}_\Gamma(y)<0$ for any $y$ for almost any $k$. 
}
\end{rem}

\subsection{The (non-minimal) strata with few zeros}
\label{subsec:fewzero}

We exclude the minimal strata from this section, since one source of ramification
divisors at the boundary, the HBB graphs, occurs only for the minimal strata and since we will
make the bound effective for them.
\par
We estimate the summands $T_1$ and $T_2$  of~\autoref{eq:gencGammaMainEstimate}.
Recall that by \autoref{le:wcoeffapp} and the subsequent remark $w_\lambda^{\app}
\to 1/2$ as $g \to \infty$ for any sequence of strata with a uniformly bounded number
of zeros and that $w_\lambda^{\app} = w_\lambda^{\mid}$ in case all the zeros are of even
order.
\par
We say that a stratum has \emph{few zeros} if
\be \label{eq:fewzerodef}
M_{-1} \leq \frac{12}{w_{\lambda}^{\app}}\frac{N}{\kappa}\left(\frac{1}{4}-\ve\right) -1
\ee
where $\ve$ was defined in \autoref{eq:epsdef}. The above condition implies that $w_{\lambda}^{\app}\to 1/2$. Indeed, since by \autoref{rem:wlambdaapp} we know that $w_{\lambda}^{\app}$ is a bounded function of $g$, by \autoref{eq:fewzerodef} we also have that $M_{-1}$ is a bounded function in the case of strata with few zeros. But this, again by \autoref{rem:wlambdaapp}, implies that the limit for $g\to \infty$ of $w_{\lambda}^{\app}$ is equal to the limit of $w_{\lambda}^{\mid}$, which is $1/2$. 
Hence, since $\kappa <N$, the condition 'few zeros' is implied by
\be\label{eq:fewzeroconseq}
M_{-1} \leq 5-24\ve'\,
\ee
where $\ve'$ is a function going to zero for $g\to \infty$.
In particular choosing $\ve' < 11/48$ implies that strata with $n \leq 10$
qualify for 'strata with few zeros'. Obviously, the condition depends on 
the distribution of zero orders. For example, if the zero type~$\mu$ does not
have simple zeros, then strata with up to~$15$ higher order zeros qualify for 'strata with
few zeros'.
\par
\begin{lemma} \label{le:fewzeroevenT2est}  
Let $y = 1/4 - \ve$ be as above. Then, for all but a finite number of strata
with few zeros, the following estimates hold:
If~$\Gamma$ is a dumbbell graph with a rational bottom vertex that carries all
marked points, then $T_2 \geq 3P_{-1}^{\rm RBT}$.
For all other $\Gamma$ we have $T_2 \geq 2P_{-1}^{\rm RBT}$.
\end{lemma}
\par
Note that the condition 'few zeros' depends on~$\ve$, i.e., the
smaller~$\ve$ is chosen the larger we need to take~$g$ for both the $T_1$-estimate
and the $T_2$-estimate to hold.
\par
\begin{proof}
Using that
\[M ^\top+n^\top+ P \= N-\kappa^\bot-1-M_{{-1}}^\bot+ P_{-1}\]
we can rewrite the $T_2$ expression above as 
\begin{align}\label{eq:T2fewzeros}
T_2&\=\left(y \frac{12}{w_{\lambda}^{\app}\kappa}-\frac{1+M_{{-1}}}{N}\right)
\kappa^\bot+\left(1-\frac{1+M_{{-1}}}{N}\right)M_{{-1}}^\bot
+\frac{1+M_{-1}}{N} P_{-1}\,.
\end{align}
The condition on few zeros ensures that with $y=\tfrac14 - \ve$
the first summand is positive for $g$ large. In the absence of rational tails each
of the terms is positive and we are done. We need to refine this in the
presence of rational tails. Note that each of the quantities
$\kappa^\bot$, $M^\bot$, $M_{{-1}}^\bot$ and $P_{-1}$ can be interpreted as a sum over
the vertices on bottom level. For each such vertex~$v$ we thus define accordingly $P_v$ and $P_{-1,v}$ etc, 
and write $\kappa_v := \kappa^\bot_v$ or
$M_v := M_v^\bot$ etc, since our focus is on bottom level anyway. 
\par
We will apply this mainly for $v$ being a rational tail vertex. In this case
$\kappa_{v} =n_{v}-1-M_{-1,v}+ 1/p_e$ 
where~$e$ is the rational tail edge
and $p_e = M_v +1$. Using the trivial estimate for non-rational tails we deduce that 
\begin{align}\label{eq:rationaltails}
 T_2 &\geq  \sum_{v \in V^{\RBT}} \left(y \frac{12}{w_{\lambda}^{\app}\kappa}-\frac{1+M_{{-1}}}{N}\right)
(n_{v}-1) \nonumber\\
& +\left(1-y \frac{12}{w_{\lambda}^{\app}\kappa}\right)M_{-1,v}
+\frac{12y}{w_{\lambda}^{\app}\kappa} \frac{1}{M_{v}+1}\,. 
\end{align}
\par
The first term is non-negative for every~$v$ and the last is positive, a negligibly 
small multiple of $1/p_e$. We use the middle summand to get the required positivity.
For this we note that 
\be \label{eq:Mmest}
M_{-1,v}\geq \frac{3}{M_{v}+1} = \frac{3}{p_e}
\ee
with equality if and only if $n_v =2$ and $M_v=2$, i.e., for rational tails
the three legs. Indeed, since the sum of reciprocals is minimized by the
equidistributed situation, we have
\bes
M_{-1,v}\= \sum_{i=1}^{n_{v}} \frac{1}{m_{i}+1}\geq  \frac{n_{v}^2}{M_{v}+n_{v}}
\ees
and 
\begin{equation}\label{eq:Minvbound}
\frac{n_{v}^2}{M_{v}+n_{v}}-\frac{3}{M_{v}+1}\
=\frac{(n_{v}^2-3)M_{v}+n_{v}(n_{v}-3)}{(M_{v}+n_{v})(M_{v}+1)}
\end{equation}
which is zero if $n_{v}=M_{v}=2$ and positive otherwise, since $n_{v}\geq 2$ and
since $M_{v}\geq n_{v}$.
\par
Suppose that $\Gamma$ is not a rational bottom dumbbell. Then we use~\autoref{eq:Mmest}
to get that 
\begin{equation}\label{eq:raitonaltailsfewzero}
\left(1-  \frac{1+M_{{-1}}}{N}\right) M_{-1,v} > \frac{2}{p_e}
\end{equation} 
for large~$g$ (since $\frac{1+M_{{-1}}}{N}\leq y \frac{12}{w_{\lambda}^{\app}\kappa}=
O(1/g)$ by condition \autoref{eq:fewzerodef}). Summing these
contributions gives the term $2P_{-1}^{\rm RBT}$ we wanted.
\par
Finally consider rational bottom dumbbell graphs with all marked points on bottom level.
We have to improve the above estimate by $1/p_e = 1/(2g-1)$. For these graphs with
$M_v=2g-2$, for any small $1>\delta>0$ we find the analogue of \autoref{eq:Minvbound} in
this situation to be
\begin{align*}
(1-\delta)M_{-1}-\frac{3}{p_e}&\geq \frac{((1-\delta)n^2-3)(2g-2)+n((1-\delta)n-3)}{(2g-2+n)(2g-1)}\,.
\end{align*}
The previous expression is positive for $n\geq 4$. For $n=2,3$, one can check that
it is positive for $\delta=1/8$ and $g$ large enough. Since  $\frac{1+M_{{-1}}}{N}\leq y
\frac{12}{w_{\lambda}^{\app}\kappa}\leq 1/8$ for $g$ large enough, we have proven the statement.
\end{proof}
\par

\subsection{The equidistributed strata}
\label{subsec:equidistributed}
We consider now  strata of type $\mu=(s^n)$. For these strata the parameters
$M_\Gamma^\top$ and $n_\Gamma^\top$ are dependent, since $M_\Gamma^\top=s\cdot n_\Gamma^\top$.
The main quantities for such strata in terms of $s$ and $n$ are
\ba \label{eq:equiexpr}
g &\= \frac{s}2 n + 1\,, &\quad N &\= (s+1)n + 1\,, \quad &\kappa_\mu
&\= n\frac{s(s+2)}{s+1}\,.
\ea
Moreover in this case 
\[\frac{\kappa^\bot}{\kappa} \=1-\frac{n^\top}{n}
+\frac{P_{-1} - P}{\kappa}.\]
We present now the analogue of \autoref{le:fewzeroevenT2est} in the case of
equidistributed strata.
\par
\begin{lemma} \label{le:equisitrT2est}  
For all but a finite number of equidistributed strata, there is some $y$ satisfying
condition \autoref{eq:yrange} such that the following estimates hold: If $\Gamma$
is a dumbbell graph with a rational bottom vertex that carries all marked points,
then $T_2 \geq 3P_{-1}^{\rm RBT}$. For all other $\Gamma$ we have $T_2 \geq 2P_{-1}^{\rm RBT}$.
\end{lemma}
\par
\begin{proof}
	In the case of equidistributed strata, we can either consider the expression
	for $T_2$ given in \autoref{eq:T2fewzeros} or an equivalent expression given by
	\begin{align}\label{eq:sn}
		T_2&\,=\, \frac{ns(s+2)}{(s+1)((s+1)n + 1)}\frac{n^\top}{n}+\frac{12}{w^{\mid}_\lambda}
		y\Big(1-\frac{n^\top}{n}\Big) \nonumber -\frac{\kappa}{N}\\
		&\quad +\left(\frac{1+M_{{-1}}}{N}-y\frac{12}{w_{\lambda}^{\mid}\kappa}\right)P+y\frac{12}{w_{\lambda}^{\mid}\kappa}P_{-1}\,.
	\end{align}
	Note that the coefficient of $\kappa^\bot$ in \autoref{eq:T2fewzeros} is
	exactly the negative of the coefficient of $P$ in the previous displayed equation. We
	first consider the range of parameters $n\leq 2(s+1)$  (we call this the range of few zeros) together with the choice $y=w_\lambda^{\app}/4$. With this choice the coefficient
	of $\kappa^\bot$ in \autoref{eq:T2fewzeros} is positive. In fact, plugging in the
	quantities from~\autoref{eq:equiexpr} yields a rational function in~$(n,s)$
	with positive denominator and a quadratic polynomial in~$n$ with $\bQ[s]$-coefficients
	in the numerator with top coefficient $-s(2+s)$. It thus suffices to check
	the positivity at the boundary values $n=2$ and $n=2(s+1)$.
	\par
	We then consider the complementary range of parameters $n>2(s+1)$ (we call this
	the range of many zeros) together with $y=w_\lambda^{\app}/6$. With this choice
	the coefficient of $P$ in \autoref{eq:sn} is positive. In fact, plugging in
	yields a rational function, which when expressed in the shifted variables~$s'=s-1 \geq0$
	and $n' = n - 2(s+1) \geq 0$ has exclusively non-negative coefficients.
	\par
	In the range of few zeros, since it is clear that the expression \autoref{eq:T2fewzeros} is positive,  if there are no rational tails then we are done.
	In the range of many zeros and in absence of rational tails, we only need to show that the first line of \autoref{eq:sn} is positive.
	Once can check that  the expression is  minimized for the maximum 
	value $n^\top=n-1$, and this bound already gives a positive expression for the first line of \autoref{eq:sn}.
	%
	\par 
	If there are rational tail edges we can consider  the contribution of each rational
	tail edge separately, as we did in the proof of \autoref{le:fewzeroevenT2est}.
	Via a numerical check given by specializing \autoref{eq:rationaltails} and using the
	fact that by stability every bottom level vertex of a rational tail has at least
	two legs, we can show that indeed in this case $T_2\geq 2P_{-1}^{\rm RBT}$.
	Similarly, we can also numerically check that in the case of a dumbbell with rational bottom and $n^\top=0$, we obtain the stronger bound $T_2\geq 3P_{-1}^{\rm RBT}$.
\end{proof}

\subsection{The minimal strata with odd spin}
\label{subsec:minimal}
For the minimal strata there is no discussion of odd order zeros nor
of rational tails, but we want to make the estimates effective. We give
again estimates for the terms in~\autoref{eq:gencGammaMainEstimate}.
\begin{lemma}\label{lem:minimalbound}
	Let $y=0.19$ and $g\geq 44$. Then we have $T_1> 0$ for all graphs apart from a banana graph or a  double banana graph with two  vertices of genus one on top level, which are also HBB graphs, for which we have $T_1> -1/\ell_\Gamma$. Moreover
	\begin{equation}\label{eq:T2minimal}
	T_2\geq 2\frac{g}{g-1}\left(1-\frac{P}{\kappa}\right).
	\end{equation}
\end{lemma}
\par
\begin{proof}
Since one can check that, for $g\geq 44$, our choice of $y$ satisfies the condition \autoref{eq:yrange}, 
we can argue as in the proof of \autoref{le:fewzeroT1est}  and effectively check
that $T_1\geq 0$ apart from the case of an HBB. 

In the case of an HBB, if there is a separating edge $e$, the additional term $-1/\ell_\Gamma$
is compensated by the $1/p_e$ contribution. Moreover, if there are at least three vertices,
then the $v^\top$-term in $T_1$ is also enough to compensate the negative term coming
from ramifications. The same is true if $v^\top=2$ and $\ell_\Gamma>1$.
Hence the first part of the statement is proved.
\par
In order to show the second part of the claim, it is enough to note that in the case
of strata with $n^\top=0$, which is the case of the minimal strata, we can
simply specialize $T_2$ and obtain the  estimate 
\begin{align*}
	T_2&\,\geq\, \frac{12}{w^{\mid}_\lambda}
y\left(1-\frac{P-P_{-}}{\kappa}\right)+(1+M_{{-1}})\frac{P}{N}
-\frac{\kappa}{N} \geq \left(\frac{12}{w^{\mid}_\lambda}
	y-\frac{\kappa}{N}\right)\left(1-\frac{P}{\kappa}\right).
\end{align*}
A a numerical check shows that for $g\geq 44$ and $y=0.19$ the coefficient
of the previous expression satisfies the desired bound.
\end{proof}

\subsection{Proofs of  \autoref{prop:numericsforstothen}
  and \autoref{prop:numericsforminimal}}
\label{sec:proofandalgo}

We are now ready to prove the main result of this section. First we show the
result for strata with few zeros and for equidistributed strata.
\par
\begin{proof}[Proof of \autoref{prop:numericsforstothen}]
In the case of  strata with \emph{few even order zeros}, we combine \autoref{le:fewzeroT1est} 
and \autoref{le:fewzeroevenT2est} to deduce from~\autoref{eq:gencGammaMainEstimate}
that $s_\Gamma(1/4-\ve) \geq T_1 + T_2 \geq 0$ for almost all strata with few zeros.
\par
For \emph{equidistributed strata with $s$ even}, we similarly
combine \autoref{le:fewzeroT1est} and \autoref{le:equisitrT2est} to obtain the
result. Indeed note that by \autoref{rem:wlambdaapp}, both $y=w_\lambda^{\app}/4$ in the  range of few zeros and $y=w_\lambda^{\app}/6$ in the range of many zeros satisfy the condition \autoref{eq:yrange}.
\par
The rest of the proof deals with the modification for strata with \emph{odd entries}.
In this case we want to improve the $T_2$-bounds of \autoref{le:fewzeroevenT2est}
and \autoref{le:equisitrT2est} by the absolute value of the
lower bound
\begin{equation}\label{eq:twozeroodd}
	y \frac{12 }{w^{\app}_{\lambda}}(w^{\app}_{\Gamma}(\mu)- {w}^{\mid}_\Gamma)\geq
	-y \frac{12 }{8w^{\app}_{\lambda}}\left(M_{-1}^\bot-\frac{\kappa^\bot}{\kappa}M_{-1}\right)
\end{equation}
coming from \autoref{le:wcoeffapp}.
\par
Consider the case of \emph{equidistributed strata with $s$ odd}. In order to
strengthen \autoref{le:equisitrT2est}, we first look at the situation without
rational tail edges. In the range of few zeros, i.e., for $n\leq 2(s+1)$, with $y=w_\lambda^{\app}/4$ (so when
the coefficient of $\kappa^\bot$ in \autoref{eq:T2fewzeros} is positive),  the
negative term \autoref{eq:twozeroodd} is compensated  by the coefficient
$\left(1-\frac{1+M_{{-1}}}{N}\right)$  of $M_{-1}^\bot$ in \autoref{eq:T2fewzeros}.
Indeed  $\frac{1+M_{{-1}}}{N}$ in the range of few zeros is $O(1/g)$, while
$y \frac{12 }{8w^{\app}_{\lambda}}=3/8$. In the range of many zeros and for
$y=w_\lambda^{\app}/6$, we can rewrite the right-hand side of the odd negative
contribution \autoref{eq:twozeroodd} as 
\[-\frac{1}{4}\left(M_{-1}^\bot-\frac{n^\bot}{n}M_{-1}\right)-\frac{P-P_{-}}{4\kappa}M_{-1}\,.\]
One can check that the first line of \autoref{eq:sn} compensates for the
first negative term of the previous displayed expression, while the second line
of \autoref{eq:sn} compensates for the second negative term.
In the presence of rational tail edges, the proof for $s$ even is a numerical check.
Further more, one checks that the additional negative odd contribution
\autoref{eq:twozeroodd}  does not spoil the required positivity. 
\par
In the case of strata with few zeros which are of odd order, we concentrate only in
the situation of strata with two odd zeros, i.e., $\mu=(m, 2g-2-m)$
with $m>0$ odd.
\par
In this case we can follow the same strategy that we used for strata with few zeros. We want to choose  $y$ such that
the coefficient of $\kappa^\bot$ in \autoref{eq:T2fewzeros} is positive, which means
\[y\geq \frac{\kappa}{N}(1+M_{-1})\frac{w_\lambda^{\app}}{12}\,.\]
We choose $y= 1/6$ (which also satisfies the condition imposed
by \autoref{prop:RamGenTypeCrit} of $y> 1/24$), so that the previous condition is
satisfied.  Since for our choice of $y$ the coefficient of $M_{-1}^\bot$
in~\autoref{eq:twozeroodd} becomes $-1/(4w^{\app}_{\lambda})$ and $w^{\app}_{\lambda}\to 1/2$,
we see that if there are no rational tails, the coefficient of $M_{-1}^\bot$
of \autoref{eq:T2fewzeros} compensates this term and hence $T_1+T_2\geq 0$.
\par
In the case of a dumbbell graph with a rational bottom vertex (and hence with
both legs on bottom level), one can numerically check that the expression 
\[T_2-3/P-\frac{1}{4w^{\app}_{\lambda}}M_{-1}\left(1-\frac{\kappa^\bot}{\kappa}\right)\]
is positive. Indeed one can express the previous expression as a rational function
depending on $g$ and $M_{-1}$ and use the bounds $1/2+1/(2g-2)\geq M_{-1}\geq 2/g$ in
order to find an expression only depending on $g$, which is then easy to check to be
positive for large $g$.
%
%
\end{proof}
\par
Now we show the effective estimate for the minimal strata.
\par
\begin{proof}[Proof of \autoref{prop:numericsforminimal}]
Thanks to \autoref{lem:minimalbound}, except for two special HBB cases, since $P/\kappa<1$ we have $s_\Gamma(0.19)\geq T_1+T_2>0$ for $g\geq 44$.
For the  two special HBB graphs, we only need to show that the bound of the
$T_2$~term dominates the negative term $-1/\ell_\Gamma$ in the bound for the $T_1$ term.
\par
In the cases of a double banana HBB graph with two vertices of genus one on top level,
we have $P=4$, hence from~\autoref{eq:T2minimal} we obtain $T_2>1$ for $g\geq 4$, which
is enough to compensate the negative term $-1/\ell_\Gamma=-1$.  
\par
In the case of an HBB banana graph, we can consider the bound $-1/\ell_\Gamma\geq -2/P$
and write $P=2g^\top$. Hence we need to check that the expression
\[T_1+T_2> -\frac{1}{g^\top} +2\frac{g}{g-1}\left(1-\frac{2g^\top}{\kappa}\right)\]
is positive. Using the bounds $1\leq g^\top\leq g-1$, we can check that the expression
is positive for $g\geq 44$.
\par
The justification for the general type statement for $g<44$ is assisted
by computer programs and proved for the ranges of~$y$ in \autoref{cap:Rangeygt}.
In this case we use the version of $D_{\NC}$ given by \autoref{prop:compensationv2}
and its refinement \autoref{prop:nc-refinement}. We remark that
for $g=14$ we need to use the $\NF_{(2g-2)}$ divisor instead of the Hurwitz
divisor $\Hur_\mu$. Moreover, for $g\leq 18$, we need to use the full shape of
the more refined $D_{\NC}$ compensation divisor given by \autoref{prop:nc-refinement}.
In order to use this refinement, we need to be able to list all possible prong
distributions on multi-banana graphs, which is not feasible for large $g$. Therefore, 
we use two different programs for $g\leq 18$  
and $18< g<44$.
Furthermore, it is still not feasible to list all two-level graphs in this range,
hence we give some explanations on how to simplify the check.  
\par
First, for all the ranges of $y$ in \autoref{cap:Rangeygt}, the compact type
contribution of the Brill--Noether divisor or the Hurwitz divisor is larger than
the non-compact type contribution, by an amount that beats the larger
$D_{\NC}$-correction for compact type, with two exceptions: EDB graphs, that
are checked separately, and tails with elliptic top in $g=13$, which require
us to run an additional extra loop over these tails, that is performed on top of the below described procedure. This check about compact-type contributions can be done by comparing the $P_{-1}^{\NCT}$-coefficient and
the $P_{-1}^{\OCT}$-coefficient appearing in~\autoref{eq:T1rewrite}, but we have to use a slightly different expression for them since we now have to use \autoref{prop:nc-refinement}.
\par
With this observation only the edges of $\Gamma$ and their prongs are relevant, not
the full graph structure. More precisely, in the minimal strata
(where $M^\top =  n^\top =0$) the expression~\autoref{eq:gencGammaMainEstimate}
depends on~$v^\top$, $P$ and $P_{-1}$ only, where the latter two implicitly
depend on $g^\top$ and~$E$. For $g\leq 18$ the program checks positivity
of~\autoref{eq:gencGammaMainEstimate} by a loop over $v^\top, g^\top, E$ and all
possible prong distributions.
\par
Second, in the range $18< g<44$, we only consider $v^\top =1$. To justify this,
we check that the coefficient of~$(v^\top-1)$ in~\autoref{eq:gencGammaMainEstimate}
is positive so that we may drop this term. Now, given that this main estimate does
no longer depend on $v^\top$, we may thus compare each graph with the graph where
all top level vertices have been merged to one. This graph is being checked in our
loop, and since we do not need to distinguish between compact type and non-compact type
edges by the first observation, we obtain valid bounds by the merging procedure.
In fact, this merging procedure might turn a graph of a ramified boundary divisor
into an unramified case, but we can check that the positivity of the $(v^\top-1)$-term
outweighs this loss.
\par
Finally, to avoid a loop over all
partitions of~$P$ to cover all prong assignments to the edges, we instead
make a case distinction on the sign of the $P_{-1}$-coefficient as a function
of~$y$. Depending on this sign, the interval of~$y$ that works for all the graphs
with fixed $(E,g^\top)$ (and thus $P= 2g^\top -2 + E$) is only constrained by
the prong distribution that is either most equidistributed or most unbalanced.
The computer program can thus be reduced to a simple loop over all possible
$(E,g^\top)$, this case distinction on the $P_{-1}$-coefficient sign, and checking
additionally the EDB graphs as well as the $D_h$-constraint.
\par
The version of $R$ given in \autoref{prop:compensationv2} has a $v^\top$-term which is not present in the version given in \autoref{eq:non-canonical-term} (which is used for the proof for  $g \geq 44$). 
This term makes the lower bound of \autoref{eq:yrange} bigger, since this bound is the one ensuring that the $v^\top$-coefficient is positive.
As a result, the range given in \autoref{cap:Rangeygt} for $g=44$ does not
include $y=0.19$, that we proved to work for all $g \geq 44$. 
\end{proof}


\printbibliography

\end{document}